\documentclass[11pt,twoside]{article}

\usepackage{jmlr2e}

\usepackage{url}
\usepackage{hyperref}
\usepackage{amstext,latexsym,amsbsy,amsfonts,amsmath,amssymb,amsthm}
\usepackage{algorithm}
\usepackage{algorithmic}
\usepackage{subfig}
\usepackage{booktabs}
\usepackage{xcolor}
\usepackage{bbm}
\usepackage{diagbox}
\usepackage{mathtools}

\DeclarePairedDelimiter{\ceil}{\lceil}{\rceil}

\newcommand{\brackets}[1]{\left[ #1 \right]}
\newcommand{\parenth}[1]{\left( #1 \right)}

\newcommand{\abss}[1]{\left| #1 \right |}

\newcommand{\Rspace}{\ensuremath{\mathbb{R}}}
\newcommand{\Ecal}{\ensuremath{\mathcal{E}}}

\newcommand{\vecnorm}[2]{\left\| #1\right\|_{#2}}
\newcommand{\enorm}[1]{\vecnorm{#1}{2}} 

\newcommand{\Exs}{\ensuremath{{\mathbb{E}}}}
\newcommand{\Prob}{\ensuremath{{\mathbb{P}}}}

\newtheoremstyle{named}{}{}{\itshape}{}{\bfseries}{.}{.5em}{\thmnote{#3's }#1}
\theoremstyle{named}

\theoremstyle{plain}

\newtheorem{theorem}{Theorem}
\newtheorem{proposition}{Proposition}
\newtheorem{lemma}{Lemma}

\newtheorem{definition}{Definition}

\newlength{\widebarargwidth}
\newlength{\widebarargheight}
\newlength{\widebarargdepth}

\makeatletter
\long\def\@makecaption#1#2{
        \vskip 0.8ex
        \setbox\@tempboxa\hbox{\small {\bf #1:} #2}
        \parindent 1.5em  
        \dimen0=\hsize
        \advance\dimen0 by -3em
        \ifdim \wd\@tempboxa >\dimen0
                \hbox to \hsize{
                        \parindent 0em
                        \hfil
                        \parbox{\dimen0}{\def\baselinestretch{0.96}\small
                                {\bf #1.} #2
                                }
                        \hfil}
        \else \hbox to \hsize{\hfil \box\@tempboxa \hfil}
        \fi
        }
\makeatother

\long\def\comment#1{}
\renewcommand\vec[1]{\ensuremath\boldsymbol{#1}}

\newcommand{\Ocal}{\ensuremath{\mathcal{O}}}

\newcommand{\rtil}{\ensuremath{\operatorname{\widetilde{r}}}}

\newcommand{\fcal}{\ensuremath{\overline{f}}}

\newcommand{\nhat}[1]{{\color{red}\textbf{#1} --NH}}

\makeatletter
\newcommand\mathcircled[1]{%
  \mathpalette\@mathcircled{#1}%
}
\newcommand\@mathcircled[2]{%
  \tikz[baseline=(math.base)] \node[draw,circle,inner sep=1pt] (math) {$\m@th#1#2$};%
}
\makeatother

\theoremstyle{plain}

\newtheorem{example}{Example}

\numberwithin{example}{section}

\numberwithin{remark}{section}

\ShortHeadings{Convergence Rates for Gaussian Mixtures of Experts}{Ho, Yang, Jordan}
\firstpageno{1}

\begin{document}

\title{Convergence Rates for Gaussian Mixtures of Experts}

\author{\name Nhat Ho\footnotemark[1] \email minhnhat@utexas.edu  \\
       \addr Division of Statistics and Data Sciences \\
       University of Texas\\
       Austin, TX 78712, USA
       \AND
       \name Chiao-Yu Yang \email chiaoyu@berkeley.edu   \\
       \addr Department of Statistics \\
       University of California\\
       Berkeley, CA 94720-1776, USA
       \AND 
       \name Michael I.\ Jordan \email jordan@cs.berkeley.edu \\
       \addr Division of Computer Science and Department of Statistics\\
       University of California\\
       Berkeley, CA 94720-1776, USA}
\editor{Francis Bach, David Blei, and Bernhard Sch{\"o}lkopf}

\maketitle


%
%
%
%
%
%
%
%
%

\begin{abstract}
We provide a theoretical treatment of over-specified Gaussian mixtures of experts with covariate-free gating networks. We establish the convergence rates of the maximum likelihood estimation MLE) for these models. Our proof technique is based on a novel notion of \emph{algebraic independence} of the expert functions. Drawing on optimal transport, we establish a connection between the algebraic independence of the expert functions and a certain class of partial differential equations (PDEs) with respect to the parameters. Exploiting this connection allows us to derive convergence rates for parameter estimation. 
\end{abstract}
\begin{keywords}
  Mixture of experts, maximum likelihood estimation, convergence rate, optimal transport, partial differential equation, algebraic geometry.
\end{keywords}


\section{Introduction}
Gaussian mixtures of experts, a class of piece-wise regression models introduced by~\citep{Jacob_Jordan-1991, Jordan-1994, Xu_Jordan-1995}, have found applications in many fields including social science~\citep{Huang-2012, Huang-2013, Kitamura-2016}, speech recognition~\citep{Jacobs-1996, Ashok_learning_nips}, natural language processing~\citep{Eigen_learning_2014, Quoc-conf-2017, Ashok_breaking_icml, Ashok_learning_nips}, and system identification~\citep{Rasmussen-2002}. Gaussian mixtures of experts differ from classical finite Gaussian mixture models in two ways. First, the mixture components (the ``experts'') are regression models, linking the location and scale of a Gaussian model of the response variable to a covariate vector $X$ via parametric models $h_{1}(X,\theta_{1})$ and $h_{2}(X,\theta_{2})$, where $\theta_{1}$, $\theta_{2}$ are parameters. Second, the mixing proportions (the ``gating network'') are also functions of the covariate vector $X$, via a parametric model $\pi(X,\gamma)$ that maps $X$ to a probability distribution over the labels of the experts. The overall model can be viewed as  a covariate-dependent finite mixture. Despite their popularity in applications, the theoretical understanding of Gaussian mixtures of experts has proved challenging and lagged behind that of finite mixture models. The inclusion of covariates $X$ in the experts and the gating networks leads to complex interactions of their parameters, which complicates the theoretical analysis of parameter estimation. 

In the setting of finite mixture models, while the early literature focused on identifiability issues~\citep{Teicher-1960, Teicher-1961, Teicher-1963, Lindsay-1995}, recent work has provided a substantive inferential theory; see for example~\citep{Rousseau-Mengersen-11, Nguyen-13, Jonas-2017, Ho-Nguyen-AOS-17, Ho-Nguyen-SIAM-18}. To the best of our knowledge,~\cite{Chen-95} set the stage for these recent developments by establishing a convergence rate of $n^{-1/4}$ for parameter estimation in the univariate setting of over-specified mixture models. Later,~\cite{Nguyen-13} used the Wasserstein metric to analyze the posterior convergence rates of parameter estimation for both finite and infinite mixtures. Recently,~\cite{Ho-Nguyen-SIAM-18} provided a unified framework to rigorously characterize the convergence rates of parameter estimation based on the singularity structures of finite mixture models. Their results demonstrated that there is a connection between the singularities of these models and the algebraic-geometric structure of the parameter space. 

Moving to Gaussian mixtures of experts, a classical line of research focused on the identifiability in these models~\citep{Jiang-1999} and on parameter estimation in the setting of exact-fitted models where the true number of components is assumed known~\citep{Tanner-conf-1999, Tanner-1999, Tanner-Neural-1999}. This assumption is, however, overly strong for most applications; the true number of components is rarely known in practice. There are two common practical approaches to deal with this issue. The first approach relies on model selection, most notably the BIC penalty~\citep{Wang-1996, Khalili-2007}.  This approach is, however, computationally expensive as we need to search for the optimal number of components over all the possible values. Furthermore, the sample size may not be large enough to support this form of inference. The second approach is to over-specify the true model, by using rough prior knowledge to specify more components than is necessary. However, theoretical analysis of parameter estimation is challenging in the over-specified setting, given the complicated interaction among the parameters of the expert functions, a phenomenon that does not occur in the exact-fitted setting of Gaussian mixtures of experts. Another challenge arises from \emph{inhomogeneity}---some parameters tend to have faster convergence rates than other parameters. This inhomogeneity makes it nontrivial to develop an appropriate distance for characterizing convergence rates. 

In the current paper we focused on a simplified setting in which the expert functions are covariate-dependent, but the gating network is not. We refer to this as the \emph{Gaussian mixture of experts with covariate-free gating functions} (GMCF) model. This model is also widely known as finite Gaussian mixture of regression~\citep{Khalili-2007}. Although simplified, this model captures the core of the mixtures-of-experts problem, which is the interactions among the different mixture components.  We believe that the general techniques that we develop here can be extended to the full mixtures-of-experts model---in particular by an appropriate generalization of the transportation distance to capture the variation of parameters from the gating networks---but we leave the development of that direction to future work. 

\subsection{Setting} We propose a general theoretical framework for analyzing the statistical performance of maximum likelihood estimation (MLE) for parameters in the setting of over-specified Gaussian mixtures of experts with covariate-free gating functions. In particular, we assume that $(X_{1},Y_{1}), \ldots, (X_{n}, Y_{n})$ are i.i.d.\ samples from a Gaussian mixture of experts with covariate-free gating functions (GMCF) of order $k_{0}$, with conditional density function $g_{G_{0}}(Y | X)$:
\begin{align}
\label{eq:true_conditional_density}
    g_{G_{0}}(Y | X) : = \sum \limits_{i=1}^{k_{0}}{\pi_{i}^0f(Y|h_{1}(X,\theta_{1i}^0), h_{2}(X,\theta_{2i}^0))},
\end{align}
where $G_{0} : = \sum_{i=1}^{k_{0}}{\pi_{i}^{0}\delta_{(\theta_{1i}^{0}, \theta_{2i}^{0})}}$ is a true but unknown probability measure (mixing measure) and $\theta_{ji}^{0} \in \Omega_{j} \subset \mathbb{R}^{q_{j}}$ for all $i, j$. Furthermore, we denote $\left\{f(\cdot|\mu,\sigma)\right\}$ as the family of location-scale univariate Gaussian distributions. We over-specify the true model by choosing $k > k_{0}$ components. 

We estimate $G_{0}$ under the over-specified GMCF model via maximum likelihood estimation (MLE). We denote the MLE as $\widehat{G}_{n}$.  Our results reveal a fundamental connection between the algebraic structure of the expert functions $h_{1}$ and $h_{2}$ and the convergence rates of the MLE through a general version of the optimal transport distance, which refers to as the \emph{generalized transportation distance}. A similar distance has been used to study the effect of algebraic singularities on parameter estimation in classical finite mixtures~\citep{Ho-Nguyen-SIAM-18}.

\subsection{Generalized transportation distance}
\label{subsection:generalized_Wasserstein_metric}
In contrast to the traditional Wasserstein metric~\citep{Villani-03}, the generalized transportation distance assigns different orders to each parameter. This special property of generalized transportation distance provides us with a tool to capture the inhomogeneity of parameter estimation in Gaussian mixtures of experts. In order to define the generalized transportation distance, we first define the semi-metric $d_{\kappa}(.,.)$ for any vector $\kappa = \parenth{ \kappa_{1},\ldots,\kappa_{q_{1}+q_{2}}} \in \mathbb{N}^{q_{1}+q_{2}}$ as follows:
\begin{eqnarray}
d_{\kappa}(\theta_{1},\theta_{2}) : = \biggr(\sum \limits_{i=1}^{q_{1} + q_{2}}{|\theta_{1}^{(i)}-\theta_{2}^{(i)}|^{\kappa_{i}}}\biggr)^{1/\| \kappa \|_{\infty}}, \nonumber
\end{eqnarray}
for any $\theta_{i} = \parenth{ \theta_{i}^{(1)},\ldots,\theta_{i}^{(q_{1}+q_{2})}} \in \mathbb{R}^{q_{1} + q_{2}}$. Generally, $d_{\kappa}(.,.)$ does not satisfy the standard triangle inequality. More precisely, when not all $\kappa_{i}$ are identical, $d_{\kappa}$ satisfies a triangle inequality only up to some positive constant less than one. When all $\kappa_{i}$ are identical, $d_{\kappa}$ becomes a metric.

Now, we let $G = \sum_{i=1}^{k} \pi_{i}\delta_{(\theta_{1i},\theta_{2i})}$ be some probability measure. The generalized transportation distance between $G$ and $G_{0}$ with respect to $\kappa = (\kappa_{1},\ldots,\kappa_{q_{1}+q_{2}}) \in \mathbb{N}^{q_{1}+q_{2}}$ is given by:
\begin{eqnarray}
\widetilde{W}_{\kappa}(G,G_{0}) := \biggr (\inf \sum_{i,j} q_{ij}
d_{\kappa}^{\| \kappa \|_{\infty}}(\eta_i ,\eta_j^0) \biggr )^{1/\| \kappa \|_{\infty}}, \label{eq:generalized_transport}
\end{eqnarray}
where the infimum is taken over all couplings $\vec{q}$ between $\vec{\pi}$ and $
\vec{\pi}^0$; i.e., where $\sum_{j} q_{ij} = \pi_{i}$ and $\sum_{i} q_{ij} = \pi_{j}^{0}$. Additionally, $\eta_{i} = (\theta_{1i}, \theta_{2i})$ and $\eta_{j}^{0} = (\theta_{1j}^{0}, \theta_{2j}^{0})$ for all $i,j$. When $\kappa = (2, \ldots, 2)$, we can check that $\widetilde{W}_{\kappa}(G,G_{0}) \equiv W_{2}(G, G_{0})$, the second order Wasserstein metric~\citep{Villani-03}.

In general, the convergence rates of mixing measures under generalized Wasserstein distance translate directly to the convergence rates of their associated atoms or parameters. More precisely, assume that there exist a sequence $\{G_{n}\}$ and a vector $\kappa = (\kappa_{1},\ldots,\kappa_{q_{1}+q_{2}}) \in \mathbb{N}^{q_{1} + q_{2}}$ such that $\widetilde{W}_{\kappa}(G_{n},G_{0}) \to 0$ at rate $\omega_{n} = o(1)$ as $n \to \infty$. Then, we can find a sub-sequence of $G_{n}$ such that each atom (support) $(\theta_{1i}^{0},\theta_{2i}^{0})$ of $G_{0}$ is the limit point of atoms of $G_{n}$. Additionally, the convergence rates for estimating $(\theta_{1i}^{0})^{(u)}$, the $u$th component of $\theta_{1i}^{0}$, are $\omega_{n}^{\|\kappa\|_{\infty}/\kappa_{u}}$ while those for estimating $(\theta_{2i}^{0})^{(v)}$ are $\omega_{n}^{\|\kappa\|_{\infty}/\kappa_{q_{1}+v}}$ for $1 \leq u \leq q_{1}$ and $1 \leq v \leq q_{2}$. Furthermore, the convergence rates for estimating the weights associated with these parameters are $\omega_{n}^{\|\kappa\|_{\infty}}$. Finally, there may exist some atoms of $G_{n}$ that converge to limit points outside the atoms of $G_{0}$. The convergence rates of these limit points are also similar to those for estimating the atoms of $G_{0}$. 
\subsection{Main contribution} The generalized transportation distance in equation~\eqref{eq:generalized_transport} allows us to introduce a notion of \emph{algebraic independence} between expert functions $h_{1}$ and $h_{2}$ that is expressed in the language of partial differential equations (PDEs). Using this notion, we are able to characterize the convergence rates of parameter estimation for several choices of expert functions $h_{1}$ and $h_{2}$ when they are either algebraically independent or not. Our overall contributions in the paper can be summarized as follows:
\begin{itemize}
\item \textbf{Algebraically independent settings:} When the expert functions $h_{1}$ and $h_{2}$ are algebraically independent, we establish the best possible convergence rate of order $n^{-1/4}$ for $\widetilde{W}_{\kappa}(\widehat{G}_{n}, G_{0})$ (up to a logarithmic factor) where $\kappa = (2, \ldots, 2)$. That result directly translates to a convergence rate of $n^{-1/ 4}$ for the support of $\widehat{G}_{n}$. 
\item \textbf{Algebraically dependent settings:} When the expert functions $h_{1}$ and $h_{2}$ are algebraically dependent, we prove that the convergence rates of parameter estimation are very slow and inhomogeneous. More precisely, the rates of convergence are either determined by the solvability of a system of polynomial equations or by the admissibility of a system of polynomial limits. The formulations of these systems depend on the PDEs that capture the interactions among the parameters for the expert functions. Furthermore, we show that the inhomogeneity of parameter estimation can be characterized based on the generalized transportation distance. 
\end{itemize}
We note in passing that while our results in the paper are specifically for the MLE, the proof technique and algebraic independence notion can also be used to analyze the convergence rate of parameter estimation from moment methods~\citep{Anima_2012, Anima-2014} with the over-specified GMCF model. 

\textbf{Organization.} The remainder of the paper is organized as follows. In Section~\ref{Section:preliminary}, we introduce the problem setup for Gaussian mixtures of experts with covariate-free gating functions. Section~\ref{Section:weak_interaction_expert} establishes convergence rates for parameter estimation under the algebraically independent setting. In Section~\ref{Section:strong_interaction_expert}, we consider various settings in which the expert functions are algebraically dependent and establish the convergence rates of parameter estimation under these settings. We provide proofs for a few key results in Section~\ref{Section:proofs} while deferring the majority of the proofs to the Appendices. Finally, we conclude in Section~\ref{Section:discussion}.

\textbf{Notation.} For any vector $x \in \mathbb{R}^{d}$, we use superscript and subscript notation interchangeably, letting $x = (x^{(1)},\ldots,x^{(d)})$ or $x = (x_{1},\ldots,x_{d})$.  Thus, either $x^{(i)}$ or $x_{i}$ is the $i$-th component of $x$. For each $x \in \mathbb{R}^{d}$, we denote $x^{\kappa} = \prod \limits_{i=1}^{d} (x^{(i)})^{\kappa_{i}}$ for any $\kappa = (\kappa^{(1)},\ldots, \kappa^{(d)}) \in \mathbb{N}^{d}$. For any two vectors $x,y \in \mathbb{R}^{d}$, we write $x \preceq y$ if $x^{(i)} \leq y^{(i)}$ for all $1 \leq i \leq d$ and $x \prec y$ if $x \preceq y$ and $x \neq y$. For any two sequences $\{a_{n}\}$ and $\{b_{n}\}$, the notation $a_{n} \precsim b_{n}$ means $a_{n} \leq C b_{n}$ for all $n \geq 1$ where $C$ is some constant independent of $n$. 

For any two density functions $p, q$ (with respect to the Lebesgue measure $\mu$), the total variation distance is given by $V(p,q)={\displaystyle \frac{1}{2} \int {|p(x)-q(x)|}\textrm{d}\mu(x)}$.
The squared Hellinger distance is defined as $h^{2}(p,q)=
{\displaystyle \frac{1}{2} \int (\sqrt{p(x)}-\sqrt{q(x)})^{2}\textrm{d}\mu(x)}$.

\section{Background} \label{Section:preliminary}
In this section, we provide the necessary background for our analysis of the convergence rates of the MLE under over-specified Gaussian mixtures of experts with covariate-free gating functions. In particular, in Section~\ref{sec:Gaussian_expert_form}, we define the over-specified Gaussian mixture of experts with covariate-free gating functions, and in Section~\ref{sec:iden_smooth}, we establish identifiability and smoothness properties for these models as well as establishing the convergence rates of density estimation.
\subsection{Problem setup}
\label{sec:Gaussian_expert_form}
Let $Y \in \mathcal{Y} \subset \mathbb{R}$ be a response variable of interest and let $X \in \mathcal{X} \subset \mathbb{R}^{d}$ be a vector of covariates believed to have an effect on $Y$. We start with a definition of identifiable expert functions.
\begin{definition}
\label{definition:identifiability}
Given $\Theta \subset \mathbb{R}^{q}$ for some $q \geq 1$. We say that an expert function $h_{1}: \mathcal{X} \times \Theta \to \mathbb{R}$ is \emph{identifiable} if for each $k \in \mathbb{N}$ the following holds:
\begin{itemize}
    \item[(I.1)] If there exist distinct parameters $(\eta_{1}, \ldots, \eta_{k})$ and $(\eta_{1}', \ldots, \eta_{k}')$ such that for almost surely $X \in \mathcal{X}$, we can find permutation function $\sigma_{X}:\{1,2,\ldots, k\} \to \{1,2,\ldots,k\}$ to satisfy $h_{1}(X, \eta_{\sigma_{X}(i)}) = h_{1}(X, \eta_{i}')$ for all $1 \leq i \leq k$, then $\{\eta_{1}, \ldots, \eta_{k}\} \equiv \{\eta_{1}', \ldots, \eta_{k}'\}$.
\end{itemize}
\end{definition}
Examples of identifiable expert function $h_{1}$ include $h_{1}(X, \eta) = g(\eta^{\top}[1,X])$ for any injective function $g: \mathbb{R} \to \mathbb{R}$ where $X \in \mathbb{R}^{d}$ and $\eta \in \mathbb{R}^{d + 1}$. Recall that we focus on Gaussian mixtures of experts~\citep{Jacob_Jordan-1991, Jordan-1994, Xu_Jordan-1995} for which the gating functions are independent of covariate $X$. We now start with the following definition of Gaussian mixtures of experts with covariate-free gating functions. 
\begin{definition} \label{definition_finite_mixture_experts}
Assume that we are given two identifiable expert functions $h_{1}:\mathcal{X} \times \Omega_{1} \to \Theta_{1} \subset \mathbb{R}$ and $h_{2}:\mathcal{X} \times \Omega_{2} \to \Theta_{2} \subset \mathbb{R}_{+}$ where $\Omega_{i} \subset \mathbb{R}^{q_{i}}$ for given dimensions $q_{i} \geq 1$ as $1 \leq i \leq 2$. Let $\{\pi_{i}\}_{i=1}^{k}$ denote $k$ weights with $\sum \limits_{i=1}^{k} \pi_{i} = 1$. We say that $(X,Y)$ follows a \emph{Gaussian mixtures of experts with covariate-free gating functions} (GMCF) of order $k$, with respect to expert functions $h_{1}$, $h_{2}$ and gating functions $\pi_{i}$, if the conditional density function of $Y$ given $X$ has the following form
\begin{align}
g_{G}(Y|X) & : = \int f\parenth{Y|h_{1}(X,\theta_{1}),h_{2}(X,\theta_{2})}dG(\theta_{1},\theta_{2}) \nonumber \\ 
& = \sum \limits_{i=1}^{k}{\pi_{i}f(Y|h_{1}(X,\theta_{1i}), h_{2}(X,\theta_{2i}))}, \nonumber
\end{align}
where $G = \sum \limits_{i=1}^{k}\pi_{i}\delta_{(\theta_{1i},\theta_{2i})}$ is a discrete probability measure that has exactly $k$ atoms on $\Omega := \Omega_{1} \times \Omega_{2}$.
\end{definition}
 As an example, when $q_{1} = q_{2} = d+1$, generalized linear expert functions take the form $h_{1}(X,\theta_{1})=\theta_{1}^{\top}[1,X]$ and $h_{2}(X,\theta_{2})=\exp\parenth{\theta^{\top}_2[1,X]}$.
\paragraph{Over-specified GMCF:} Let $(X_{1},Y_{1}), \ldots, (X_{n}, Y_{n})$ be i.i.d.\ draws from a GMCF of order $k_{0}$ with conditional density function $g_{G_{0}}(Y|X)$ where $G_{0} : = \sum_{i=1}^{k_{0}}{\pi_{i}^{0}\delta_{(\theta_{1i}^{0}, \theta_{2i}^{0})}}$ is a true but unknown probability measure (mixing measure). Since $k_{0}$ is generally unknown in practice, one popular approach to estimate the mixing measure $G_{0}$ is based on over-specifying the true number of components $k_{0}$. In particular, we fit the true model with $k > k_{0}$ number of components where $k$ is a given threshold that is chosen based on prior domain knowledge. We refer to this setting as the \textit{over-specified GMCF}. 
\paragraph{Maximum likelihood estimation (MLE):} To obtain an estimate of $G_{0}$, we define the MLE as follows:
\begin{eqnarray}
\widehat{G}_{n} : = \mathop {\arg \max}\limits_{G \in \mathcal{G}}{\sum \limits_{i=1}^{n}{\log(g_{G}(Y_{i}|X_{i})}}), \label{eqn:MLE_formulation}
\end{eqnarray}
where $\mathcal{G}$ is some subset of $\Ocal_{k}(\Omega) : = \{G = \sum_{i=1}^{l} \pi_{i}\delta_{(\theta_{1i},\theta_{2i})}: 1 \leq l \leq k \}$, namely, the set of all discrete probability measures with at most $k$ components. Detailed formulations of $\mathcal{G}$ will be given later based on the specific structures of expert functions $h_{1}$ and $h_{2}$. 
\paragraph{Universal assumptions and notation:} Throughout this paper, we assume that $\Omega_{1}$ and $\Omega_{2}$ are compact subsets of $\Rspace^{q_{1}}$ and $\Rspace^{q_{2}}$ respectively. Additionally, $\Omega := \Omega_{1} \times \Omega_{2}$ and $X$ is a random vector and has a given prior density function $\fcal(X)$, which is independent of the choices of expert functions $h_{1}$, $h_{2}$. Furthermore, $\mathcal{X}$ is a fixed compact set of $\Rspace^{d}$. Finally we denote
\begin{align*}
p_{G}(X,Y) : = g_{G} (Y|X) \overline{f}(X)
\end{align*} 
as the joint distribution (or equivalently mixing density) of $X$ and $Y$ for any $G \in \Ocal_{k}(\Omega)$.
\subsection{General identifiability, smoothness condition, and density estimation}
\label{sec:iden_smooth}
In order to establish the convergence rates of $\widehat{G}_{n}$, our analysis relies on three main ingredients: general identifiability of the GMCF, H\"{o}lder continuity of the GMCF up to any order $r \geq 1$, and parametric convergence rates for density estimation under the over-specified GMCF. We begin with the following result regarding the identifiability of GMCF.
\begin{proposition}\label{proposition:identifiability_fixed_weights_setting}
For given identifiable expert functions $h_{1}$ and $h_{2}$, the GMCF is identifiable with respect to $h_{1}$ and $h_{2}$, namely, whenever there are finite discrete probability measures $G$ and $G'$ on $\Omega$ such that $p_{G}(X,Y)=p_{G'}(X,Y)$ almost surely $(X,Y) \in \mathcal{X} \times \mathcal{Y}$, then it follows that $G \equiv G'$.
\end{proposition}
The proof of Proposition~\ref{proposition:identifiability_fixed_weights_setting} is in Appendix~\ref{Section:proof_identifiable}. A second result that plays a central role in analyzing convergence of the MLE in over-specified GMCF is the uniform H\"{o}lder continuity, formulated as follows:
\begin{proposition} \label{proposition:Lipschitz_continuity}
For any $r \geq 1$, the GMCF admits the uniform H\"{o}lder continuity up to the $r$th order, with respect to the expert functions $h_{1}$, $h_{2}$ and prior density function $\overline{f}$:
\begin{eqnarray}
& & \hspace{-3 em} \sum \limits_{|\kappa|=r} \overline{f}(x) \biggr| \biggr(\dfrac{\partial^{|\kappa|}f}{\partial{\theta_{1}^{\kappa_{1}}} \partial{\theta_{2}^{\kappa_{2}}}} \parenth{y|h_{1}(x,\theta_{1}), h_{2}(x,\theta_{2})} \nonumber \\
& & \hspace{1 em } - \dfrac{\partial^{|\kappa|}f}{\partial{\theta_{1}^{\kappa_{1}}} \partial{\theta_{2}^{\kappa_{2}}}} f \parenth{y|h_{1}(x,\theta_{1}'), h_{2}(x,\theta_{2}')} \biggr) \gamma^{\kappa} \biggr| \leq C \|(\theta_{1}, \theta_{2}) - (\theta_{1}', \theta_{2}')\|^{\delta}\|\gamma\|^{r}, \nonumber
\end{eqnarray}
for any $\gamma \in \mathbb{R}^{q_{1}+q_{2}}$ and for some positive constants $\delta$ and $C$ that are independent of $x, y$ and $(\theta_{1}, \theta_{2}), (\theta_{1}', \theta_{2}') \in \Omega$. Here, $\kappa = (\kappa_{1}, \kappa_{2}) \in \mathbb{N}^{q_{1}+ q_{2}}$ where $\kappa_{i} \in \mathbb{N}^{q_{i}}$ for any $1 \leq i \leq 2$.
\end{proposition}
Finally, when the expert functions $h_{1}$ and $h_{2}$ are sufficiently smooth in terms of their parameters, we can guarantee the parametric convergence rate of density estimation.
\begin{proposition} \label{proposition:convergence_rates_density_estimation_MECFG}
Assume that the expert functions $h_{1}$ and $h_{2}$ are twice differentiable with respect to their parameters. Additionally, assume that there exist positive constants $a, \underline{\gamma}, \overline{\gamma}$ such that $h_{1}(X,\theta_{1}) \in [-a,a]$, $h_{2}(X,\theta_{2}) \in [\underline{\gamma},\overline{\gamma}]$ for all $X \in \mathcal{X}, \theta_{1} \in \Omega_{1}, \theta_{2}\in \Omega_{2}$. Then, the following holds:
\begin{align}
\Prob(h(p_{\widehat{G}_{n}},p_{G_{0}}) > C(\log n/n)^{1/2}) \precsim \exp(-c \log n) \label{eqn:density_estimation_rootn_rate}
\end{align}
for universal positive constants $C$ and $c$ that depend only on $\Omega$.
\end{proposition}
\noindent
The proof of Proposition~\ref{proposition:convergence_rates_density_estimation_MECFG} is in Appendix~\ref{Section:Appendix_C}.
\section{Algebraically independent expert functions}
\label{Section:weak_interaction_expert}
In this section, we consider the MLE in equation~\eqref{eqn:MLE_formulation} over the entire parameter space $\Ocal_{k}(\Omega)$.  That is, we let $\mathcal{G} = \Ocal_{k}(\Omega)$. To analyze the convergence rates of MLE under over-specified GMCF we capture the algebraic interaction among the expert functions $h_{1}$ and $h_{2}$ via the following definition. 
\begin{definition} \label{definition:second_order_MECFG}
We say that the expert functions $h_{1}, h_{2}$ are \emph{algebraically independent} if they are twice differentiable with respect to their parameters $\theta_{1}$ and $\theta_{2}$ and the following holds:
\begin{itemize}
\item[(O.1)] For any $(\theta_{1}, \theta_{2})$, if we have $\alpha_{i}, \beta_{uv} \in \mathbb{R}$ (for $1 \leq i \leq q_{2}$, and $1 \leq u \leq v \leq q_{1}$) such that
\begin{align}
\sum \limits_{i = 1}^{q_{2}} \alpha_{i}\dfrac{\partial{h_{2}^{2}}}{\partial{\theta_{2}^{(i)}}}(X,\theta_{2}) + \sum \limits_{1 \leq u \leq v \leq q_{1}} \beta_{uv} \dfrac{\partial{h_{1}}}{\partial{\theta_{1}^{(u)}}}(X,\theta_{1})\dfrac{\partial{h_{1}}}{\partial{\theta_{1}^{(v)}}}(X,\theta_{1}) = 0, \nonumber
\end{align}
almost surely in $X$, then we must also have $\alpha_{i}= \beta_{uv} = 0$ for all $1 \leq i \leq q_{2}$ and $1 \leq u \leq v \leq q_{1}$.
\item[(O.2)] For any $\theta_{2}$, if we have $\gamma_{uv} \in \mathbb{R}$ (for $1 \leq u \leq v \leq q_{2}$) such that
\begin{align*}
    \sum_{1 \leq u \leq v \leq q_{2}} \gamma_{uv} \dfrac{\partial{h_{2}^2}}{\partial{\theta_{2}^{(u)}}}(X,\theta_{2})\dfrac{\partial{h_{2}^2}}{\partial{\theta_{2}^{(v)}}}(X,\theta_{2}) = 0,
\end{align*}
almost surely in $X$, then we have $\gamma_{uv} = 0$ for all $1 \leq u \leq v \leq q_{2}$.
\item[(O.3)] For any $(\theta_{1}, \theta_{2})$, if we have $\eta_{uv} \in \mathbb{R}$ (for $1 \leq u \leq q_{1}$ and $1 \leq v \leq q_{2}$) such that
\begin{align*}
    \sum_{u = 1}^{q_{1}} \sum_{v = 1}^{q_{2}} \eta_{uv} \dfrac{\partial{h_{1}}}{\partial{\theta_{1}^{(u)}}}(X,\theta_{1})\dfrac{\partial{h_{2}^2}}{\partial{\theta_{2}^{(v)}}}(X,\theta_{2}) = 0,
\end{align*}
almost surely in $X$, then we have $\eta_{uv} = 0$ for all $1 \leq u \leq q_{1}$ and $1 \leq v \leq q_{2}$.
\end{itemize}
\end{definition}
Note that in this definition we use the convention that if $\dfrac{\partial{h_{2}^{2}}}{\partial{\theta_{2}^{(i)}}}(X,\theta_{2}) = 0$ almost surely for some $1 \leq i \leq q_{2}$, then we have $\alpha_{i} = 0$. The same convention goes for other derivatives in Conditions (O.1), (O.2), and (O.3). An equivalent way to express the Condition (O.1) in Definition~\ref{definition:second_order_MECFG} is that the elements in a set of partial derivatives,
\begin{align}
\left\{\dfrac{\partial{h_{1}}}{\partial{\theta_{1}^{(u)}}}(X,\theta_{1})\dfrac{\partial{h_{1}}}{\partial{\theta_{1}^{(v)}}}(X,\theta_{1}), \dfrac{\partial{h_{2}^{2}}}{\partial{\theta_{2}^{(i)}}}(X,\theta_{2}): \ 1 \leq i \leq q_{2}, 1 \leq u \leq v \leq q_{1} \right\}, \nonumber
\end{align}
are linearly independent with respect to $X$. Similarly, the Conditions (O.2) and (O.3) indicate that the elements of the following sets of partial derivatives
\begin{align*}
    \left\{\dfrac{\partial{h_{2}^2}}{\partial{\theta_{2}^{(u)}}}(X,\theta_{2})\dfrac{\partial{h_{2}^2}}{\partial{\theta_{2}^{(v)}}}(X,\theta_{2}): \ 1 \leq u \leq v \leq q_{2} \right\}, \\
    \left\{\dfrac{\partial{h_{1}}}{\partial{\theta_{1}^{(u)}}}(X,\theta_{1})\dfrac{\partial{h_{2}^2}}{\partial{\theta_{2}^{(v)}}}(X,\theta_{2}): \ 1 \leq u \leq q_{1}, \ 1 \leq v \leq q_{2} \right\}
\end{align*}
are linearly independent with respect to $X$. To exemplify the algebraic independence notion in Definition~\ref{definition:second_order_MECFG}, we consider the following simple examples of expert functions $h_{1}$ and $h_{2}$ that are algebraically independent.

\begin{example} \label{example:not_second_order_MECFG}
(a) Let $\mathcal{X} \subset \Rspace^{d}$. If we choose expert functions $h_{1}(X,\theta_{1}) = \theta_{1}^{\top} X$ and $h_{2}^{2}(X,\theta_{2}) = \theta_{2}$ for all $\theta_{1} \in \Omega_{1} \subset \Rspace^{d}$ and $\theta_{2} \in \Omega_{2} \subset \Rspace_{+}$, then $h_{1}$ and $h_{2}$ are algebraically independent.

(b) Let $\mathcal{X} \subset \Rspace_{+}$. If we choose expert functions $h_{1}(X, \theta_{1}) = (\theta_{1}^{(1)}+\theta_{1}^{(2)}X)^{m}$ for all $\theta_{1} = (\theta_{1}^{(1)},\theta_{1}^{(2)}) \in \Omega_{1} \subset \mathbb{R}^{2}$, where $m>1$ and $h_{2}^{2}(X, \theta_{2}) = \theta_{2}X$ for all $\theta_{2} \in \Omega_{2} \subset \mathbb{R}_{+}$, then $h_{1}, h_{2}$ are algebraically independent.

\end{example}
There are also several standard settings of GMCF  where the expert functions $h_{1}$ and $h_{2}$ are algebraically dependent (See Section~\ref{Section:strong_interaction_expert} for these examples). For instance, for the standard univariate mixtures of Gaussian distributions, namely, when the expert functions $h_{1}(X, \theta_{1}) = \theta_{1}$ and $h_{2}^2(X, \theta_{2}) = \theta_{2}$ for all $\theta_{1}, \theta_{2}$, then these expert functions are algebraically dependent as they violate Condition (O.1) in Definition~\ref{definition:second_order_MECFG}, which is due to following PDE: $\parenth{\frac{\partial{h_{1}}}{\partial{\theta_{1}}}(X,\theta_{1})}^{2} = \dfrac{\partial{h_{2}^{2}}}{\partial{\theta_{2}}}(X,\theta_{2})$ for all $\theta_{1}, \theta_{2}$. The convergence rate of MLE for the over-specified Gaussian mixtures had been established in~\citep{Ho-Nguyen-AOS-17}. Another example of GMCF when the expert functions are algebraically dependent is the Gaussian mixture of regression~\citep{Khalili-2007}, namely, $h_{1}(X, \theta_{1}) = \theta_{1}^{\top}[1, X]$ and $h_{2}^2(X, \theta_{2})$ for all $\theta_{1} = (\theta_{1}^{(1)}, \ldots, \theta_{1}^{(d + 1)})$ and $\theta_{2}$. These expert functions also violate Condition (O.1) due to the PDE: $\parenth{\frac{\partial{h_{1}}}{\partial{\theta_{1}}^{(1)}}(X,\theta_{1})}^{2} = \dfrac{\partial{h_{2}^{2}}}{\partial{\theta_{2}}}(X,\theta_{2})$. That relation between $\theta_{1}^{(1)}$ and $\theta_{2}$ leads to distinct behaviors of the elements of $\theta_{1}$ in the convergence rates of MLE (see Theorem~\ref{theorem:lower_bound_Gaussian_family} for a detailed statement). Further examples of algebraically dependent functions as well as convergence rates of their parameter estimation are provided in Section~\ref{Section:strong_interaction_expert}.

Going back to the algebraic independence condition for the expert functions $h_{1}$ and $h_{2}$, we have the following result regarding the convergence rates of parameter estimation $\widehat{G}_{n}$ under the over-specified GMCF model.
\begin{theorem} \label{theorem:total_variation_bound_over-fitted_MECFG}
Assume that expert functions $h_{1}$ and $h_{2}$ are algebraically independent and are twice differentiable with respect to their parameters. Additionally, assume that there exist positive constants $a, \underline{\gamma}, \overline{\gamma}$ such that $h_{1}(X,\theta_{1}) \in [-a,a]$, $h_{2}(X,\theta_{2}) \in [\underline{\gamma},\overline{\gamma}]$ for all $X \in \mathcal{X}, \theta_{1} \in \Omega_{1}, \theta_{2}\in \Omega_{2}$. Then, the following holds:
\begin{itemize}
\item[(a)] (Convergence rate of MLE) There exists a positive constant $C_{0}$ depending on $G_{0}$ and $\Omega$ such that
\begin{eqnarray}
\mathbb{P}(\widetilde{W}_{\kappa}(\widehat{G}_{n},G_{0}) > C_{0} (\log n/n)^{1/4}) \precsim \exp(-c\log n), \nonumber
\end{eqnarray}
where $\kappa = (2,\ldots,2) \in \mathbb{N}^{q_{1}+q_{2}}$ and $c$ is a positive constant depending only on $\Omega$.
\item[(b)] (Minimax lower bound) For any $\kappa'$such that $(1,\ldots,1) \preceq \kappa' \prec \kappa = (2, \ldots, 2)$, 
\begin{align}
\inf \limits_{\overline{G}_{n}} \sup \limits_{G \in \Ocal_{k}(\Omega) \backslash \Ocal_{k_{0}-1}(\Omega)} \Exs_{p_{G}} \parenth{\widetilde{W}_{\kappa'}(\overline{G}_{n},G)} \geq c'n^{-1/(2\|\kappa'\|_{\infty})}. \nonumber
\end{align}
Here, the infimum is taken over all sequences of estimates $\overline{G}_{n} \in \Ocal_{k}(\Omega)$. Furthermore, $\Exs_{p_{G}}$ denotes the expectation taken with respect to the product measure with mixture density $p_{G}^{n}$, and $c'$ stands for a universal constant depending only on $\Omega$. 
\end{itemize}
\end{theorem}
\noindent
The proof of Theorem~\ref{theorem:total_variation_bound_over-fitted_MECFG} is in Section~\ref{Section:proof_strong_identifiable}.

A few comments are in order. First, part (a) of Theorem~\ref{theorem:total_variation_bound_over-fitted_MECFG} establishes a best possible convergence rate of $n^{-1/4}$ (up to a logarithmic factor) of $\widehat{G}_{n}$ to $G_{0}$ under the generalized transportation distance $\widetilde{W}_{\kappa}$ while part (b) of that theorem demonstrates that this rate is sharp. The convergence rate $n^{-1/4}$ of $\widehat{G}_{n}$ suggests that the rate of estimating individual components $(\theta_{1i}^{0})^{(u)}$ and $(\theta_{2i}^{0})^{(v)}$ is $n^{-1/4}$ for $1 \leq u \leq q_{1}$ and $1 \leq v \leq q_{2}$. The main reason for these slow convergence rates is the singularity of Fisher information matrix for these components. Such a singularity phenomenon is caused by the effect of fitting the true model by larger model, a phenomenon which has been observed previously in traditional mixture models settings under strong identifiability~\citep{Chen-95, Nguyen-13}. 

Second, we would like to emphasize that Theorem~\ref{theorem:total_variation_bound_over-fitted_MECFG} is not only of theoretical interest. Indeed, it provides insight into the choice of expert functions that are likely to have favorable convergence in practice. When the expert functions are not algebraically independent, we demonstrate in the next section that the convergence rates of parameter estimation in over-specified GMCF are very slow and depend on a notion of complexity level of over-specification. 

\section{Algebraically dependent expert functions}
\label{Section:strong_interaction_expert}
In the previous section we established a convergence rate $n^{-1/4}$ for the MLE when the expert functions $h_{1}$ and $h_{2}$ are algebraically independent. In many scenarios, however, the expert functions are taken to be \textit{algebraically dependent}. Here we show that in some of these settings the convergence rates of the MLE can be much slower than $n^{-1/4}$.

In order to simplify our proofs in the algebraically-dependent cases, in this section we only consider single covariate settings, i.e., $X \in \mathcal{X} \subset \mathbb{R}$. With the similar proof technique, most of the results in this section can be generalized to their corresponding multivariate covariate settings, i.e., $X \in \mathcal{X} \subset \mathbb{R}^{d}$. Furthermore, we focus on the case in which the MLE is restrained to a parameter space $\mathcal{G}$ that has the following structure:
\begin{align}
\mathcal{G} = \Ocal_{k, \bar{c}_{0}}(\Omega) = \{G = \sum_{i=1}^{l} \pi_{i}\delta_{(\theta_{1i},\theta_{2i})}: 1 \leq l \leq k \ \text{and} \ \pi_{i} \geq \overline{c}_{0} \ \forall i \}. \nonumber
\end{align}
That is, we consider the set of discrete probability measures with at most $k$ components such that their weights are lower bounded by $\overline{c}_{0}$ for some given sufficiently small positive number $\overline{c}_{0}$. Under this assumption, the true but unknown mixing measure $G_{0} = \sum \limits_{i=1}^{k_{0}}{\pi_{i}^{0}\delta_{(\theta_{1i}^{0},\theta_{2i}^{0})} \in \Ecal_{k_{0}}(\Omega)}$ is assumed to have $\pi_{i}^{0} \geq \overline{c}_{0}$ for $1 \leq i \leq k_{0}$.

We first study a few specific settings when the expert functions $h_{1}$ and $h_{2}$ do not satisfy Condition (O.1) in Sections~\ref{Section:linear_expert_functions} and~\ref{Section:non-linear_expert}. Then, we study a few other representative settings when the expert functions $h_{1}$ and $h_{2}$ do not satisfy either Conditions (O.2) or (O.3) in Section~\ref{sec:other_algebraically_dependent_expert}.
\subsection{Beyond Condition (O.1): Linear expert functions and uniform convergence rates of the MLE}
\label{Section:linear_expert_functions}
In this section, we consider a few representative examples involving expert functions $h_{1}$ and $h_{2}$ that are algebraically dependent. We establish the corresponding convergence rates of the MLE for these examples. Our analysis will be divided into two distinct choices for $h_{2}$: when $h_{2}$ is covariate independent and when $h_{2}$ depends on the covariate.
\subsubsection{Covariate-independent expert function $h_{2}$} We first consider an algebraic dependence setting where the expert function $h_{2}$ is independent of the covariate $X$.
\begin{example} \label{example:not_second_order_MECFG}
Let the expert functions be $h_{1}(X|\theta_{1}) = \theta_{1}^{(1)}+\theta_{1}^{(2)}X$ for all $\theta_{1} = (\theta_{1}^{(1)},\theta_{1}^{(2)}) \in \Omega_{1} \subset \mathbb{R}^{2}$ and $h_{2}^{2}(X|\theta_{2}) = \theta_{2}$ for all $\theta_{2} \in \Omega_{2} \subset \mathbb{R}_{+}$. These expert functions $h_{1}$ and $h_{2}$ are algebraically dependent, as characterized via the following PDE relating $h_{1}$ and $h_{2}$
\begin{align}
\parenth{\frac{\partial{h_{1}}}{\partial{\theta_{1}^{(1)}}}(X,\theta_{1})}^{2} = \dfrac{\partial{h_{2}^{2}}}{\partial{\theta_{2}}}(X,\theta_{2}), \label{eqn:PDE_structure_cova_independent}
\end{align}
for all $\theta_{1}$ and $\theta_{2}$. 
\end{example}
It is clear that the expert functions $h_{1}$ and $h_{2}$ in Example~\ref{example:not_second_order_MECFG} satisfy Conditions (O.2) and (O.3) in Definition~\ref{definition:identifiability}. Therefore, these expert functions only do not satisfy condition (O.1) in that definition. Now, let $\overline{r}: = \overline{r}(k-k_{0}+1)$ be the minimum value of $r$ such that the following system of polynomial equations:
\begin{eqnarray}
\sum \limits_{j=1}^{k-k_{0}+1} \sum \limits_{n_{1}, n_{2}} \dfrac{c_{j}^{2}a_{j}^{n_{1}}b_{j}^{n_{2}}}{n_{1}!n_{2}!} = 0 \ \text{for each} \ \alpha=1,\ldots,r, \label{eqn:system_polynomial_Gaussian_first}
\end{eqnarray}
does not have any nontrivial solution for the unknown variables $(a_{j},b_{j},c_{j})_{j=1}^{k-k_{0}+1}$. The ranges of $n_{1}$ and $n_{2}$ in the second sum consist of all natural pairs satisfying the equation $n_{1}+2n_{2}=\alpha$. A solution to the above system is considered \textit{nontrivial} if all of variables $c_{j}$ are non-zeroes, while at least one of the $a_{j}$ is non-zero. 

Our use of the parameter $\overline{r}$ builds on earlier work by~\cite{Ho-Nguyen-AOS-17} who used it to establish convergence rates in the setting of over-specified location-scale Gaussian mixtures, which is a special case of over-specified GMCF when the expert are identity functions. As demonstrated in Proposition 2.1 in~\cite{Ho-Nguyen-AOS-17}, when $k - k_{0} = 1$, then $\bar{r} = 4$. When $k - k_{0} =2$, we have $\bar{r} = 6$. When $k - k_{0} \geq 3$, then $\bar{r} \geq 7$. As the authors indicated, the actual value of $\bar{r}$ when $k - k_{0} \geq 3$ is non-trivial to determine as we need to use computational algebra tools, such as Groebner bases, and these tools become computationally expensive to use when $k - k_{0} \geq 3$.  The following theorem shows that $\overline{r}$ plays a role in the sharp convergence rate of the MLE under the setting of expert functions in Example~\ref{example:not_second_order_MECFG}. 
\begin{theorem} \label{theorem:lower_bound_Gaussian_family}
Assume that expert functions $h_{1}(X|\theta_{1}) = \theta_{1}^{(1)} + \theta_{1}^{(2)} X$ for $\theta_{1} = (\theta_{1}^{(1)}, \theta_{1}^{(2)}) \in \Omega_{1} \subset \mathbb{R}^{2}$ and $h_{2}^{2}(X|\theta_{2}) = \theta_{2}$ for $\theta_{2} \in \Omega_{2} \subset \mathbb{R}_{+}$. Then, we have the following results:
\begin{itemize}
\item[(a)] (Convergence rate of MLE) There exists a positive constant $C_{0}$ depending only on $G_{0}$ and $\Omega$ such that
\begin{eqnarray}
\mathbb{P}\parenth{\widetilde{W}_{\kappa}(\widehat{G}_{n},G_{0}) > C_{0}(\log n/n)^{1/2\overline{r}}} \precsim \exp(-c\log n), \nonumber
\end{eqnarray}
where $\kappa = (\overline{r},2,\ceil{\overline{r}/2})$ and $\overline{r}$ is defined in equation~\eqref{eqn:system_polynomial_Gaussian_first}. Here, $c$ is a positive constant depending only on $\Omega$.
\item[(b)] (Minimax lower bound) For any $\kappa'$such that $(1,1,1) \preceq \kappa' \prec \kappa = (\overline{r},2,\ceil{\overline{r}/2})$, 
\begin{align}
\inf \limits_{\overline{G}_{n}} \sup \limits_{G \in \Ocal_{k}(\Omega) \backslash \Ocal_{k_{0}-1}(\Omega)} \Exs_{p_{G}} \parenth{\widetilde{W}_{\kappa'}(\overline{G}_{n},G)} \succsim n^{-1/(2\|\kappa'\|_{\infty})}. \nonumber
\end{align}
Here, the infimum is taken over all sequences of estimates $\overline{G}_{n} \in \Ocal_{k}(\Omega)$. Furthermore, $\Exs_{p_{G}}$ denotes the expectation taken with respect to the product measure with mixture density $p_{G}^{n}$. 
\end{itemize}
\end{theorem}
\noindent
The proof of Theorem~\ref{theorem:lower_bound_Gaussian_family} is in Section~\ref{Section:proof_weakly_identifiable_covariate_independent}. 

The sharp convergence rates of MLE in Theorem~\ref{theorem:lower_bound_Gaussian_family} demonstrate that the best 
possible convergence rates of estimating $(\theta_{1i}^{0})^{(1)}$, $
(\theta_{1i}^{0})^{(2)}$, and $\theta_{2i}^{0}$ are not uniform. In 
particular, the rates for estimating $(\theta_{1i}^{0})^{(1)}$ and $
(\theta_{1i}^{0})^{(2)}$ are $n^{-1/2\overline{r}}$ and 
$n^{-1/4}$, respectively, while the rate for estimating $\theta_{2i}^{0}$ is $n^{-1/
2\ceil{\overline{r}/2}}$ (up to a logarithmic factor) for all $1 \leq i \leq 
k_{0}$. Therefore, estimation of the second component of $\theta_{1i}
^{0}$ is generally much faster than estimation of the first component of $
\theta_{1i}^{0}$ and $\theta_{2i}^{0}$. As is seen in the proof,
the slow convergence of $(\theta_{1i}^{0})^{(1)}$ and $\theta_{2i}
^{0}$ arises from the way in which the structure of the PDE~\eqref{eqn:PDE_structure_cova_independent} captures the statistically relevant dependence of the expert functions 
$h_{1}$ and $h_{2}$. In particular, the PDE shows that $(\theta_{1i}
^{0})^{(1)}$ and $\theta_{2i}^{0}$ are linearly dependent, but, since the second component of $\theta_{1i}^{0}$ is 
associated with the covariate $X$, it does not have any interaction with 
$\theta_{2i}^{0}$, which explains why it enjoys a much faster 
convergence rate than the other parameters.

Second, if we choose expert functions $h_{1}(X,\theta_{1}) = 
\theta_{1}^{(1)}+\theta_{1}^{(2)}X+\ldots+\theta_{1}^{(q_{1})}
X^{q_{1}}$ for any $q_{1} \geq 2$ and $h_{2}^{2}(X,\theta_{2}) = 
\theta_{2}$ where $\theta_{1} = (\theta_{1}^{(1)},\ldots,\theta_{1}
^{(q_{1})})$, then with a similar argument we obtain that the best possible convergence 
rates for estimating $(\theta_{1i}^{0})^{(j)}$ for $j \neq 1$ are 
$n^{-1/4}$ for all $1 \leq i \leq k_{0}$ while those for $(\theta_{1i}
^{0})^{(1)}$ and $\theta_{2i}^{0}$ are 
$n^{-1/2\overline{r}}$ and $n^{-1/2\ceil{\overline{r}/2}}$, respectively (up to a logarithmic factor).

\subsubsection{Covariate-dependent expert function $h_{2}$}
We now turn to the setting of algebraic dependence between the parameters associated with covariate $X$ in $h_{1}$ and the parameters of $h_{2}$.
\begin{example} \label{example:not_second_order_MECFG_second}
Define expert functions $h_{1}(X,\theta_{1}) = \theta_{1}^{(1)}+\theta_{1}^{(2)}X$ for all $\theta_{1} = (\theta_{1}^{(1)},\theta_{1}^{(2)}) \in \Omega_{1} \subset \mathbb{R}^{2}$ and $h_{2}^{2}(X,\theta_{2}) = \theta_{2}^{(1)} + \theta_{2}^{(2)}X^{2}$, for all $\theta_{2}=(\theta_{2}^{(1)},\theta_{2}^{(2)}) \in \Omega_{2} \subset \mathbb{R}^{2}$ such that $\theta_{2}^{(1)}, \theta_{2}^{(2)} \geq 0$ and $\theta_{2}^{(1)}+\theta_{2}^{(2)} \geq \overline{\gamma}$ for some positive constant $\overline{\gamma}$. We have the following PDE for these expert functions:
\begin{align}
\parenth{\frac{\partial{h_{1}}}{\partial{\theta_{1}^{(1)}}}(X,\theta_{1})}^{2} = \dfrac{\partial{h_{2}^{2}}}{\partial{\theta_{2}^{(1)}}}(X,\theta_{2}), \label{eqn:PDE_structure_cova_dependent_first} \\ \parenth{\dfrac{\partial h_{1}}{\partial{\theta_{1}^{(2)}}}(X, \theta_{1})}^{2} = \dfrac{\partial{h_{2}^{2}}}{\partial{\theta_{2}^{(2)}}}(X, \theta_{2}), \label{eqn:PDE_structure_cova_dependent_second}
\end{align} 
which shows that $h_{1}$ and $h_{2}$ are algebraically dependent.
\end{example}
We can check that the expert functions $h_{1}$ and $h_{2}$ in Example~\ref{example:not_second_order_MECFG_second} satisfy Conditions (O.2) and (O.3) in Definition~\ref{definition:identifiability}; therefore, they only do not satisfy Condition (O.1) in that definition. The main distinction between Example~\ref{example:not_second_order_MECFG_second} and Example~\ref{example:not_second_order_MECFG} is that we have the covariate $X^2$ in the formulation of the expert function $h_{2}$ in Example~\ref{example:not_second_order_MECFG_second}. This inclusion leads to a rather rich spectrum of convergence rates for the MLE. To illustrate these convergence rates, we consider two distinct cases for the expert function $h_{2}$: 
\begin{itemize}
\item \textbf{without offset}: $\theta_{2}^{(1)}=0$, i.e., $h_{2}^{2}(X,\theta_{2}) = \theta_{2}^{(2)}X^{2}$. 
\item \textbf{with offset}: $\theta_{2}^{(1)}$ is taken into account; i.e., $h_{2}^{2}(X,\theta_{2}) = \theta_{2}^{(1)} + \theta_{2}^{(2)}X^{2}$.
\end{itemize}
\begin{theorem} (Without offset) \label{theorem:lower_bound_Gaussian_family_first_first}
Let $\overline{r}$ be defined as in equation~\eqref{eqn:system_polynomial_Gaussian_first}. Assume that expert functions $h_{1}(X,\theta_{1}) = \theta_{1}^{(1)} + \theta_{1}^{(2)} X$ for $\theta_{1} = (\theta_{1}^{(1)}, \theta_{1}^{(2)}) \in \Omega_{1} \subset \mathbb{R}^{2}$ and $h_{2}^{2}(X,\theta_{2}) = \theta_{2}X^{2}$ for $\theta_{2} \in \Omega_{2} \subset \mathbb{R}_{+}$. Then, the following holds:
\begin{itemize}
\item[(a)] (Convergence rate of MLE) There exists a positive constant $C_{0}$ depending only on $G_{0}$ and $\Omega$ such that
\begin{eqnarray}
\mathbb{P}\parenth{\widetilde{W}_{\kappa}(\widehat{G}_{n},G_{0}) > C_{0}(\log n/n)^{1/2\overline{r}}} \precsim \exp(-c\log n), \nonumber
\end{eqnarray}
where $\kappa = (2,\overline{r},\ceil{\overline{r}/2})$. Here, $c$ is a positive constant depending only on $\Omega$.
\item[(b)] (Minimax lower bound) For any $\kappa'$ such that $(1,1,1) \preceq \kappa' \prec \kappa = (2,\overline{r},\ceil{\overline{r}/2})$, 
\begin{align}
\inf \limits_{\overline{G}_{n}} \sup \limits_{G \in \Ocal_{k}(\Omega) \backslash \Ocal_{k_{0}-1}(\Omega)} \Exs_{p_{G}} \parenth{\widetilde{W}_{\kappa'}(\overline{G}_{n},G)} \succsim n^{-1/(2\|\kappa'\|_{\infty})}. \nonumber
\end{align}
Here, the infimum is taken over all sequences of estimates $\overline{G}_{n} \in \Ocal_{k}(\Omega)$. Furthermore, $\Exs_{p_{G}}$ denotes the expectation taken with respect to the product measure with mixture density $p_{G}^{n}$. 
\end{itemize}
\end{theorem}
The proof of Theorem~\ref{theorem:lower_bound_Gaussian_family_first_first} is in Appendix~\ref{Section:proof_weakly_identifiable_covariate_dependent_first}. 

In contrast to the setting of Theorem~\ref{theorem:lower_bound_Gaussian_family}, the expert function $h_{2}$ is now a function of $X^{2}$. The sharp convergence rate of $\widehat{G}_{n}$ in Theorem \ref{theorem:lower_bound_Gaussian_family_first_first} demonstrates that the convergence rates for estimating $(\theta_{1i}^{0})^{(1)}$, $(\theta_{1i}^{0})^{(2)}$, and $\theta_{2i}^{0}$ are $n^{-1/4}$, $n^{-1/2\overline{r}}$, and $n^{-1/2\ceil{\overline{r}/2}}$, respectively, for all $1 \leq i \leq k_{0}$. Therefore, with the formulation of expert functions given in Theorem~\ref{theorem:lower_bound_Gaussian_family_first_first}, estimation of the first component of $\theta_{1i}^{0}$ is much faster than estimation of the second component of $\theta_{1i}^{0}$. This is in contrast to the results in Theorem~\ref{theorem:lower_bound_Gaussian_family}. A high-level explanation for this phenomenon is again obtained by considering the PDE structure, which in this case is given by equation~\eqref{eqn:PDE_structure_cova_dependent_second}:
\begin{align}
\parenth{\dfrac{\partial h_{1}}{\partial{\theta_{1}^{(2)}}}(X, \theta_{1})}^{2} = \dfrac{\partial{h_{2}^{2}}}{\partial{\theta_{2}}}(X, \theta_{2}). \nonumber
\end{align} 
Such a structure implies the dependence of the second component of $\theta_{1i}^{0}$ and $\theta_{2i}^{0}$; therefore, there exists a strong interaction between $(\theta_{1i}^{0})^{(2)}$ and $\theta_{2i}^{0}$ in terms of their convergence rates. On the other hand, the first component of $\theta_{1i}^{0}$ and $\theta_{2}^{0}$ are linearly independent, which implies that there is virtually no interaction between these two terms. As a consequence, $(\theta_{1i}^{0})^{(1)}$ will enjoy much faster convergence rates than $(\theta_{1i}^{0})^{(2)}$ and $\theta_{2i}^{0}$.

In contrast to the setting without an offset term in the expert function $h_{2}$, 
the convergence rate of the MLE under the setting with the offset term in $h_{2}$ suffers from two ways: one which is captured by the PDE structure with 
respect to $\theta_{1}^{(1)}$ and $\theta_{2}^{(1)}$ in equation~\eqref{eqn:PDE_structure_cova_dependent_first} and another from the PDE 
structure with respect to $\theta_{1}^{(2)}$ and $\theta_{2}^{(2)}$ 
in equation~\eqref{eqn:PDE_structure_cova_dependent_second}.
\begin{theorem} (With offset)  \label{theorem:lower_bound_Gaussian_family_first_second}
Let $\overline{r}$ be defined as in equation~\eqref{eqn:system_polynomial_Gaussian_first}. Assume that expert functions $h_{1}(X|\theta_{1}) = \theta_{1}^{(1)} + \theta_{1}^{(2)} X$ for $\theta_{1} = (\theta_{1}^{(1)}, \theta_{2}^{(2)}) \in \Omega_{1} \subset \mathbb{R}^{2}$ and $h_{2}^{2}(X|\theta_{2}) = \theta_{2}^{(1)} + \theta_{2}^{(2)}X^{2}$ for $\theta_{2}=(\theta_{2}^{(1)},\theta_{2}^{(2)}) \in \Omega_{2} \subset \mathbb{R}^{2}$ such that $\theta_{2}^{(1)}, \theta_{2}^{(2)} \geq 0$ and $\theta_{2}^{(1)}+\theta_{2}^{(2)} \geq \overline{\gamma}$ for some given positive $\overline{\gamma}$. Then, the following holds:
\begin{itemize}
\item[(a)] (Convergence rate of MLE) There exists a positive constant $C_{0}$ depending only on $G_{0}$ and $\Omega$ such that
\begin{eqnarray}
\mathbb{P}\parenth{\widetilde{W}_{\kappa}(\widehat{G}_{n},G_{0}) > C_{0}(\log n/n)^{1/2\overline{r}}} \precsim \exp(-c\log n),
\end{eqnarray}
where $\kappa = (\overline{r},\overline{r},\ceil{\overline{r}/2},\ceil{\overline{r}/2})$. Here, $c$ is a positive constant depending only on $\Omega$.
\item[(b)] (Minimax lower bound) For any $\kappa'$ such that $(1,1,1,1) \preceq \kappa' \prec \kappa = (\overline{r},\overline{r},\ceil{\overline{r}/2},\ceil{\overline{r}/2})$, 
\begin{align}
\inf \limits_{\overline{G}_{n}} \sup \limits_{G \in \Ocal_{k}(\Omega) \backslash \Ocal_{k_{0}-1}(\Omega)} \Exs_{p_{G}} \parenth{\widetilde{W}_{\kappa'}(\overline{G}_{n},G)} \succsim n^{-1/(2\|\kappa'\|_{\infty})}. \nonumber
\end{align}
Here, the infimum is taken over all sequences of estimates $\overline{G}_{n} \in \Ocal_{k}(\Omega)$. Furthermore, $\Exs_{p_{G}}$ denotes the expectation taken with respect to the product measure with mixture density $p_{G}^{n}$. 
\end{itemize}
\end{theorem}
\noindent
The proof of Theorem~\ref{theorem:lower_bound_Gaussian_family_first_second} is in Appendix~\ref{Section:proof_weakly_identifiable_covariate_dependent_second}. 

Note that when there is an offset term in the expert function $h_{2}$, the convergence rate 
of $\widehat{G}_{n}$ suggests that the convergence rates for 
estimating $(\theta_{1i}^{0})^{(1)}$, $(\theta_{1i}^{0})^{(2)}$, $
(\theta_{2i}^{0})^{(1)}$, and $(\theta_{2i}^{0})^{(2)}$ are 
 $n^{-1/2\overline{r}}$, $n^{-1/2\overline{r}}$, $n^{-1/
2\ceil{\overline{r}/2}}$,  and $n^{-1/2\ceil{\overline{r}/2}}$, respectively, for all $1 \leq i \leq k_{0}$. In comparison to the convergence 
rate $n^{-1/4}$ for estimating $(\theta_{1i}^{0})^{(2)}$ under the  
setting without covariate $X^{2}$ in $h_{2}$ in Theorem~\ref{theorem:lower_bound_Gaussian_family}, the convergence rate 
$n^{-1/2\overline{r}}$ for estimating $(\theta_{1i}^{0})^{(2)}$ under 
the setting of Theorem~\ref{theorem:lower_bound_Gaussian_family_first_second} is much 
slower. Furthermore, the convergence rate $n^{-1/2\overline{r}}$ for 
estimating $(\theta_{1i}^{0})^{(2)}$ in the setting of 
Theorem~\ref{theorem:lower_bound_Gaussian_family_first_second} is 
much slower than the corresponding rate $n^{-1/4}$ for estimating $
(\theta_{1i}^{0})^{(2)}$ in the setting of Theorem~\ref{theorem:lower_bound_Gaussian_family_first_first}. 

Note also that if we choose more general expert functions, $h_{1}(X,\theta_{1}) = 
\theta_{1}^{(1)}+\theta_{1}^{(2)}X+\ldots+\theta_{1}^{(q_{1})}
X^{q_{1}}$, for any $q_{1} \geq 1$ and $h_{2}^{2}(X,\theta_{2}) = 
\theta_{2}^{(1)}+\theta_{2}^{(2)}X^2+\ldots+\theta_{2}^{(q_{1})}
X^{2q_{1}}$, where $\theta_{1} = (\theta_{1}^{(1)},\ldots,\theta_{1}
^{(q_{1})})$ and $\theta_{2} = (\theta_{2}^{(1)},\ldots,\theta_{2}
^{(q_{1})})$, i.e., letting $q_{2} = q_{1}$, then we also obtain that the best 
possible convergence rates for estimating $(\theta_{1i}^{0})^{(j)}$ are 
$n^{-1/2\overline{r}}$ while those for estimating $(\theta_{2i}
^{0})^{(j)}$ are $n^{-1/2\ceil{\overline{r}/2}}$ for all $1 \leq i \leq k_{0}$ and 
$1 \leq j \leq q_{1}$. Such results can be explained by the
following system of PDEs characterizing the dependence between $
\theta_{1}^{(i)}$ and $\theta_{2}^{(i)}$ for $1 \leq i \leq q_{1}$:
\begin{align}
\parenth{\frac{\partial{h_{1}}}{\partial{\theta_{1}^{(i)}}}(X,\theta_{1})}^{2} = \dfrac{\partial{h_{2}^{2}}}{\partial{\theta_{2}^{(i)}}}(X,\theta_{2}), \ \text{for all} \ 1 \leq i \leq q_{1}, \nonumber 
\end{align}
for any $(\theta_{1}, \theta_{2})$. 

\subsection{Beyond Condition (O.1): Nonlinear expert functions and non-uniform convergence rates of MLE}
\label{Section:non-linear_expert}
Thus far we have considered various algebraic dependence settings for linear expert functions $h_{1}$ and $h_{2}$ with respect to their parameters. Under these settings, the convergence rates of the MLE are uniform; i.e., they are independent of the values of the true mixing measure $G_{0}$. In this section, we demonstrate that in the case of nonlinear expert functions $h_{1}$ and $h_{2}$ that are algebraically dependent and do not satisfy Condition (O.1), the convergence rates of $\widehat{G}_{n}$ strongly depend on the values of $G_{0}$. 

The specific setting that we consider is when $h_{1}$ is nonlinear in terms of its parameter $\theta_{1}$ while $h_{2}$ is independent of the covariate $X$. In that setting, we have the following simple example of algebraically dependent expert functions:
\begin{align}
h_{1}(X,\theta_{1}) = \parenth{\theta_{1}^{(1)}+\theta_{1}^{(2)}X}^2, \ h_{2}^{2}(X|\theta_{2}) = \theta_{2}, \label{eqn:non_linear_expert}
\end{align}
for all $\theta_{1} = (\theta_{1}^{(1)},\theta_{1}^{(2)}) \in \Omega_{1} = [0,\overline{\tau}_{1}] \times [0,\overline{\tau}_{2}]$ and $\theta_{2} \in \Omega_{2} \subset \mathbb{R}_{+}$ where $\overline{\tau}_{1}$, $\overline{\tau}_{2}$ are given positive numbers. Here, the choice regarding the ranges of $\theta_{1}$ is to ensure that the expert function $h_{1}$ is identifiable with respect to its parameter $\theta_{1}$. The following result shows that the expert functions $h_{1}$ and $h_{2}$ are algebraically dependent as they do not satisfy Condition (O.1).
\begin{proposition} \label{proposition:not_second_order_MECFG_third}
Assume that the expert functions $h_{1}$ and $h_{2}$ take the forms in equation~\eqref{eqn:non_linear_expert}. Then the expert functions $h_{1}$ and $h_{2}$ are algebraically dependent, as captured in the following PDE that relates $h_{1}$ and $h_{2}$:
\begin{align}
\parenth{\dfrac{\partial h_{1}}{\partial{\theta_{1}^{(1)}}}(X, \theta_{1})}^{2} = 4(\theta_{1}^{(1)})^{2}\dfrac{\partial{h_{2}^{2}}}{\partial{\theta_{2}}}(X, \theta_{2}), \label{eqn:PDE_structure_nonlinear_expert}
\end{align}
for all $\theta_{1} = (\theta_{1}^{(1)}, 0)$ and $\theta_{2}$. 
\end{proposition}
Unlike the previous PDEs in equations~\eqref{eqn:PDE_structure_cova_independent},~\eqref{eqn:PDE_structure_cova_dependent_first}, and~\eqref{eqn:PDE_structure_cova_dependent_second}, which hold 
for all $(\theta_{1},\theta_{2})$, the PDE in equation~\eqref{eqn:PDE_structure_nonlinear_expert} holds only under a special 
structure for $\theta_{1}$; namely, $\theta_{1} = (\theta^{(1)}, 0)$, where the second component of $\theta_{1}$ needs to be zero. Such a special 
structure of the PDE leads to an interesting phase transition regarding the 
convergence rates of the MLE under specific values of true mixing measure 
$G_{0}$. In order to capture this phase transition precisely, we distinguish two separate settings of $G_{0}$:
\begin{itemize}
\item[•] \textbf{Nonlinearity setting I}: As long as there exists $(\theta_{1i}^{0})^{(2)} = 0$ for some $1 \leq i \leq k_{0}$, we have $(\theta_{1i}^{0})^{(1)} = 0$.
\item[•] \textbf{Nonlinearity setting II}: There exists $\theta_{1i}^{0}$ such that $(\theta_{1i}^{0})^{(1)} \neq 0$ and $(\theta_{1i}^{0})^{(2)} = 0$ for some index $1 \leq i \leq k_{0}$. 
\end{itemize}
\subsubsection{Nonlinearity setting I}
Under the nonlinearity setting I for the true mixing measure $G_{0}$, we have the following result regarding the sharp convergence rate of the MLE.
\begin{theorem} \label{theorem:lower_bound_Gaussian_family_third}
Let the expert functions $h_{1}$ and $h_{2}$ be defined as in equation~\eqref{eqn:non_linear_expert}. Under the nonlinearity setting I for $G_{0}$, the following holds:
\begin{itemize}
    \item[(a)] (Convergence rate of MLE) There exists a positive constant $C_{0}$ depending only on $G_{0}$ and $\Omega$ such that
\begin{eqnarray}
\mathbb{P}(\widetilde{W}_{\kappa}(\widehat{G}_{n},G_{0}) > C_{0}(\log n/n)^{1/4}) \precsim \exp(-c\log n), \nonumber
\end{eqnarray}
where $\kappa = (2,2,2)$. Here, $c$ is a positive constant depending only on $\Omega$.
\item[(b)] (Minimax lower bound) For any $\kappa'$ such that $(1,1,1) \preceq \kappa' \prec \kappa = (2,2,2)$, 
\begin{align}
\inf \limits_{\overline{G}_{n} \in \mathcal{G}_{1}} \sup \limits_{G \in \mathcal{G}_{1}} \Exs_{p_{G}} \parenth{\widetilde{W}_{\kappa'}(\overline{G}_{n},G)} \succsim n^{-1/(2\|\kappa'\|_{\infty})}, \nonumber
\end{align}
where the structure of the parameter space $\mathcal{G}_{1} \subset \Ocal_{k}(\Omega) \backslash \Ocal_{k_{0}-1}(\Omega)$ is given by
\begin{align}
\mathcal{G}_{1} = \biggr\{G = \sum_{i=1}^{k'} \pi_{i}\delta_{(\theta_{1i},\theta_{2i})}: k_{0} \leq k' \leq k \ \text{and as long as} \ \theta_{1i}^{(2)} = 0 \ \text{for some} \ 1 \leq i \leq k', \nonumber \\
\text{then} \ \theta_{1i}^{(1)} = 0 \biggr\}. \nonumber
\end{align}
Here, $\Exs_{p_{G}}$ denotes the expectation taken with respect to the product measure with mixture density $p_{G}^{n}$. 
\end{itemize}
\end{theorem}
\noindent
The proof of Theorem~\ref{theorem:lower_bound_Gaussian_family_third} is in Appendix~\ref{Section:proof_non_linear_covariate_dependent_first}. 

Under the nonlinearity setting I for $G_{0}$, the results of Theorem~\ref{theorem:lower_bound_Gaussian_family_third} suggest that the convergence rates for estimating $(\theta_{1i}^{0})^{(1)}$, $(\theta_{1i}^{0})^{(2)}$, and $\theta_{2i}^{0}$ are $n^{-1/4}$ and these convergence rates are sharp. Furthermore, these convergence rates match those under the settings in which the expert functions $h_{1}$ and $h_{2}$ are algebraically independent. This phenomenon arises because there is no linkage between $\theta_{1i}^{0}$ and $\theta_{2i}^{0}$ in the PDE for the nonlinearity setting I. 

\subsubsection{Nonlinearity setting II}
Unlike the nonlinearity setting I of $G_{0}$, the convergence rate of MLE under nonlinearity setting II is more complicated to analyze due to the existence of the zero-valued coefficient $(\theta_{1i}^{0})^{(2)}$ for some $1 \leq i \leq k_{0}$. To simplify the presentation, we first start with a result regarding the structure of the partial derivatives of $f$ when the second component of $\theta_{1}$ is zero. We then define an \emph{inhomogeneous system of polynomial limits} based on this structural assumption to analyze the behavior of the MLE. Finally, we state a formal convergence rate result of the MLE under the general nonlinearity setting II for $G_{0}$. 
\paragraph{Partial derivative structures:} 
Since there exists a zero-valued coefficient $(\theta_{1i}^{0})^{(2)}$ for some $1 \leq i \leq k_{0}$ under the nonlinearity setting II of $G_{0}$, we will focus on understanding the partial derivatives of $f$ when the second component of $\theta_{1}$ is 0, i.e., $\theta_{1}^{(2)} = 0$. To facilitate the discussion, we firstly consider a few specific simple examples of these derivatives: 
\begin{align}
\dfrac{\partial{f}}{\partial{\theta_{1}^{(1)}}} = 2\theta_{1}^{(1)}\dfrac{\partial{f}}{\partial{h_{1}}}, \ \quad \ \dfrac{\partial{f}}{\partial{\theta_{1}^{(2)}}} = 2\theta_{1}^{(1)}X \dfrac{\partial{f}}{\partial{h_{1}}}, \ \quad \ \dfrac{\partial{f}}{\partial{\theta_{2}}} = \dfrac{\partial{f}}{\partial{h_{2}^{2}}} = \dfrac{1}{2}\dfrac{\partial^{2}{f}}{\partial{h_{1}^{2}}}, \nonumber \\
\dfrac{\partial^{2}{f}}{\partial{(\theta_{1}^{(1)})^{2}}} =  2\dfrac{\partial{f}}{\partial{h_{1}}} + 4(\theta_{1}^{(1)})^{2}\dfrac{\partial^{2}{f}}{\partial{h_{1}^{2}}}, \ \quad \ \dfrac{\partial{f}}{\partial{(\theta_{1}^{(2)})^{2}}} = 2 X^{2}\dfrac{\partial{f}}{\partial{h_{1}}} + 4(\theta_{1}^{(1)})^{2}X^{2}\dfrac{\partial^{2}{f}}{\partial{h_{1}^{2}}}, \nonumber \\
\dfrac{\partial^{2}{f}}{\partial{\theta_{2}^{2}}} = \dfrac{\partial^{2}{f}}{\partial{h_{2}^{4}}} = \dfrac{1}{4}\dfrac{\partial^{4}{f}}{\partial{h_{1}^{4}}}, \ \quad \ \dfrac{\partial^{3}{f}}{\partial{(\theta_{1}^{(1)})^{3}}} =  12\theta_{1}^{(1)}\dfrac{\partial^{2}{f}}{\partial{h_{1}^{2}}} + 8(\theta_{1}^{(1)})^{3}\dfrac{\partial^{3}{f}}{\partial{h_{1}^{3}}}. \nonumber
\end{align}
Here, we suppress the condition on  $h_{1}(X,\theta_{1})$ and $h_{2}(X,\theta_{2})$ in the notation to simplify the presentation. From this computation, it is clear that $\frac{\partial{f}}{\partial{\theta_{1}^{(1)}}}$, $\frac{\partial{f}}{\partial{\theta_{2}}}$, and $\frac{\partial^{2}{f}}{\partial{(\theta_{1}^{(1)})^{2}}}$ are not linearly independent with respect to $X$ and $Y$. This dependence among these partial derivatives underlies the complex behavior of the MLE in this setting.

By iterating this computation of partial derivatives of $f$ up to a high order, we obtain the following key lemma generalizing the structure of partial derivatives of $f$ with respect to $\theta_{1}^{(1)}$ and $\theta_{2}$.
\begin{lemma} \label{lemma:representation_partial_derivative}
Assume that $\theta_{1}^{(2)}=0$. For any value of $\theta_{1}^{(1)} \neq 0$, $\theta_{2}$, and $\gamma = (\gamma_{1},\gamma_{2}) \in \mathbb{N}^{2}$, the following holds:
\begin{itemize}
\item[(a)] When $\gamma_{1}$ is an odd number, we have: 
\begin{eqnarray}
& & \dfrac{\partial^{|\gamma|}f}{\partial{(\theta_{1}^{(1)})^{\gamma_{1}}\partial{\theta_{2}^{\gamma_{2}}}}}(Y|h_{1}(X,\theta_{1}),h_{2}(X,\theta_{2})) \nonumber \\
& & \hspace{7 em} = \dfrac{1}{2^{\gamma_{2}}} \biggr(\sum \limits_{u = 0} ^{(\gamma_{1}-1)/2} P_{u}^{(\gamma_{1})}(\theta_{1}^{(1)}) \dfrac{\partial^{\frac{\gamma_{1}+1}{2}+ u+2\gamma_{2}}f}{\partial{h_{1}^{\frac{\gamma_{1}+1}{2}+ u+2\gamma_{2}}}}(Y|h_{1}(X, \theta_{1}),h_{2}(X, \theta_{2}))\biggr). \nonumber
\end{eqnarray}
\item[(b)] When $\gamma_{1}$ is an even number, then:
\begin{eqnarray}
& & \dfrac{\partial^{|\gamma|}f}{\partial{(\theta_{1}^{(1)})^{\gamma_{1}}\partial{\theta_{2}^{\gamma_{2}}}}}(Y|h_{1}(X,\theta_{1}),h_{2}(X,\theta_{2})) \nonumber \\
& &  \hspace{ 7 em}  = \dfrac{1}{2^{\gamma_{2}}} \biggr(\sum \limits_{u = 0} ^{\gamma_{1}/2} P_{u}^{(\gamma_{1})}(\theta_{1}^{(1)}) \dfrac{\partial^{\frac{\gamma_{1}}{2}+ u +2\gamma_{2}}f}{\partial{h_{1}^{\frac{\gamma_{1}}{2}+ u +2\gamma_{2}}}}(Y|h_{1}(X, \theta_{1}),h_{2}(X, \theta_{2}))\biggr). \nonumber
\end{eqnarray}
\end{itemize}
Here, $P_{u}^{(\gamma_{1})}(\theta_{1}^{(1)})$ are polynomials in terms of $\theta_{1}^{(1)}$ that satisfy the following iterative equations:
\begin{eqnarray}
P_{0}^{(1)} (\theta_{1}^{(1)}) : = 2\theta_{1}^{(1)}, \ P_{0}^{(\gamma_{1}+1)}(\theta_{1}^{(1)}) : = \dfrac{\partial P_{0}^{(\gamma_{1})}}{\partial{\theta_{1}^{(1)}}}(\theta_{1}^{(1)}), \nonumber \\
P_{\tau}^{(\gamma_{1}+1)}(\theta_{1}^{(1)}) : = 2\theta_{1}^{(1)} P_{\tau-1}^{(\gamma_{1})}(\theta_{1}^{(1)}) + \dfrac{\partial P_{\tau-1}^{(\gamma_{1})}}{\partial{\theta_{1}^{(1)}}}(\theta_{1}^{(1)}), \nonumber
\end{eqnarray}
for any $1 \leq u \leq (\gamma_{1}-1)/2$ when $\gamma_{1}$ is an 
odd number or for any $1 \leq u \leq (\gamma_{1}-2)/2$ such that $\gamma_{1}$ is an even number. Additionally, $P_{(\gamma_{1}+1)/
2}^{(\gamma_{1}+1)}(\theta_{1}^{(1)}) = 2\theta_{1}^{(1)} 
P_{(\gamma_{1}-1)/2}^{(\gamma_{1})}(\theta_{1}^{(1)})$ if $
\gamma_{1}$ is an odd number while $P_{\gamma_{1}/2}
^{(\gamma_{1}+1)}(\theta_{1}^{(1)}) = 2\theta_{1}^{(1)} 
P_{(\gamma_{1}-2)/2}^{(\gamma_{1})}(\theta_{1}^{(1)})$ when $
\gamma_{1} \geq 2$ is an even number. 
\end{lemma}
\paragraph{Inhomogeneous system of polynomial limits:} Given the 
specifications of the polynomials $P_{\tau}^{(\gamma_{1})}(\theta_{1}
^{(1)})$ in Lemma~\ref{lemma:representation_partial_derivative}, we define a system of polynomial limits that is useful for studying convergence rates under the nonlinearity setting II as follows.
Assume that we are given $s \in \mathbb{N}$ and $3s$ sequences $\left\{a_{i,n}\right\}_{n 
\geq 1}$, $\left\{b_{i,n}\right\}_{n \geq 1}$, and $\left\{c_{i,n}\right
\}_{n \geq 1}$ such that $a_{i,n} \to 0, \ b_{i,n} \to 0$ as $n \to 
\infty$ for $1 \leq i \leq s$, while $c_{i,n} \geq 0$ as $1 \leq i \leq s$ 
and $\sum \limits_{i=1}^{s} c_{i,n} \leq \overline{c}$ for some given $
\overline{c} > 0$. For each $\theta_{1}^{(1)}$ and $r \in \mathbb{N}$, we 
denote the following inhomogeneous system of polynomial limits:
\begin{eqnarray}
\dfrac{\sum \limits_{\gamma_{1},\gamma_{2}, u} \dfrac{P_{u}^{(\gamma_{1})}(\theta_{1}^{(1)})}{2^{\gamma_{2}}}\biggr( \sum \limits_{i=1}^{s} c_{i,n}\dfrac{a_{i,n}^{\gamma_{1}}b_{i,n}^{\gamma_{2}}}{\gamma_{1}!\gamma_{2}!}\biggr)}{\sum \limits_{i=1}^{s} c_{i,n}\biggr(|a_{i,n}|^{r}+|b_{i,n}|^{\ceil{r/2}}\biggr)} \to 0, \label{eqn:system_complex_polynomial_limits}
\end{eqnarray}
as $n \to \infty$ for all $1 \leq l \leq 2r$ where the summation with 
respect to $\gamma_{1},\gamma_{2}, u$ in the numerator satisfies $
\gamma_{1}/2+ u +2\gamma_{2}=l$, $u \leq \gamma_{1}/2$ when $
\gamma_{1}$ is an even number while $(\gamma_{1}+1)/2+ u
+2\gamma_{2} = l$, $u \leq (\gamma_{1}-1)/2$ when $\gamma_{1}$ 
is an odd number. Additionally, $\gamma_{1}+\gamma_{2} \leq r$. 

From these conditions, it is clear that the system of polynomial 
limits \eqref{eqn:system_complex_polynomial_limits} contains exactly 
$2r$ polynomial limits. For example, when $r=2$, the  system of 
polynomial limits contains four polynomial limits, which take the following form:
\begin{eqnarray}
\biggr(\sum \limits_{i=1}^{s} c_{i,n} a_{i,n}^{2} + 2\theta_{1}^{(1)} \sum \limits_{i=1}^{s} c_{i,n} a_{i,n} \biggr)\bigg/\biggr(\sum \limits_{i=1}^{s} c_{i,n} \biggr(|a_{i,n}|^{2}+|b_{i,n}| \biggr)\biggr) \to 0, \nonumber \\
\biggr(4(\theta_{1}^{(1)})^{2}\biggr(\sum \limits_{i=1}^{s} c_{i,n} a_{i,n}^{2}\biggr) + \sum \limits_{i=1}^{s} c_{i,n} b_{i,n} \biggr)\bigg/\biggr(\sum \limits_{i=1}^{s} c_{i,n} \biggr(|a_{i,n}|^{2}+|b_{i,n}| \biggr)\biggr) \to 0, \nonumber \\
\theta_{1}^{(1)}\biggr(\sum \limits_{i=1}^{s} c_{i,n} a_{i,n}b_{i,n} \biggr)\bigg/\biggr(\sum \limits_{i=1}^{s} c_{i,n} \biggr(|a_{i,n}|^{2}+|b_{i,n}| \biggr)\biggr) \to 0, \nonumber \\
\biggr(\sum \limits_{i=1}^{s} c_{i,n} b_{i,n}^{2} \biggr)\bigg/\biggr(\sum \limits_{i=1}^{s} c_{i,n} \biggr(|a_{i,n}|^{2}+|b_{i,n}| \biggr)\biggr) \to 0. \nonumber
\end{eqnarray} 
\paragraph{Studying system of polynomial limits:} In general, when $r$ is large, the system of polynomial limits~\eqref{eqn:system_complex_polynomial_limits} does not have a solution; i.e., not all 
the polynomial limits go to zero. We can therefore find a smallest 
value of $r$ such that this system of polynomial limits has no solution. This motivates the following definition that plays a key role in 
obtaining a convergence rate for the MLE.
\begin{definition} \label{definition:smallest_number}
For any $s \geq 1$ and $\theta_{1}^{(1)}$,  define $\rtil(\theta_{1}^{(1)},s)$ as the smallest positive integer $r$ such that system of polynomial limits \eqref{eqn:system_complex_polynomial_limits} does not hold for any choices of sequences $\left\{a_{i,n}\right\}_{n \geq 1}$, $\left\{b_{i,n}\right\}_{n \geq 1}$, and $\left\{c_{i,n}\right\}_{n \geq 1}$.
\end{definition} 
In general, determining the exact value of $\rtil(\theta_{1}^{(1)},s)$ is difficult as the system of polynomial limits~\eqref{eqn:system_complex_polynomial_limits} is intricate. In the following lemma, we demonstrate that we can obtain an upper bound of $\rtil(\theta_{1}^{(1)},s)$ based on the system of polynomial equations~\eqref{eqn:system_polynomial_Gaussian_first}, for any $s \geq 1$ and $\theta_{1}^{(1)} \neq 0$. 
\begin{lemma} \label{lemma:upper_bound_polynomial_limits}
For a general value of $s \geq 2$ and for all $\theta_{1}^{(1)} \neq 0$, we have
\begin{align}
3 \leq \rtil(\theta_{1}^{(1)},s) \leq \overline{r}(s), \nonumber
\end{align} 
where $\overline{r}(s)$ is defined as in~\eqref{eqn:system_polynomial_Gaussian_first}.
\end{lemma}
\paragraph{Convergence rates of MLE:} Equipped with the definition of $\rtil(\theta_{1}^{(1)},s)$, we have the following result for the convergence rate of the MLE under the nonlinearity setting II.
\begin{theorem} \label{theorem:lower_bound_Gaussian_family_second_singularity_case}
Given the nonlinearity setting II for $G_{0}$ and the expert functions $h_{1}$ and $h_{2}$  in equation~\eqref{eqn:non_linear_expert}, we define $\mathcal{A} : = \{i \in [k_{0}]: \ (\theta_{1i}^{0})^{(1)} \neq 0 \ \text{and} \ (\theta_{1i}^{0})^{(2)} = 0 \}$ and
\begin{align}
i_{\text{max}} : = \mathop{ \arg \max} \limits_{i \in \mathcal{A}} \rtil((\theta_{1i}^{0})^{(1)},k-k_{0}+1). \nonumber
\end{align}
Additionally, we denote $\rtil_{\text{sin}} := \rtil((\theta_{1i_{\text{max}}}^{0})^{(1)},k-k_{0}+1)$. 
Then, there exists a positive constant $C_{0}$ depending only on $G_{0}$ and $\Omega$ such that
\begin{eqnarray}
\mathbb{P}(\widetilde{W}_{\kappa}(\widehat{G}_{n},G_{0}) > C_{0}(\log n/n)^{1/2\rtil_{\text{sin}}}) \precsim \exp(-c\log n), \nonumber
\end{eqnarray}
where $\kappa = \parenth{\rtil_{\text{sin}},2,\ceil{\rtil_{\text{sin}}/2}}$.
\end{theorem}
\noindent
The proof of Theorem~\ref{theorem:lower_bound_Gaussian_family_second_singularity_case} is in Appendix~\ref{Section:proof_non_linear_covariate_dependent_second}. 

A few comments are in order. First, the result of Theorem~\ref{theorem:lower_bound_Gaussian_family_second_singularity_case} indicates that the convergence rates for estimating $(\theta_{1i}^{0})^{(1)}, (\theta_{1i}^{0})^{(2)}, \theta_{2i}^{0}$ are  $n^{-1/2\rtil_{\text{sin}}}, n^{-1/4}$, and $n^{-1/2\ceil{\rtil_{\text{sin}}/2}}$, respectively, for $1 \leq i \leq k_{0}$. The slow convergence rates of estimating $(\theta_{1}^{0})^{(1)}$ and $\theta_{2}^{0}$ under nonlinearity setting II is captured by the PDE~\eqref{eqn:PDE_structure_nonlinear_expert}, which indicates that $(\theta_{1i}^{0})^{(1)}$ and $\theta_{2i}^{0}$ are linearly dependent when the second component of $\theta_{1i}^{0}$ is zero. 

Second, since $\rtil_{\text{sin}} \leq \overline{r} = \overline{r}( k - k_{0} + 1)$, the convergence rates for estimating $(\theta_{1i}^{0})^{(1)}$ and $\theta_{2i}^{0}$ under the settings of expert functions $h_{1}$ and $h_{2}$ in equation~\eqref{eqn:non_linear_expert} may be faster than those of $(\theta_{1i}^{0})^{(1)}$ and $\theta_{2i}^{0}$ under the choice of expert functions $h_{1}$ and $h_{2}$ in Example~\ref{example:not_second_order_MECFG}, i.e., $h_{1}(X|\theta_{1}) = \theta_{1}^{(1)}+\theta_{1}^{(2)}X$ and $h_{2}^{2}(X|\theta_{2}) = \theta_{2}$. Therefore, parameter estimation when $h_{1}$ is quadratic in terms of $\theta_{1}^{(1)}+\theta_{1}^{(2)}X$ is generally easier than when $h_{1}$ is linear in terms of $\theta_{1}^{(1)}+\theta_{1}^{(2)}X$. Finally, we wish to remark that it is unclear whether the convergence rate of MLE in Theorem~\ref{theorem:lower_bound_Gaussian_family_second_singularity_case} is sharp due to the complex behaviors of the system of limits~\eqref{eqn:system_complex_polynomial_limits}. We leave the sharpness of that rate for the future work.
\paragraph{General picture:} In general, if we have an expert function $h_{1}(X|\theta_{1}) = (\theta_{1}^{(1)}+\theta_{1}^{(2)}X)^{m}$ for some positive integer $m \geq 1$, and expert function $h_{2}$ is independent of covariate $X$ as in equation~\eqref{eqn:non_linear_expert}, then we also have that $h_{1}$ and $h_{2}$ are algebraically dependent. The corresponding PDE strucure is the following:
\begin{align}
\parenth{\dfrac{\partial h_{1}}{\partial{\theta_{1}^{(1)}}}(X, \theta_{1})}^{2} = m^2(\theta_{1}^{(1)})^{2( m-1)}\dfrac{\partial{h_{2}^{2}}}{\partial{\theta_{2}}}(X, \theta_{2}), \label{eqn:PDE_structure_nonlinear_expert_general}
\end{align}
for all $\theta_{1} = (\theta_{1}^{(1)}, 0)$ and $\theta_{2}$. This PDE structure captures a phase transition between nonlinearity setting I and nonlinearity setting II. More precisely, we can check that the convergence rate of $\widehat{G}_{n}$ will be $n^{-1/4}$, which is similar to that in Theorem~\ref{theorem:lower_bound_Gaussian_family_third} under the nonlinearity setting I. Under the nonlinearity setting II, the convergence rates of $\widehat{G}_{n}$ are again determined by a system of polynomial limits, which is dependent on $m$ and much more complicated than that in equation~\eqref{eqn:system_complex_polynomial_limits}. A useful insight that arises from these systems is that the convergence rates of $(\theta_{1i}^{0})^{(1)}$ and $\theta_{2i}^{0}$ are better than $n^{-1/2\overline{r}}$ and $n^{-1/2\ceil{\overline{r}/2}}$ respectively while that of $(\theta_{1i}^{0})^{(2)}$ is $n^{-1/4}$ for any $1 \leq i \leq k_{0}$. As a consequence, the convergence rates for parameter estimation when $m \geq 2$ are always better than $m=1$.
\subsection{Beyond Condition (O.2) or (O.3): Other algebraically dependent expert functions}
\label{sec:other_algebraically_dependent_expert}
Thus far, we have studied the convergence rates of MLE under a few representative settings when the expert functions are algebraically dependent and do not satisfy Condition (O.1) in Definition~\ref{definition:identifiability}. In this section, we study these convergence rates under some specific settings when the expert functions do not satisfy either Condition (O.2) or (O.3). 
\subsubsection{Beyond Only One Condition}
\label{subsec:beyond_O2}
We first study a setting when the expert functions $h_{1}$ and $h_{2}$ do not satisfy Condition (O.2) while they satisfy Conditions (O.1) and (O.3). 
\begin{example} \label{example:condition_O2_without_offset}
For $X \in \mathcal{X} \subset \mathbb{R}_{+}$, we define expert functions $h_{1}(X,\theta_{1}) = \theta_{1} X^2$ for all $\theta_{1} \in \Omega_{1} \subset \mathbb{R}$ and $h_{2}^{2}(X,\theta_{2}) = \theta_{2}^{(1)} + \theta_{2}^{(2)}X + \theta_{2}^{(3)}X^{2}$ for all $\theta_{2}=(\theta_{2}^{(1)}, \theta_{2}^{(2)}, \theta_{2}^{(3)}) \in \Omega_{2} \subset \mathbb{R}_{+}^{3}$ such that $\theta_{2}^{(1)}+\theta_{2}^{(2)} + \theta_{2}^{(3)} \geq \overline{\gamma}$ for some positive constant $\overline{\gamma}$. Then, we have the following PDE for these expert functions:
\begin{align}
    \parenth{\frac{\partial{h_{2}^2}}{\partial{\theta_{2}^{(2)}}}(X,\theta_{1})}^{2} = \dfrac{\partial{h_{2}^{2}}}{\partial{\theta_{2}^{(1)}}}(X,\theta_{2})\dfrac{\partial{h_{2}^{2}}}{\partial{\theta_{2}^{(3)}}}(X,\theta_{2}), \label{eqn:PDE_structure_O2}
\end{align} 
which shows that $h_{1}$ and $h_{2}$ are algebraically dependent.
\end{example}
Equation~\eqref{eqn:PDE_structure_O2} indicates that the expert functions in Example~\ref{example:condition_O2_without_offset} do not satisfy Condition (O.2). We can verify that these expert functions satisfy Conditions (O.1) and (O.3). Hence, these expert functions only do not satisfy Condition (O.2) in Definition~\ref{definition:identifiability}. The following theorem establishes the sharp convergence rate of MLE under this setting of expert functions. 
\begin{theorem} \label{theorem:lower_bound_Gaussian_family_without_offset_O2}
Assume that for $X \in \mathcal{X} \subset \mathbb{R}_{+}$ the expert functions $h_{1}(X,\theta_{1}) = \theta_{1} X^2$ for all $\theta_{1} \in \Omega_{1} \subset \mathbb{R}$ and $h_{2}^{2}(X,\theta_{2}) = \theta_{2}^{(1)} + \theta_{2}^{(2)}X + \theta_{2}^{(3)}X^{2}$ for all $\theta_{2}=(\theta_{2}^{(1)}, \theta_{2}^{(2)}, \theta_{2}^{(3)}) \in \Omega_{2} \subset \mathbb{R}_{+}^{3}$ such that $\theta_{2}^{(1)}+\theta_{2}^{(2)} + \theta_{2}^{(3)} \geq \overline{\gamma}$ for some positive constant $\overline{\gamma}$. Then, the following holds:
\begin{itemize}
\item[(a)] (Convergence rate of MLE) There exists a positive constant $C_{0}$ depending only on $G_{0}$ and $\Omega$ such that
\begin{eqnarray}
\mathbb{P}\parenth{\widetilde{W}_{\kappa}(\widehat{G}_{n},G_{0}) > C_{0}(\log n/n)^{1/4}} \precsim \exp(-c\log n), \nonumber
\end{eqnarray}
where $\kappa = (2, 2, 2, 2)$. Here, $c$ is a positive constant depending only on $\Omega$.
\item[(b)] (Minimax lower bound) For any $\kappa'$ such that $(1,1,1,1) \preceq \kappa' \prec \kappa = (2, 2, 2, 2)$, 
\begin{align}
\inf \limits_{\overline{G}_{n}} \sup \limits_{G \in \Ocal_{k}(\Omega) \backslash \Ocal_{k_{0}-1}(\Omega)} \Exs_{p_{G}} \parenth{\widetilde{W}_{\kappa'}(\overline{G}_{n},G)} \succsim n^{-1/(2\|\kappa'\|_{\infty})}. \nonumber
\end{align}
Here, the infimum is taken over all sequences of estimates $\overline{G}_{n} \in \Ocal_{k}(\Omega)$. Furthermore, $\Exs_{p_{G}}$ denotes the expectation taken with respect to the product measure with mixture density $p_{G}^{n}$. 
\end{itemize}
\end{theorem}
\noindent
Proof of Theorem~\ref{theorem:lower_bound_Gaussian_family_without_offset_O2} is in Appendix~\ref{subsec:proof:theorem:lower_bound_Gaussian_family_without_offset_O2}.

The result of Theorem~\ref{theorem:lower_bound_Gaussian_family_without_offset_O2} entails that even though the expert functions are algebraically dependent due to their violation of Condition (O.2), the convergence rate of MLE under this setting is still $n^{-1/4}$ (up to some logarithmic factor), which is similar that of MLE when the expert functions are algebraically independent in Theorem~\ref{theorem:total_variation_bound_over-fitted_MECFG}. The convergence rate $n^{-1/4}$ of the MLE under the generalized transportation distance also leads the the uniform convergence rates $n^{-1/4}$ for estimating the individual components $\theta_{1i}^{0}$ and $(\theta_{2i}^{0})^{(u)}$ for $1 \leq u \leq 3$ and $1 \leq i \leq k_{0}$.

We now move to another example of algebraically dependent expert functions $h_{1}$, $h_{2}$ when they do not satisfy Condition (O.3) while they satisfy Conditions (O.1) and (O.2).
\begin{example} \label{example:condition_O3}
For $X \in \mathcal{X} \subset \mathbb{R}_{+}$, we define expert functions $h_{1}(X,\theta_{1}) = \theta_{1}^{(1)} + \theta_{1}^{(2)} X^2$ for all $\theta_{1} = (\theta_{1}^{(1)}, \theta_{1}^{(2)}) \in \Omega_{1} \subset \mathbb{R}^2$ and $h_{2}^{2}(X,\theta_{2}) = \theta_{2}^{(1)} X + \theta_{2}^{(2)}X^{3}$ for all $\theta_{2}=(\theta_{2}^{(1)}, \theta_{2}^{(2)}) \in \Omega_{2} \subset \mathbb{R}_{+}^{2}$ such that $\theta_{2}^{(1)}+\theta_{2}^{(2)} \geq \overline{\gamma}$ for some positive constant $\overline{\gamma}$. Then, we have the following PDE for these expert functions:
\begin{align}
    \frac{\partial{h_{1}}}{\partial{\theta_{1}^{(1)}}}(X,\theta_{1}) \frac{\partial{h_{2}^2}}{\partial{\theta_{2}^{(2)}}}(X,\theta_{2}) = \frac{\partial{h_{1}}}{\partial{\theta_{1}^{(2)}}}(X,\theta_{1}) \frac{\partial{h_{2}^2}}{\partial{\theta_{2}^{(1)}}}(X,\theta_{2}), \label{eqn:PDE_structure_O3}
\end{align} 
which shows that $h_{1}$ and $h_{2}$ are algebraically dependent.
\end{example}
The PDE~\eqref{eqn:PDE_structure_O3} indicates that the expert functions $h_{1}$ and $h_{2}$ in Example~\ref{example:condition_O3} do not satisfy Condition (O.3). We can check that these expert functions still satisfy Conditions (O.1) and (O.2). The sharp convergence rate of MLE under this setting of expert functions is established in the following theorem.
\begin{theorem} \label{theorem:lower_bound_Gaussian_family_O3}
Assume that for $X \in \mathcal{X} \subset \mathbb{R}_{+}$ the expert functions $h_{1}(X,\theta_{1}) = \theta_{1}^{(1)} + \theta_{1}^{(2)} X^2$ for all $\theta_{1} = (\theta_{1}^{(1)}, \theta_{1}^{(2)}) \in \Omega_{1} \subset \mathbb{R}^2$ and $h_{2}^{2}(X,\theta_{2}) = \theta_{2}^{(1)} X + \theta_{2}^{(2)}X^{3}$ for all $\theta_{2}=(\theta_{2}^{(1)}, \theta_{2}^{(2)}) \in \Omega_{2} \subset \mathbb{R}_{+}^{2}$ such that $\theta_{2}^{(1)}+\theta_{2}^{(2)} \geq \overline{\gamma}$ for some positive constant $\overline{\gamma}$. Then, the following holds:
\begin{itemize}
\item[(a)] (Convergence rate of MLE) There exists a positive constant $C_{0}$ depending only on $G_{0}$ and $\Omega$ such that
\begin{eqnarray}
\mathbb{P}\parenth{\widetilde{W}_{\kappa}(\widehat{G}_{n},G_{0}) > C_{0}(\log n/n)^{1/4}} \precsim \exp(-c\log n), \nonumber
\end{eqnarray}
where $\kappa = (2, 2, 2, 2)$. Here, $c$ is a positive constant depending only on $\Omega$.
\item[(b)] (Minimax lower bound) For any $\kappa'$ such that $(1,1,1,1) \preceq \kappa' \prec \kappa = (2, 2, 2, 2)$, 
\begin{align}
\inf \limits_{\overline{G}_{n}} \sup \limits_{G \in \Ocal_{k}(\Omega) \backslash \Ocal_{k_{0}-1}(\Omega)} \Exs_{p_{G}} \parenth{\widetilde{W}_{\kappa'}(\overline{G}_{n},G)} \succsim n^{-1/(2\|\kappa'\|_{\infty})}. \nonumber
\end{align}
Here, the infimum is taken over all sequences of estimates $\overline{G}_{n} \in \Ocal_{k}(\Omega)$. Furthermore, $\Exs_{p_{G}}$ denotes the expectation taken with respect to the product measure with mixture density $p_{G}^{n}$. 
\end{itemize}
\end{theorem}
\noindent
The proof of Theorem~\ref{theorem:lower_bound_Gaussian_family_O3} is in Appendix~\ref{subsec:proof:theorem:lower_bound_Gaussian_family_O3}. 

Interestingly, similar to the result of Theorem~\ref{theorem:lower_bound_Gaussian_family_without_offset_O2}, the convergence rate of parameter estimation for the particular setting in Example~\ref{example:condition_O3} that the expert functions do not satisfy Condition (O.3) but still satisfy Conditions (O.1) and (O.2) is still $n^{-1/4}$ (up to some logarithmic factor). That convergence rate is identical to the convergence rate of MLE when the expert functions are algebraically independent. Furthermore, that convergence rate of MLE directly leads to the uniformly convergence rates $n^{-1/4}$ of individual components $(\theta_{1i}^0)^{(u)}$ and $(\theta_{1i}^0)^{(v)}$ for $1 \leq u, v \leq 2$ and $1 \leq i \leq k_{0}$. 
\subsubsection{Beyond More Than One Conditions}
\label{subsec:beyond_O3}
We now discuss a few specific settings of the expert functions $h_{1}$ and $h_{2}$ when they do not satisfy more than one conditions in Definition~\ref{definition:identifiability}. The first example is when the expert functions $h_{1}$ and $h_{2}$ satisfy Condition (O.2) but do not satisfy Conditions (O.1) and (O.3).
\begin{example} \label{example:condition_O3_O1}
For $X \in \mathcal{X} \subset \mathbb{R}_{+}$, we define expert functions $h_{1}(X,\theta_{1}) = \theta_{1}^{(1)} + \theta_{1}^{(2)} X^2$ for all $\theta_{1} = (\theta_{1}^{(1)}, \theta_{1}^{(2)}) \in \Omega_{1} \subset \mathbb{R}^2$ and $h_{2}^{2}(X,\theta_{2}) = \theta_{2}^{(1)} + \theta_{2}^{(2)}X^{2}$ for all $\theta_{2}=(\theta_{2}^{(1)}, \theta_{2}^{(2)}) \in \Omega_{2} \subset \mathbb{R}_{+}^{2}$ such that $\theta_{2}^{(1)}+\theta_{2}^{(2)} \geq \overline{\gamma}$ for some positive constant $\overline{\gamma}$. Then, we have the following PDEs for these expert functions:
\begin{align}
    \frac{\partial{h_{2}^2}}{\partial{\theta_{2}^{(1)}}}(X,\theta_{2}) = \parenth{\dfrac{\partial{h_{1}}}{\partial{\theta_{1}^{(1)}}}(X,\theta_{1})}^{2}, \label{eqn:PDE_structure_O3_O1_first} \\
    \frac{\partial{h_{2}^2}}{\partial{\theta_{2}^{(2)}}}(X,\theta_{2}) = \dfrac{\partial{h_{1}}}{\partial{\theta_{1}^{(1)}}}(X,\theta_{1}) \dfrac{\partial{h_{1}}}{\partial{\theta_{1}^{(2)}}}(X,\theta_{1}), \label{eqn:PDE_structure_O3_O1_second} \\
    \frac{\partial{h_{1}}}{\partial{\theta_{1}^{(1)}}}(X,\theta_{1}) \frac{\partial{h_{2}^2}}{\partial{\theta_{2}^{(2)}}}(X,\theta_{2}) = \frac{\partial{h_{1}}}{\partial{\theta_{1}^{(2)}}}(X,\theta_{1}) \frac{\partial{h_{2}^2}}{\partial{\theta_{2}^{(1)}}}(X,\theta_{2}), \label{eqn:PDE_structure_O3_O1_third}
\end{align} 
which shows that $h_{1}$ and $h_{2}$ are algebraically dependent.
\end{example}
Equations~\eqref{eqn:PDE_structure_O3_O1_first} and~\eqref{eqn:PDE_structure_O3_O1_second} demonstrate that the expert functions in Example~\ref{example:condition_O3_O1} do not satisfy Condition (O.1) while equation~\eqref{eqn:PDE_structure_O3_O1_third} proves that these expert functions do not satisfy Condition (O.3). We can check that these expert functions still satisfy Condition (O.2). The following result establishes the sharp convergence rate of MLE under this setting of these expert functions.
\begin{theorem} \label{theorem:lower_bound_Gaussian_family_O3_O1}
Assume that for $X \in \mathcal{X} \subset \mathbb{R}_{+}$ the expert functions $h_{1}(X,\theta_{1}) = \theta_{1}^{(1)} + \theta_{1}^{(2)} X^2$ for all $\theta_{1} = (\theta_{1}^{(1)}, \theta_{1}^{(2)}) \in \Omega_{1} \subset \mathbb{R}^2$ and $h_{2}^{2}(X,\theta_{2}) = \theta_{2}^{(1)} + \theta_{2}^{(2)}X^{2}$ for all $\theta_{2}=(\theta_{2}^{(1)}, \theta_{2}^{(2)}) \in \Omega_{2} \subset \mathbb{R}_{+}^{2}$ such that $\theta_{2}^{(1)}+\theta_{2}^{(2)} \geq \overline{\gamma}$ for some positive constant $\overline{\gamma}$. Then, the following holds:
\begin{itemize}
\item[(a)] (Convergence rate of MLE) There exists a positive constant $C_{0}$ depending only on $G_{0}$ and $\Omega$ such that
\begin{eqnarray}
\mathbb{P}\parenth{\widetilde{W}_{\kappa}(\widehat{G}_{n},G_{0}) > C_{0}(\log n/n)^{1/2\overline{r}}} \precsim \exp(-c\log n), \nonumber
\end{eqnarray}
where $\kappa = (\bar{r}, 2, \ceil{\overline{r}/2}, 2)$. Here, $c$ is a positive constant depending only on $\Omega$.
\item[(b)] (Minimax lower bound) For any $\kappa'$ such that $(1,1,1,1) \preceq \kappa' \prec \kappa = (\bar{r}, 2, \ceil{\overline{r}/2}, 2)$, 
\begin{align}
\inf \limits_{\overline{G}_{n}} \sup \limits_{G \in \Ocal_{k}(\Omega) \backslash \Ocal_{k_{0}-1}(\Omega)} \Exs_{p_{G}} \parenth{\widetilde{W}_{\kappa'}(\overline{G}_{n},G)} \succsim n^{-1/(2\|\kappa'\|_{\infty})}. \nonumber
\end{align}
Here, the infimum is taken over all sequences of estimates $\overline{G}_{n} \in \Ocal_{k}(\Omega)$. Furthermore, $\Exs_{p_{G}}$ denotes the expectation taken with respect to the product measure with mixture density $p_{G}^{n}$. 
\end{itemize}
\end{theorem}
\noindent 
The proof of Theorem~\ref{theorem:lower_bound_Gaussian_family_O3_O1} is in Appendix~\ref{subsec:proof:theorem:lower_bound_Gaussian_family_O3_O1}. 

The results of Theorem~\ref{theorem:lower_bound_Gaussian_family_O3_O1} indicate that the convergence rate of MLE is $n^{-1/2 \bar{r}}$ and sharp when the expert functions $h_{1}(X,\theta_{1}) = \theta_{1}^{(1)} + \theta_{1}^{(2)} X^2$ and $h_{2}^2(X, \theta_{2}) = \theta_{2}^{(1)} + \theta_{2}^{(2)} X^2$. This convergence rate indicates that the convergence rates for estimating $(\theta_{1i}^{0})^{(1)}, (\theta_{1i}^{0})^{(2)}, (\theta_{2i}^{0})^{(1)}, (\theta_{2i}^{0})^{(2)}$ are respectively $n^{-1/2\bar{r}}, n^{-1/4}, n^{-1/2\ceil{\overline{r}/2}}$, and $n^{-1/4}$ for all $1 \leq i \leq k_{0}$. The slow convergence rates for estimating $(\theta_{1i}^{0})^{(1)}$ and $(\theta_{2i}^{0})^{(1)}$ are due to the linear dependence of these components in equation~\eqref{eqn:PDE_structure_O3_O1_first}. Interestingly, even though we have linear dependence of $(\theta_{1i}^{0})^{(2)}$ and $(\theta_{2i}^{0})^{(2)}$ in equations~\eqref{eqn:PDE_structure_O3_O1_second} and~\eqref{eqn:PDE_structure_O3_O1_third}, the convergence rates for estimating these components are still $n^{-1/4}$ and align with those when the algebraic independence is satisfied.

Our final example is when the expert functions $h_{1}$ and $h_{2}$ do not satisfy all Conditions (O.1), (O.2), and (O.3).
\begin{example} \label{example:condition_O3_O2_O1}
For $X \in \mathcal{X} \subset \mathbb{R}_{+}$, we define expert functions $h_{1}(X,\theta_{1}) = \theta_{1}^{(1)} + \theta_{1}^{(2)} X^2$ for all $\theta_{1} = (\theta_{1}^{(1)}, \theta_{1}^{(2)}) \in \Omega_{1} \subset \mathbb{R}^2$ and $h_{2}^{2}(X,\theta_{2}) = \theta_{2}^{(1)} + \theta_{2}^{(2)}X + \theta_{2}^{(3)}X^2$ for all $\theta_{2}=(\theta_{2}^{(1)}, \theta_{2}^{(2)}, \theta_{2}^{(3)}) \in \Omega_{2} \subset \mathbb{R}_{+}^{3}$ such that $\theta_{2}^{(1)}+\theta_{2}^{(2)} + \theta_{2}^{(3)} \geq \overline{\gamma}$ for some positive constant $\overline{\gamma}$. Then, we have the following PDEs for these expert functions:
\begin{align}
    \frac{\partial{h_{2}^2}}{\partial{\theta_{2}^{(1)}}}(X,\theta_{2}) = \parenth{\dfrac{\partial{h_{1}}}{\partial{\theta_{1}^{(1)}}}(X,\theta_{1})}^{2}, \label{eqn:PDE_structure_O3_O2_O1_first} \\
    \frac{\partial{h_{2}^2}}{\partial{\theta_{2}^{(2)}}}(X,\theta_{2}) = \dfrac{\partial{h_{1}}}{\partial{\theta_{1}^{(1)}}}(X,\theta_{1}) \dfrac{\partial{h_{1}}}{\partial{\theta_{1}^{(2)}}}(X,\theta_{1}), \label{eqn:PDE_structure_O3_O2_O1_second} \\
    \parenth{\frac{\partial{h_{2}^2}}{\partial{\theta_{2}^{(2)}}}(X,\theta_{2})}^{2} = \dfrac{\partial{h_{2}^{2}}}{\partial{\theta_{2}^{(1)}}}(X,\theta_{2})\dfrac{\partial{h_{2}^{2}}}{\partial{\theta_{2}^{(3)}}}(X,\theta_{2}), \label{eqn:PDE_structure_O3_O2_O1_third} \\
    \frac{\partial{h_{1}}}{\partial{\theta_{1}^{(1)}}}(X,\theta_{1}) \frac{\partial{h_{2}^2}}{\partial{\theta_{2}^{(3)}}}(X,\theta_{2}) = \frac{\partial{h_{1}}}{\partial{\theta_{1}^{(2)}}}(X,\theta_{1}) \frac{\partial{h_{2}^2}}{\partial{\theta_{2}^{(1)}}}(X,\theta_{2}), \label{eqn:PDE_structure_O3_O2_O1_fourth}
\end{align} 
which shows that $h_{1}$ and $h_{2}$ are algebraically dependent.
\end{example}
Equations~\eqref{eqn:PDE_structure_O3_O2_O1_first} and~\eqref{eqn:PDE_structure_O3_O2_O1_second} demonstrate that the expert functions in Example~\ref{example:condition_O3_O2_O1} do not satisfy Condition (O.1). Equation~\eqref{eqn:PDE_structure_O3_O2_O1_third} proves that these expert functions do not satisfy Condition (O.2) and equation~\eqref{eqn:PDE_structure_O3_O2_O1_fourth} indicates that these expert functions do not satisfy Condition (O.3). The following theorem establishes the optimal convergence rate of parameter estimation $\widehat{G}_{n}$ of the over-specified GMCF model under the setting of these expert functions.
\begin{theorem} \label{theorem:lower_bound_Gaussian_family_O3_O2_O1}
Assume that for $X \in \mathcal{X} \subset \mathbb{R}_{+}$ the expert functions $h_{1}(X,\theta_{1}) = \theta_{1}^{(1)} + \theta_{1}^{(2)} X^2$ for all $\theta_{1} = (\theta_{1}^{(1)}, \theta_{1}^{(2)}) \in \Omega_{1} \subset \mathbb{R}^2$ and $h_{2}^{2}(X,\theta_{2}) = \theta_{2}^{(1)} + \theta_{2}^{(2)}X + \theta_{2}^{(3)}X^2$ for all $\theta_{2}=(\theta_{2}^{(1)}, \theta_{2}^{(2)}, \theta_{2}^{(3)}) \in \Omega_{2} \subset \mathbb{R}_{+}^{3}$ such that $\theta_{2}^{(1)}+\theta_{2}^{(2)} + \theta_{2}^{(3)} \geq \overline{\gamma}$ for some positive constant $\overline{\gamma}$. Then, the following holds:
\begin{itemize}
\item[(a)] (Convergence rate of MLE) There exists a positive constant $C_{0}$ depending only on $G_{0}$ and $\Omega$ such that
\begin{eqnarray}
\mathbb{P}\parenth{\widetilde{W}_{\kappa}(\widehat{G}_{n},G_{0}) > C_{0}(\log n/n)^{1/2\overline{r}}} \precsim \exp(-c\log n), \nonumber
\end{eqnarray}
where $\kappa = (\bar{r}, 2, \ceil{\overline{r}/2}, 2, 2)$. Here, $c$ is a positive constant depending only on $\Omega$.
\item[(b)] (Minimax lower bound) For any $\kappa'$ such that $(1,1,1,1) \preceq \kappa' \prec \kappa = (\bar{r}, 2, \ceil{\overline{r}/2}, 2, 2)$, 
\begin{align}
\inf \limits_{\overline{G}_{n}} \sup \limits_{G \in \Ocal_{k}(\Omega) \backslash \Ocal_{k_{0}-1}(\Omega)} \Exs_{p_{G}} \parenth{\widetilde{W}_{\kappa'}(\overline{G}_{n},G)} \succsim n^{-1/(2\|\kappa'\|_{\infty})}. \nonumber
\end{align}
Here, the infimum is taken over all sequences of estimates $\overline{G}_{n} \in \Ocal_{k}(\Omega)$. Furthermore, $\Exs_{p_{G}}$ denotes the expectation taken with respect to the product measure with mixture density $p_{G}^{n}$. 
\end{itemize}
\end{theorem}
\noindent 
The proof of Theorem~\ref{theorem:lower_bound_Gaussian_family_O3_O2_O1} is in Appendix~\ref{subsection:proof_theorem:lower_bound_Gaussian_family_O3_O2_O1}.

A few comments with the results of Theorem~\ref{theorem:lower_bound_Gaussian_family_O3_O2_O1} are in order. First, the convergence rate of the MLE under generalized transportation distance indicates that the convergence rates for estimating $(\theta_{1i}^{0})^{(1)}, (\theta_{1i}^{0})^{(2)}, (\theta_{2i}^{0})^{(1)}, (\theta_{2i}^{0})^{(2)}, (\theta_{2i}^{0})^{(3)}$ are respectively $n^{-1/2 \bar{r}}, n^{-1/4}, n^{-1/2\ceil{\overline{r}/2}}, n^{-1/4}$, and $n^{-1/4}$ (up to some logarithmic factors) for all $1 \leq i \leq k_{0}$. The slow convergence rates for estimating $(\theta_{1i}^{0})^{(1)}$ and $(\theta_{2i}^{0})^{(1)}$ are due to the PDE~\eqref{eqn:PDE_structure_O3_O2_O1_first}, which entails a dependency among these parameters. On the other hand, despite of the dependency of $(\theta_{1i}^{0})^{(2)}$, $(\theta_{2i}^{0})^{(2)}$, and $(\theta_{2i}^{0})^{(3)}$ via PDEs~\eqref{eqn:PDE_structure_O3_O2_O1_second}-\eqref{eqn:PDE_structure_O3_O2_O1_fourth}, the convergence rates for estimating these elements are still comparable to those when the algebraically independent assumption holds. This result again suggests that the Condition (O.1) seems to be the key condition to slow down the convergence rates of the MLE $\widehat{G}_{n}$ and its individual components. 

\section{Proofs of key results} \label{Section:proofs}
In this section, we provide the proofs of the key theoretical results in the paper while deferring the rest to the Appendices.
Our proof techniques build on previous work for establishing the sharp convergence rates for parameter estimation under traditional finite mixture models~\citep{Chen-95, Jonas-2017, Ho-Nguyen-SIAM-18}
and are based on using a generalized transportation distance to provide controls on various Taylor expansions.  We begin with a lemma that presents a general strategy for obtaining sharp convergence rates. 
\begin{lemma} \label{lemma:convergence_rates_MLE}
(a) (MLE estimation) Assume that there exists some $\kappa \in \mathbb{N}^{q_{1}+q_{2}}$ such that
\begin{eqnarray}
\inf \limits_{G \in \mathcal{G}} h(p_{G},p_{G_{0}})/ \widetilde{W}_{\kappa}^{\|\kappa\|_{\infty}}(G,G_{0}) > 0, \label{eqn:lower_bound_Hellinger_generalized_Wasserstein}
\end{eqnarray}
where $\mathcal{G}$ is a subset of $\Ocal_{k}(\Omega)$ for the over-fitted setting of the GMCF model. Then there exists some positive constant $C_{0}$ depending only on $G_{0}$ and $\Omega$ such that
\begin{eqnarray}
\mathbb{P}(\widetilde{W}_{\kappa}(\widehat{G}_{n},G_{0}) > C_{0}(\log n/n)^{1/2\|\kappa\|_{\infty}}) \precsim \exp(-c\log n), \nonumber
\end{eqnarray}
where $c$ is a positive constant depending only on $\Omega$.

\noindent
(b) (Minimax lower bound) Assume that inequality \eqref{eqn:lower_bound_Hellinger_generalized_Wasserstein} holds for any $G_{0} \in \mathcal{G}$. Furthermore, as long as $G_{0} \in \mathcal{G}$, the following holds
\begin{eqnarray}
\inf \limits_{G \in \mathcal{G}} h(p_{G},p_{G_{0}})/W_{\kappa'}^{\|\kappa'\|_{\infty}}(G,G_{0}) = 0 \label{eqn:tight_lower_bound}
\end{eqnarray}
for all $\kappa' \prec \kappa$. Then, for any $\kappa'$such that $(1,\ldots,1) \preceq \kappa' \prec \kappa$, 
\begin{align}
\inf \limits_{\overline{G}_{n} \in \mathcal{G}} \sup \limits_{G \in \mathcal{G} \backslash \Ocal_{k_{0}-1}(\Omega)} \Exs_{p_{G}} \parenth{\widetilde{W}_{\kappa'}(\overline{G}_{n},G)} \geq c'n^{-1/(2\|\kappa'\|_{\infty})}. \nonumber
\end{align}
Here, $\Exs_{p_{G}}$ denotes the expectation taken with respect to product measure with mixture density $p_{G}^{n}$, and $c'$ stands for a universal constant depending on $\Omega$.
\end{lemma}
\noindent
Proof of Lemma~\ref{lemma:convergence_rates_MLE} is in Appendix~\ref{subsec:proof:lemma:convergence_rates_MLE}.
\subsection{Proof of Theorem \ref{theorem:total_variation_bound_over-fitted_MECFG}}
\label{Section:proof_strong_identifiable}
Given Lemma~\ref{lemma:convergence_rates_MLE}, we obtain the conclusion of Theorem~\ref{theorem:total_variation_bound_over-fitted_MECFG}, by demonstrating the following results:
\begin{eqnarray}
\inf \limits_{G \in \Ocal_{k}(\Omega)} h(p_{G},p_{G_{0}})/\widetilde{W}_{\kappa}^{\|\kappa\|_{\infty}}(G,G_{0}) & >0, \label{eqn:general_overfit_first} \\
\inf \limits_{G \in \Ocal_{k}(\Omega)} h(p_{G},p_{G_{0}})/\widetilde{W}_{\kappa'}^{\|\kappa'\|_{\infty}}(G,G_{0}) & = 0 \label{eqn:general_overfit_second}
\end{eqnarray}
for any $\kappa' \prec \kappa$ where $\kappa=(2,\ldots,2)$. The proof of inequality~\eqref{eqn:general_overfit_first} is in Section~\ref{subsection:key_inequality_algebraic_independent} while the proof of equality~\eqref{eqn:general_overfit_second} is in Section~\ref{subsection:equality_overfit}.
\subsubsection{Proof for inequality \eqref{eqn:general_overfit_first}}
\label{subsection:key_inequality_algebraic_independent}
The proof of inequality~\eqref{eqn:general_overfit_first} is divided into two parts: local structure and global structure. 
\paragraph{Local structure:} We first demonstrate that inequality~\eqref{eqn:general_overfit_first} holds when $\widetilde{W}_{\kappa}(G,G_{0})$ is sufficiently small. In particular, we will prove that
\begin{align}
\lim \limits_{\epsilon \to 0} \inf \limits_{G \in \Ocal_{k}(\Omega): \widetilde{W}_{\kappa}(G,G_{0}) \leq \epsilon} h(p_{G},p_{G_{0}})/\widetilde{W}_{\kappa}^{\|\kappa\|_{\infty}}(G,G_{0})>0. \nonumber
\end{align}
Due to the standard lower bound $h \geq V$, it is sufficient to show that
\begin{align}
\lim \limits_{\epsilon \to 0} \inf \limits_{G \in \Ocal_{k}(\Omega): \widetilde{W}_{\kappa}(G,G_{0}) \leq \epsilon} V(p_{G},p_{G_{0}})/\widetilde{W}_{\kappa}^{\|\kappa\|_{\infty}}(G,G_{0})>0. \nonumber
\end{align}
Assume that the above statement does not hold. This implies that we can find a sequence $G_{n} \in \Ocal_{k}(\Omega)$ such that $V(p_{G_{n}},p_{G_{0}})/\widetilde{W}_{\kappa}^{\|\kappa\|_{\infty}}(G_{n},G_{0}) \to 0$ and $\widetilde{W}_{\kappa}(G_{n},G_{0}) \to 0$ as $n \to \infty$. As being demonstrated in Lemma~\ref{lemma:relabel_sequence} in Appendix~\ref{sec:auxi_result}, we can assume the sequence $G_{n}$ has exactly $\bar{k}$ atoms, where $k_{0} \leq \bar{k} \leq k$, and can be represented as follows:
\begin{align}
G_{n} = \sum \limits_{i=1}^{k_{0}+\overline{l}} \sum \limits_{j=1}^{s_{i}} p_{ij}^{n} \delta_{(\theta_{1ij}^{n}, \theta_{2ij}^{n})}, \nonumber
\end{align}
where $\overline{l} \geq 0$ is some nonnegative integer and $s_{i} \geq 1$ for $1 \leq i \leq k_{0}+\overline{l}$ such that $\sum \limits_{i=1}^{k_{0}+\overline{l}} s_{i} = \overline{k}$. Additionally, $(\theta_{1ij}^{n}, \theta_{2ij}^{n}) \to (\theta_{1i}^{0}, \theta_{2i}^{0})$ and $\sum \limits_{j=1}^{s_{i}} p_{ij}^{n} \to \pi_{i}^{0}$ for all $1 \leq i \leq k_{0}+\overline{l}$. Here, $\pi_{i}^{0} = 0$ as $k_{0} +1 \leq i \leq \overline{k}$ while $(\theta_{1i}^{0},\theta_{2i}^{0})$ are possible extra limit points from the convergence of components of $G_{n}$ as $k_{0} +1 \leq i \leq \overline{k}$. 

Now, according to Lemma~\ref{lemma:generalized_Wasserstein_distance_polynomials} in Appendix~\ref{sec:auxi_result}, we have
\begin{align}
\widetilde{W}_{\kappa}^{\|\kappa\|_{\infty}}(G,G_{0}) & \precsim \sum \limits_{i=1}^{k_{0} + \overline{l}}\sum \limits_{j=1}^{s_{i}} p_{ij}^{n} d_{\kappa}^{\|\kappa\|_{\infty}}\parenth{\eta_{ij}^{n},\eta_{i}^{0}} + \sum \limits_{i=1}^{k_{0} + \overline{l}} \abss{\sum \limits_{j=1}^{s_{i}}p_{ij}^{n} - \pi_{i}^{0}} \nonumber \\
& = \sum \limits_{i=1}^{k_{0} + \overline{l}}\sum \limits_{j=1}^{s_{i}} p_{ij}^{n} \parenth{\enorm{\theta_{1ij}^{n} - \theta_{1i}^{0}}^{2} + \enorm{\theta_{2ij}^{n} - \theta_{2i}^{0}}^{2}} + \sum \limits_{i=1}^{k_{0} + \overline{l}}\abss{\sum \limits_{j=1}^{s_{i}}p_{ij}^{n} - \pi_{i}^{0}} \nonumber \\
& : = D_{\kappa}(G_{n},G_{0}), \nonumber
\end{align}
where $\kappa=(2,\ldots,2)$, $\eta_{i}^{0} = \parenth{\theta_{1i}^{0},\theta_{2i}^{0}}$ and $\eta_{ij}^{n} = \parenth{\theta_{1ij}^{n}, \theta_{2ij}^{n}}$ for $1 \leq i \leq k_{0}$ and $1 \leq j \leq s_{i}$. For the simplicity of presentation, we introduce the following notation: $\Delta \theta_{1ij}^{n} : =  \theta_{1ij}^{n} - \theta_{1i}^{0}$, $\Delta \theta_{2ij}^{n} : =  \theta_{2ij}^{n} - \theta_{2i}^{0}$ for $1 \leq i \leq k_{0} + \overline{l}$ and $1 \leq j \leq s_{i}$. Additionally, we denote $\Delta \theta_{1ij}^{n} : = \parenth{\parenth{\Delta \theta_{1ij}^{n}}^{(1)}, \ldots, \parenth{\Delta \theta_{1ij}^{n}}^{(q_{1})}}$ and $\Delta \theta_{2ij}^{n} : = \parenth{\parenth{\Delta \theta_{2ij}^{n}}^{(1)}, \ldots, \parenth{\Delta \theta_{2ij}^{n}}^{(q_{2})}}$ for all $i,j$. 

Since $V(p_{G_{n}},p_{G_{0}})/\widetilde{W}_{\kappa}^{\|\kappa\|_{\infty}}(G_{n},G_{0}) \to 0$ as $n \to \infty$, we obtain that $V(p_{G_{n}},p_{G_{0}})/D_{\kappa}(G_{n},G_{0})$\\$ \to 0$. To facilitate the proof argument, we divide it into several steps.
\paragraph{Step 1 - Structure of Taylor expansion:} By means of a Taylor expansion up to the second order, for any $1 \leq i \leq k_{0} + \overline{l}$ and $1 \leq j \leq s_{i}$, the following holds:
\begin{align}
f\parenth{Y| h_{1}(X, \theta_{1ij}^{n}), h_{2}(X, \theta_{2ij}^{n})} - f\parenth{Y| h_{1}(X, \theta_{1i}^{0}), h_{2}(X, \theta_{2i}^{0})} & \nonumber \\
& \hspace{- 24 em} = \sum \limits_{1 \leq |\alpha|+|\beta| \leq 2} \dfrac{1}{\alpha!\beta!}\prod \limits_{u=1}^{q_{1}} \biggr\{(\Delta \theta_{1ij}^{n})^{(u)}\biggr\}^{\alpha_{u}}\prod \limits_{v=1}^{q_{2}}\biggr\{(\Delta \theta_{2ij}^{n})^{(v)}\biggr\}^{\beta_{v}} \dfrac{\partial^{|\alpha| + |\beta|}{f}}{\partial{\theta_{1}^{\alpha}}\partial{\theta_{2}^{\beta}}}\parenth{Y|h_{1}(X,\theta_{1i}^{0}),h_{2}(X,\theta_{2i}^{0})} \nonumber \\
& + R_{ij}(X,Y), \nonumber
\end{align}
where $\alpha = \parenth{\alpha_{1}, \ldots, \alpha_{q_{1}}}$, $\beta = \parenth{\beta_{1}, \ldots, \beta_{q_{2}}}$, $\abss{\alpha} = \alpha_{1} + \ldots + \alpha_{q_{1}}$, and $\abss{\beta} = \beta_{1} + \ldots + \beta_{q_{2}}$. $R_{ij}(X,Y)$ is the remainder from the Taylor expansion and it satisfies
\begin{align}
R_{ij}(X,Y)\overline{f}(X) = \mathcal{O}\parenth{\enorm{\Delta \theta_{1ij}}^{2+\gamma} + \enorm{\Delta \theta_{2ij}}^{2+\gamma}}, \nonumber
\end{align}
for some universal constant $\gamma > 0$ for all $1 \leq i \leq k_{0}$ and $1 \leq j \leq s_{i}$. We thus have:
\begin{align}
& \hspace{-3 em} p_{G_{n}}(X,Y)-p_{G_{0}}(X,Y) \nonumber \\
& =  \sum \limits_{i=1}^{k_{0} + \overline{l}} \sum \limits_{j = 1}^{s_{i}} p_{ij}^{n} \brackets{f(Y| h_{1}(X, \theta_{1ij}^{n}), h_{2}(X, \theta_{2ij}^{n}) - f(Y| h_{1}(X, \theta_{1i}^{0}), h_{2}(X, \theta_{2i}^{0})}\overline{f}(X) \nonumber \\
& \hspace{ 2 em} + \sum \limits_{i=1}^{k_{0}+\overline{l}}\biggr(\sum \limits_{j=1}^{s_{i}}p_{ij}^{n} - \pi_{i}^{0}\biggr)f(Y|h_{1}(X,\theta_{1i}^{0}),h_{2}(X,\theta_{2i}^{0}))\overline{f}(X) \nonumber
\end{align}
\begin{align}
& = \sum \limits_{i=1}^{k_{0} + \overline{l}}\sum \limits_{j=1}^{s_{i}}p_{ij}^{n}\sum \limits_{1 \leq |\alpha|+|\beta| \leq 2} \dfrac{1}{\alpha!\beta!}\prod \limits_{u=1}^{q_{1}} \biggr\{(\Delta \theta_{1ij}^{n})^{(u)}\biggr\}^{\alpha_{u}}\prod \limits_{v=1}^{q_{2}}\biggr\{(\Delta \theta_{2ij}^{n})^{(v)}\biggr\}^{\beta_{v}} \nonumber \\
&  \hspace{ 10 em} \times \dfrac{\partial^{|\alpha| + |\beta|}{f}}{\partial{\theta_{1}^{\alpha}}\partial{\theta_{2}^{\beta}}}\parenth{Y|h_{1}(X,\theta_{1i}^{0}),h_{2}(X,\theta_{2i}^{0})}\overline{f}(X) \nonumber \\
&  \hspace{ 2 em} + \sum \limits_{i=1}^{k_{0}+\overline{l}}\biggr(\sum \limits_{j=1}^{s_{i}}p_{ij}^{n} - \pi_{i}^{0}\biggr)f(Y|h_{1}(X,\theta_{1i}^{0}),h_{2}(X,\theta_{2i}^{0}))\overline{f}(X) + R(X,Y)
\nonumber \\
&   : = A_{n} + B_{n} + R(X,Y), \nonumber
\end{align}
where $R(X,Y) = \parenth{ \sum \limits_{i=1}^{k_{0} + \overline{l}} \sum \limits_{j = 1}^{s_{i}} R_{ij}(X,Y)} \overline{f}(X) = \mathcal{O}\parenth{\sum \limits_{i = 1}^{k_{0} + \overline{l}} \sum \limits_{j = 1}^{s_{i}} p_{ij}^{n}\brackets{\enorm{\Delta \theta_{1ij}}^{2+\gamma} + \enorm{\Delta \theta_{2ij}}^{2+\gamma}}}$. From the formulation of $D_{\kappa}(G_{n},G_{0})$, it is clear that
\begin{align}
R(X,Y) / D_{\kappa}(G_{n},G_{0}) \precsim \sum \limits_{i = 1}^{k_{0} + \overline{l}} \sum \limits_{j = 1}^{s_{i}} \brackets{\enorm{\Delta \theta_{1ij}}^{\gamma} + \enorm{\Delta \theta_{2ij}}^{\gamma}} \to 0 \label{eqn:Taylor_remainder_limit}
\end{align}
as $n \to \infty$. For the univariate location-scale Gaussian distribution, we have the following characteristic PDE:
\begin{align}
\dfrac{\partial^{2}{f}}{\partial{\mu^2}}(x,\mu,\sigma) = 2\dfrac{\partial{f}}{\partial{\sigma^{2}}}(x,\mu,\sigma), \label{eqn:PDE_Gaussian_family}
\end{align}
where $\mu$ and $\sigma$ respectively stand for the location and scale parameter in a location-scale Gaussian distribution. Governed by that PDE, we find that
\begin{align}
\frac{\partial^{2}{f}}{\partial{h_{1}^{2}}}\parenth{Y|h_{1}(X,\theta_{1}),h_{2}(X,\theta_{2})} = 2\frac{\partial{f}}{\partial{h_{2}^{2}}}\parenth{Y|h_{1}(X,\theta_{1}),h_{2}(X,\theta_{2})},
\end{align}
for all $(\theta_{1},\theta_{2})$. Therefore, for any $(\theta_{1}, \theta_{2})$, a straightforward calculation yields the following:
\begin{align}
\frac{\partial{f}}{\partial{\theta_{1}^{(u)}}}\parenth{Y|h_{1}(X,\theta_{1}),h_{2}(X,\theta_{2})} & = \frac{\partial{h_{1}}}{\partial{\theta_{1}^{(u)}}}(X,\theta_{1})\frac{\partial{f}}{\partial{h_{1}}}\parenth{Y|h_{1}(X,\theta_{1}),h_{2}(X,\theta_{2})}, \nonumber \\
\frac{\partial{f}}{\partial{\theta_{2}^{(v)}}}\parenth{Y|h_{1}(X,\theta_{1}),h_{2}(X,\theta_{2})} & = \frac{\partial{h_{2}^{2}}}{\partial{\theta_{2}^{(v)}}}(X,\theta_{2})\frac{\partial{f}}{\partial{h_{2}^{2}}}\parenth{Y|h_{1}(X,\theta_{1}),h_{2}(X,\theta_{2})} \nonumber \\
& = \dfrac{1}{2}\frac{\partial{h_{2}^{2}}}{\partial{\theta_{2}^{(v)}}}(X,\theta_{2})\frac{\partial^{2}{f}}{\partial{h_{1}^{2}}}\parenth{Y|h_{1}(X,\theta_{1}),h_{2}(X,\theta_{2})}, \nonumber
\end{align}
for all $1 \leq u \leq q_{1}$ and $1 \leq v \leq q_{2}$.
Similarly, the PDE structure~\eqref{eqn:PDE_Gaussian_family} leads to
\begin{align}
\frac{\partial^{2}{f}}{\partial{\theta_{1}^{(u)}}\partial{\theta_{1}^{(v)}}}\parenth{Y|h_{1}(X,\theta_{1}),h_{2}(X,\theta_{2})} = \frac{\partial^{2}{h_{1}}}{\partial{\theta_{1}^{(u)}}\partial{\theta_{1}^{(v)}}}(X,\theta_{1})\frac{\partial{f}}{\partial{h_{1}}}\parenth{Y|h_{1}(X,\theta_{1}),h_{2}(X,\theta_{2})} \nonumber \\
+ \frac{\partial{h_{1}}}{\partial{\theta_{1}^{(u)}}}(X,\theta_{1})\frac{\partial{h_{1}}}{\partial{\theta_{1}^{(v)}}}(X,\theta_{1})\frac{\partial^{2}{f}}{\partial{h_{1}^{2}}}\parenth{Y|h_{1}(X,\theta_{1}),h_{2}(X,\theta_{2})}, \nonumber \\
\frac{\partial^{2}{f}}{\partial{\theta_{2}^{(u)}}\partial{\theta_{2}^{(v)}}}\parenth{Y|h_{1}(X,\theta_{1}),h_{2}(X,\theta_{2})} = \frac{\partial^{2}{h_{2}^{2}}}{\partial{\theta_{2}^{(u)}}\partial{\theta_{2}^{(v)}}}(X,\theta_{2})\frac{\partial{f}}{\partial{h_{2}^{2}}}\parenth{Y|h_{1}(X,\theta_{1}),h_{2}(X,\theta_{2})} \nonumber \\
+ \frac{\partial{h_{2}^{2}}}{\partial{\theta_{2}^{(u)}}}(X,\theta_{2})\frac{\partial{h_{2}^{2}}}{\partial{\theta_{2}^{(v)}}}(X,\theta_{2})\frac{\partial^{2}{f}}{\partial{h_{2}^{4}}}\parenth{Y|h_{1}(X,\theta_{1}),h_{2}(X,\theta_{2})} \nonumber \\
= \dfrac{1}{2}\frac{\partial^{2}{h_{2}^{2}}}{\partial{\theta_{2}^{(u)}}\partial{\theta_{2}^{(v)}}}(X,\theta_{2})\frac{\partial^{2}{f}}{\partial{h_{1}^{2}}}\parenth{Y|h_{1}(X,\theta_{1}),h_{2}(X,\theta_{2})} \nonumber \\
+ \dfrac{1}{4}\frac{\partial{h_{2}^{2}}}{\partial{\theta_{2}^{(u)}}}(X,\theta_{2})\frac{\partial{h_{2}^{2}}}{\partial{\theta_{2}^{(v)}}}(X,\theta_{2})\frac{\partial^{4}{f}}{\partial{h_{1}^{4}}}\parenth{Y|h_{1}(X,\theta_{1}),h_{2}(X,\theta_{2})}, \nonumber \\
\frac{\partial^{2}{f}}{\partial{\theta_{1}^{(u)}}\partial{\theta_{2}^{(v)}}}\parenth{Y|h_{1}(X,\theta_{1}),h_{2}(X,\theta_{2})} = \frac{\partial{h_{1}}}{\partial{\theta_{1}^{(u)}}}(X,\theta_{1})\frac{\partial{h_{2}^{2}}}{\partial{\theta_{2}^{(v)}}}(X,\theta_{2})\frac{\partial^{2}{f}}{\partial{h_{1}}\partial{h_{2}^{2}}}\parenth{Y|h_{1}(X,\theta_{1}),h_{2}(X,\theta_{2})} \nonumber \\
= \frac{1}{2} \frac{\partial{h_{1}}}{\partial{\theta_{1}^{(u)}}}(X,\theta_{1})\frac{\partial{h_{2}^{2}}}{\partial{\theta_{2}^{(v)}}}(X,\theta_{2})\frac{\partial^{3}{f}}{\partial{h_{1}^{3}}}\parenth{Y|h_{1}(X,\theta_{1}),h_{2}(X,\theta_{2})}, \nonumber
\end{align} 
for all $u,v$. Equipped with the above equations, we can rewrite $A_{n}$ as follows
\begin{align}
A_{n}  = \sum \limits_{i = 1}^{k_{0} + \overline{l}} \sum \limits_{\tau = 1}^{4} A_{n,\tau}^{(i)}(X) \frac{\partial^{\tau}{f}}{\partial{h_{1}^{\tau}}}\parenth{Y|h_{1}(X,\theta_{1i}^{0}),h_{2}(X,\theta_{2i}^{0})}\overline{f}(X) : = \sum \limits_{i=1}^{k_{0} + \overline{l}} \overline{A}_{n,\tau}^{(i)}(X,Y), \nonumber
\end{align}
where the explicit forms of $A_{n,\tau}^{(i)}(X)$ are
\begin{align}
& A_{n,1}^{(i)}(X) : = \sum \limits_{j=1}^{s_{i}} p_{ij}^{n}\biggr(\sum \limits_{u=1}^{q_{1}} \parenth{\Delta \theta_{1ij}^{n}}^{(u)}\frac{\partial{h_{1}}}{\partial{\theta_{1}^{(u)}}}(X, \theta_{1i}^{0}) \nonumber \\
& \hspace{ 16 em} + \sum \limits_{1 \leq u,v \leq q_{1}} \frac{\parenth{\Delta \theta_{1ij}^{n}}^{(u)}\parenth{\Delta \theta_{1ij}^{n}}^{(v)}}{1+1_{\{u = v\}}} \frac{\partial^{2}{h_{1}}}{\partial{\theta_{1}^{(u)}}\partial{\theta_{1}^{(v)}}}(X, \theta_{1i}^{0})\biggr), \nonumber \\
& A_{n,2}^{(i)}(X) : =  \sum \limits_{j=1}^{s_{i}} p_{ij}^{n}\biggr\{\frac{1}{2} \sum \limits_{u=1}^{q_{2}} \parenth{\Delta \theta_{2ij}^{n}}^{(u)}\frac{\partial{h_{2}^{2}}}{\partial{\theta_{2}^{(u)}}}(X, \theta_{2i}^{0}) \nonumber \\
& \hspace{ 15 em} + \sum \limits_{1 \leq u,v \leq q_{1}} \frac{\parenth{\Delta \theta_{1ij}^{n}}^{(u)}\parenth{\Delta \theta_{1ij}^{n}}^{(v)}}{1+1_{\{u = v\}}} \frac{\partial{h_{1}}}{\partial{\theta_{1}^{(u)}}}(X, \theta_{1i}^{0})\frac{\partial{h_{1}}}{\partial{\theta_{1}^{(v)}}}(X, \theta_{1i}^{0})  \nonumber \\
& \hspace{ 15 em} + \frac{1}{2} \sum \limits_{1 \leq u,v \leq q_{2}} \frac{\parenth{\Delta \theta_{2ij}^{n}}^{(u)}\parenth{\Delta \theta_{2ij}^{n}}^{(v)}}{1+1_{\{u = v\}}} \frac{\partial^{2}{h_{2}^{2}}}{\partial{\theta_{2}^{(u)}}\partial{\theta_{2}^{(v)}}}(X, \theta_{2i}^{0}) \biggr\}, \nonumber
\end{align}
\begin{align}
& A_{n,3}^{(i)}(X) : =  \frac{1}{2} \sum \limits_{j=1}^{s_{i}} p_{ij}^{n}  \sum \limits_{u=1}^{q_{1}}\sum \limits_{v=1}^{q_{2}} \parenth{\Delta \theta_{1ij}^{n}}^{(u)}\parenth{\Delta \theta_{2ij}^{n}}^{(v)} \frac{\partial{h_{1}}}{\partial{\theta_{1}^{(u)}}}(X,\theta_{1i}^{0})\frac{\partial{h_{2}^{2}}}{\partial{\theta_{2}^{(v)}}}(X,\theta_{2i}^{0}), \nonumber \\
& A_{n,4}^{(i)}(X) : =  \frac{1}{4} \sum \limits_{j=1}^{s_{i}} p_{ij}^{n}  \sum \limits_{1 \leq u,v \leq q_{2}} \frac{\parenth{\Delta \theta_{2ij}^{n}}^{(u)}\parenth{\Delta \theta_{2ij}^{n}}^{(v)}}{1 + 1_{\{u = v\}}} \frac{\partial{h_{2}^{2}}}{\partial{\theta_{2}^{(u)}}}(X,\theta_{2i}^{0})\frac{\partial{h_{2}^{2}}}{\partial{\theta_{2}^{(v)}}}(X,\theta_{2i}^{0}). \nonumber 
\end{align}

In view of the above computations, we can treat $\overline{A}_{n,\tau}^{(i)}(X, Y)/ D_{\kappa}(G_{n},G_{0})$ as a linear combinations of elements from $\mathcal{F}_{\tau}(i)$ for $1 \leq \tau \leq 4$, which can be defined as follows:
\begin{align}
\mathcal{F}_{1}(i) &: = \left\{ \frac{\partial{h_{1}}}{\partial{\theta_{1}^{(u)}}}(X, \theta_{1i}^{0})\frac{\partial{f}}{\partial{h_{1}}}\parenth{Y|h_{1}(X,\theta_{1i}^{0}),h_{2}(X,\theta_{2i}^{0})}\overline{f}(X): \ 1 \leq u \leq q_{1} \right\} \nonumber \\
& \cup \left\{ \frac{\partial^{2}{h_{1}}}{\partial{\theta_{1}^{(u)}}\partial{\theta_{1}^{(v)}}}(X, \theta_{1i}^{0})\frac{\partial^{2}{f}}{\partial{h_{1}^{2}}}\parenth{Y|h_{1}(X,\theta_{1i}^{0}),h_{2}(X,\theta_{2i}^{0})}\overline{f}(X): \ 1 \leq u,v \leq q_{1} \right\}, \nonumber \\
\mathcal{F}_{2}(i) & : = \left\{ \frac{\partial{h_{2}^{2}}}{\partial{\theta_{2}^{(u)}}}(X, \theta_{2i}^{0})\frac{\partial^{2}{f}}{\partial{h_{1}^{2}}}\parenth{Y|h_{1}(X,\theta_{1i}^{0}),h_{2}(X,\theta_{2i}^{0})}\overline{f}(X): \ 1 \leq u \leq q_{2} \right\} \nonumber \\
& \cup \left\{ \frac{\partial{h_{1}}}{\partial{\theta_{1}^{(u)}}}(X, \theta_{1i}^{0})\frac{\partial{h_{1}}}{\partial{\theta_{1}^{(v)}}}(X, \theta_{1i}^{0})\frac{\partial^{2}{f}}{\partial{h_{1}^{2}}}\parenth{Y|h_{1}(X,\theta_{1i}^{0}),h_{2}(X,\theta_{2i}^{0})}\overline{f}(X): \ 1 \leq u ,v\leq q_{1} \right\} \nonumber \\
& \cup \left\{\frac{\partial^{2}{h_{2}^{2}}}{\partial{\theta_{2}^{(u)}}\partial{\theta_{2}^{(v)}}}(X, \theta_{2i}^{0})\frac{\partial^{2}{f}}{\partial{h_{1}^{2}}}\parenth{Y|h_{1}(X,\theta_{1i}^{0}),h_{2}(X,\theta_{2i}^{0})}\overline{f}(X): \ 1 \leq u,v \leq q_{2} \right\}, \nonumber 
\end{align}
\begin{align}
\mathcal{F}_{3}(i) &: = \biggr\{ \frac{\partial{h_{1}}}{\partial{\theta_{1}^{(u)}}}(X,\theta_{1i}^{0})\frac{\partial{h_{2}^{2}}}{\partial{\theta_{2}^{(v)}}}(X,\theta_{2i}^{0})\frac{\partial^{3}{f}}{\partial{h_{1}^{3}}}\parenth{Y|h_{1}(X,\theta_{1i}^{0}),h_{2}(X,\theta_{2i}^{0})}\overline{f}(X): \nonumber \\
& \hspace{ 24 em} 1 \leq u \leq q_{1}, \ 1 \leq v \leq q_{2} \biggr\}, \nonumber \\
\mathcal{F}_{4}(i) &: = \left\{ \frac{\partial{h_{2}^{2}}}{\partial{\theta_{2}^{(u)}}}(X,\theta_{2i}^{0})\frac{\partial{h_{2}^{2}}}{\partial{\theta_{2}^{(v)}}}(X,\theta_{2i}^{0}) \frac{\partial^{4}{f}}{\partial{h_{1}^{4}}}\parenth{Y|h_{1}(X,\theta_{1i}^{0}),h_{2}(X,\theta_{2i}^{0})}\overline{f}(X): \ 1 \leq u, v \leq q_{2} \right\}. \nonumber
\end{align}
Therefore, we can view $A_{n}/ D_{\kappa}(G_{n},G_{0})$ as a linear combination of elements from $\mathcal{F} : = \cup_{i=1}^{k_{0} + \overline{l}} \cup_{j = 1}^{4} \mathcal{F}_{j}(i)$. 
Similarly, we can view $B_{n}/ D_{\kappa}(G_{n},G_{0})$ as a linear combination of elements of the form $f(Y|h_{1}(X,\theta_{1i}^{0}),h_{2}(X,\theta_{2i}^{0}))\overline{f}(X)$ for $1 \leq i \leq k_{0} + \overline{l}$. 
\paragraph{Step 2 - Non-vanishing coefficients:} Assume that all of the coefficients in the representation of $A_{n}/D_{\kappa}(G_{n},G_{0})$ and $B_{n}/D_{\kappa}(G_{n},G_{0})$ go to 0 as $n \to \infty$, namely, the coefficients of the elements from $\mathcal{F}$ and of elements of the form $f(Y|h_{1}(X,\theta_{1i}^{0}),h_{2}(X,\theta_{2i}^{0}))\overline{f}(X)$ for $1 \leq i \leq k_{0} + \overline{l}$ go to 0. By taking the summation of the absolute values of the coefficients of $B_{n} / D_{\kappa}(G_{n},G_{0})$, the following limit holds
\begin{align}
\sum \limits_{i=1}^{k_{0}+\overline{l}}\abss{\sum \limits_{j=1}^{s_{i}}p_{ij}^{n} - \pi_{i}^{0}}\bigg/D_{\kappa}(G_{n},G_{0}) \to 0. \nonumber
\end{align}
From the expression for $D_{\kappa}(G_{n},G_{0})$, this yields:
\begin{align}
\sum \limits_{i=1}^{k_{0} + \overline{l}}\sum \limits_{j=1}^{s_{i}} p_{ij}^{n} \parenth{\enorm{\Delta \theta_{1ij}^{n}}^{2} + \enorm{\Delta \theta_{2ij}^{n}}^{2}}\bigg/D_{\kappa}(G_{n},G_{0}) \to 1. \label{eqn:first_limit}
\end{align}
On the other hand, according to the formulation of $A_{n,4}^{(i)}(X)$, the coefficients associated with the elements $\parenth{\frac{\partial{h_{2}^{2}}}{\partial{\theta_{2}^{(u)}}}(X,\theta_{2i}^{0})}^{2} \frac{\partial^{4}{f}}{\partial{h_{1}^{4}}}$ in $\mathcal{F}_{4}(i)$ are $\sum \limits_{j=1}^{s_{i}} p_{ij}^{n} \biggr\{\parenth{\Delta \theta_{2ij}^{n}}^{(u)}\biggr\}^{2}/\brackets{8 D_{\kappa}(G_{n},G_{0})}$ as $1 \leq i \leq k_{0} + \overline{l}$ and $1 \leq u \leq q_{2}$. According to the hypothesis, these coefficients go to zero; therefore, by taking the summation of all of these coefficients, we obtain that
\begin{align}
\sum \limits_{i=1}^{k_{0} + \overline{l}}\sum \limits_{j=1}^{s_{i}} p_{ij}^{n} \enorm{\Delta \theta_{2ij}^{n}}^{2} \bigg/D_{\kappa}(G_{n},G_{0}) \to 0.  \label{eqn:first_sublimit}
\end{align}
Furthermore, from the formulation of $A_{n.2}^{(i)}(X)$, we can check that the coefficients attached to the elements $\parenth{\frac{\partial{h_{1}}}{\partial{\theta_{1}^{(u)}}}(X, \theta_{1i}^{0})}^{2} \frac{\partial^{2}{f}}{\partial{h_{1}^{2}}}\parenth{Y|h_{1}(X,\theta_{1i}^{0}),h_{2}(X,\theta_{2i}^{0})}\overline{f}(X)$ in $\mathcal{F}_{2}(i)$ are 
\begin{align}
\sum \limits_{j=1}^{s_{i}} p_{ij}^{n} \biggr\{\parenth{\Delta \theta_{1ij}^{n}}^{(u)}\biggr\}^{2}\bigg/\brackets{2 D_{\kappa}(G_{n},G_{0})}, \nonumber
\end{align}
as $1 \leq i \leq k_{0} + \overline{l}$ and $1 \leq u \leq q_{1}$. As all of these coefficients go to zero, by taking the summation of these coefficients, we obtain the following limit:
\begin{align}
\sum \limits_{i=1}^{k_{0} + \overline{l}}\sum \limits_{j=1}^{s_{i}} p_{ij}^{n} \enorm{\Delta \theta_{1ij}^{n}}^{2} \bigg/D_{\kappa}(G_{n},G_{0}) \to 0.  \label{eqn:second_sublimit}
\end{align}
Combining the results from equations~\eqref{eqn:first_sublimit} and~\eqref{eqn:second_sublimit}, the following limit holds:
\begin{align}
\sum \limits_{i=1}^{k_{0} + \overline{l}}\sum \limits_{j=1}^{s_{i}} p_{ij}^{n} \parenth{\enorm{\Delta \theta_{1ij}^{n}}^{2} + \enorm{\Delta \theta_{2ij}^{n}}^{2}} \bigg/D_{\kappa}(G_{n},G_{0}) \to 0, \nonumber
\end{align}
which is a contradiction to equation~\eqref{eqn:first_limit}. Therefore, not all the coefficients in the representation of $A_{n}/ D_{\kappa}(G_{n},G_{0})$ and $B_{n}/ D_{\kappa}(G_{n},G_{0})$ go to zero as $n \to \infty$. 
\paragraph{Step 3 - Fatou's argument:} We denote $m_{n}$ as the maximum of the absolute values of the coefficients in the representation of $A_{n}/ D_{\kappa}(G_{n},G_{0})$ and $B_{n}/ D_{\kappa}(G_{n},G_{0})$. From here, we define $d_{n} : = 1/m_{n}$. Since not all the coefficients of $A_{n}/ D_{\kappa}(G_{n},G_{0})$ and $B_{n}/ D_{\kappa}(G_{n},G_{0})$ vanish, we have $d_{n} \not \to \infty$ as $n \to \infty$. From the definition of $m_{n}$, we denote 
\begin{align}
\parenth{\sum \limits_{j=1}^{s_{i}}p_{ij}^{n} - \pi_{i}^{0}}/ m_{n} \to \alpha(i); \ 
\parenth{\sum \limits_{j=1}^{s_{i}} p_{ij}^{n} \parenth{\Delta \theta_{\tau ij}^{n}}^{(u)}}/m_{n} \to \beta_{\tau u}(i), \nonumber \\
\parenth{\sum \limits_{j=1}^{s_{i}} p_{ij}^{n} \parenth{\Delta \theta_{\tau ij}^{n}}^{(u)} \parenth{\Delta \theta_{\tau ij}^{n}}^{(v)}} / m_{n} \to \gamma_{\tau uv}(i), \nonumber \\
\parenth{\sum \limits_{j=1}^{s_{i}} p_{ij}^{n} \parenth{\Delta \theta_{1ij}^{n}}^{(u)} \parenth{\Delta \theta_{2ij}^{n}}^{(v)}} / m_{n} \to \eta_{uv}(i), \nonumber
\end{align}
as $n \to \infty$ for all $1 \leq i \leq k_{0} + \overline{l}$ and all $u, v$. Here, at least one among $\alpha(i), \beta_{\tau u}(i), \gamma_{\tau uv}(i)$, and $\eta_{uv}(i)$ is different from zero for all $i, u, v$. Invoking Fatou's lemma, we have: 
\begin{align}
0 = \lim \limits_{n \to \infty} d_{n}\frac{V(p_{G_{n}},p_{G_{0}})}{D_{\kappa}(G_{n},G_{0})} \geq \int \mathop{\lim \inf} \limits_{n \to \infty} d_{n}\frac{\abss{p_{G_{n}}(X,Y)-p_{G_{0}}(X,Y)}}{D_{\kappa}(G_{n},G_{0})}d(X,Y). \label{eqn:Fatou_argument}
\end{align}
From the definition of $\alpha(i), \beta_{\tau u}(i), \gamma_{\tau uv}(i), \eta_{uv}(i)$, the following holds:
\begin{align}
d_{n}\frac{p_{G_{n}}(X,Y)-p_{G_{0}}(X,Y)}{D_{\kappa}(G_{n},G_{0})} & \to \sum \limits_{i = 1}^{k_{0} + \overline{l}} \sum \limits_{\tau = 0}^{4} E_{\tau}^{(i)}(X) \frac{\partial^{\tau}{f}}{\partial{h_{1}^{\tau}}}\parenth{Y|h_{1}(X,\theta_{1i}^{0}),h_{2}(X,\theta_{2i}^{0})}\overline{f}(X), \label{eqn:second_limit}
\end{align}
for all $(X,Y)$ where the expressions for $E_{\tau}^{(i)}(X)$ are:
\begin{align}
E_{0}^{(i)}(X) & : = \alpha(i), \ E_{1}^{(i)}(X) = \sum \limits_{u=1}^{q_{1}} \beta_{1u}(i) \frac{\partial{h_{1}}}{\partial{\theta_{1}^{(u)}}}(X, \theta_{1i}^{0}) + \sum \limits_{1 \leq u,v \leq q_{1}} \frac{\gamma_{1uv}(i)}{1+1_{\{u = v\}}} \frac{\partial^{2}{h_{1}}}{\partial{\theta_{1}^{(u)}}\partial{\theta_{1}^{(v)}}}(X, \theta_{1i}^{0}), \nonumber \\
E_{2}^{(i)}(X) & : = \frac{1}{2} \sum \limits_{u=1}^{q_{2}} \beta_{2u}(i)\frac{\partial{h_{2}^{2}}}{\partial{\theta_{2}^{(u)}}}(X, \theta_{2i}^{0}) + \sum \limits_{1 \leq u,v \leq q_{1}} \frac{\gamma_{1uv}(i)}{1+1_{\{u = v\}}} \frac{\partial{h_{1}}}{\partial{\theta_{1}^{(u)}}}(X, \theta_{1i}^{0})\frac{\partial{h_{1}}}{\partial{\theta_{1}^{(v)}}}(X, \theta_{1i}^{0})  \nonumber \\
& + \frac{1}{2} \sum \limits_{1 \leq u,v \leq q_{2}} \frac{\gamma_{2uv}(i)}{1+1_{\{u = v\}}} \frac{\partial^{2}{h_{2}^{2}}}{\partial{\theta_{2}^{(u)}}\partial{\theta_{2}^{(v)}}}(X, \theta_{2i}^{0}), \nonumber \\
E_{3}^{(i)}(X) & : =  \frac{1}{2} \sum \limits_{u=1}^{q_{1}}\sum \limits_{v=1}^{q_{2}} \eta_{uv}(i) \frac{\partial{h_{1}}}{\partial{\theta_{1}^{(u)}}}(X,\theta_{1i}^{0})\frac{\partial{h_{2}^{2}}}{\partial{\theta_{2}^{(v)}}}(X,\theta_{2i}^{0}), \nonumber \\
E_{4}^{(i)}(X) & : =  \frac{1}{4} \sum \limits_{1 \leq u,v \leq q_{2}} \frac{\gamma_{2uv}(i)}{1 + 1_{\{u = v\}}} \frac{\partial{h_{2}^{2}}}{\partial{\theta_{2}^{(u)}}}(X,\theta_{2i}^{0})\frac{\partial{h_{2}^{2}}}{\partial{\theta_{2}^{(v)}}}(X,\theta_{2i}^{0}). \nonumber
\end{align}
Combining the results from equations~\eqref{eqn:Fatou_argument} and~\eqref{eqn:second_limit}, the following equation holds
\begin{align}
\sum \limits_{i = 1}^{k_{0} + \overline{l}} \sum \limits_{\tau = 0}^{4} E_{\tau}^{(i)}(X) \frac{\partial^{\tau}{f}}{\partial{h_{1}^{\tau}}}\parenth{Y|h_{1}(X,\theta_{1i}^{0}),h_{2}(X,\theta_{2i}^{0})}\overline{f}(X) = 0, \nonumber
\end{align}
almost surely $(X,Y)$. For almost surely $X$, the set
\begin{align}
\left\{\frac{\partial^{\tau}{f}}{\partial{h_{1}^{\tau}}}\parenth{Y|h_{1}(X,\theta_{1i}^{0}),h_{2}(X,\theta_{2i}^{0})}: 0 \leq \tau \leq 4 \right\} \nonumber
\end{align}
is linearly independent with respect to $Y$. Therefore, the above equation eventually leads to $E_{\tau}^{(i)}(X) = 0$ almost surely $X$ for $0 \leq \tau \leq 4$ and $1 \leq i \leq k_{0} + \overline{l}$. 

When $\tau = 0$, it is clear that the equation $E_{\tau}^{(i)}(X) = 0$ almost surely $X$ demonstrates that $\alpha(i) = 0$ for all $i$. When $\tau \geq 3$, since the expert functions $h_{1}$ and $h_{2}$ are algebraically independent, the equations $E_{\tau}^{(i)}(X) = 0$ almost surely $X$ lead to $\gamma_{2uv}(i) = 0$ and $\eta_{uv}(i) = 0$ for all $(u,v)$ and $i$. Furthermore, invoking the fact that the expert functions $h_{1}$ and $h_{2}$ are algebraically independent and the result that $\gamma_{2uv}(i) = 0$ for all $(u,v)$, the equation $E_{2}^{(i)}(X) = 0$ almost surely $X$ implies that $\beta_{2u}(i) = 0$ and $\gamma_{1uv}(i) = 0$ for all $i$ and $(u,v)$. Collecting the previous results, the equation $E_{1}^{(i)}(X) = 0$ almost surely $X$ leads to $\beta_{1u}(i) = 0$ for all $i$ and $u$. Therefore, all the coefficients $\alpha(i), \beta_{\tau u}(i), \gamma_{\tau uv}(i)$, and $\eta_{uv}(i)$ are equal to zero for all $i$ and $u,v$, which is a contradiction.

As a consequence, we can find some $\epsilon_{0}>0$ such that 
\begin{align}
\inf \limits_{G \in \Ocal_{k}(\Omega): \widetilde{W}_{\kappa}(G,G_{0}) \leq \epsilon_{0}} h(p_{G},p_{G_{0}})/\widetilde{W}_{\kappa}^{\|\kappa\|_{\infty}}(G,G_{0})>0. \nonumber
\end{align}
\paragraph{Global structure:} Given the local bound that we have just established, to obtain the conclusion of inequality~\eqref{eqn:general_overfit_first}, it is sufficient to demonstrate that
\begin{align}
\inf \limits_{G \in \Ocal_{k}(\Omega): \widetilde{W}_{\kappa}(G,G_{0}) > \epsilon_{0}} h(p_{G},p_{G_{0}})/\widetilde{W}_{\kappa}^{\|\kappa\|_{\infty}}(G,G_{0})>0. \nonumber
\end{align}
Assume that the above result does not hold. This indicates that we can find a sequence $\overline{G}_{n} \in \Ocal_{k}(\Omega)$ such that $h(p_{\overline{G}_{n}},p_{G_{0}})/\widetilde{W}_{\kappa}^{\|\kappa\|_{\infty}}(\overline{G}_{n},G_{0}) \to 0$ as $n \to \infty$ while $\widetilde{W}_{\kappa}(\overline{G}_{n},G_{0}) > \epsilon_{0}$ for all $n \geq 1$. Since the set $\Omega$ is bounded, there exists a subsequence of $G_{n}$ such that $G_{n} \to G'$ for some mixing measure $G' \in \Ocal_{k}(\Omega)$. To facilitate the discussion, we replace this subsequence by the whole sequence of $G_{n}$. Then, as $\widetilde{W}_{\kappa}(\overline{G}_{n},G_{0}) > \epsilon_{0}$ for all $n \geq 1$, this implies that $\widetilde{W}_{\kappa}(G',G_{0}) \geq \epsilon_{0}$. Combining the previous bound with $h(p_{\overline{G}_{n}},p_{G_{0}})/\widetilde{W}_{\kappa}^{\|\kappa\|_{\infty}}(\overline{G}_{n},G_{0}) \to 0$, we obtain that $h(p_{\overline{G}_{n}},p_{G_{0}}) \to 0$ as $n \to \infty$. Invoking Fatou's lemma, the following inequality holds:
\begin{align}
0 = \lim \limits_{n \to \infty} h^{2}(p_{\overline{G}_{n}},p_{G_{0}}) & \geq \frac{1}{2}\int \mathop {\lim \inf} \limits_{n \to \infty} \parenth{\sqrt{p_{\overline{G}_{n}}(X,Y)} - \sqrt{p_{G_{0}}(X,Y)}}^{2}d(X,Y) \nonumber \\
& = \frac{1}{2} \int \parenth{\sqrt{p_{G'}(X,Y)} - \sqrt{p_{G_{0}}(X,Y)}}^{2}d(X,Y). \nonumber
\end{align}
This inequality leads to $p_{G'}(X,Y) = p_{G_{0}}(X,Y)$ for almost surely $X,Y$. Due to the identifiability of GMCF, this leads to $G' \equiv G_{0}$, which is a contradiction to the result that $\widetilde{W}_{\kappa}(G',G_{0}) \geq \epsilon_{0} > 0$. Hence, we achieve the conclusion of inequality~\eqref{eqn:general_overfit_first}.
\subsubsection{Proof for equality \eqref{eqn:general_overfit_second}}
\label{subsection:equality_overfit}
To achieve the conclusion of equality~\eqref{eqn:general_overfit_second}, it is equivalent to find a sequence $G_{n} \in \Ocal_{k}(\Omega)$ such that $h(p_{G_{n}},p_{G_{0}})/\widetilde{W}_{\kappa'}^{\|\kappa'\|_{\infty}}(G_{n},G_{0}) \to 0$ as $n \to \infty$ for every $\kappa' \prec \kappa$. In fact, for any $\kappa' \prec \kappa$, we have $\min \limits_{1 \leq i \leq q_{1} + q_{2}} \kappa'^{(i)} < 2$. Without loss of generality, we assume $\kappa'^{(1)} = \min \limits_{1 \leq i \leq q_{1} + q_{2}} \kappa'^{(i)} < 2$. Now, we construct a sequence of mixing measures, $G_{n} = \sum_{i=1}^{k_{0}+1} \pi_{i}^{n}\delta_{\parenth{\theta_{1i}^{n},\theta_{2i}^{n}}}$, with $k_{0}+1$ components as follows: $(\pi_{i}^{n}, \theta_{1i}^{n},\theta_{2i}^{n}) \equiv (\pi_{i-1}^{0}, \theta_{1(i-1)}^{0}, \theta_{2(i-1)}^{0})$ for $3 \leq i \leq k_{0} + 1$. Additionally, $\pi_{1}^{n} = \pi_{2}^{n} = 1/2$, $(\theta_{11}^{n},\theta_{21}^{n}) \equiv (\theta_{11}^{0} - \vec{1}_{q_{1}}/n, \theta_{21}^{0} - \vec{1}_{q_{2}}/n)$, and $(\theta_{12}^{n},\theta_{22}^{n}) \equiv (\theta_{11}^{0} + \vec{1}_{q_{1}}/n, \theta_{21}^{0} + \vec{1}_{q_{2}}/n)$. Now, by means of Taylor expansion up to the first order, we have
\begin{align}
p_{G_{n}}(X,Y) - p_{G_{0}}(X,Y) & = \sum \limits_{i=1}^{2} \pi_{i}^{n}\parenth{f(Y|h_{1}(X,\theta_{1i}^{n},\theta_{2i}^{n}) - f(Y|h_{1}(X,\theta_{11}^{0},\theta_{21}^{0})}\overline{f}(X) \nonumber \\
& = \sum \limits_{i=1}^{2} \pi_{i}^{n} \sum \limits_{|\alpha|+|\beta| = 1} \dfrac{1}{\alpha!\beta!}\prod \limits_{u=1}^{q_{1}} \biggr\{(\Delta \theta_{1i}^{n})^{(u)}\biggr\}^{\alpha_{u}}\prod \limits_{v=1}^{q_{2}}\biggr\{(\Delta \theta_{2i}^{n})^{(v)}\biggr\}^{\beta_{v}} \nonumber \\
& \times \dfrac{\partial{f}}{\partial{\theta_{1}^{\alpha}}\partial{\theta_{2}^{\beta}}}\parenth{Y|h_{1}(X,\theta_{11}^{0}),h_{2}(X,\theta_{21}^{0})}\overline{f}(X) + \overline{R}(X,Y), \nonumber
\end{align}
where $\Delta \theta_{1i}^{n} = \theta_{1i}^{n} - \theta_{11}^{0}$ and $\Delta \theta_{2i}^{n} = \theta_{2i}^{n} - \theta_{21}^{0}$ for $1 \leq i \leq 2$. Here $\overline{R}(X,Y)$ is a Taylor remainder from the above expansion. With the choice of $\pi_{i}^{n}, \theta_{1i}^{n}$, and $\theta_{2i}^{n}$ for $1 \leq i \leq 2$, we can verify that:
\begin{align}
\sum \limits_{i=1}^{2} \pi_{i}^{n}\prod \limits_{u=1}^{q_{1}} \biggr\{(\Delta \theta_{1i}^{n})^{(u)}\biggr\}^{\alpha_{u}}\prod \limits_{v=1}^{q_{2}}\biggr\{(\Delta \theta_{2i}^{n})^{(v)}\biggr\}^{\beta_{v}} = 0, \nonumber
\end{align}
for all $|\alpha| + |\beta| = 1$. Therefore, we have the following representation
\begin{align}
p_{G_{n}}(X,Y) - p_{G_{0}}(X,Y) = \overline{R}(X,Y), \nonumber
\end{align}
where the explicit form of the Taylor remainder $\overline{R}(X,Y)$ is as follows:
\begin{align}
\overline{R}(X,Y) & = \sum \limits_{i=1}^{2} \pi_{i}^{n} \sum \limits_{|\alpha| + |\beta| = 2} \dfrac{2}{\alpha!\beta!}\prod \limits_{u=1}^{q_{1}} \biggr\{(\Delta \theta_{1i}^{n})^{(u)}\biggr\}^{\alpha_{u}}\prod \limits_{v=1}^{q_{2}}\biggr\{(\Delta \theta_{2i}^{n})^{(v)}\biggr\}^{\beta_{v}} \nonumber \\
& \times \int \limits_{0}^{1} (1 - t) \dfrac{\partial^{2}{f}}{\partial{\theta_{1}^{\alpha}}\partial{\theta_{2}^{\beta}}}\parenth{Y|h_{1}(X,\theta_{11}^{0} + t\Delta \theta_{1i}^{n}),h_{2}(X,\theta_{21}^{0} + t\Delta \theta_{2i}^{n})}\overline{f}(X) dt. \nonumber
\end{align}
From the properties of a univariate location-scale Gaussian distribution, we can verify that
\begin{align}
T_{\alpha,\beta} = \sup \limits_{t \in [0,1]} \int \frac{\parenth{\dfrac{\partial^{2}{f}}{\partial{\theta_{1}^{\alpha}}\partial{\theta_{2}^{\beta}}}\parenth{Y|h_{1}(X,\theta_{11}^{0} + t\Delta \theta_{1i}^{n}),h_{2}(X,\theta_{21}^{0} + t\Delta \theta_{2i}^{n})}}^{2}}{f(Y|h_{1}(X,\theta_{11}^{0}),h_{2}(X,\theta_{21}^{0})} d(X,Y) < \infty, \label{eqn:finite_Taylor_first}
\end{align}
for all $|\alpha| + |\beta| = 2$. Additionally, the expressions for $\theta_{1i}^{n}$ and $\theta_{2i}^{n}$ indicate that
\begin{align} 
F_{\alpha,\beta} = \sum \limits_{i=1}^{2} \pi_{i}^{n} \frac{2}{\alpha!\beta!}\prod \limits_{u=1}^{q_{1}} \biggr\{(\Delta \theta_{1i}^{n})^{(u)}\biggr\}^{\alpha_{u}}\prod \limits_{v=1}^{q_{2}}\biggr\{(\Delta \theta_{2i}^{n})^{(v)}\biggr\}^{\beta_{v}} = \mathcal{O}(n^{-2}), \label{eqn:finite_Taylor_second}
\end{align}
for all $|\alpha| + |\beta| = 2$. Now a direct computation yields that
\begin{align}
\frac{h^{2}(p_{G_{n}}, p_{G_{0}})}{\widetilde{W}_{\kappa'}^{2\|\kappa'\|_{\infty}}(G_{n},G_{0})} & = \frac{1}{2}\int \frac{\parenth{p_{G_{n}}(X,Y) - p_{G_{0}}(X,Y)}^{2}}{\parenth{\sqrt{p_{G_{n}}(X,Y)} + \sqrt{p_{G_{0}}(X,Y)}}^{2}\widetilde{W}_{\kappa'}^{2\|\kappa'\|_{\infty}}(G_{n},G_{0})}d(X,Y) \nonumber \\
& \leq \frac{1}{2} \int \frac{\overline{R}^{2}(X,Y)}{p_{G_{0}}(X,Y)\widetilde{W}_{\kappa'}^{2\|\kappa'\|_{\infty}}(G_{n},G_{0})}d(X,Y). \nonumber 
\end{align}
The Cauchy-Schwartz inequality implies that the following inequality holds:
\begin{align}
\int \frac{\overline{R}^{2}(X,Y)}{p_{G_{0}}(X,Y)} d(X,Y) \precsim \sum \limits_{|\alpha| + |\beta| = 2} T_{\alpha,\beta}F_{\alpha,\beta}^{2} = \mathcal{O}(n^{-4}), \nonumber
\end{align}
where the final bound comes from the bounds on $T_{\alpha,\beta}$ and $F_{\alpha,\beta}$ in equations~\eqref{eqn:finite_Taylor_first} and~\eqref{eqn:finite_Taylor_second}. On the other hand, the choice of $G_{n}$ guarantees that $\widetilde{W}_{\kappa'}^{2\|\kappa'\|_{\infty}}(G_{n},G_{0}) = \mathcal{O}(n^{-2\kappa'^{(1)}})$ as $\kappa'^{(1)} = \min \limits_{1 \leq i \leq q_{1} + q_{2}} \kappa'^{(i)}$. Since $\kappa'^{(1)} < 2$, it is clear that 
\begin{align}
\int \frac{\overline{R}^{2}(X,Y)}{p_{G_{0}}(X,Y)\widetilde{W}_{\kappa'}^{2\|\kappa'\|_{\infty}}(G_{n},G_{0})}d(X,Y) \to 0, \nonumber
\end{align} 
as $n \to \infty$. Therefore, $h^{2}(p_{G_{n}}, p_{G_{0}})\bigg/\widetilde{W}_{\kappa'}^{2\|\kappa'\|_{\infty}}(G_{n},G_{0}) \to 0$. As a consequence, we obtain the conclusion of equality~\eqref{eqn:general_overfit_second}.
\subsection{Proof of Theorem \ref{theorem:lower_bound_Gaussian_family}}
\label{Section:proof_weakly_identifiable_covariate_independent}
By means of Lemma~\ref{lemma:convergence_rates_MLE}, we prove Theorem~\ref{theorem:lower_bound_Gaussian_family} by establishing the following results:
\begin{eqnarray}
\inf \limits_{G \in \Ocal_{k,\overline{c}_{0}} (\Omega)} h(p_{G},p_{G_{0}})/\widetilde{W}_{\kappa}^{\|\kappa\|_{\infty}}(G,G_{0})>0, \label{eqn:proof_Gaussian_firstcase_first} \\
\inf \limits_{G \in \Ocal_{k} (\Omega) } h(p_{G},p_{G_{0}})/\widetilde{W}_{\kappa'}^{\|\kappa'\|_{\infty}}(G,G_{0}) = 0, \label{eqn:proof_Gaussian_firstcase_second}
\end{eqnarray}
for any $\kappa' \prec \kappa$ where $\kappa = (\overline{r},2,\ceil{\overline{r}/2})$. To simplify the presentation, we assume that $\overline{r}$ is an even number throughout this proof, which leads to $\kappa = (\overline{r}, 2, \overline{r}/2)$. The proof when $\overline{r}$ is an odd number can be obtained in a similar fashion.
\subsubsection{Proof for inequality \eqref{eqn:proof_Gaussian_firstcase_first}} 
\label{subsection:key_inequality_algebraic_dependent_linear_first}
To streamline the argument, we provide a proof only for the local structural inequality:
\begin{eqnarray}
\lim \limits_{\epsilon \to 0} \inf \limits_{G \in \Ocal_{k,\overline{c}_{0}} (\Omega): \widetilde{W}_{\kappa}(G,G_{0}) \leq \epsilon} V(p_{G},p_{G_{0}})/\widetilde{W}_{\kappa}^{\|\kappa\|_{\infty}}(G,G_{0}) > 0; \nonumber
\end{eqnarray}
the global structural result, for  inequality~\eqref{eqn:proof_Gaussian_firstcase_first}, can be argued in a similar fashion as in the proof of Theorem~\ref{theorem:total_variation_bound_over-fitted_MECFG}. Assume now that the local structure inequality does not hold. This implies that we can find a sequence $G_{n} \in \Ocal_{k,\overline{c}_{0}} (\Omega)$ such that $V(p_{G_{n}},p_{G_{0}})/\widetilde{W}_{\kappa}^{\|\kappa\|_{\infty}}(G_{n},G_{0}) \to 0$ and $\widetilde{W}_{\kappa}(G_{n},G_{0}) \to 0$ as $n \to \infty$. Employing the similar argument as in Theorem~\ref{theorem:total_variation_bound_over-fitted_MECFG} in Section~\ref{Section:proof_strong_identifiable}, we can represent the sequence $G_{n}$ as follows:
\begin{eqnarray}
G_{n} = \sum \limits_{i=1}^{k_{0}}\sum \limits_{j=1}^{s_{i}}p_{ij}^{n}\delta_{(\theta_{1ij}^{n},\theta_{2ij}^{n})}, \label{eqn:proof_Gaussian_firstcase_notation}
\end{eqnarray}
where $(\theta_{1ij}^{n},\theta_{2ij}^{n}) \to (\theta_{1i}^{0},\theta_{2i}^{0})$ for all $1 \leq i \leq k_{0}, 1 \leq j \leq s_{i}$ and $\sum \limits_{j=1}^{s_{i}}p_{ij}^{n} \to \pi_{i}^{0}$ for all $1 \leq i \leq k_{0}$. Note that we do not have $\overline{l}$ in the representation of $G_{n} \in \Ocal_{k,\overline{c}_{0}} (\Omega)$, in contrast to the result in Section~\ref{Section:proof_strong_identifiable}. The reason is that the weights of $G_{n}$ are lower bounded by a positive number $\overline{c}_{0}$, which entails that there exists no extra components $(\theta_{1i}^{0},\theta_{2i}^{0})$ as the limit points of the components of $G_{n}$. In this proof, for the simplicity of presentation, we denote $\Delta \theta_{1ij}^{n} = \theta_{1ij}^{n} - \theta_{1i}^{0}$ and $\Delta \theta_{2ij}^{n} = \theta_{2ij}^{n} - \theta_{2i}^{0}$ for all $1 \leq i \leq k_{0}, 1 \leq j \leq s_{i}$. Additionally, $\Delta \theta_{1ij}^{n}=((\Delta \theta_{1ij}^{n})^{(1)}, (\Delta \theta_{1ij}^{n})^{(2)})$ for all $1 \leq i \leq k_{0}, 1 \leq j \leq s_{i}$. Now, according to Lemma~\ref{lemma:generalized_Wasserstein_distance_polynomials} in Appendix~\ref{sec:auxi_result}, we have:
\begin{eqnarray}
\widetilde{W}_{\kappa}^{\|\kappa\|_{\infty}}(G_{n},G_{0}) \precsim \sum \limits_{i=1}^{k_{0}}\sum \limits_{j=1}^{s_{i}}{p_{ij}^{n}\biggr(\biggr|(\Delta \theta_{1ij}^{n})^{(1)}\biggr|^{\overline{r}}+\biggr|(\Delta \theta_{1ij}^{n})^{(2)}\biggr|^{2}+|\Delta \theta_{2ij}^{n}|^{\overline{r}/2}\biggr)} \nonumber \\
+ \sum \limits_{i=1}^{k_{0}} |\sum \limits_{j=1}^{s_{i}}p_{ij}^{n} - \pi_{i}^{0}| := D_{\kappa}(G_{n},G_{0}), \nonumber
\end{eqnarray}
where $\kappa = (\overline{r}, 2, \overline{r}/2)$. Since $V(p_{G_{n}},p_{G_{0}})/\widetilde{W}_{\kappa}^{\|\kappa\|_{\infty}}(G,G_{0}) \to 0$, we have $V(p_{G_{n}},p_{G_{0}})/D_{\kappa}(G_{n},G_{0})$\\$\to 0$. We again divide our proof argument into several steps.
\paragraph{Step 1 - Structure of Taylor expansion:} Using the decomposition $p_{G_{n}}(X,Y)-p_{G_{0}}(X,Y)$, as in the proof of inequality~\eqref{eqn:general_overfit_first} in Section~\ref{Section:proof_strong_identifiable}, we carry out a Taylor expansion up to the order $\overline{r}$:
\begin{eqnarray}
p_{G_{n}}(X,Y)-p_{G_{0}}(X,Y)& = & \sum \limits_{i=1}^{k_{0}}\sum \limits_{j=1}^{s_{i}}p_{ij}^{n}\sum \limits_{1 \leq |\alpha| \leq \overline{r}} \dfrac{1}{\alpha!}\biggr\{(\Delta \theta_{1ij}^{n})^{(1)}\biggr\}^{\alpha_{1}}\biggr\{(\Delta \theta_{1ij}^{n})^{(2)}\biggr\}^{\alpha_{2}}(\Delta \theta_{2ij}^{n})^{\alpha_{3}} \nonumber \\
& & \times \dfrac{\partial^{|\alpha|}{f}}{\partial{(\theta_{1}^{(1)})^{\alpha_{1}}}\partial{(\theta_{1}^{(2)})^{\alpha_{2}}}\partial{\theta_{2}^{\alpha_{3}}}}\parenth{Y|h_{1}(X,\theta_{1i}^{0}),h_{2}(X,\theta_{2i}^{0})}\overline{f}(X) \nonumber \\
& & \hspace{-2 em} + \sum \limits_{i=1}^{k_{0}}\biggr(\sum \limits_{j=1}^{s_{i}}p_{ij}^{n} - \pi_{i}^{0}\biggr)f(Y|h_{1}(X,\theta_{1i}^{0}),h_{2}(X,\theta_{2i}^{0}))\overline{f}(X) + R(X,Y)
\nonumber \\
& &  : = A_{n} + B_{n} + R(X,Y), \label{eqn:proof_Gaussian_firstcase_notation_second}
\end{eqnarray}
where $R(X,Y)$ is a remainder term.  This remainder term is such that $R(X,Y)/D_{\kappa}(G_{n},G_{0}) \to 0$ as $n \to \infty$, is due to the uniform H\"{o}lder continuity of a location-scale Gaussian family with respect to expert functions $h_{1}$, $h_{2}$, and prior density $\overline{f}$ (cf.\ Proposition~\ref{proposition:Lipschitz_continuity}). 

From the formulation of the expert functions $h_{1}$, $h_{2}$ as well as the structural form of the PDE for location-scale Gaussian kernel, we obtain the following:
\begin{eqnarray}
\dfrac{\partial^{|\alpha|}{f}}{\partial{(\theta_{1}^{(1)})^{\alpha_{1}}}\partial{(\theta_{1}^{(2)})^{\alpha_{2}}}\partial{\theta_{2}^{\alpha_{3}}}}(Y|h_{1}(X,\theta_{1}),h_{2}(X,\theta_{2})) = \dfrac{X^{\alpha_{2}}}{2^{\alpha_{3}}}\dfrac{\partial^{\alpha_{1}+\alpha_{2}+2\alpha_{3}}{f}}{\partial{h_{1}^{\alpha_{1}+\alpha_{2}+2\alpha_{3}}}}(Y|h_{1}(X,\theta_{1}),h_{2}(X,\theta_{2})), \nonumber
\end{eqnarray}
for any $\alpha_{1},\alpha_{2},\alpha_{3} \in \mathbb{N}$, $\theta_{1} \in \Omega_{1}$, and $\theta_{2} \in \Omega_{2}$. From this equation, we can rewrite $A_{n}$ as follows:
\begin{eqnarray}
A_{n} & = & \sum \limits_{i=1}^{k_{0}}\sum \limits_{j=1}^{s_{i}}p_{ij}^{n}\sum \limits_{1 \leq |\alpha| \leq \overline{r}} \dfrac{1}{\alpha!}\biggr\{(\Delta \theta_{1ij}^{n})^{(1)}\biggr\}^{\alpha_{1}}\biggr\{(\Delta \theta_{1ij}^{n})^{(2)}\biggr\}^{\alpha_{2}}(\Delta \theta_{2ij}^{n})^{\alpha_{3}} \nonumber \\
& & \times \dfrac{X^{\alpha_{2}}}{2^{\alpha_{3}}}\dfrac{\partial^{\alpha_{1}+\alpha_{2}+2\alpha_{3}}{f}}{\partial{h_{1}^{\alpha_{1}+\alpha_{2}+2\alpha_{3}}}}(Y|h_{1}(X,\theta_{1i}^{0}),h_{2}(X,\theta_{2i}^{0}))\overline{f}(X) \nonumber \\
& = & \sum \limits_{i=1}^{k_{0}}\sum \limits_{j=1}^{s_{i}}p_{ij}^{n}\sum \limits_{\alpha_{2}=0}^{\overline{r}}\sum \limits_{l=0}^{2(\overline{r}-\alpha_{2})}\sum \limits_{\alpha_{1},\alpha_{3}}\dfrac{1}{2^{\alpha_{3}}\alpha!}\biggr\{(\Delta \theta_{1ij}^{n})^{(1)}\biggr\}^{\alpha_{1}}\biggr\{(\Delta \theta_{1ij}^{n})^{(2)}\biggr\}^{\alpha_{2}}(\Delta \theta_{2ij}^{n})^{\alpha_{3}} \nonumber \\
& & \times  X^{\alpha_{2}}\dfrac{\partial^{l+\alpha_{2}}{f}}{\partial{h_{1}^{l+\alpha_{2}}}}(Y|h_{1}(X,\theta_{1i}^{0}),h_{2}(X,\theta_{2i}^{0}))\overline{f}(X), \label{eqn:proof_Gaussian_firstcase_third}
\end{eqnarray}
where $\alpha_{1}, \alpha_{3} \in \mathbb{N}$ in the sum of the second equation satisfies $\alpha_{1}+2\alpha_{3}=l$ and $1 - \alpha_{2} \leq \alpha_{1}+\alpha_{3} \leq \overline{r}-\alpha_{2}$. We define
\begin{eqnarray}
\mathcal{F} : = \biggr\{X^{\alpha_{2}}\dfrac{\partial^{l+\alpha_{2}}{f}}{\partial{h_{1}^{l+\alpha_{2}}}}(Y|h_{1}(X,\theta_{1i}^{0}),h_{2}(X,\theta_{2i}^{0}))\overline{f}(X): \ 0 \leq \alpha_{2} \leq \overline{r}, \ 0 \leq l \leq 2(\overline{r}-\alpha_{2}), \ 1 \leq i \leq k_{0} \biggr\}. \nonumber
\end{eqnarray}
We claim that the elements of $\mathcal{F}$ are linearly independent with respect to $X$ and $Y$. We prove this claim at the end of this proof. Assume that this claim is given at the moment. Inspecting the explicit form of $\mathcal{F}$, we can treat $A_{n}/D_{\kappa}(G_{n},G_{0})$, $B_{n}/D_{\kappa}(G_{n},G_{0})$ as a linear combination of elements of $\mathcal{F}$.
\paragraph{Step 2 - Non-vanishing coefficients:} To simplify the proof, we denote $E_{\alpha_{2},l}(\theta_{1i}^{0},\theta_{2i}^{0})$ as the coefficient of $X^{\alpha_{2}}\dfrac{\partial^{l+\alpha_{2}}{f}}{\partial{h_{1}^{l+\alpha_{2}}}}(Y|h_{1}(X|\theta_{1i}^{0}),h_{2}(X|\theta_{2i}^{0}))\overline{f}(X)$ in $A_{n}$ and $B_{n}$ for any $0 \leq \alpha_{2} \leq \overline{r}$, $0 \leq l \leq 2(\overline{r}-\alpha_{2})$, and $1 \leq i \leq k_{0}$. Then, the coefficients associated with $X^{\alpha_{2}}\dfrac{\partial^{l+\alpha_{2}}{f}}{\partial{h_{1}^{l+\alpha_{2}}}}(Y|h_{1}(X,\theta_{1i}^{0}),h_{2}(X,\theta_{2i}^{0}))\overline{f}(X)$ in $A_{n}/D_{\kappa}(G_{n},G_{0})$ and $B_{n}/D_{\kappa}(G_{n},G_{0})$ take the form $E_{\alpha_{2},l}(\theta_{1i}^{0},\theta_{2i}^{0})/D_{\kappa}(G_{n},G_{0})$. 

Assume that all of the coefficients in the representation of $A_{n}/D_{\kappa}(G_{n},G_{0})$, $B_{n}/D_{\kappa}(G_{n},G_{0})$ go to zero as $n \to \infty$. By taking the summation of $|E_{0,0}(\theta_{1i}^{0},\theta_{2i}^{0})/D_{\kappa}(G_{n},G_{0})|$ for all $1 \leq i \leq k_{0}$, we obtain that
\begin{eqnarray}
\biggr(\sum \limits_{i=1}^{k_{0}} |\sum \limits_{j=1}^{s_{i}}{p_{ij}^{n}} - \pi_{i}^{0}|\biggr)/D_{\kappa}(G_{n},G_{0}) \to 0. \nonumber
\end{eqnarray}
Additionally, according to equation \eqref{eqn:proof_Gaussian_firstcase_third}, we can verify that
\begin{eqnarray}
\dfrac{E_{2,0}(\theta_{1i}^{0},\theta_{2i}^{0})}{D_{\kappa}(G_{n},G_{0})} = \dfrac{\sum \limits_{j=1}^{s_{i}}p_{ij}^{n}\biggr|(\Delta \theta_{1ij}^{n})^{(2)}\biggr|^{2}}{D_{\kappa}(G_{n},G_{0})} \to 0, \nonumber
\end{eqnarray}
for all $1 \leq i \leq k_{0}$. From the formulation of $D_{\kappa}(G_{n},G_{0})$, the above limits lead to
\begin{eqnarray}
\biggr\{\sum \limits_{i=1}^{k_{0}}\sum \limits_{j=1}^{s_{i}}{p_{ij}^{n}\biggr(\biggr|(\Delta \theta_{1ij}^{n})^{(1)}\biggr|^{\overline{r}}+|\Delta \theta_{2ij}^{n}|^{\overline{r}/2}\biggr)}\biggr\}/D_{\kappa}(G_{n},G_{0}) \to 1 \ \text{as} \ n \to \infty. \nonumber
\end{eqnarray}
Therefore, we can find an index $i^{*} \in \left\{1,\ldots,k_{0} \right\}$ such that
\begin{eqnarray}
L = \biggr\{\sum \limits_{j=1}^{s_{i}}{p_{i^{*}j}^{n}\biggr(\biggr|(\Delta \theta_{1i^{*}j}^{n})^{(1)}\biggr|^{\overline{r}}+|\Delta \theta_{2i^{*}j}^{n}|^{\overline{r}/2}\biggr)}\biggr\}/D_{\kappa}(G_{n},G_{0}) \not \to 0, \nonumber
\end{eqnarray}
as $n \to \infty$. Given the definition of $L$, we find that
\begin{align}
	\dfrac{\sum_{i  = 1}^{k_{0}} \sum \limits_{j=1}^{s_{i}}p_{ij}^{n}\biggr|(\Delta \theta_{1ij}^{n})^{(2)}\biggr|^{2}}{\sum \limits_{j=1}^{s_{i}}{p_{i^{*}j}^{n}\biggr(\biggr|(\Delta \theta_{1i^{*}j}^{n})^{(1)}\biggr|^{\overline{r}}+|\Delta \theta_{2i^{*}j}^{n}|^{\overline{r}/2}\biggr)}} =  \dfrac{\sum_{i  = 1}^{k_{0}} \sum \limits_{j=1}^{s_{i}}p_{ij}^{n}\biggr|(\Delta \theta_{1ij}^{n})^{(2)}\biggr|^{2}}{D_{\kappa}(G_{n},G_{0})} \cdot \dfrac{1}{L} \to 0.  \label{eq:explain_limit_first} 
\end{align}
Now, since we have $E_{\alpha_{2},l}(\theta_{1i}^{0},\theta_{2i}^{0})/D_{\kappa}(G_{n},G_{0}) \to 0$ for all values of $\alpha_{2}, l, i$ from the hypothesis, we obtain that
\begin{eqnarray}
M_{\alpha_{2},l}(\theta_{1i}^{0},\theta_{2i}^{0}) = \dfrac{E_{\alpha_{2},l}(\theta_{1i}^{0},\theta_{2i}^{0})}{\sum \limits_{j=1}^{s_{i}}{p_{1j}^{n}\biggr(\biggr|(\Delta \theta_{1i^{*}j}^{n})^{(1)}\biggr|^{\overline{r}}+|\Delta \theta_{2i^{*}j}^{n}|^{\overline{r}/2}\biggr)}} = \dfrac{1}{L} \dfrac{E_{\alpha_{2},l}(\theta_{1i}^{0},\theta_{2i}^{0})}{D_{\kappa}(G_{n},G_{0})} \to 0, \nonumber
\end{eqnarray}
for any $0 \leq \alpha_{2} \leq \overline{r}$, $0 \leq l \leq 2(\overline{r}-\alpha_{2})$, and $1 \leq i \leq k_{0}$. Note that, given the limit~\eqref{eq:explain_limit_first}, when $\alpha_{2} \geq 2$ we find that 
\begin{align*}
 M_{\alpha_{2}, l}(\theta_{1i}^{0}, \theta_{2i}^{0}) \precsim \dfrac{\sum \limits_{j=1}^{s_{i}}p_{ij}^{n}\biggr|(\Delta \theta_{1ij}^{n})^{(2)}\biggr|^{2}}{\sum \limits_{j=1}^{s_{i}}{p_{i^{*}j}^{n}\biggr(\biggr|(\Delta \theta_{1i^{*}j}^{n})^{(1)}\biggr|^{\overline{r}}+|\Delta \theta_{2i^{*}j}^{n}|^{\overline{r}/2}\biggr)}} \to 0.
\end{align*}
Therefore, to obtain a contradiction, it is sufficient to consider only when $\alpha_{2} \leq 1$. When $\alpha_{2} = 1$, direct computation leads to 
\begin{align*}
	M_{1,l}(\theta_{1i}^{0}, \theta_{2i}^{0}) = \dfrac{\sum \limits_{j=1}^{s_{i}}p_{ij}^{n}\sum \limits_{\substack{\alpha_{1}+2\alpha_{3}=l \\ \alpha_{1}+\alpha_{3} \leq \overline{r}}}{\dfrac{\biggr\{(\Delta \theta_{1ij}^{n})^{(1)}\biggr\}^{\alpha_{1}} (\Delta \theta_{1ij}^{n})^{(2)}(\Delta \theta_{2ij}^{n})^{\alpha_{3}}}{2^{\alpha_{3}}\alpha_{1}!\alpha_{3}!}}}{\sum \limits_{j=1}^{s_{i}}{p_{1j}^{n}\biggr(\biggr|(\Delta \theta_{1i^{*}j}^{n})^{(1)}\biggr|^{\overline{r}}+|\Delta \theta_{2i^{*}j}^{n}|^{\overline{r}/2}\biggr)}},
\end{align*}
for any $0 \leq l \leq 2(\bar{r} - 1)$. It is clear that if we have $(\Delta \theta_{1ij}^{n})^{(2)} = 0$ for all $1 \leq i \leq k_{0}$ and $1 \leq j \leq s_{i}$, it is possible that $M_{1, l}(\theta_{1i}^{0}, \theta_{2i}^{0}) \to 0$ for all $1 \leq i \leq k_{0}$. As that existence of the sequence $(\Delta \theta_{1ij}^{n})^{(2)}$ does not violate any of the previous limits, it indicates that when $\alpha_{2} = 1$, it is possible that $M_{\alpha_{2}, l}(\theta_{1i}^{0}, \theta_{2i}^{0})$ go to 0 for all $i$ and $j$. Hence, to obtain a contradiction with the system of limits from $M_{\alpha_{2}, l}(\theta_{1i}^{0}, \theta_{2i}^{0})$, we only need to consider $\alpha_{2} = 0$.

By the representation of $A_{n}$ in equation~\eqref{eqn:proof_Gaussian_firstcase_third}, we can verify that
\begin{eqnarray}
M_{0,l}(\theta_{1i^{*}}^{0},\theta_{2i^{*}}^{0}) = \dfrac{\sum \limits_{j=1}^{s_{i^{*}}}p_{i^{*}j}^{n}\sum \limits_{\substack{\alpha_{1}+2\alpha_{3}=l \\ \alpha_{1}+\alpha_{3} \leq \overline{r}}}{\dfrac{\biggr\{(\Delta \theta_{1i^{*}j}^{n})^{(1)}\biggr\}^{\alpha_{1}}(\Delta \theta_{2i^{*}j}^{n})^{\alpha_{3}}}{2^{\alpha_{3}}\alpha_{1}!\alpha_{3}!}}}{\sum \limits_{j=1}^{s_{i}}{p_{1j}^{n}\biggr(\biggr|(\Delta \theta_{1i^{*}j}^{n})^{(1)}\biggr|^{\overline{r}}+|\Delta \theta_{2i^{*}j}^{n}|^{\overline{r}/2}\biggr)}} \to 0. \nonumber
\end{eqnarray}
\paragraph{Step 3 - Understanding the system of polynomial limits:} The technique for studying the above system of polynomial limits is similar to that of Step 1 in the proof of Proposition 3.3 in~\citep{Ho-Nguyen-SIAM-18}. Here, we briefly sketch the proof for completeness. We denote $\overline{M} = \max \limits_{1 \leq j \leq s_{i^{*}}} \left\{|(\Delta \theta_{1i^{*}j}^{n})^{(1)}|, |\Delta \theta_{2i^{*}j}^{n}|^{1/2}\right\}$ and $\overline{p}=\max \limits_{1 \leq j \leq s_{i^{*}}}\left\{p_{j}\right\}$. Given this notation, let $(\Delta \theta_{1i^{*}j}^{n})^{(1)}/\overline{M} \to a_{j}$, $\Delta \theta_{2i^{*}j}^{n}/\overline{M}^{2} \to b_{j}$, and $p_{i^{*}j}^{n}/\overline{p} \to c_{j}^{2}$ for all $1 \leq j \leq s_{i^{*}}$. Since $p_{i^{*}j}^{n} \geq \overline{c}_{0}$, we will have $c_{j}>0$ for all $1 \leq j \leq s_{i^{*}}$. By dividing both the numerators and the denominators of $M_{0,l}(\theta_{1i^{*}}^{0},\theta_{2i^{*}}^{0})$ by $\overline{M}^{l}$, we obtain the following system of polynomial equations:
\begin{eqnarray}
\sum \limits_{j=1}^{s_{1}}\sum \limits_{\substack{\alpha_{1}+2\alpha_{3}=l \\ \alpha_{1}+\alpha_{3} \leq \overline{r}}}{\dfrac{c_{j}^{2}a_{j}^{\alpha_{1}}b_{j}^{\alpha_{3}}}{2^{\alpha_{3}}\alpha_{1}!\alpha_{3}!}} = 0 \nonumber
\end{eqnarray}
for all $1 \leq l \leq \overline{r}$. Since $s_{1} \leq k-k_{0}+1$ (as $s_{i} \geq 1$ for all $1 \leq i \leq k_{0}$), this system of polynomial equations will not admit any nontrivial solutions $(a_{j},b_{j},c_{j})_{j=1}^{s_{1}}$ according to the definition of $\overline{r}$.  This is a contradiction. As a consequence, not all the coefficients of $A_{n}/D_{\kappa}(G_{n},G_{0})$ and $B_{n}/D_{\kappa}(G_{n},G_{0})$ go to zero as $n \to \infty$. 
\paragraph{Step 4 - Fatou's argument:} Equipped with the above result, we utilize Fatou's argument in Step 3 of the proof of inequality~\eqref{eqn:general_overfit_first} to obtain a contradiction. We denote 
\begin{align}
m_{n} = \max \limits_{0 \leq \alpha_{2} \leq \overline{r}, \ 0 \leq l \leq 2(\overline{r} - \alpha_{2}), \ 1 \leq i \leq k_{0}} \abss{E_{\alpha_{2},l}(\theta_{1i}^{0},\theta_{2i}^{0})}/D_{\kappa}(G_{n},G_{0}); \nonumber
\end{align}
i.e., $m_{n}$ is the maximum of the absolute values of the coefficients in the representation of $A_{n}/D_{\kappa}(G_{n},G_{0})$ and $B_{n}/D_{\kappa}(G_{n},G_{0})$. We now define $E_{\alpha_{2},l}(\theta_{1i}^{0},\theta_{2i}^{0})/ m_{n} \to \tau_{\alpha_{2},l}(i)$ as $n \to \infty$ for all $1 \leq i \leq k_{0}$, $0 \leq \alpha_{2} \leq \overline{r}$, and $0 \leq l \leq 2(\overline{r} - \alpha_{2})$. Here, at least one among $\tau_{\alpha_{2},l}(i)$ is different from zero. Armed with Fatou's lemma as in the proof of inequality~\eqref{eqn:general_overfit_first}, we  obtain the following equation:
\begin{align}
\sum \limits_{i, \alpha_{2}, l} \tau_{\alpha_{2},l}(i)X^{\alpha_{2}}\dfrac{\partial^{l+\alpha_{2}}{f}}{\partial{h_{1}^{l+\alpha_{2}}}}(Y|h_{1}(X,\theta_{1i}^{0}),h_{2}(X,\theta_{2i}^{0}))\overline{f}(X) = 0, \label{eqn:linear_independence_covariate_inde}
\end{align}
almost surely $(X,Y)$ where the ranges of $(i,\alpha_{2},l)$ in the sum satisfy $1 \leq i \leq k_{0}$, $0 \leq \alpha_{2} \leq \overline{r}$, and $0 \leq l \leq 2(\overline{r} - \alpha_{2})$. According to the claim that the elements of $\mathcal{F}$ are linearly independent with respect to $X$ and $Y$,  equation~\eqref{eqn:linear_independence_covariate_inde} indicates that $\tau_{\alpha_{2},l}(i) = 0$ for all $i, \alpha_{2}, l$, which is a contradiction. As a consequence, we prove inequality~\eqref{eqn:proof_Gaussian_firstcase_first}. 
\paragraph{Proof for claim that the elements of $\mathcal{F}$ are linearly independent:} To facilitate the presentation, we reuse the notation from Step 4. In particular, assume that we can find $\tau_{\alpha_{2},l}(i) \in \mathbb{R}$ ($1 \leq i \leq k_{0}$, $0 \leq \alpha_{2} \leq \overline{r}$, and $0 \leq l \leq 2(\overline{r} - \alpha_{2})$) such that equation~\eqref{eqn:linear_independence_covariate_inde} holds
 almost surely $X$ and $Y$. This equation is equivalent to
\begin{eqnarray}
\sum \limits_{i=1}^{k_{0}}\sum \limits_{u=0}^{2\overline{r}}\biggr(\sum \limits_{\alpha_{2} + l=u}\tau_{\alpha_{2}, l}(i) X^{\alpha_{2}}\biggr)\dfrac{\partial^{u}{f}}{\partial{h_{1}^{u}}}(Y|h_{1}(X|\theta_{1i}^{0}),h_{2}(X|\theta_{2i}^{0})) = 0 \label{eqn:lemma_linear_independence_Gaussian_firstcase_first}
\end{eqnarray}
for almost surely $X$ and $Y$. Since $(\theta_{11}^{0},\theta_{21}^{0}), \ldots, (\theta_{1k_{0}}^{0},\theta_{2k_{0}}^{0})$ are $k_{0}$ distinct pairs, we also obtain that $(h_{1}(X|\theta_{11}^{0}),h_{2}(X|\theta_{21}^{0})), \ldots, (h_{1}(X|\theta_{1k_{0}}^{0}),h_{2}(X|\theta_{2k_{0}}^{0}))$ are $k_{0}$ distinct pairs for almost surely $X$. With that result, for $X$ almost surely, we  have that $\dfrac{\partial^{u}{f}}{\partial{h_{1}^{u}}}(Y|h_{1}(X|\theta_{1i}^{0}),h_{2}(X|\theta_{2i}^{0}))$ are linearly independent with respect to $Y$ for $0 \leq u \leq 2\overline{r}$. Therefore, equation \eqref{eqn:lemma_linear_independence_Gaussian_firstcase_first} implies that $\sum \limits_{j+l=u} \tau_{\alpha_{2},l}(i) X^{j}=0$ for all $1 \leq i \leq k_{0}$ and $0 \leq u \leq 2\overline{r}$. As it is a polynomial of $X \in \mathcal{X}$, which is a bounded subset of $\mathbb{R}$, equation~\eqref{eqn:lemma_linear_independence_Gaussian_firstcase_first} only holds when all the coefficients are zero; i.e., $\tau_{\alpha_{2},l}(i) = 0$ for all $\alpha_{2} + l = u$, $1 \leq i \leq k_{0}$ and $0 \leq u \leq 2\overline{r}$. Hence, we establish the claim.
\subsubsection{Proof for equality \eqref{eqn:proof_Gaussian_firstcase_second}} 
\label{subsection:key_equality_algebraic_dependent_linear_first}
In a manner similar to the proof strategy in Theorem \ref{theorem:total_variation_bound_over-fitted_MECFG}, to obtain the conclusion for~\eqref{eqn:proof_Gaussian_firstcase_second}, it is sufficient to construct some sequence $G_{n} \in \Ocal_{k}(\Omega)$ such that
\begin{align}
h(p_{G_{n}},p_{G_{0}})/\widetilde{W}_{\kappa'}^{\|\kappa'\|_{\infty}}(G_{n},G_{0}) \to 0, \nonumber
\end{align}  
for any $\kappa' \prec \kappa = (\overline{r},2, \overline{r}/2)$. The construction for $G_{n}$ will be carried out under two particular settings of $\kappa'$. 
\paragraph{Case 1:} $\kappa' = (\kappa'^{(1)},\kappa'^{(2)},\kappa'^{(3)})$ where $\kappa'^{(2)}<2$. Under this setting, we construct $G_{n} = \sum_{i = 1}^{k_{0}+1} \pi_{i}^{n}\delta_{(\theta_{1i}^{n},\theta_{2i}^{n})}$ such that $(\pi_{i}^{n},\theta_{1i}^{n}, \theta_{2i}^{n}) \equiv (\pi_{i-1}^{0}, \theta_{1(i-1)}^{0}, \theta_{2(i-1)}^{0})$ for $3 \leq i \leq k_{0}+1$. Additionally, $\pi_{1}^{n} = \pi_{2}^{n} = \pi_{1}^{0}/2$, $\parenth{(\theta_{1i}^{n})^{(1)},\theta_{2i}^{n}} = \parenth{(\theta_{11}^{0})^{(1)},\theta_{21}^{0}}$ for $1 \leq i \leq 2$, and $(\theta_{11}^{n})^{(2)} = (\theta_{11}^{0})^{(2)} - 1/n$, $(\theta_{12}^{n})^{(2)} = (\theta_{11}^{0})^{(2)} + 1/n$. From this construction for $G_{n}$, we can verify that $\widetilde{W}_{\kappa'}^{\|\kappa'\|_{\infty}}(G_{n},G_{0}) \asymp n^{-\kappa'^{(2)}}$.  Denote $\Delta \theta_{1i}^{n} = \theta_{1i}^{n} - \theta_{1i}^{0}$ for $1 \leq i \leq 2$. Now, by means of Taylor expansion up to the first order around $(\theta_{11}^{0})^{(2)}$, we have
\begin{align}
p_{G_{n}}(X,Y) - p_{G_{0}}(X,Y) & = \sum \limits_{i=1}^{2} \pi_{i}^{n}\parenth{f(Y|h_{1}(X,\theta_{1i}^{n}),h_{2}(X,\theta_{2i}^{n})) - f(Y|h_{1}(X,\theta_{11}^{0}),h_{2}(X,\theta_{21}^{0}))}\overline{f}(X) \nonumber \\
& = \sum \limits_{i=1}^{2} \pi_{i}^{n}\parenth{\Delta \theta_{1i}^{n}}^{(2)} \dfrac{\partial{f}}{\partial{\theta_{1}^{(2)}}}\parenth{Y|h_{1}(X,\theta_{11}^{0}),h_{2}(X,\theta_{21}^{0})}\overline{f}(X) + \overline{R}_{1}(X,Y), \nonumber
\end{align}
where $\parenth{\Delta \theta_{1i}^{n}}^{(2)} = (\theta_{1i}^{n})^{(2)} - (\theta_{1i}^{0})^{(2)}$ for $1 \leq i \leq 2$ and $\overline{R}_{1}(X,Y)$ is Taylor remainder such that
\begin{align}
\overline{R}_{1}(X,Y) & = \sum \limits_{i=1}^{2} \pi_{i}^{n} \biggr\{\parenth{\Delta \theta_{1i}^{n}}^{(2)}\biggr\}^{2} \int \limits_{0}^{1} (1 - t) \dfrac{\partial^{2}{f}}{\partial{(\theta_{1}^{(2)})^{2}}}(Y|h_{1}(X,\theta_{11}^{0} + t\Delta \theta_{1i}^{n}), \nonumber \\
& \hspace{ 18 em} h_{2}(X,\theta_{21}^{0} + t\Delta \theta_{2i}^{n}))\overline{f}(X) dt. \nonumber
\end{align}
It is not hard to check that $\sum \limits_{i=1}^{2} \pi_{i}^{n}\biggr\{\parenth{\Delta \theta_{1i}^{n}}^{(2)}\biggr\}^{2} = \mathcal{O}(n^{-2})$ and
\begin{align}
\sup \limits_{t \in [0,1]} \int \frac{\parenth{\dfrac{\partial^{2}{f}}{\partial{(\theta_{1}^{(2)})^{2}}}\parenth{Y|h_{1}(X,\theta_{11}^{0} + t\Delta \theta_{1i}^{n}),h_{2}(X,\theta_{21}^{0} + t\Delta \theta_{2i}^{n})}}^{2}}{f(Y|h_{1}(X,\theta_{11}^{0}),h_{2}(X,\theta_{21}^{0})} d(X,Y) < \infty. \nonumber
\end{align}
Therefore, using the same argument as in the proof of equality~\eqref{eqn:general_overfit_second}, the following holds:
\begin{align}
\frac{h^{2}(p_{G_{n}}, p_{G_{0}})}{\widetilde{W}_{\kappa'}^{2\|\kappa'\|_{\infty}}(G_{n},G_{0})} \precsim \int \frac{\overline{R}_{1}^{2}(X,Y)}{p_{G_{0}}(X,Y)\widetilde{W}_{\kappa'}^{2\|\kappa'\|_{\infty}}(G_{n},G_{0})}d(X,Y) \precsim \frac{\mathcal{O}(n^{-4})}{n^{-2\kappa'^{(2)}}} \to 0 \nonumber
\end{align}
as $n \to \infty$. Therefore, we achieve the conclusion of equality~\eqref{eqn:proof_Gaussian_firstcase_second} under Case 1. 
\paragraph{Case 2:} $\kappa' = (\kappa'^{(1)},2,\kappa'^{(3)})$ where $(\kappa'^{(1)},\kappa'^{(3)}) \prec \parenth{\overline{r},\overline{r}/2}$. Under this setting, we construct $G_{n} = \sum_{i=1}^{k}\pi_{i}^{n}\delta_{(\theta_{1i}^{n},\theta_{2i}^{n})}$ such that $(\pi_{i+k-k_{0}}^{n}, \theta_{1(i+k-k_{0})}^{n}, \theta_{2(i+k-k_{0})}^{n}) = (\pi_{i}^{0},\theta_{1i}^{0},\theta_{2i}^{0})$ for $2 \leq i \leq k_{0}$. For $1 \leq j \leq k - k_{0} + 1$, we choose $(\theta_{1j}^{n})^{(2)} = (\theta_{11}^{0})^{(2)}$ and
\begin{align}
(\theta_{1j}^{n})^{(1)} = (\theta_{11}^{0})^{(1)} + \frac{a_{j}^{*}}{n}, \ \theta_{2j}^{n} = \theta_{21}^{0} + \frac{2b_{j}^{*}}{n^{2}}, \ \pi_{j}^{n} = \frac{\pi_{1}^{0}(c_{j}^{*})^{2}}{\sum_{i=1}^{k-k_{0}+1} (c_{j}^{*})^{2}}, \nonumber
\end{align}
where $(c_{i}^{*},a_{i}^{*},b_{i}^{*})_{i=1}^{k-k_{0}+1}$ are the nontrivial solution of the system of polynomial equations~\eqref{eqn:system_polynomial_Gaussian_first} when $r = \overline{r} - 1$. With this formulation of $G_{n}$, it is clear that
\begin{align}
& \hspace{- 3 em} p_{G_{n}}(X,Y) - p_{G_{0}}(X,Y) \nonumber \\
& = \sum \limits_{i=1}^{k - k_{0} +1} \pi_{i}^{n}\parenth{f(Y|h_{1}(X,\theta_{1i}^{n}),h_{2}(X,\theta_{2i}^{n})) - f(Y|h_{1}(X,\theta_{11}^{0}),h_{2}(X,\theta_{21}^{0}))}\overline{f}(X). \nonumber
\end{align}
By means of a Taylor expansion up to the $(\overline{r} - 1$)th order around $\parenth{(\theta_{11}^{0})^{(1)},\theta_{21}^{0}}$, i.e., along the direction of the first component of $\theta_{11}^{0}$ and $\theta_{21}^{0}$, the following equation holds:
\begin{align}
& [f(Y|h_{1}(X,\theta_{1i}^{n}),h_{2}(X,\theta_{2i}^{n})) - f(Y|h_{1}(X,\theta_{11}^{0}),h_{2}(X,\theta_{21}^{0}))] \overline{f}(X) \nonumber \\
& = \sum \limits_{1 \leq |\alpha| \leq \overline{r} - 1} \dfrac{1}{\alpha!}\biggr\{(\Delta \theta_{1i}^{n})^{(1)}\biggr\}^{\alpha_{1}}(\Delta \theta_{2i}^{n})^{\alpha_{2}} \dfrac{\partial^{|\alpha|}{f}}{\partial{(\theta_{1}^{(1)})^{\alpha_{1}}}\partial{\theta_{2}^{\alpha_{2}}}}\parenth{Y|h_{1}(X,\theta_{11}^{0}),h_{2}(X,\theta_{21}^{0})}\overline{f}(X) + \overline{R}_{2i}(X,Y) \nonumber \\
& = \sum \limits_{1 \leq |\alpha| \leq \overline{r} - 1} \dfrac{1}{\alpha!}\biggr\{(\Delta \theta_{1i}^{n})^{(1)}\biggr\}^{\alpha_{1}}(\Delta \theta_{2i}^{n})^{\alpha_{2}}\dfrac{\partial^{\alpha_{1}+2\alpha_{2}}{f}}{\partial{h_{1}^{\alpha_{1}+2\alpha_{2}}}}(Y|h_{1}(X,\theta_{1i}^{0}),h_{2}(X,\theta_{2i}^{0}))\overline{f}(X) + \overline{R}_{2i}(X,Y), \nonumber
\end{align}
where $\alpha = (\alpha_{1}, \alpha_{2})$ in the sum and $\overline{R}_{2i}(X,Y)$ is a remainder. Equipped with this equation, we can rewrite $p_{G_{n}}(X,Y) - p_{G_{0}}(X,Y)$ as 
\begin{align}
p_{G_{n}}(X,Y) - p_{G_{0}}(X,Y) & = \sum \limits_{i = 1}^{k - k_{0} + 1} \pi_{i}^{n} \sum \limits_{1 \leq |\alpha| \leq \overline{r} - 1} \dfrac{1}{\alpha!}\biggr\{(\Delta \theta_{1i}^{n})^{(1)}\biggr\}^{\alpha_{1}}(\Delta \theta_{2i}^{n})^{\alpha_{2}} \nonumber \\
& \times \dfrac{\partial^{\alpha_{1}+2\alpha_{2}}{f}}{\partial{h_{1}^{\alpha_{1}+2\alpha_{2}}}}(Y|h_{1}(X,\theta_{1i}^{0}),h_{2}(X,\theta_{2i}^{0}))\overline{f}(X) + \overline{R}_{2}(X,Y) \nonumber \\
& = \sum \limits_{l = 1}^{2(\overline{r}-1)} \brackets{\sum \limits_{\substack{\alpha_{1}+2\alpha_{2} = l \\ \alpha_{1} + \alpha_{2} \leq \overline{r}-1}} \frac{1}{\alpha !}\sum \limits_{i=1}^{k - k_{0}+1} \pi_{i}^{n}\biggr\{(\Delta \theta_{1i}^{n})^{(1)}\biggr\}^{\alpha_{1}}(\Delta \theta_{2i}^{n})^{\alpha_{2}}} \nonumber \\
& \times \dfrac{\partial^{l}{f}}{\partial{h_{1}^{l}}}(Y|h_{1}(X,\theta_{1i}^{0}),h_{2}(X,\theta_{2i}^{0}))\overline{f}(X) + \overline{R}_{2}(X,Y), \nonumber
\end{align}
where $\overline{R}_{2}(X,Y) = \sum_{i=1}^{k-k_{0}+1} \pi_{i}^{n}\overline{R}_{2i}(X,Y)$ and the range of $\alpha$ in the second equality satisfies $\alpha_{1} + 2\alpha_{2} = l$ and $\alpha_{1} + \alpha_{2} \leq \overline{r}-1$. From the formulations of $\pi_{i}^{n}$, $\theta_{1i}^{n}$, and $\theta_{2i}^{n}$ as $1 \leq i \leq k - k_{0} + 1$, we can check that 
\begin{align}
\sum \limits_{\substack{\alpha_{1}+2\alpha_{2} = l \\ \alpha_{1} + \alpha_{2} \leq \overline{r}-1}} \frac{1}{\alpha !}\sum \limits_{i=1}^{k - k_{0}+1} \pi_{i}^{n}\biggr\{(\Delta \theta_{1i}^{n})^{(1)}\biggr\}^{\alpha_{1}}(\Delta \theta_{2i}^{n})^{\alpha_{2}} = 0, \nonumber
\end{align}
when $1 \leq l \leq \overline{r} - 1$. Additionally, we also have 
\begin{align}
L_{l} = \sum \limits_{\substack{\alpha_{1}+2\alpha_{2} = l \\ \alpha_{1} + \alpha_{2} \leq \overline{r}-1}} \frac{1}{\alpha !}\sum \limits_{i=1}^{k - k_{0}+1} \pi_{i}^{n}\biggr\{(\Delta \theta_{1i}^{n})^{(1)}\biggr\}^{\alpha_{1}}(\Delta \theta_{2i}^{n})^{\alpha_{2}} = \mathcal{O}(n^{-\overline{r}}), \nonumber
\end{align}
when $\overline{r} \leq l \leq 2(\overline{r} - 1)$. Furthermore, the explicit form of $\overline{R}_{2}(X,Y)$ is as follows:
\begin{align}
\overline{R}_{2}(X,Y) & = \sum_{i=1}^{k-k_{0}+1} \pi_{i}^{n}\sum \limits_{|\alpha| = \overline{r}} \dfrac{\overline{r}}{\alpha!}\biggr\{(\Delta \theta_{1i}^{n})^{(1)}\biggr\}^{\alpha_{1}}(\Delta \theta_{2i}^{n})^{\alpha_{2}} \nonumber \\
& \times \int \limits_{0}^{1} (1-t)^{\overline{r}-1} \dfrac{\partial^{\overline{r}}{f}}{\partial{(\theta_{1}^{(1)})^{\alpha_{1}}}\partial{\theta_{2}^{\alpha_{2}}}}\parenth{Y|h_{1}(X,\theta_{11}^{0} + t\Delta \theta_{1i}^{n}),h_{2}(X,\theta_{21}^{0} + t\Delta \theta_{2i}^{n})}\overline{f}(X)dt. \nonumber
\end{align}
It is not hard to check that $\sum \limits_{i=1}^{k-k_{0}+1} \pi_{i}^{n} \biggr\{(\Delta \theta_{1i}^{n})^{(1)}\biggr\}^{\alpha_{1}}(\Delta \theta_{2i}^{n})^{\alpha_{2}} = \mathcal{O}(n^{-\overline{r}})$ and
\begin{align}
\sup \limits_{t \in [0,1]} \int \frac{\parenth{\dfrac{\partial^{\overline{r}}{f}}{\partial{(\theta_{1}^{(1)})^{\alpha_{1}}}\partial{\theta_{2}^{\alpha_{2}}}}\parenth{Y|h_{1}(X,\theta_{11}^{0} + t\Delta \theta_{1i}^{n}),h_{2}(X,\theta_{21}^{0} + t\Delta \theta_{2i}^{n})}}^{2}}{f(Y|h_{1}(X,\theta_{11}^{0}),h_{2}(X,\theta_{21}^{0})} d(X,Y) < \infty, \nonumber
\end{align}
for any $|\alpha| = r$. By the Cauchy-Schwartz inequality, the following inequality holds:
\begin{align}
\frac{h^{2}(p_{G_{n}}, p_{G_{0}})}{\widetilde{W}_{\kappa'}^{2\|\kappa'\|_{\infty}}(G_{n},G_{0})} & \precsim \sum_{l=\overline{r}}^{2(\overline{r}-1)} \frac{L_{l}^{2}}{\widetilde{W}_{\kappa'}^{2\|\kappa'\|_{\infty}}(G_{n},G_{0})}\int \frac{\parenth{\frac{\partial^{l}{f}}{\partial{h_{1}^{l}}}(Y|h_{1}(X,\theta_{1i}^{0}),h_{2}(X,\theta_{2i}^{0}))\overline{f}(X)}^{2}}{p_{G_{0}}(X,Y)}d(X,Y) \nonumber \\
& + \frac{\overline{R}_{2}^{2}(X,Y)}{\widetilde{W}_{\kappa'}^{2\|\kappa'\|_{\infty}}(G_{n},G_{0})}. \nonumber
\end{align}
From a property of the location-scale Gaussian distribution, we have
\begin{align}
\int \frac{\parenth{\frac{\partial^{l}{f}}{\partial{h_{1}^{l}}}(Y|h_{1}(X,\theta_{1i}^{0}),h_{2}(X,\theta_{2i}^{0}))\overline{f}(X)}^{2}}{p_{G_{0}}(X,Y)}d(X,Y) < \infty, \nonumber
\end{align}
for any $\overline{r} \leq l \leq 2(\overline{r}-1)$. 
Furthermore, by means of a similar argument as in the proof of equality~\eqref{eqn:general_overfit_second}, we can argue that 
\begin{align}
\frac{\overline{R}_{2}^{2}(X,Y)}{\widetilde{W}_{\kappa'}^{2\|\kappa'\|_{\infty}}(G_{n},G_{0})} \precsim \frac{\mathcal{O}(n^{-2\overline{r}})}{n^{-2\min\{\kappa'^{(1)},\kappa'^{(3)}\}}}. \nonumber
\end{align}
Putting these results together, we have
\begin{align}
\frac{h^{2}(p_{G_{n}}, p_{G_{0}})}{\widetilde{W}_{\kappa'}^{2\|\kappa'\|_{\infty}}(G_{n},G_{0})} \precsim \frac{\mathcal{O}(n^{-2\overline{r}})}{n^{-2\min\{\kappa'^{(1)},\kappa'^{(3)}\}}} \to 0, 
\end{align}
as $n \to \infty$. Therefore, we obtain the conclusion of equality~\eqref{eqn:proof_Gaussian_firstcase_second} under Case 2.
\section{Discussion} \label{Section:discussion}
We have provided a systematic theoretical understanding of the convergence rates of parameter estimation under over-specified Gaussian mixtures of experts based on an analysis of an underlying algebraic structure. In particular, we have introduced a new theoretical tool, which we refer to as algebraic independence, and we have established a connection between this algebraic structure and a certain family of PDEs.  This connection allows us to determine convergence rates of the MLE under various choices of expert functions $h_{1}$ and $h_{2}$.

There are several directions for future research. First, the current convergence rates of the MLE are established under the assumptions that the parameter spaces are bounded; it would be important to remove this assumption for wider practical applicability. Second, the results of the paper demonstrate that the convergence rates of MLE are only very slow when the expert functions are algebraically dependent. When we indeed fit the models with algebraically independent expert functions while the true expert functions are algebraically dependent, i.e., we misspecify the expert functions, the convergence rates of MLE become $n^{- 1/ 4}$. However, the MLE will not converge to the true mixing measure. This raises an interesting challenge of how to characterize the difference between the limiting mixing measure and the true mixing measure in terms of the generalized transportation distance. Finally, since the log-likelihood function of over-specified Gaussian mixtures of experts is nonconcave, the MLE does not have a closed form in practice. Therefore, heuristic optimization algorithms, such as Expectation-Maximization (EM) algorithm, are generally used to approximate MLE. The convergence rates of EM algorithm and other optimization algorithms in standard mixture models had been studied in~\citep{dwivedi2020, Raaz_Ho_Koulik_2018_second}. Recently,~\cite{Kwon_minimax_2021} established the minimax convergence rates of EM algorithm and ~\cite{Tongzheng_2022} studied the convergence rate of Polyak step size gradient descent algorithm for symmetric two-component Gaussian mixed linear regression, which is a special case of Gaussian mixture of experts. It is of practical importance to investigate the computational errors arising from the updates of the optimization algorithms, such as EM algorithm, on the convergence rates of MLE under general Gaussian mixtures of experts.

\section{Acknowledgements} \label{Section:acknowledgement}
This work was supported in part by the Mathematical Data Science program of the Office of Naval Research under grant number N00014-18-1-2764.

\appendix 
\section{Appendix} \label{Section:Appendix_A}
In this appendix, we provide proofs for remaining results in the paper. 
\subsection{Proof of Proposition~\ref{proposition:identifiability_fixed_weights_setting}}
\label{Section:proof_identifiable}
Assume that there exist $G = \sum_{i = 1}^{k} \pi_{i} \delta_{(\theta_{1i}, \theta_{2i})}$ and $G' = \sum_{i = 1}^{k'} \pi_{i}' \delta_{(\theta_{1i}', \theta_{2i}')}$ such that $p_{G}(X, Y) = p_{G'}(X, Y)$ for almost surely $(X, Y) \in \mathcal{X} \times \mathcal{Y}$. It is equivalent to 
\begin{align}
    \sum_{i = 1}^{k} \pi_{i} f(Y|h_{1}(X, \theta_{1i}), h_{2}(X, \theta_{2i})) = \sum_{i = 1}^{k'} \pi_{i}' f(Y|h_{1}(X, \theta_{1i}'), h_{2}(X, \theta_{2i}')) \label{eq:identifiability_first_equation}
\end{align}
for almost surely $(X, Y)$. Due to the identifiability of location-scale Gaussian mixtures~\cite{Teicher-1960, Teicher-1961}, we have $k = k'$ and $\{\pi_{1}, \pi_{2}, \ldots, \pi_{k}\} \equiv \{\pi_{1}', \pi_{2}', \ldots, \pi_{k}'\}$. Without loss of generality, we assume that $\pi_{i} = \pi_{i}'$ for all $1 \leq i \leq k$. We denote by $J_{1}, \ldots, J_{l}$ the partition of $\{1, 2,\ldots, k\}$ for some $l \leq k$ such that $\pi_{i} = \pi_{i}'$ for any $i, i' \in J_{j}$ and $1 \leq j \leq l$. Furthermore, $\pi_{i} \neq \pi_{i}'$ when $i$ and $i'$ do not belong to the same set $J_{j}$ for any $1 \leq j \leq l$. Therefore, we can rewrite equation~\eqref{eq:identifiability_first_equation} as follows:
\begin{align*}
    \sum_{j = 1}^{l} \sum_{i \in J_{j}} \pi_{i} f(Y|h_{1}(X,\theta_{1i}), h_{2}(X, \theta_{2i})) = \sum_{j = 1}^{l} \sum_{i \in J_{j}} \pi_{i}' f(Y|h_{1}(X,\theta_{1i}'), h_{2}(X, \theta_{2i}')).
\end{align*}
From these results, for almost surely $X$, for each $1 \leq j \leq l$ there exist permutation functions $\sigma_{X}^{j}: J_{j} \to J_{j}$ such that $(h_{1}(X,\theta_{1\sigma_{X}^{j}(i)}), h_{2}(X,\theta_{2\sigma_{X}^{j}(i)})) \equiv (h_{1}(X,\theta_{1i}'), h_{2}(X,\theta_{2i}'))$ for all $i \in J_{j}$. Since the expert functions $h_{1}$ and $h_{2}$ are identifiable, from Definition~\ref{definition:identifiability} these equations indicate that $\{(\theta_{1i}, \theta_{2i}): \ i \in J_{j} \} \equiv \{(\theta_{1i}', \theta_{2i}'): \ i \in J_{j} \}$ for all $1 \leq j \leq l$. As a consequence,
\begin{align*}
    G = \sum_{j = 1}^{l} \sum_{i \in J_{j}} \pi_{i} \delta_{(\theta_{1i}, \theta_{2i})} = \sum_{j = 1}^{l} \sum_{i \in J_{j}} \pi_{i}' \delta_{(\theta_{1i}', \theta_{2i}')} = G'.
\end{align*}
We obtain the conclusion of the proposition.
\subsection{Proof of Lemma~\ref{lemma:convergence_rates_MLE}}
\label{subsec:proof:lemma:convergence_rates_MLE}
The proof of part (a) of Lemma~\ref{lemma:convergence_rates_MLE} is straightforward from the parametric convergence rate of $h(p_{\widehat{G}_{n}},p_{G_{0}})$ established in Proposition~\ref{proposition:convergence_rates_density_estimation_MECFG}. We therefore omit the proof of part (a) of Lemma~\ref{lemma:convergence_rates_MLE} for the brevity of presentation.

We now provide proof of part (b) of Lemma~\ref{lemma:convergence_rates_MLE}. It follows the same
argument as that of Lemma 1 in~\citep{Yu-97}. Fix $(1, \ldots, 1) \preceq \kappa' \precsim \kappa$ and the true mixing measure $G_{0}$. Let $C_{0} > 0$ be any fixed constant. From the hypothesis of part (b), for any sufficiently small $\varepsilon > 0$, we can find $G_{0}' \in \mathcal{G}$ such that $\widetilde{W}_{\kappa'}(G_{0}', G_{0}) = 2 \varepsilon$ and $h(p_{G_{0}'}, p_{G_{0}}) \leq C_{0} \varepsilon^{\|\kappa'\|_{\infty}}$. Now, by taking any sequence of estimates $\bar{G}_{n} \in \mathcal{G}$, we obtain that
\begin{align*}
    2 \max_{G \in \{G_{0}, G_{0}'\}} \Exs_{p_{G}} \brackets{\widetilde{W}_{\kappa'}(\bar{G}_{n}, G)} \geq \Exs_{p_{G_{0}}} \brackets{\widetilde{W}_{\kappa'}(\bar{G}_{n}, G_{0})} + \Exs_{p_{G_{0}'}} \brackets{\widetilde{W}_{\kappa'}(\bar{G}_{n}, G_{0}')}.
\end{align*}
Here, $\Exs_{p_{G}}$ denotes the expectation taken with respect to the product measure with mixture density $p_{G}^{n}$. Since $\widetilde{W}_{\kappa'}$ satisfies the weak triangle inequality, we have a positive constant $C_{1}$ depending on $\kappa'$ such that
\begin{align*}
    \widetilde{W}_{\kappa'}(\bar{G}_{n}, G_{0}) + \widetilde{W}_{\kappa'}(\bar{G}_{n}, G_{0}') \geq C_{1} \widetilde{W}_{\kappa'}(G_{0}, G_{0}') = 2 C_{1} \varepsilon.
\end{align*}
Therefore, we find that
\begin{align*}
    \Exs_{p_{G_{0}}} \brackets{\widetilde{W}_{\kappa'}(\bar{G}_{n}, G_{0})} + \Exs_{p_{G_{0}'}} \brackets{\widetilde{W}_{\kappa'}(\bar{G}_{n}, G_{0}')} \geq 2 C_{1} \varepsilon \inf_{f_{1}, f_{2}} \parenth{ \Exs_{p_{G_{0}} } \brackets{f_{1}} + \Exs_{p_{G_{0}}'} \brackets{f_{2}}},
\end{align*}
where the infimum is taken over non-negative measurable functions $f_{1}$ and $f_{2}$ defined in terms of $X_{1}, \ldots, X_{n}$ such that $f_{1} + f_{2} = 1$. The definition of total variation distance indicates that we can rewrite the right-hand-side of the above inequality as $1 - V(p_{G_{0}}^{n}, p_{G_{0}'}^{n})$. Putting the above results together, we find that
\begin{align*}
    \max_{G \in \{G_{0}, G_{0}'\}} \Exs_{p_{G}} \brackets{\widetilde{W}_{\kappa'}(\bar{G}_{n}, G)} \geq C_{1} \varepsilon \parenth{1 - V(p_{G_{0}}^{n}, p_{G_{0}'}^{n})}. 
\end{align*}
Since $V \leq h$, we obtain that
\begin{align*}
    V(p_{G_{0}}^{n}, p_{G_{0}'}^{n}) \leq h(p_{G_{0}}^{n}, p_{G_{0}'}^{n}) = \sqrt{1 - \parenth{1 - h^{2}(p_{G_{0}}, p_{G_{0}'})}^{n}} \leq \sqrt{1 - \parenth{1 - C_{0}^2 \varepsilon^{2 \|\kappa'\|_{\infty}}}^{n}}.
\end{align*}
Therefore, we arrive at
\begin{align*}
    \max_{G \in \{G_{0}, G_{0}'\}} \Exs_{p_{G}} \brackets{\widetilde{W}_{\kappa'}(\bar{G}_{n}, G)} \geq C_{1} \varepsilon \parenth{1 - \sqrt{1 - \parenth{1 - C_{0}^2 \varepsilon^{2 \|\kappa'\|_{\infty}}}^{n}}}. 
\end{align*}
By choosing $C_{0}^2 \varepsilon^{2 \|\kappa'\|_{\infty}} = \frac{1}{n}$, the bound in the above display becomes 
\begin{align*}
\max_{G \in \{G_{0}, G_{0}'\}} \Exs_{p_{G}} \brackets{\widetilde{W}_{\kappa'}(\bar{G}_{n}, G)} \geq c_{1} n^{-1/2\|\kappa'\|_{\infty}}.
\end{align*}
As $\sup \limits_{G \in \mathcal{G} \backslash \Ocal_{k_{0}-1}(\Omega)} \Exs_{p_{G}} \parenth{\widetilde{W}_{\kappa'}(\overline{G}_{n},G)} \geq \max_{G \in \{G_{0}, G_{0}'\}} \Exs_{p_{G}} \brackets{\widetilde{W}_{\kappa'}(\bar{G}_{n}, G)}$, we reach the conclusion of part (b) of the lemma. 
\subsection{Proof of Theorem \ref{theorem:lower_bound_Gaussian_family_first_first}}
\label{Section:proof_weakly_identifiable_covariate_dependent_first}
Similar to previous proofs in Section~\ref{Section:proofs}, it is sufficient to demonstrate the following results:
\begin{eqnarray}
\lim \limits_{\epsilon \to 0} \inf \limits_{G \in \Ocal_{k,\overline{c}_{0}} (\Omega): \widetilde{W}_{\kappa}(G,G_{0}) \leq \epsilon} V(p_{G},p_{G_{0}})/\widetilde{W}_{\kappa}^{\|\kappa\|_{\infty}}(G,G_{0}) & > 0, \label{eqn:proof_Gaussian_secondcase_first} \\
\inf \limits_{G \in \Ocal_{k} (\Omega)} h(p_{G},p_{G_{0}})/\widetilde{W}_{\kappa'}^{\|\kappa'\|_{\infty}}(G,G_{0}) & = 0, \label{eqn:proof_Gaussian_secondcase_second}
\end{eqnarray}
for any $\kappa' \prec \kappa$ where $\kappa=(2,\overline{r},\ceil{\overline{r}/2})$. Without loss of generality, we assume that $\overline{r}$ is an even number throughout this proof, i.e., $\kappa = (2, \overline{r}, \overline{r}/2)$. Proof of inequality~\eqref{eqn:proof_Gaussian_secondcase_first} is in Appendix~\ref{subsection:key_inequality_algebraic_dependent_linear_second} while proof of equality~\eqref{eqn:proof_Gaussian_secondcase_second} is in Appendix~\ref{subsec:proof_Gaussian_secondcase_second}.
\subsubsection{Proof for inequality \eqref{eqn:proof_Gaussian_secondcase_first}} 
\label{subsection:key_inequality_algebraic_dependent_linear_second}
Assume the inequality~\eqref{eqn:proof_Gaussian_secondcase_first} does not hold. It indicates that there exists a sequence $G_{n} \in \Ocal_{k,\overline{c}_{0}} (\Omega)$ such that $V(p_{G},p_{G_{0}})/\widetilde{W}_{\kappa}^{\|\kappa\|_{\infty}}(G_{n},G_{0}) \to 0$ and $\widetilde{W}_{\kappa}^{\|\kappa\|_{\infty}}(G_{n},G_{0}) \to 0$. To simplify the presentation, we reuse the notation of $G_{n}$ as in equation~\eqref{eqn:proof_Gaussian_firstcase_notation} in the proof of Theorem~\ref{theorem:lower_bound_Gaussian_family} in Section~\ref{Section:proof_weakly_identifiable_covariate_independent}. Since $\kappa = (2,\overline{r},\overline{r}/2)$, we have
\begin{eqnarray}
\widetilde{W}_{\kappa}^{\|\kappa\|_{\infty}}(G_{n},G_{0}) \precsim \sum \limits_{i=1}^{k_{0}}\sum \limits_{j=1}^{s_{i}}{p_{ij}^{n}\biggr(\biggr|(\Delta \theta_{1ij}^{n})^{(1)}\biggr|^{2}+\biggr|(\Delta \theta_{1ij}^{n})^{(2)}\biggr|^{\overline{r}}+|\Delta \theta_{2ij}^{n}|^{\overline{r}/2}\biggr)} \nonumber \\
+ \sum \limits_{i=1}^{k_{0}} |\sum \limits_{j=1}^{s_{i}}p_{ij}^{n} - \pi_{i}^{0}| := D_{\kappa}(G_{n},G_{0}). \nonumber
\end{eqnarray}
Similar to the proof of Theorem \ref{theorem:lower_bound_Gaussian_family}, by means of Taylor expansion up to the $\overline{r}$ order, we can represent
\begin{eqnarray}
p_{G_{n}}(X,Y)-p_{G_{0}}(X,Y) = A_{n} + B_{n} + R(X,Y), \nonumber
\end{eqnarray}
where $A_{n}$, $B_{n}$, and $R(X,Y)$ are identifical to those in equation~\eqref{eqn:proof_Gaussian_firstcase_notation_second} such that $R(X,Y)/D_{\kappa}(G_{n},G_{0})$\\$\to 0$ as $n \to \infty$. Given the formulation of expert functions $h_{1}, h_{2}$, we have the following key equation:
\begin{eqnarray}
\dfrac{\partial^{|\alpha|}{f}}{\partial{(\theta_{1}^{(1)})^{\alpha_{1}}}\partial{(\theta_{1}^{(2)})^{\alpha_{2}}}\partial{\theta_{2}^{\alpha_{3}}}}(Y|h_{1}(X|\theta_{1}),h_{2}(X|\theta_{2})) = \dfrac{X^{\alpha_{2}+2\alpha_{3}}}{2^{\alpha_{3}}}\dfrac{\partial^{\alpha_{1}+\alpha_{2}+2\alpha_{3}}{f}}{\partial{h_{1}^{\alpha_{1}+\alpha_{2}+2\alpha_{3}}}}(Y|h_{1}(X|\theta_{1}),h_{2}(X|\theta_{2})), \nonumber
\end{eqnarray}
for any $\alpha_{1}, \alpha_{2}, \alpha_{3} \in \mathbb{N}$. Equipped with the above equation, $A_{n}$ can be rewritten as
\begin{eqnarray}
A_{n} & = & \sum \limits_{i=1}^{k_{0}}\sum \limits_{j=1}^{s_{i}} p_{ij}^{n}\sum \limits_{\alpha_{1}=0}^{\overline{r}}\sum \limits_{l=0}^{2(\overline{r}-\alpha_{1})}\sum \limits_{\alpha_{1},\alpha_{3}}\dfrac{1}{2^{\alpha_{3}}\alpha!}\biggr\{(\Delta \theta_{1ij}^{n})^{(1)}\biggr\}^{\alpha_{1}}\biggr\{(\Delta \theta_{1ij}^{n})^{(2)}\biggr\}^{\alpha_{2}}(\Delta \theta_{2ij}^{n})^{\alpha_{3}} \nonumber \\
& & \hspace{12 em} \times  X^{l}\dfrac{\partial^{\alpha_{1}+l}{f}}{\partial{h_{1}^{\alpha_{1}+l}}}(Y|h_{1}(X|\theta_{1i}^{0}),h_{2}(X|\theta_{2i}^{0}))\overline{f}(X), \label{eqn:proof_Gaussian_secondcase_third}
\end{eqnarray}
where $\alpha_{2}, \alpha_{3} \in \mathbb{N}$ in the above sum satisfies $\alpha_{2}+2\alpha_{3}=l$ and $1-\alpha_{1} \leq \alpha_{2}+\alpha_{3} \leq \overline{r}-\alpha_{1}$. If we define
\begin{eqnarray}
\mathcal{F} = \biggr\{X^{l}\dfrac{\partial^{\alpha_{1}+l}{f}}{\partial{h_{1}^{\alpha_{1}+l}}}(Y|h_{1}(X|\theta_{1i}^{0}),h_{2}(X|\theta_{2i}^{0}))\overline{f}(X): \ 0 \leq \alpha_{1} \leq \overline{r}, \ 0 \leq l \leq 2(\overline{r}-\alpha_{1}), \ 1 \leq i \leq k_{0} \biggr\}, \nonumber
\end{eqnarray}
then the elements of $\mathcal{F}$ are linearly independent with respect to $X$ and $Y$. The proof argument of this claim is similar to that in equation~\eqref{eqn:lemma_linear_independence_Gaussian_firstcase_first} in Section~\ref{Section:proof_weakly_identifiable_covariate_independent}. Therefore, we can treat $A_{n}/D_{\kappa}(G_{n},G_{0})$, $B_{n}/D_{\kappa}(G_{n},G_{0})$ as a linear combination of elements of $\mathcal{F}$. 

Similar to the proof of Theorem~\ref{theorem:lower_bound_Gaussian_family} in Section~\ref{Section:proof_weakly_identifiable_covariate_independent}, we denote $F_{\alpha_{1},l}(\theta_{1i}^{0},\theta_{2i}^{0})$ as the coefficient of $X^{l}\dfrac{\partial^{\alpha_{1}+l}{f}}{\partial{h_{1}^{\alpha_{1}+l}}}(Y|h_{1}(X|\theta_{1i}^{0}),h_{2}(X|\theta_{2i}^{0}))\overline{f}(X)$ in $A_{n}$ and $B_{n}$ for any $0 \leq \alpha_{1} \leq \overline{r}$, $0 \leq l \leq 2(\overline{r}-\alpha_{1})$, and $1 \leq i \leq k_{0}$. Then, the coefficients of $X^{l}\dfrac{\partial^{\alpha_{1}+l}{f}}{\partial{h_{1}^{\alpha_{1}+l}}}(Y|h_{1}(X|\theta_{1i}^{0}),h_{2}(X|\theta_{2i}^{0}))\overline{f}(X)$ in $A_{n}/D_{\kappa}(G_{n},G_{0})$ and $B_{n}/D_{\kappa}(G_{n},G_{0})$ will be $F_{\alpha_{1},l}(\theta_{1i}^{0},\theta_{2i}^{0})/D_{\kappa}(G_{n},G_{0})$. 

Assume that all of these coefficients go to 0 as $n \to \infty$. By taking the summation of $|F_{0,0}(\theta_{1i}^{0},\theta_{2i}^{0})/D_{\kappa}(G_{n},G_{0})|$ for all $1 \leq i \leq k_{0}$, we obtain that
\begin{eqnarray}
\biggr(\sum \limits_{i=1}^{k_{0}} |\sum \limits_{j=1}^{s_{i}}{p_{ij}^{n}} - \pi_{i}^{0}|\biggr)/D_{\kappa}(G_{n},G_{0}) \to 0. \nonumber
\end{eqnarray}
Additionally, according to equation \eqref{eqn:proof_Gaussian_secondcase_third}, we can verify that
\begin{eqnarray}
F_{2,0}(\theta_{1i}^{0},\theta_{2i}^{0})/D_{\kappa}(G_{n},G_{0}) = \biggr(\sum \limits_{j=1}^{s_{i}}p_{ij}^{n}\biggr|(\Delta \theta_{1ij}^{n})^{(1)}\biggr|^{2}\biggr)/D_{\kappa}(G_{n},G_{0}) \to 0 \nonumber
\end{eqnarray}
for all $1 \leq i \leq k_{0}$. From the formulation of $D_{\kappa}(G_{n},G_{0})$, the above limits lead to
\begin{eqnarray}
\biggr\{\sum \limits_{i=1}^{k_{0}}\sum \limits_{j=1}^{s_{i}}{p_{ij}^{n}\biggr(\biggr|(\Delta \theta_{1ij}^{n})^{(2)}\biggr|^{\overline{r}}+|\Delta \theta_{2ij}^{n}|^{\overline{r}/2}\biggr)}\biggr\}/D_{\kappa}(G_{n},G_{0}) \to 1 \ \text{as} \ n \to \infty. \nonumber
\end{eqnarray}
Thus, we can find an index $i^{*} \in \left\{1,\ldots,k_{0} \right\}$ such that
\begin{eqnarray}
J = \biggr\{\sum \limits_{j=1}^{s_{i}}{p_{i^{*}j}^{n}\biggr(\biggr|(\Delta \theta_{1i^{*}j}^{n})^{(2)}\biggr|^{\overline{r}}+|\Delta \theta_{2i^{*}j}^{n}|^{\overline{r}/2}\biggr)}\biggr\}/D_{\kappa}(G_{n},G_{0}) \not \to 0 \nonumber
\end{eqnarray}
as $n \to \infty$. Without loss of generality, we assume that $i^{*}=1$. Now, since we have $F_{\alpha_{1},l}(\theta_{1i}^{0},\theta_{2i}^{0})/D_{\kappa}(G_{n},G_{0}) \to 0$ for all values of $\alpha_{1}, l, i$, we obtain that
\begin{eqnarray}
M_{\alpha_{1},l}(\theta_{1i}^{0},\theta_{2i}^{0}) = \dfrac{F_{\alpha_{1},l}(\theta_{1i}^{0},\theta_{2i}^{0})}{\sum \limits_{j=1}^{s_{i}}{p_{1j}^{n}\biggr(\biggr|(\Delta \theta_{11j}^{n})^{(2)}\biggr|^{\overline{r}}+|\Delta \theta_{21j}^{n}|^{\overline{r}/2}\biggr)}} = \dfrac{1}{J} \dfrac{F_{\alpha_{1},l}(\theta_{1i}^{0},\theta_{2i}^{0})}{D_{\kappa}(G_{n},G_{0})} \to 0, \nonumber
\end{eqnarray}
for any $0 \leq \alpha_{1} \leq \overline{r}$, $0 \leq l \leq 2(\overline{r}-\alpha_{1})$, and $1 \leq i \leq k_{0}$. From the representation of $A_{n}$ in equation~\eqref{eqn:proof_Gaussian_secondcase_third}, we can verify that
\begin{eqnarray}
M_{0,l}(\theta_{11}^{0},\theta_{21}^{0}) = \dfrac{\sum \limits_{j=1}^{s_{1}}p_{1j}^{n}\sum \limits_{\substack{\alpha_{2}+2\alpha_{3}=l \\ \alpha_{2}+\alpha_{3} \leq \overline{r}}}{\dfrac{\biggr\{(\Delta \theta_{11j}^{n})^{(2)}\biggr\}^{\alpha_{2}}(\Delta \theta_{21j}^{n})^{\alpha_{3}}}{2^{\alpha_{3}}\alpha_{2}!\alpha_{3}!}}}{\sum \limits_{j=1}^{s_{i}}{p_{1j}^{n}\biggr(\biggr|(\Delta \theta_{11j}^{n})^{(2)}\biggr|^{\overline{r}}+|\Delta \theta_{21j}^{n}|^{\overline{r}/2}\biggr)}} \to 0. \nonumber
\end{eqnarray}
for $0 \leq l \leq 2\overline{r}$. Using the same argument as that in Step 3 of the proof of Theorem~\ref{theorem:lower_bound_Gaussian_family} in Section~\ref{Section:proof_weakly_identifiable_covariate_independent}, the above system of polynomial limits does not hold. As a consequence, not all the coefficients in the linear combinations of $A_{n}/D_{\kappa}(G_{n},G_{0})$ and $B_{n}/D_{\kappa}(G_{n},G_{0})$ go to 0 as $n \to \infty$. From here, using the Fatou's argument in Step 4 of the proof of Theorem~\ref{theorem:lower_bound_Gaussian_family} and the fact that the elements of $\mathcal{F}$ are linearly independent with respect to $X$ and $Y$, we achieve the conclusion of claim~\eqref{eqn:proof_Gaussian_secondcase_first}.
\subsubsection{Proof for equality~\eqref{eqn:proof_Gaussian_secondcase_second}} 
\label{subsec:proof_Gaussian_secondcase_second}
To alleviate the presentation, we will only provide a proof sketch of equality~\eqref{eqn:proof_Gaussian_secondcase_second}. We also divide the proof into two settings of $\kappa' \prec \kappa = (2,\overline{r}, \overline{r}/2)$.
\paragraph{Case 1:} $\kappa' = \parenth{\kappa'^{(1)}, \kappa'^{(2)}, \kappa'^{(3)}}$ when $\kappa'^{(1)} < 2$. Under this setting, we construct $G_{n} = \sum_{i = 1}^{k_{0}+1} \pi_{i}^{n}\delta_{(\theta_{1i}^{n},\theta_{2i}^{n})}$ such that $(\pi_{i}^{n},\theta_{1i}^{n}, \theta_{2i}^{n}) \equiv (\pi_{i-1}^{0}, \theta_{1(i-1)}^{0}, \theta_{2(i-1)}^{0})$ for $3 \leq i \leq k_{0}+1$. Additionally, $\pi_{1}^{n} = \pi_{2}^{n} = \pi_{1}^{0}/2$, $\parenth{(\theta_{1i}^{n})^{(2)},\theta_{2i}^{n}} = \parenth{(\theta_{11}^{0})^{(2)},\theta_{21}^{0}}$ for $1 \leq i \leq 2$, and $(\theta_{11}^{n})^{(1)} = (\theta_{11}^{0})^{(1)} - 1/n$, $(\theta_{12}^{n})^{(1)} = (\theta_{11}^{0})^{(1)} + 1/n$. From this construction of $G_{n}$, we can verify that $\widetilde{W}_{\kappa'}^{\|\kappa'\|_{\infty}}(G_{n},G_{0}) \asymp n^{-\kappa'^{(1)}}$. Given that formulation of $G_{n}$, when we perform Taylor expansion up to the first order around $(\theta_{11}^{0})^{(1)}$, the following equation holds
\begin{align}
p_{G_{n}}(X,Y) - p_{G_{0}}(X,Y) = \overline{R}_{1}(X,Y), \nonumber
\end{align}
where $\overline{R}_{1}(X,Y)$ is Taylor remainder such that
\begin{align}
& \overline{R}_{1}(X,Y) \nonumber \\
& = \sum \limits_{i=1}^{2} \pi_{i}^{n} \biggr\{\parenth{\Delta \theta_{1i}^{n}}^{(1)}\biggr\}^{2} \int \limits_{0}^{1} (1 - t) \dfrac{\partial^{2}{f}}{\partial{(\theta_{1}^{(1)})^{2}}}\parenth{Y|h_{1}(X,\theta_{11}^{0} + t\Delta \theta_{1i}^{n}),h_{2}(X,\theta_{21}^{0} + t\Delta \theta_{2i}^{n})}\overline{f}(X) dt. \nonumber
\end{align}
Using the same argument as that in Case 1 in the proof of equality~\eqref{eqn:proof_Gaussian_firstcase_second} in Section~\ref{subsection:key_equality_algebraic_dependent_linear_first}, the following holds
\begin{align}
\frac{h^{2}(p_{G_{n}}, p_{G_{0}})}{\widetilde{W}_{\kappa'}^{2\|\kappa'\|_{\infty}}(G_{n},G_{0})} \precsim \int \frac{\overline{R}_{1}^{2}(X,Y)}{p_{G_{0}}(X,Y)\widetilde{W}_{\kappa'}^{2\|\kappa'\|_{\infty}}(G_{n},G_{0})}d(X,Y) \precsim \frac{\mathcal{O}(n^{-4})}{n^{-2\kappa'^{(1)}}} \to 0 \nonumber
\end{align}
as $n \to \infty$. Therefore, we achieve the conclusion of equality~\eqref{eqn:proof_Gaussian_secondcase_second} under Case 1. 
\paragraph{Case 2:} $\kappa' = \parenth{2, \kappa'^{(2)}, \kappa'^{(3)}}$ when $(\kappa_{1}'^{(2)},\kappa'^{(3)}) \prec (\overline{r},\overline{r}/2)$. Under this setting of $\kappa'$, we construct $G_{n} = \sum_{i=1}^{k}\pi_{i}^{n}\delta_{(\theta_{1i}^{n},\theta_{2i}^{n})}$ such that $(\pi_{i+k-k_{0}}^{n}, \theta_{1(i+k-k_{0})}^{n}, \theta_{2(i+k-k_{0})}^{n}) = (\pi_{i}^{0},\theta_{1i}^{0},\theta_{2i}^{0})$ for $2 \leq i \leq k_{0}$. For $1 \leq j \leq k - k_{0} + 1$, we choose $(\theta_{1j}^{n})^{(1)} = (\theta_{11}^{0})^{(1)}$ and
\begin{align}
(\theta_{1j}^{n})^{(2)} = (\theta_{11}^{0})^{(2)} + \frac{a_{j}^{*}}{n}, \ \theta_{2j}^{n} = \theta_{21}^{0} + \frac{2b_{j}^{*}}{n^{2}}, \ \pi_{j}^{n} = \frac{\pi_{1}^{0}(c_{j}^{*})^{2}}{\sum_{i=1}^{k-k_{0}+1} (c_{j}^{*})^{2}}, \nonumber
\end{align}
where $(c_{i}^{*},a_{i}^{*},b_{i}^{*})_{i=1}^{k-k_{0}+1}$ are the non-trivial solution of system of polynomial equations~\eqref{eqn:system_polynomial_Gaussian_first} when $r = \overline{r} - 1$. From here, by performing Taylor expansion around $((\theta_{11}^{0})^{(2)},\theta_{21}^{0})$, i.e., along the direction of the second component of $\theta_{11}^{0}$ and $\theta_{21}^{0}$ and arguing similarly as Case 2 in the proof of equality~\eqref{eqn:proof_Gaussian_firstcase_second} in Section~\ref{subsection:key_equality_algebraic_dependent_linear_first}, we obtain that
\begin{align}
p_{G_{n}}(X,Y) - p_{G_{0}}(X,Y) = \sum \limits_{l = \overline{r}}^{2\overline{r}-2} \mathcal{O}(n^{-\overline{r}})\dfrac{\partial^{l}{f}}{\partial{h_{1}^{l}}}(Y|h_{1}(X,\theta_{1i}^{0}),h_{2}(X,\theta_{2i}^{0}))\overline{f}(X) + \overline{R}_{2}(X,Y), \nonumber 
\end{align}
where $\overline{R}_{2}(X,Y)$ is Taylor remainder such that the following limit holds
\begin{align}
\int \frac{\overline{R}_{2}^{2}(X,Y)}{p_{G_{0}}(X,Y)\widetilde{W}_{\kappa'}^{2\|\kappa'\|_{\infty}}(G_{n},G_{0})}d(X,Y) \precsim \frac{\mathcal{O}(n^{-2\overline{r}})}{n^{-2\min\{\kappa'^{(2)},\kappa'^{(3)}\}}} \to 0. \nonumber
\end{align}
Therefore, we achieve that $h^{2}(p_{G_{n}}, p_{G_{0}})\bigg/\widetilde{W}_{\kappa'}^{2\|\kappa'\|_{\infty}}(G_{n},G_{0}) \to 0$
as $n \to \infty$. As a consequence, we reach the conclusion of equality~\eqref{eqn:proof_Gaussian_secondcase_second} under Case 2. 
\subsection{Proof of Theorem \ref{theorem:lower_bound_Gaussian_family_first_second}}
\label{Section:proof_weakly_identifiable_covariate_dependent_second}
Similar to the previous proofs, it is sufficient to demonstrate the following results:
\begin{eqnarray}
\lim \limits_{\epsilon \to 0} \inf \limits_{G \in \Ocal_{k,\overline{c}_{0}} (\Omega): \widetilde{W}_{\kappa}(G,G_{0}) \leq \epsilon} V(p_{G},p_{G_{0}})/\widetilde{W}_{\kappa}^{\|\kappa\|_{\infty}}(G,G_{0}) > 0, \label{eqn:proof_Gaussian_thirdcase_first} \\
\inf \limits_{G \in \Ocal_{k} (\Omega)} h(p_{G},p_{G_{0}})/\widetilde{W}_{\kappa'}^{\|\kappa'\|_{\infty}}(G,G_{0}) = 0, \label{eqn:proof_Gaussian_thirdcase_second}
\end{eqnarray}
for any $\kappa' \prec \kappa$ where $\kappa=(\overline{r},\overline{r},\ceil{\overline{r}/2},\ceil{\overline{r}/2})$. Without loss of generality, we assume that $\overline{r}$ is even, i.e., $\kappa = (\overline{r},\overline{r},\overline{r}/2,\overline{r}/2)$. Proof of inequality~\eqref{eqn:proof_Gaussian_thirdcase_first} is in Appendix~\ref{subsec:proof_Gaussian_thirdcase_first} while proof of equality~\eqref{eqn:proof_Gaussian_thirdcase_second} is in Appendix~\ref{subsec:proof_Gaussian_thirdcase_second}.
\subsubsection{Proof for inequality~\eqref{eqn:proof_Gaussian_thirdcase_first}}
\label{subsec:proof_Gaussian_thirdcase_first}
Assume that the conclusion of claim~\eqref{eqn:proof_Gaussian_thirdcase_first} does not hold. By using the same notations of $G_{n}$ as in the proof of Theorem~\ref{theorem:lower_bound_Gaussian_family}, we can find a sequence $G_{n}$ that has representation~\eqref{eqn:proof_Gaussian_firstcase_notation} such that $V(p_{G_{n}},p_{G_{0}})/\widetilde{W}_{\kappa}^{\|\kappa\|_{\infty}}(G_{n},G_{0}) \to 0$ and $\widetilde{W}_{\kappa}(G_{n},G_{0}) \to 0$. Here, since $\theta_{2ij}^{n}$ and $\theta_{2i}^{(0)}$ have 2 dimensions, we denote $\Delta \theta_{2ij}^{n} = ((\Delta \theta_{2ij}^{n})^{(1)}, (\Delta \theta_{2ij}^{n})^{(2)})$ for all $1 \leq i \leq k_{0}$ and $1 \leq j \leq s_{i}$ throughout this proof. According to Lemma \ref{lemma:generalized_Wasserstein_distance_polynomials}, we have
\begin{eqnarray}
\widetilde{W}_{\kappa}^{\|\kappa\|_{\infty}}(G_{n},G_{0}) \precsim \sum \limits_{i=1}^{k_{0}}\sum \limits_{j=1}^{s_{i}}{p_{ij}^{n}\biggr(\biggr|(\Delta \theta_{1ij}^{n})^{(1)}\biggr|^{\overline{r}}+\biggr|(\Delta \theta_{1ij}^{n})^{(2)}\biggr|^{\overline{r}}+\biggr|(\Delta \theta_{2ij}^{n})^{(1)}\biggr|^{\overline{r}/2}+\biggr|(\Delta \theta_{2ij}^{n})^{(2)}\biggr|^{\overline{r}/2}\biggr)} \nonumber \\
+ \sum \limits_{i=1}^{k_{0}} |\sum \limits_{j=1}^{s_{i}}p_{ij}^{n} - \pi_{i}^{0}| := D_{\kappa}(G_{n},G_{0}). \nonumber
\end{eqnarray}
Invoking Taylor expansion up to the order $\overline{r}$, we obtain that
\begin{eqnarray}
p_{G_{n}}(X,Y)-p_{G_{0}}(X,Y)& = & \sum \limits_{i=1}^{k_{0}}\sum \limits_{j=1}^{s_{i}}p_{ij}^{n}\sum \limits_{1 \leq |\alpha| \leq \overline{r}} \dfrac{1}{\alpha!}\biggr\{(\Delta \theta_{1ij}^{n})^{(1)}\biggr\}^{\alpha_{1}}\biggr\{(\Delta \theta_{1ij}^{n})^{(2)}\biggr\}^{\alpha_{2}}\biggr\{(\Delta \theta_{2ij}^{n})^{(1)}\biggr\}^{\alpha_{3}} \nonumber \\
& & \hspace{-8 em} \times \biggr\{(\Delta \theta_{2ij}^{n})^{(2)}\biggr\}^{\alpha_{4}}\dfrac{\partial^{|\alpha|}{f}}{\partial{(\theta_{1}^{(1)})^{\alpha_{1}}}\partial{(\theta_{1}^{(2)})^{\alpha_{2}}}\partial{(\theta_{2}^{(1)})^{\alpha_{3}}}\partial{(\theta_{2}^{(2)})^{\alpha_{4}}}}(Y|h_{1}(X|\theta_{1i}^{0}),h_{2}(X|\theta_{2i}^{0}))\overline{f}(X) \nonumber \\
& & \hspace{-4 em} + \sum \limits_{i=1}^{k_{0}}\biggr(\sum \limits_{j=1}^{s_{i}}p_{ij}^{n} - \pi_{i}^{0}\biggr)f(Y|h_{1}(X|\theta_{1i}^{0},h_{2}(X|\theta_{2i}^{0}))\overline{f}(X) + R(X,Y) \nonumber \\
& : = & A_{n} + B_{n} + R(X,Y), \nonumber
\end{eqnarray}
where $\alpha = (\alpha_{1},\alpha_{2},\alpha_{3},\alpha_{4})$ and $R(X,Y)$ is a Taylor remainder such that $R(X,Y)/D_{\kappa}(G_{n},G_{0}) \to 0$ as $n \to \infty$ for all $(X,Y)$. The formulation of expert functions $h_{1}, h_{2}$ and the PDE structure of Gaussian kernel lead to
\begin{eqnarray}
& & \hspace{-6 em} \dfrac{\partial^{|\alpha|}{f}}{\partial{(\theta_{1}^{(1)})^{\alpha_{1}}}\partial{(\theta_{1}^{(2)})^{\alpha_{2}}}\partial{(\theta_{2}^{(1)})^{\alpha_{3}}}\partial{(\theta_{2}^{(2)})^{\alpha_{4}}}}(Y|h_{1}(X,\theta_{1}),h_{2}(X,\theta_{2})) \nonumber \\
& & \hspace{4 em} = \dfrac{X^{\alpha_{2}+2\alpha_{4}}}{2^{\alpha_{3}+\alpha_{4}}}\dfrac{\partial^{\alpha_{1}+\alpha_{2}+2\alpha_{3}+2\alpha_{4}}{f}}{\partial{h_{1}^{\alpha_{1}+\alpha_{2}+2\alpha_{3}+2\alpha_{4}}}}(Y|h_{1}(X,\theta_{1}),h_{2}(X,\theta_{2})), \nonumber
\end{eqnarray}
for any $\alpha_{1}, \alpha_{2}, \alpha_{3}, \alpha_{4} \in \mathbb{N}$, $\theta_{1} \in \Omega_{1}$, and $\theta_{2} \in \Omega_{2}$. With the above equation, we can rewrite $A_{n}$ as follows
\begin{eqnarray}
A_{n} & = & \sum \limits_{i=1}^{k_{0}}\sum \limits_{j=1}^{s_{i}}p_{ij}^{n}\sum \limits_{1 \leq |\alpha| \leq \overline{r}} \dfrac{1}{\alpha!}\biggr\{(\Delta \theta_{1ij}^{n})^{(1)}\biggr\}^{\alpha_{1}}\biggr\{(\Delta \theta_{1ij}^{n})^{(2)}\biggr\}^{\alpha_{2}}\biggr\{(\Delta \theta_{2ij}^{n})^{(1)}\biggr\}^{\alpha_{3}} \biggr\{(\Delta \theta_{2ij}^{n})^{(2)}\biggr\}^{\alpha_{4}} \nonumber \\
& & \times \dfrac{X^{\alpha_{2}+2\alpha_{4}}}{2^{\alpha_{3}}}\dfrac{\partial^{\alpha_{1}+\alpha_{2}+2\alpha_{3}+2\alpha_{4}}{f}}{\partial{h_{1}^{\alpha_{1}+\alpha_{2}+2\alpha_{3}+2\alpha_{4}}}}(Y|h_{1}(X|\theta_{1i}^{0}),h_{2}(X|\theta_{2i}^{0}))\overline{f}(X) \nonumber \\
& = & \sum \limits_{i=1}^{k_{0}}\sum \limits_{j=1}^{s_{i}}p_{ij}^{n}\sum \limits_{1 \leq l_{1}+l_{2} \leq 2\overline{r}} \biggr(\sum \limits_{\alpha_{1},\alpha_{2},\alpha_{3},\alpha_{4}}\dfrac{1}{2^{\alpha_{3}+\alpha_{4}}\alpha!}\biggr\{(\Delta \theta_{1ij}^{n})^{(1)}\biggr\}^{\alpha_{1}}\biggr\{(\Delta \theta_{1ij}^{n})^{(2)}\biggr\}^{\alpha_{2}}\biggr\{(\Delta \theta_{2ij}^{n})^{(1)}\biggr\}^{\alpha_{3}} \nonumber \\
& & \times \biggr\{(\Delta \theta_{2ij}^{n})^{(2)}\biggr\}^{\alpha_{4}}\biggr)
X^{l_{2}}\dfrac{\partial^{l_{1}+l_{2}}{f}}{\partial{h_{1}^{l_{1}+l_{2}}}}(Y|h_{1}(X|\theta_{1i}^{0}),h_{2}(X|\theta_{2i}^{0}))\overline{f}(X), \label{eqn:proof_Gaussian_thirdcase_third}
\end{eqnarray}
where $\alpha_{1},\alpha_{2},\alpha_{3},\alpha_{4} \in \mathbb{N}$ in the sum of second equation satisfies $\alpha_{1}+2\alpha_{3}=l_{1}$, $\alpha_{2}+2\alpha_{4}=l_{2}$, and $1 \leq \alpha_{1} + \alpha_{2} + \alpha_{3} + \alpha_{4} \leq \overline{r}$. 

As demonstrated in the earlier proofs, we can treat $A_{n}/D_{\kappa}(G_{n},G_{0})$ and $B_{n}/D_{\kappa}(G_{n},G_{0})$ as a linear combination of $X^{l_{2}}\dfrac{\partial^{l_{1}+l_{2}}{f}}{\partial{h_{1}^{l_{1}+l_{2}}}}(Y|h_{1}(X|\theta_{1i}^{0}),h_{2}(X|\theta_{2i}^{0}))\overline{f}(X)$ for $0 \leq l_{1}+l_{2} \leq 2\overline{r}$ and $1 \leq i \leq k_{0}$, which are linearly independent with respect to $X$ and $Y$. For the simplicity of presentation, we denote $E_{l_{1},l_{2}}(\theta_{1i}^{0},\theta_{2i}^{0})$ as the coefficient of $X^{l_{2}}\dfrac{\partial^{l_{1}+l_{2}}{f}}{\partial{h_{1}^{l_{1}+l_{2}}}}(Y|h_{1}(X|\theta_{1i}^{0}),h_{2}(X|\theta_{2i}^{0}))\overline{f}(X)$ in $A_{n}, B_{n}$. From the equation~\eqref{eqn:proof_Gaussian_thirdcase_third}, we can check that
\begin{eqnarray}
E_{0,l}(\theta_{1i}^{0},\theta_{2i}^{0}) = \biggr(\sum \limits_{j=1}^{s_{i}}p_{ij}^{n}\sum \limits_{\substack{\alpha_{2}+2\alpha_{4}=l \\ \alpha_{2}+\alpha_{4} \leq \overline{r}}}\biggr\{(\Delta \theta_{1ij}^{n})^{(2)}\biggr\}^{\alpha_{2}}\biggr\{(\Delta \theta_{2ij}^{n})^{(2)}\biggr\}^{\alpha_{4}}\biggr)/(2^{\alpha_{4}}\alpha_{2}!\alpha_{4}!), \nonumber \\
E_{l,0}(\theta_{1i}^{0},\theta_{2i}^{0}) = \biggr(\sum \limits_{j=1}^{s_{i}}p_{ij}^{n} \sum \limits_{\substack{\alpha_{1}+2\alpha_{3}=l \\ \alpha_{1}+\alpha_{3} \leq \overline{r}}} \biggr\{(\Delta \theta_{1ij}^{n})^{(1)}\biggr\}^{\alpha_{1}}\biggr\{(\Delta \theta_{2ij}^{n})^{(1)}\biggr\}^{\alpha_{3}}\biggr)/(2^{\alpha_{3}}\alpha_{1}!\alpha_{3}!), \nonumber
\end{eqnarray}
for any $1 \leq l \leq 2\overline{r}$. 

Assume that all of the coefficients of $A_{n}/D_{\kappa}(G_{n},G_{0})$ and $B_{n}/D_{\kappa}(G_{n},G_{0})$ go to 0 as $n \to \infty$. The summation of $|E_{0,0}(\theta_{1i}^{0},\theta_{2i}^{0})/D_{\kappa}(G_{n},G_{0})|$ for all $1 \leq i \leq k_{0}$ leads to
\begin{eqnarray}
\biggr(\sum \limits_{i=1}^{k_{0}} |\sum \limits_{j=1}^{s_{i}}{p_{ij}^{n}} - \pi_{i}^{0}|\biggr)/D_{\kappa}(G_{n},G_{0}) \to 0. \nonumber
\end{eqnarray}
From the formulation of $D_{\kappa}(G_{n},G_{0})$, the above limit implies that
\begin{eqnarray}
\biggr\{\sum \limits_{i=1}^{k_{0}}\sum \limits_{j=1}^{s_{i}}{p_{ij}^{n}\biggr(\biggr|(\Delta \theta_{1ij}^{n})^{(1)}\biggr|^{\overline{r}}+\biggr|(\Delta \theta_{1ij}^{n})^{(2)}\biggr|^{\overline{r}}+\biggr|(\Delta \theta_{2ij}^{n})^{(1)}\biggr|^{\overline{r}/2}+\biggr|(\Delta \theta_{2ij}^{n})^{(2)}\biggr|^{\overline{r}/2}\biggr)}\biggr\}/D_{\kappa}(G_{n},G_{0}) \to 1. \nonumber
\end{eqnarray}
Therefore, we can find an index $i^{*} \in \left\{1,\ldots,k_{0}\right\}$ such that
\begin{eqnarray}
\biggr\{\sum \limits_{j=1}^{s_{i^{*}}}{p_{i^{*}j}^{n}\biggr(\biggr|(\Delta \theta_{1i^{*}j}^{n})^{(1)}\biggr|^{\overline{r}}+\biggr|(\Delta \theta_{1i^{*}j}^{n})^{(2)}\biggr|^{\overline{r}}+\biggr|(\Delta \theta_{2i^{*}j}^{n})^{(1)}\biggr|^{\overline{r}/2}+\biggr|(\Delta \theta_{2i^{*}j}^{n})^{(2)}\biggr|^{\overline{r}/2}\biggr)}\biggr\}/D_{\kappa}(G_{n},G_{0}) \not \to 0. \nonumber
\end{eqnarray}
The above result leads to two distinct cases.
\paragraph{Case 1:} $\biggr\{\biggr(\sum \limits_{j=1}^{s_{i^{*}}}{p_{i^{*}j}^{n}\biggr(\biggr|(\Delta \theta_{1i^{*}j}^{n})^{(1)}\biggr|^{\overline{r}}+\biggr|(\Delta \theta_{2i^{*}j}^{n})^{(1)}\biggr|^{\overline{r}/2}\biggr)}\biggr\}/D_{\kappa}(G_{n},G_{0}) \not \to 0$. By taking the product between the inverse of the previous ratio and $E_{l,0}(\theta_{1i^{*}}^{0},\theta_{2i^{*}}^{0})/D_{\kappa}(G_{n},G_{0})$, we achieve the following system of limits
\begin{eqnarray}
\dfrac{\biggr(\sum \limits_{j=1}^{s_{i^{*}}}p_{i^{*}j}^{n}\sum \limits_{\substack{\alpha_{1}+2\alpha_{3}=l \\ \alpha_{1}+\alpha_{3} \leq \overline{r}}} \biggr\{(\Delta \theta_{1i^{*}j}^{n})^{(1)}\biggr\}^{\alpha_{1}}\biggr\{(\Delta \theta_{2i^{*}j}^{n})^{(1)}\biggr\}^{\alpha_{3}}\biggr)/(2^{\alpha_{3}}\alpha_{1}!\alpha_{3}!)}{\biggr(\sum \limits_{j=1}^{s_{i^{*}}}{p_{i^{*}j}^{n}\biggr(\biggr|(\Delta \theta_{1i^{*}j}^{n})^{(1)}\biggr|^{\overline{r}}+\biggr|(\Delta \theta_{2i^{*}j}^{n})^{(1)}\biggr|^{\overline{r}/2}\biggr)}} \to 0, \nonumber
\end{eqnarray}
for all $1 \leq l \leq 2\overline{r}$, which does not hold according to the argument of the proof of Theorem \ref{theorem:lower_bound_Gaussian_family}. Therefore, Case 1 can not hold.
\paragraph{Case 2:} $\biggr\{\biggr(\sum \limits_{j=1}^{s_{i^{*}}}{p_{i^{*}j}^{n}\biggr(\biggr|(\Delta \theta_{1i^{*}j}^{n})^{(2)}\biggr|^{\overline{r}}+\biggr|(\Delta \theta_{2i^{*}j}^{n})^{(2)}\biggr|^{\overline{r}/2}\biggr)}\biggr\}/D_{\kappa}(G_{n},G_{0}) \not \to 0$. By taking the product between the inverse of the previous ratio with $E_{0,l}(\theta_{1i^{*}}^{0},\theta_{2i^{*}}^{0})/D_{\kappa}(G_{n},G_{0})$, we obtain that following system of limits
\begin{eqnarray}
\dfrac{\biggr(\sum \limits_{j=1}^{s_{i^{*}}}p_{i^{*}j}^{n}\sum \limits_{\substack{\alpha_{2}+2\alpha_{4}=l \\ \alpha_{2}+\alpha_{4} \leq \overline{r}}} \biggr\{(\Delta \theta_{1i^{*}j}^{n})^{(2)}\biggr\}^{\alpha_{2}}\biggr\{(\Delta \theta_{2i^{*}j}^{n})^{(2)}\biggr\}^{\alpha_{4}}\biggr)/(2^{\alpha_{4}}\alpha_{2}!\alpha_{4}!)}{\biggr(\sum \limits_{j=1}^{s_{i^{*}}}{p_{i^{*}j}^{n}\biggr(\biggr|(\Delta \theta_{1i^{*}j}^{n})^{(2)}\biggr|^{\overline{r}}+\biggr|(\Delta \theta_{2i^{*}j}^{n})^{(2)}\biggr|^{\overline{r}/2}\biggr)}} \to 0, \nonumber
\end{eqnarray}
for all $1 \leq l \leq 2\overline{r}$, which does not hold. Thus, Case 2 can not happen.

As a consequence, not all the coefficients of $A_{n}/D_{\kappa}(G_{n},G_{0})$ and $B_{n}/D_{\kappa}(G_{n},G_{0})$ go to 0 as $n \to \infty$. From here, using the same argument as the Fatou's argument in Step 4 of the proof of Theorem~\ref{theorem:lower_bound_Gaussian_family}, we achieve the conclusion of inequality~\eqref{eqn:proof_Gaussian_thirdcase_first}.
\subsubsection{Proof for equality~\eqref{eqn:proof_Gaussian_thirdcase_second}} 
\label{subsec:proof_Gaussian_thirdcase_second}
To avoid unnecessary repetition, we only sketch the proof for equality~\eqref{eqn:proof_Gaussian_thirdcase_second}. Since $\kappa' = (\kappa'^{(1)}, \kappa'^{(2)}, \kappa'^{(3)}, \kappa'^{(4)}) \prec (\overline{r}, \overline{r}, \overline{r}/2 , \overline{r}/2)$, one of the two pairs $(\kappa'^{(1)},\kappa'^{(3)})$, $(\kappa'^{(2)},\kappa'^{(4)})$ is strictly dominated by $(\overline{r},\overline{r}/2)$. Without loss of generality, we assume that $(\kappa'^{(1)},\kappa'^{(3)}) \prec (\overline{r}, \overline{r}/2)$. Under this setting of $\kappa'$, we construct a sequence of mixing measures $G_{n} = \sum_{i=1}^{k}\pi_{i}^{n}\delta_{(\theta_{1i}^{n},\theta_{2i}^{n})}$ as follows. We choose $(\pi_{i+k-k_{0}}^{n}, \theta_{1(i+k-k_{0})}^{n}, \theta_{2(i+k-k_{0})}^{n}) = (\pi_{i}^{0},\theta_{1i}^{0},\theta_{2i}^{0})$ for $2 \leq i \leq k_{0}$. For $1 \leq j \leq k - k_{0} + 1$, we choose $\parenth{(\theta_{1j}^{n})^{(2)},(\theta_{2j}^{n})^{(2)}} \equiv \parenth{(\theta_{11}^{0})^{(2)},(\theta_{21}^{0})^{(2)}}$ and
\begin{align}
(\theta_{1j}^{n})^{(1)} = (\theta_{11}^{0})^{(1)} + \frac{a_{j}^{*}}{n}, \ (\theta_{2j}^{n})^{(1)} = (\theta_{21}^{0})^{(1)} + \frac{2b_{j}^{*}}{n^{2}}, \ \pi_{j}^{n} = \frac{\pi_{1}^{0}(c_{j}^{*})^{2}}{\sum_{i=1}^{k-k_{0}+1} (c_{j}^{*})^{2}}, \nonumber
\end{align}
where $(c_{i}^{*},a_{i}^{*},b_{i}^{*})_{i=1}^{k-k_{0}+1}$ are the non-trivial solution of system of polynomial equations~\eqref{eqn:system_polynomial_Gaussian_first} when $r = \overline{r} - 1$. From here, by performing Taylor expansion around $((\theta_{11}^{0})^{(1)},(\theta_{21}^{0})^{(1)})$, i.e., along the direction of the first component of $\theta_{11}^{0}$ and $\theta_{21}^{0}$ and arguing similarly as Case 2 in the proof of equality~\eqref{eqn:proof_Gaussian_firstcase_second} in Section~\ref{subsection:key_equality_algebraic_dependent_linear_first}, we obtain that
\begin{align}
p_{G_{n}}(X,Y) - p_{G_{0}}(X,Y) = \sum \limits_{l = \overline{r}}^{2\overline{r}-2} \mathcal{O}(n^{-\overline{r}})\dfrac{\partial^{l}{f}}{\partial{h_{1}^{l}}}(Y|h_{1}(X,\theta_{1i}^{0}),h_{2}(X,\theta_{2i}^{0}))\overline{f}(X) + \overline{R}(X,Y), \nonumber 
\end{align}
where $\overline{R}(X,Y)$ is a Taylor remainder such that the following limit holds
\begin{align}
\int \frac{\overline{R}^{2}(X,Y)}{p_{G_{0}}(X,Y)\widetilde{W}_{\kappa'}^{2\|\kappa'\|_{\infty}}(G_{n},G_{0})}d(X,Y) \precsim \frac{\mathcal{O}(n^{-2\overline{r}})}{n^{-2\min\{\kappa'^{(1)},\kappa'^{(3)}\}}} \to 0. \nonumber
\end{align}
As a consequence, we eventually achieve that
\begin{align}
h^{2}(p_{G_{n}}, p_{G_{0}})\bigg/\widetilde{W}_{\kappa'}^{2\|\kappa'\|_{\infty}}(G_{n},G_{0}) \to 0 \nonumber
\end{align}
as $n \to \infty$, which leads to the conclusion of equality~\eqref{eqn:proof_Gaussian_thirdcase_second}.
\subsection{Proof of Theorem \ref{theorem:lower_bound_Gaussian_family_third}}
\label{Section:proof_non_linear_covariate_dependent_first}
Similar to the proof of Theorem \ref{theorem:lower_bound_Gaussian_family}, to obtain the conclusion of Theorem~\ref{theorem:lower_bound_Gaussian_family_third}, it is sufficient to demonstrate that
\begin{eqnarray}
\lim \limits_{\epsilon \to 0} \inf \limits_{G \in \Ocal_{k,\overline{c}_{0}} (\Omega): W_{\kappa}(G,G_{0}) \leq \epsilon} V(p_{G},p_{G_{0}})/\widetilde{W}_{\kappa}^{\|\kappa\|_{\infty}}(G,G_{0}) > 0, \label{eqn:proof_Gaussian_non_linear_settingI_first} \\
\lim \limits_{\epsilon \to 0} \inf \limits_{G \in \mathcal{G}_{1}: \ W_{\kappa'}(G,G_{0}) \leq \epsilon} h(p_{G},p_{G_{0}})/\widetilde{W}_{\kappa'}^{\|\kappa'\|_{\infty}}(G,G_{0}) = 0, \label{eqn:proof_equality_Gaussian_non_linear_settingI_first}
\end{eqnarray}
for all $\kappa' \prec \kappa=(2,2,2)$ where $\kappa : = (2,2,2)$ and $\mathcal{G}_{1}$ is defined in Theorem~\ref{theorem:lower_bound_Gaussian_family_third}.
\subsubsection{Proof of inequality~\eqref{eqn:proof_Gaussian_non_linear_settingI_first}}
\label{subsec:proof_Gaussian_non_linear_settingI_first}
Assume that the above result does not hold, which leads to the existence of sequence $G_{n} = \sum \limits_{i=1}^{k_{0}}\sum \limits_{j=1}^{s_{i}}p_{ij}^{n}\delta_{(\theta_{1ij}^{n},\theta_{2ij}^{n})}$ such that  $V\parenth{p_{G_{n}},p_{G_{0}}}/\widetilde{W}_{\kappa}^{\|\kappa\|_{\infty}}(G_{n},G_{0}) \to 0$ and $\widetilde{W}_{\kappa}(G_{n},G_{0}) \to 0$. Here, $(\theta_{1ij}^{n},\theta_{2ij}^{n}) \to (\theta_{1i}^{0},\theta_{2i}^{0})$ for all $1 \leq i \leq k_{0}, 1 \leq j \leq s_{i}$ and $\sum \limits_{j=1}^{s_{i}}p_{ij}^{n} \to \pi_{i}^{0}$ for all $1 \leq i \leq k_{0}$. In this proof, we denote $\Delta \theta_{1ij}^{n} : = ((\Delta \theta_{1ij}^{n})^{(1)}, (\Delta \theta_{1ij}^{n})^{(2)})$ for all $1 \leq i \leq k_{0}$ and $1 \leq j \leq s_{i}$. According to Lemma \ref{lemma:generalized_Wasserstein_distance_polynomials} in Appendix~\ref{sec:auxi_result}, we have
\begin{eqnarray}
\widetilde{W}_{\kappa}^{\|\kappa\|_{\infty}}(G_{n},G_{0}) \precsim \sum \limits_{i=1}^{k_{0}}\sum \limits_{j=1}^{s_{i}}{p_{ij}^{n}\biggr(\biggr|(\Delta \theta_{1ij}^{n})^{(1)}\biggr|^{\kappa^{(1)}}+\biggr|(\Delta \theta_{1ij}^{n})^{(2)}\biggr|^{\kappa^{(2)}}+\biggr|\Delta \theta_{2ij}^{n}\biggr|^{\kappa^{(3)}} \biggr)} \nonumber \\
+ \sum \limits_{i=1}^{k_{0}} \abss{\sum \limits_{j=1}^{s_{i}}p_{ij}^{n} - \pi_{i}^{0}} := D_{\kappa}(G_{n},G_{0}). \nonumber
\end{eqnarray}
Since the proof argument for claim~\eqref{eqn:proof_Gaussian_non_linear_settingI_first} is rather intricate, we divide this argument into several steps. 
\paragraph{Step 1 - Structure of Taylor expansion:} By means of Taylor expansion up to the order $\|\kappa\|_{\infty} = 2$, we obtain that
\begin{eqnarray}
p_{G_{n}}(X,Y)-p_{G_{0}}(X,Y)& = & \sum \limits_{i=1}^{k_{0}}\sum \limits_{j=1}^{s_{i}}p_{ij}^{n}\sum \limits_{1 \leq |\alpha| \leq \|\kappa\|_{\infty}} \dfrac{1}{\alpha!}\biggr\{(\Delta \theta_{1ij}^{n})^{(1)}\biggr\}^{\alpha_{1}}\biggr\{(\Delta \theta_{1ij}^{n})^{(2)}\biggr\}^{\alpha_{2}}(\Delta \theta_{2ij}^{n})^{\alpha_{3}} \nonumber \\
& \times & \dfrac{\partial^{|\alpha|}{f}}{\partial{(\theta_{1}^{(1)})^{\alpha_{1}}}\partial{(\theta_{1}^{(2)})^{\alpha_{2}}}\partial{\theta_{2}^{\alpha_{3}}}}(Y|h_{1}(X|\theta_{1i}^{0}),h_{2}(X|\theta_{2i}^{0}))\overline{f}(X) \nonumber \\
& + & \sum \limits_{i=1}^{k_{0}}\biggr(\sum \limits_{j=1}^{s_{i}}p_{ij}^{n} - p_{i}^{0}\biggr)f(Y|h_{1}(X|\theta_{1i}^{0},h_{2}(X|\theta_{2i}^{0}))\overline{f}(X) + R(X,Y) \nonumber \\
& : = & A_{n} + B_{n} + R(X,Y), \nonumber
\end{eqnarray}
where $\alpha = (\alpha_{1},\alpha_{2},\alpha_{3})$ and $R(X,Y)$ is a Taylor remainder such that $R(X,Y)/D_{\kappa}(G_{n},G_{0}) \to 0$ as $n \to \infty$. 

With the formation of expert functions $h_{1}$ and $h_{2}$, we can check that $\dfrac{\partial^{|\alpha|}{f}}{\partial{(\theta_{1}^{(1)})^{\alpha_{1}}}\partial{(\theta_{1}^{(2)})^{\alpha_{2}}}\partial{\theta_{2}^{\alpha_{3}}}}$\\$(Y|h_{1}(X|\theta_{1i}^{0}),h_{2}(X|\theta_{2i}^{0}))\overline{f}(X)$ are not linearly independent with respect to $X$ and $Y$. Therefore, as being argued in the proof of Theorem \ref{theorem:lower_bound_Gaussian_family}, we can not consider $A_{n}$ as a linear combinations of these derivatives. To see clearly the influence of non-linearity setting I of $G_{0}$ on the set of linear independent elements of $A_{n}$, we will provide the detail formulations of key partial derivatives of $f$ with respect to $\theta_{1}$ and $\theta_{2}$ up to the second order. 
\paragraph{Key partial derivatives up to the second order:} In particular, for any $\theta_{1}$ and $\theta_{2}$, by means of direct computation and the PDE equation $\dfrac{\partial^{2}{f}}{\partial{h_{1}^{2}}} = 2 \dfrac{\partial{f}}{\partial{h_{2}^{2}}}$, we can verify that
\begin{align}
\dfrac{\partial{f}}{\partial{\theta_{1}^{(1)}}} = 2\parenth{\theta_{1}^{(1)}+\theta_{1}^{(2)}X}\dfrac{\partial{f}}{\partial{h_{1}}}, \ \dfrac{\partial{f}}{\partial{\theta_{1}^{(2)}}} = 2X\parenth{\theta_{1}^{(1)}+\theta_{1}^{(2)}X}\dfrac{\partial{f}}{\partial{h_{1}}}, \nonumber \\
\dfrac{\partial^{2}{f}}{\partial{(\theta_{1}^{(1)})^{2}}} = 2\dfrac{\partial{f}}{\partial{h_{1}}} +  4\parenth{\theta_{1}^{(1)}+\theta_{1}^{(2)}X}^{2}\dfrac{\partial^{2}{f}}{\partial{h_{1}^{2}}}, \nonumber
\end{align}
\begin{align}
\dfrac{\partial^{2}{f}}{\partial{(\theta_{1}^{(2)})^{2}}} = 2X^2\dfrac{\partial{f}}{\partial{h_{1}}} +  4X^2\parenth{\theta_{1}^{(1)}+\theta_{1}^{(2)}X}^{2}\dfrac{\partial^{2}{f}}{\partial{h_{1}^{2}}}, \nonumber \\
\dfrac{\partial^{2}{f}}{\partial{\theta_{1}^{(1)}}\partial{\theta_{1}^{(2)}}} = 2X\dfrac{\partial{f}}{\partial{h_{1}}} +  4X\parenth{\theta_{1}^{(1)}+\theta_{1}^{(2)}X}^{2}\dfrac{\partial^{2}{f}}{\partial{h_{1}^{2}}},  \nonumber \\
\dfrac{\partial{f}}{\partial{\theta_{2}}} = \dfrac{\partial{f}}{\partial{h_{2}^{2}}} = \dfrac{1}{2}\dfrac{\partial^{2}{f}}{\partial{h_{1}^{2}}}, \ \dfrac{\partial^{2}{f}}{\partial{\theta_{2}^{2}}} = \dfrac{\partial^{2}{f}}{\partial{h_{2}^{4}}} = \dfrac{1}{4}\dfrac{\partial^{4}{f}}{\partial{h_{1}^{4}}}. \label{eqn:partial_derivatives_second_order}
\end{align}
Here, we suppress the condition on  $h_{1}(X,\theta_{1})$ and $h_{2}(X,\theta_{2})$ in the notation to simplify the presentation. 
\paragraph{Set of linear independent elements:} We define
\begin{align}
\mathcal{F} : = \left\{X^{l_{1}}\dfrac{\partial^{l_{2}}{f}}{\partial{h_{1}^{l_{2}}}}(Y|h_{1}(X,\theta_{1i}^{0}),h_{2}(X,\theta_{2i}^{0}))\overline{f}(X): \ (l_{1},l_{2}) \in \mathcal{B}, \ 1 \leq i \leq k_{0} \right\}, \nonumber
\end{align}
where $\mathcal{B} = \left\{(0,0), (0,1) (0,2) , (0,3), (0,4) , (1,1), (1,2), (1,3), (2,1), (2,2), (2,3), (3,2), (4,2) \right\}$. According to the key partial derivatives of $f$ up to the second order given by equation~\eqref{eqn:partial_derivatives_second_order}, we can validate that the elements of $\mathcal{F}$ are linearly independent with respect to $X$ and $Y$. Therefore, we can treat $A_{n}/D_{\kappa}(G_{n},G_{0})$, $B_{n}/D_{\kappa}(G_{n},G_{0})$ as a linear combination of linear independent elements of $\mathcal{F}$. 
\paragraph{Step 2 - Non-vanishing coefficients:} Assume that all the coefficients in the representation of $A_{n}/D_{\kappa}(G_{n},G_{0})$ and $B_{n}/D_{\kappa}(G_{n},G_{0})$ go to 0 as $n \to \infty$. By taking the summation of the absolute value of coefficients in $B_{n}/D_{\kappa}(G_{n},G_{0})$, it implies that
\begin{align}
\parenth{\sum \limits_{i=1}^{k_{0}} |\sum \limits_{j=1}^{s_{i}} p_{ij}^{n} - \pi_{i}^{0}|}/D_{\kappa}(G_{n},G_{0}) \to 0. \nonumber
\end{align}
Furthermore, from the formulations of key partial derivatives in equation~\eqref{eqn:partial_derivatives_second_order}, the vanishing of coefficients of $\dfrac{\partial^{4}{f}}{\partial{h_{1}^{4}}}(Y|h_{1}(X,\theta_{1i}^{0}),h_{2}(X,\theta_{2i}^{0}))$ to 0 as $1 \leq i \leq k_{0}$ leads to
\begin{align}
\parenth{\sum \limits_{i=1}^{k_{0}} \sum \limits_{j=1}^{s_{i}}p_{ij}^{n}\abss{\Delta \theta_{2ij}^{n}}^{2}}/D_{\kappa}(G_{n},G_{0}) \to 0. \nonumber
\end{align}
Given the above results, the following holds
\begin{align}
\parenth{\sum \limits_{i=1}^{k_{0}} \left\{\sum \limits_{j=1}^{s_{i}} p_{ij}^{n} \abss{\Delta \theta_{2ij}^{n}}^{2} + |\sum \limits_{j=1}^{s_{i}} p_{ij}^{n} - \pi_{i}^{0}|\right\}}/D_{\kappa}(G_{n},G_{0}) \to 0. \label{eqn:limit_first_uniform_non-linearity_first}
\end{align}
On the other hand, the hypothesis that the coefficients of $X^{l_{1}}\dfrac{\partial{f}}{\partial{h_{1}}}(Y|h_{1}(X,\theta_{1i}^{0}),h_{2}(X,\theta_{2i}^{0}))$ as $l_{1} \in \left\{0,2\right\}$ and $X^{l_{2}}\dfrac{\partial^{2}{f}}{\partial{h_{1}^{2}}}(Y|h_{1}(X,\theta_{1i}^{0}),h_{2}(X,\theta_{2i}^{0}))$ as $l_{2} \in \left\{0,2\right\}$ as $1 \leq l_{2} \leq 4$ go to 0 as $1 \leq i \leq k_{0}$ respectively lead to the following system of polynomial limits:
\begin{align}
2(\theta_{1i}^{0})^{(1)}\parenth{\sum \limits_{j=1}^{s_{i}} p_{ij}^{n}(\Delta \theta_{1ij}^{n})^{(1)}}/D_{\kappa}(G_{n},G_{0}) + 2 I_{n,i} \to 0, \nonumber \\
2(\theta_{1i}^{0})^{(2)}\parenth{\sum \limits_{j=1}^{s_{i}} p_{ij}^{n}(\Delta \theta_{1ij}^{n})^{(2)}}/D_{\kappa}(G_{n},G_{0}) + 2 K_{n,i} \to 0, \nonumber \\
2(\theta_{1i}^{0})^{(1)}(\theta_{1i}^{0})^{(2)}I_{n,i} + \left\{(\theta_{1i}^{0})^{(1)}\right\}^{2}J_{n,i} \to 0, \nonumber \\
\left\{(\theta_{1i}^{0})^{(2)}\right\}^{2}I_{n,i} + 2(\theta_{1i}^{0})^{(1)}(\theta_{1i}^{0})^{(2)}J_{n,i} + \left\{(\theta_{1i}^{0})^{(1)}\right\}^{2}K_{n,i} \to 0, \nonumber \\
2(\theta_{1i}^{0})^{(1)}(\theta_{1i}^{0})^{(2)}K_{n,i} + \left\{(\theta_{1i}^{0})^{(2)}\right\}^{2}J_{n,i} \to 0, \nonumber \\
\left\{(\theta_{1i}^{0})^{(2)}\right\} K_{n,i} \to 0, \label{eqn:simple_system_polynomial_limits}
\end{align}
for all $1 \leq i \leq k_{0}$ where the explicit forms of $I_{n,i}$, $J_{n,i}$, and $K_{n,i}$ are as follows:
\begin{align}
I_{n,i} : = \parenth{\sum \limits_{j=1}^{s_{i}}p_{ij}^{n}\abss{(\Delta \theta_{1ij}^{n})^{(1)}}^{2}}/D_{\kappa}(G_{n},G_{0}), \nonumber \\ J_{n,i} : = \parenth{\sum \limits_{j=1}^{s_{i}}p_{ij}^{n}(\Delta \theta_{1ij}^{n})^{(1)}(\Delta \theta_{1ij}^{n})^{(2)}}/D_{\kappa}(G_{n},G_{0}), \nonumber \\
K_{n,i} : = \parenth{\sum \limits_{j=1}^{s_{i}}p_{ij}^{n}\abss{(\Delta \theta_{1ij}^{n})^{(2)}}^{2}}/D_{\kappa}(G_{n},G_{0}). \nonumber
\end{align}
According to the formulation of non-linearity setting I of $G_{0}$, we only have two possible cases to consider with respect to a pair $\parenth{(\theta_{1i}^{0})^{(1)},(\theta_{1i}^{0})^{(2)}}$:
\paragraph{Case 1:} $(\theta_{1i}^{0})^{(2)} \neq 0$. Under this case, the final limit in system of limits~\eqref{eqn:simple_system_polynomial_limits} indicates that $K_{n,i} \to 0$ as $n \to \infty$. Plugging this result into the fifth limit in this system, we achieve that $J_{n,i} \to 0$ as $n \to \infty$. Putting the previous results together, the third limit in the system of limits leads to $I_{n,i} \to 0$ as $n \to \infty$. 
\paragraph{Case 2:} $(\theta_{1i}^{0})^{(1)} = (\theta_{1i}^{0})^{(2)} = 0$. Under this case, the final four limits in system of polynomial limits~\eqref{eqn:simple_system_polynomial_limits} always hold. On the other hand, the first two limits of this system leads to $I_{n,i} \to 0$ and $K_{n,i} \to 0$ as $n \to \infty$. 

Given the results from Case 1 and Case 2, the following limit holds
\begin{align}
\dfrac{\sum \limits_{i=1}^{k_{0}} \sum \limits_{j=1}^{s_{i}} p_{ij}^{n}\parenth{\biggr|(\Delta \theta_{1ij}^{n})^{(1)}\biggr|^{2}+\biggr|(\Delta \theta_{1ij}^{n})^{(2)}\biggr|^{2}}}{D_{\kappa}(G_{n},G_{0})} \to 0. \label{eqn:limit_first_uniform_non-linearity_second}
\end{align}
Putting the results from~\eqref{eqn:limit_first_uniform_non-linearity_first} and~\eqref{eqn:limit_first_uniform_non-linearity_second} together, we obtain that
\begin{align}
1 = D_{\kappa}(G_{n},G_{0}) / D_{\kappa}(G_{n},G_{0}) \to 0, \nonumber
\end{align}
which is a contradiction. Therefore, not all the coefficients of $A_{n}/D_{\kappa}(G_{n},G_{0})$ and $B_{n}/D_{\kappa}(G_{n},G_{0})$ go to 0 as $n \to \infty$. From here, using the same argument as that of using the Fatou's argument in Step 4 of the proof of Theorem~\ref{theorem:lower_bound_Gaussian_family}, we achieve the conclusion of inequality~\eqref{eqn:proof_Gaussian_non_linear_settingI_first} under non-linearity setting I of $G_{0}$. 
\subsubsection{Proof of equality~\eqref{eqn:proof_equality_Gaussian_non_linear_settingI_first}} 
We will construct a similar sequence of mixing measures $G_{n}$ as that in the proof of~\eqref{eqn:general_overfit_second} in Section~\ref{subsection:equality_overfit}. More precisely, we define $G_{n} = \sum_{i=1}^{k_{0}+1} \pi_{i}^{n}\delta_{\parenth{\theta_{1i}^{n},\theta_{2i}^{n}}}$ with $k_{0}+1$ components as follows: $(\pi_{i}^{n}, \theta_{1i}^{n},\theta_{2i}^{n}) \equiv (\pi_{i-1}^{0}, \theta_{1(i-1)}^{0}, \theta_{2(i-1)}^{0})$ for $3 \leq i \leq k_{0} + 1$. Additionally, $\pi_{1}^{n} = \pi_{2}^{n} = 1/2$, $(\theta_{11}^{n},\theta_{21}^{n}) \equiv (\theta_{11}^{0} - \vec{1}_{2}/n, \theta_{21}^{0} - 1/n)$, and $(\theta_{12}^{n},\theta_{22}^{n}) \equiv (\theta_{11}^{0} + \vec{1}_{2}/n, \theta_{21}^{0} + 1/n)$. Now, by means of Taylor expansion up to the first order, the detail formulations of first order derivatives in~\eqref{eqn:partial_derivatives_second_order}, and the choice of $G_{n}$, we have
\begin{align}
p_{G_{n}}(X,Y) - p_{G_{0}}(X,Y) & = \sum \limits_{i=1}^{2} \pi_{i}^{n}\parenth{f(Y|h_{1}(X,\theta_{1i}^{n},\theta_{2i}^{n}) - f(Y|h_{1}(X,\theta_{11}^{0},\theta_{21}^{0})}\overline{f}(X) \nonumber \\
& = \sum \limits_{i=1}^{2} \pi_{i}^{n} \sum \limits_{|\alpha|+|\beta| = 1} \dfrac{1}{\alpha!\beta!}\prod \limits_{u=1}^{q_{1}} \biggr\{(\Delta \theta_{1i}^{n})^{(u)}\biggr\}^{\alpha_{u}}\prod \limits_{u=1}^{q_{1}}\biggr\{(\Delta \theta_{2i}^{n})^{(v)}\biggr\}^{\beta_{v}} \nonumber \\
& \times \dfrac{\partial{f}}{\partial{\theta_{1}^{\alpha}}\partial{\theta_{2}^{\beta}}}\parenth{Y|h_{1}(X,\theta_{11}^{0}),h_{2}(X,\theta_{21}^{0})}\overline{f}(X) + \overline{R}(X,Y) \nonumber \\
& = \overline{R}(X,Y), \nonumber
\end{align}
where $\Delta \theta_{1i}^{n} = \theta_{1i}^{n} - \theta_{11}^{0}$ and $\Delta \theta_{2i}^{n} = \theta_{2i}^{n} - \theta_{21}^{0}$ for $1 \leq i \leq 2$. Using the similar argument as that in the proof of~\eqref{eqn:general_overfit_second} in Section~\ref{subsection:equality_overfit}, $\overline{R}(X,Y)$ is a Taylor remainder from the above expansion such that
\begin{align}
\int \frac{\overline{R}_{1}^{2}(X,Y)}{p_{G_{0}}(X,Y)\widetilde{W}_{\kappa'}^{2\|\kappa'\|_{\infty}}(G_{n},G_{0})}d(X,Y) \precsim \frac{\mathcal{O}(n^{-4})}{n^{-2\min\{\kappa'^{(1)},\kappa'^{(2)},\kappa'^{(3)}\}}} \to 0 \nonumber
\end{align}
as $n \to \infty$. As a consequence, we achieve the conclusion of equality~\eqref{eqn:proof_equality_Gaussian_non_linear_settingI_first}.
\subsection{Proof of Theorem~\ref{theorem:lower_bound_Gaussian_family_second_singularity_case}}
\label{Section:proof_non_linear_covariate_dependent_second}
To achieve the conclusion of the theorem, it is sufficient to demonstrate that
\begin{eqnarray}
\lim \limits_{\epsilon \to 0} \inf \limits_{G \in \Ocal_{k,\overline{c}_{0}(\Omega)}: \widetilde{W}_{\widetilde{\kappa}_{\text{sin}}}(G,G_{0}) \leq \epsilon} V(p_{G},p_{G_{0}})/\widetilde{W}_{\widetilde{\kappa}_{\text{sin}}}^{\rtil_{\text{sin}}}(G,G_{0}) > 0 \nonumber
\end{eqnarray}
where $\rtil_{\text{sin}} = \rtil\parenth{(\theta_{1i_{\text{max}}}^{0})^{(1)},k-k_{0}+1}$ and $\widetilde{\kappa}_{\text{sin}} = \parenth{\rtil_{\text{sin}}, 2, \ceil{\rtil_{\text{sin}}/2}}$. Assume that the above result does not hold. It implies that we can find sequence $G_{n}$ such that 
\begin{align}
V\parenth{p_{G_{n}},p_{G_{0}}}/\widetilde{W}_{\widetilde{\kappa}_{\text{sin}}}^{\rtil_{\text{sin}}}(G_{n},G_{0}) \to 0, \nonumber
\end{align}
and $\widetilde{W}_{\widetilde{\kappa}_{\text{sin}}}(G_{n},G_{0}) \to 0$. To avoid unnecessary repetition, we utilize the same notation of $G_{n}$ as in the proof of Theorem~\ref{theorem:lower_bound_Gaussian_family_third} in Appendix~\ref{Section:proof_non_linear_covariate_dependent_first}. 
\paragraph{Step 1 - Structure of Taylor expansion:} Similar to the proof of Theorem~\ref{theorem:lower_bound_Gaussian_family_third}, we have the following representation when we perform Taylor expansion up to the order $\rtil_{\text{sin}}$:
\begin{align}
p_{G_{n}}(X,Y)-p_{G_{0}}(X,Y) : = A_{n} + B_{n} +R(X,Y), \nonumber
\end{align} 
where $R(X,Y)$ is a Taylor remainder such that $R(X,Y)/D_{\widetilde{\kappa}_{\text{sin}}}(G_{n},G_{0}) \to 0$ as $n \to \infty$. The forms of $B_{n}$ and $D_{\widetilde{\kappa}_{\text{sin}}}(G_{n},G_{0})$ are similar to that in Step 1 of Theorem~\ref{theorem:lower_bound_Gaussian_family_third} except that we use $\widetilde{\kappa}_{\text{sin}}$ instead of $\kappa$. Furthermore, $A_{n}$ has the following form:
\begin{align}
A_{n} & : = \sum \limits_{i=1}^{k_{0}} A_{n}(i) = \sum \limits_{i \in \mathcal{A}} A_{n}(i) + \sum \limits_{i \in \mathcal{A}^{c}} A_{n}(i), \nonumber
\end{align}
where $\mathcal{A} : = \{i \in [k_{0}]: \ (\theta_{1i}^{0})^{(1)} \neq 0 \ \text{and} \ (\theta_{1i}^{0})^{(2)} = 0 \}$ and 
\begin{align} 
A_{n}(i) & : = \sum \limits_{j=1}^{s_{i}}p_{ij}^{n}\sum \limits_{1 \leq |\alpha| \leq \rtil_{\text{sin}}} \dfrac{1}{\alpha!}\biggr\{(\Delta \theta_{1ij}^{n})^{(1)}\biggr\}^{\alpha_{1}}\biggr\{(\Delta \theta_{1ij}^{n})^{(2)}\biggr\}^{\alpha_{2}}(\Delta \theta_{2ij}^{n})^{\alpha_{3}} \nonumber \\
& \times \dfrac{\partial^{|\alpha|}{f}}{\partial{(\theta_{1}^{(1)})^{\alpha_{1}}}\partial{(\theta_{1}^{(2)})^{\alpha_{2}}}\partial{\theta_{2}^{\alpha_{3}}}}\parenth{Y|h_{1}(X, \theta_{1i}^{0}),h_{2}(X, \theta_{2i}^{0})}\overline{f}(X) \nonumber 
\end{align}
for $i \in [k_{0}]$. Under the non-linearity setting II of $G_{0}$, there exists an index $i$ such that $(\theta_{1i}^{0})^{(1)} \neq 0$ and $(\theta_{1i}^{0})^{(2)} = 0$. Therefore, we have $\abss{\mathcal{A}} \geq 1$. To analyze the structure of $A_{n}(i)$, we consider two settings of index $i$: $i \in \mathcal{A}$ and $i \in \mathcal{A}^{c}$. 
\paragraph{Index $i \in \mathcal{A}$:} For any $i \in \mathcal{A}$, the collection of full partial derivatives $\dfrac{\partial^{|\alpha|}{f}}{\partial{(\theta_{1}^{(1)})^{\alpha_{1}}}\partial{(\theta_{1}^{(2)})^{\alpha_{2}}}\partial{\theta_{2}^{\alpha_{3}}}}$\\$(Y|h_{1}(X|\theta_{1i}^{0}),h_{2}(X|\theta_{2i}^{0}))\overline{f}(X)$ up to order $\rtil_{\text{sin}} \geq 2$ is not linearly independent with respect to $X$ and $Y$. Therefore, we cannot treat $A_{n}(i)$ as a linear combination of these derivatives as long as $i \in \mathcal{A}$. Our strategy is to reduce this collection of full partial derivatives into a collection of linearly independent terms of the forms $X^{l_{1}}\dfrac{\partial^{l_{2}}{f}}{\partial{h_{1}^{l_{2}}}}(Y|h_{1}(X,\theta_{1i}^{0}),h_{2}(X,\theta_{2i}^{0}))\overline{f}(X)$ as those in the previous proofs for some $(l_{1},l_{2})$. Given that idea, we define
\begin{align}
\mathcal{F}(i) = \left\{X^{l_{1}}\dfrac{\partial^{l_{2}}{f}}{\partial{h_{1}^{l_{2}}}}(Y|h_{1}(X,\theta_{1i}^{0}),h_{2}(X,\theta_{2i}^{0}))\overline{f}(X): (l_{1},l_{2}) \in \mathcal{B}(i) \right\}, \nonumber
\end{align}
the set of all linear independent terms deriving from computing the partial derivatives of $f$ up to order $\rtil_{\text{sin}}$ with respect to $\theta_{1}$ and $\theta_{2}$. In general, the exact form of $\mathcal{B}(i)$ is very difficult to obtain. For the purpose of this proof, we only need to focus on a subset of $\mathcal{B}(i)$ in which we have a closed form. In particular, we denote a set $\mathcal{B}_{\text{sub}}$ as follows:
\begin{align}
\mathcal{B}_{\text{sub}} : = \{(2,0)\} \cup \{(0,l_{2}): \ 1 \leq l_{2} \leq 2\rtil_{\text{sin}}\}. \nonumber
\end{align}
We claim that $\mathcal{B}_{\text{sub}}$ is a subset of $\mathcal{B}(i)$ for any $i \in \mathcal{A}$. We prove this claim at the end of this proof. From now on, we assume that this claim is given.
\paragraph{Index $i \in \mathcal{A}^{c}$:} For any $i \in \mathcal{A}^{c}$, we also have the linear dependence of the set of full partial derivatives $\dfrac{\partial^{|\alpha|}{f}}{\partial{(\theta_{1}^{(1)})^{\alpha_{1}}}\partial{(\theta_{1}^{(2)})^{\alpha_{2}}}\partial{\theta_{2}^{\alpha_{3}}}}(Y|h_{1}(X|\theta_{1i}^{0}),h_{2}(X|\theta_{2i}^{0}))\overline{f}(X)$ up to order $\rtil_{\text{sin}} \geq 2$. Similar to the strategy of case $i \in \mathcal{A}$, we also reduce the previous set into a collection of linearly independent terms of the forms $X^{l_{1}}\dfrac{\partial^{l_{2}}{f}}{\partial{h_{1}^{l_{2}}}}(Y|h_{1}(X,\theta_{1i}^{0}),h_{2}(X,\theta_{2i}^{0}))\overline{f}(X)$, which can be defined as:
\begin{align}
\overline{\mathcal{F}}(i) = \left\{X^{l_{1}}\dfrac{\partial^{l_{2}}{f}}{\partial{h_{1}^{l_{2}}}}(Y|h_{1}(X,\theta_{1i}^{0}),h_{2}(X,\theta_{2i}^{0}))\overline{f}(X): (l_{1},l_{2}) \in \overline{\mathcal{B}}(i) \right\}. \nonumber
\end{align}
However, the structure of $\overline{\mathcal{B}}(i)$ is also very complicated. For the purpose of this proof, we only consider its subset $\overline{\mathcal{B}}_{\text{sub}}$, which has the following form:
\begin{align}
\overline{\mathcal{B}}_{\text{sub}} : = \left\{(0,1), (0,4) , (1,2), (2,1), (2,2), (2,3), (3,2), (4,2) \right\}. \nonumber
\end{align}
The proof for the claim that $\overline{\mathcal{B}}_{\text{sub}} \subset \overline{\mathcal{B}}(i)$ for any $i \in \mathcal{A}^{c}$ is similar to that from claim $\mathcal{B}_{\text{sub}} \subset \mathcal{B}(i)$ as $i \in \mathcal{A}$; therefore, it is omitted. From now on, we also assume that the above claim is true.

Given the formulations of $\mathcal{F}(i)$ and $\overline{\mathcal{F}}(i)$, we can treat $A_{n}(i)/D_{\widetilde{\kappa}_{\text{sin}}}(G_{n},G_{0})$ as the linear combinations of elements from $\mathcal{F}(i)$ divided by $D_{\widetilde{\kappa}_{\text{sin}}}(G_{n},G_{0})$ for $i \in \mathcal{A}$ and from $\overline{\mathcal{F}}(i)$ divided by $D_{\widetilde{\kappa}_{\text{sin}}}(G_{n},G_{0})$ for $i \in \mathcal{A}^{c}$. 
\paragraph{Non-vanishing coefficients:} Similar to the previous proofs, we assume that all the coefficients in the representation of $A_{n}(i)/D_{\widetilde{\kappa}_{\text{sin}}}(G_{n},G_{0})$ and $B_{n}/D_{\widetilde{\kappa}_{\text{sin}}}(G_{n},G_{0})$ go to 0 as $n \to \infty$ for all $i \in [k_{0}]$. From the definitions of $\mathcal{B}_{\text{sub}}$ and $\overline{\mathcal{B}}_{\text{sub}}$, we have the coefficients associated with $X^{l_{1}}\dfrac{\partial^{l_{2}}{f}}{\partial{h_{1}^{l_{2}}}}(Y|h_{1}(X,\theta_{1i}^{0}),h_{2}(X,\theta_{2i}^{0}))\overline{f}(X)$ go to 0 when $(l_{1},l_{2}) \in \mathcal{B}_{\text{sub}}$ for $i \in \mathcal{A}$ or $(l_{1},l_{2}) \in \overline{\mathcal{B}}_{\text{sub}}$ for $i \in \mathcal{A}^{c}$. 

For the simplicity of the presentation, we denote $E_{(l_{1},l_{2})}(i)$ the coefficients of the element $X^{l_{1}}\dfrac{\partial^{l_{2}}{f}}{\partial{h_{1}^{l_{2}}}}(Y|h_{1}(X,\theta_{1i}^{0}),h_{2}(X,\theta_{2i}^{0}))\overline{f}(X)$ when $(l_{1},l_{2}) \in \mathcal{B}_{\text{sub}}$ and $i \in \mathcal{A}$. Similarly, $\overline{E}_{(l_{1},l_{2})}(i)$ are the coefficients of $X^{l_{1}}\dfrac{\partial^{l_{2}}{f}}{\partial{h_{1}^{l_{2}}}}(Y|h_{1}(X,\theta_{1i}^{0}),h_{2}(X,\theta_{2i}^{0}))\overline{f}(X)$ when $(l_{1},l_{2}) \in \overline{\mathcal{B}}_{\text{sub}}$ and $i \in \mathcal{A}^{c}$.

For $(l_{1},l_{2}) = (0,l)$ as $1 \leq l \leq 2\rtil_{\text{sin}}$, the exact formulation of $E_{(l_{1},l_{2})}(i)$ can be derived from determining the coefficient of $X^{l_{1}}\dfrac{\partial^{l_{2}}{f}}{\partial{h_{1}^{l_{2}}}}(Y|h_{1}(X,\theta_{1i}^{0}),h_{2}(X,\theta_{2i}^{0}))\overline{f}(X)$ in the following term:
\begin{align}
\sum \limits_{j=1}^{s_{i}}p_{ij}^{n}\sum \limits_{1 \leq |\gamma| \leq \rtil_{\text{sin}}} \dfrac{1}{\gamma!}\biggr\{(\Delta \theta_{1ij}^{n})^{(1)}\biggr\}^{\gamma_{1}}(\Delta \theta_{2ij}^{n})^{\gamma_{2}} \dfrac{\partial^{|\gamma|}{f}}{\partial{(\theta_{1}^{(1)})^{\gamma_{1}}}\partial{\theta_{2}^{\gamma_{2}}}}\parenth{Y|h_{1}(X, \theta_{1i}^{0}),h_{2}(X, \theta_{2i}^{0})}\overline{f}(X), \nonumber
\end{align}
for $\gamma = (\gamma_{1},\gamma_{2})$. Equipped with the result of Lemma~\ref{lemma:representation_partial_derivative}, we can verify that
\begin{align}
E_{(0,l)}(i) = \brackets{\sum \limits_{\gamma_{1},\gamma_{2},\tau} \dfrac{P_{\tau}^{(\gamma_{1})}\parenth{(\theta_{1i}^{0})^{(1)}}}{2^{\gamma_{2}}}\biggr(\sum \limits_{j=1}^{s_{i}} p_{ij}^{n}\dfrac{\biggr\{(\Delta \theta_{1ij}^{n})^{(1)}\biggr\}^{\gamma_{1}}(\Delta \theta_{2ij}^{n})^{\gamma_{2}}}{\gamma_{1}!\gamma_{2}!} \biggr)}\bigg/D_{\widetilde{\kappa}_{\text{sin}}}(G_{n},G_{0}), \label{eqn:formulation_proof_non_linear_first}
\end{align}
where the summation with respect to $\gamma_{1},\gamma_{2},\tau$ in the numerator satisfies $\gamma_{1}/2+\tau+2\gamma_{2}=l$, $\tau \leq \gamma_{1}/2$ when $\gamma_{1}$ is an even number while $(\gamma_{1}+1)/2+\tau+2\gamma_{2} = l$, $\tau \leq (\gamma_{1}-1)/2$ when $\gamma_{1}$ is an odd number. Furthermore, $\gamma_{1}+\gamma_{2} \leq \rtil_{\text{sin}}$. 

By taking the summation of the absolute value of coefficients in $B_{n}/D_{\kappa}(G_{n},G_{0})$, it implies that
\begin{align}
\parenth{\sum \limits_{i=1}^{k_{0}} |\sum \limits_{j=1}^{s_{i}} p_{ij}^{n} - \pi_{i}^{0}|}/D_{\widetilde{\kappa}_{\text{sin}}}(G_{n},G_{0}) \to 0. \nonumber
\end{align}
From the definition of $D_{\widetilde{\kappa}_{\text{sin}}}(G_{n},G_{0})$, it leads to
\begin{align}
\brackets{\sum \limits_{i=1}^{k_{0}}\sum \limits_{j=1}^{s_{i}}{p_{ij}^{n}\biggr(\biggr|(\Delta \theta_{1ij}^{n})^{(1)}\biggr|^{\rtil_{\text{sin}}}+\biggr|(\Delta \theta_{1ij}^{n})^{(2)}\biggr|^{2}+\biggr|\Delta \theta_{2ij}^{n}\biggr|^{\ceil{\rtil_{\text{sin}}/2}} \biggr)}}\bigg/D_{\widetilde{\kappa}_{\text{sin}}}(G_{n},G_{0}) \to 1. \nonumber
\end{align}
Therefore, there exists an index $i^{*} \in [k_{0}]$ such that 
\begin{align}
\brackets{\sum \limits_{j=1}^{s_{i^{*}}}{p_{i^{*}j}^{n}\biggr(\biggr|(\Delta \theta_{1i^{*}j}^{n})^{(1)}\biggr|^{\rtil_{\text{sin}}}+\biggr|(\Delta \theta_{1i^{*}j}^{n})^{(2)}\biggr|^{2}+\biggr|\Delta \theta_{2i^{*}j}^{n}\biggr|^{\ceil{\rtil_{\text{sin}}/2}} \biggr)}}\bigg/D_{\widetilde{\kappa}_{\text{sin}}}(G_{n},G_{0}) \not \to 0. \nonumber
\end{align}
We denote
\begin{align}
\overline{D}_{\widetilde{\kappa}_{\text{sin}}}(G_{n},G_{0}) := \sum \limits_{j=1}^{s_{i^{*}}}{p_{i^{*}j}^{n}\biggr(\biggr|(\Delta \theta_{1i^{*}j}^{n})^{(1)}\biggr|^{\rtil_{\text{sin}}}+\biggr|(\Delta \theta_{1i^{*}j}^{n})^{(2)}\biggr|^{2}+\biggr|\Delta \theta_{2i^{*}j}^{n}\biggr|^{\ceil{\rtil_{\text{sin}}/2}} \biggr)}. \nonumber
\end{align}
As $E_{(l_{1},l_{2})}(i) \to 0$ and $\overline{E}_{(l_{1}',l_{2}')}(j) \to 0$ for $i \in \mathcal{A}$, $j \in \mathcal{A}^{c}$, $(l_{1},l_{2}) \in \mathcal{B}_{\text{sub}}$, and $(l_{1}',l_{2}') \in \overline{\mathcal{B}}_{\text{sub}}$, the following holds:
\begin{align}
K_{(l_{1},l_{2})}(i) : = \frac{D_{\widetilde{\kappa}_{\text{sin}}}(G_{n},G_{0})}{\overline{D}_{\widetilde{\kappa}_{\text{sin}}}(G_{n},G_{0})} E_{(l_{1},l_{2})}(i) \to 0, \nonumber \\
\overline{K}_{(l_{1}',l_{2}')}(j) : = \frac{D_{\widetilde{\kappa}_{\text{sin}}}(G_{n},G_{0})}{\overline{D}_{\widetilde{\kappa}_{\text{sin}}}(G_{n},G_{0})} \overline{E}_{(l_{1}',l_{2}')}(j) \to 0, \nonumber
\end{align}
for all $i \in \mathcal{A}$, $j \in \mathcal{A}^{c}$, $(l_{1},l_{2}) \in \mathcal{B}_{\text{sub}}$, and $(l_{1}',l_{2}') \in \overline{\mathcal{B}}_{\text{sub}}$. Now, we consider two possible settings of $i^{*}$.
\paragraph{Setting 1 - $i^{*} \in \mathcal{A}$:} By direct computation, the vanishing of $K_{(2,0)}(i^{*})$ to 0 is equivalent to
\begin{align}
\parenth{\sum \limits_{j=1}^{s_{i^{*}}}{p_{i^{*}j}^{n} \biggr|(\Delta \theta_{1i^{*}j}^{n})^{(2)}\biggr|^{2}}} \bigg/\overline{D}_{\widetilde{\kappa}_{\text{sin}}} \to 0. \nonumber
\end{align}
From the definition of $\widetilde{\kappa}_{\text{sin}}$, the above result leads to
\begin{align}
L_{n} : = \sum \limits_{j=1}^{s_{i^{*}}}{p_{i^{*}j}^{n}\biggr(\biggr|(\Delta \theta_{1i^{*}j}^{n})^{(1)}\biggr|^{\rtil_{\text{sin}}}+\biggr|\Delta \theta_{2i^{*}j}^{n}\biggr|^{\ceil{\rtil_{\text{sin}}/2}} \biggr)}\bigg/\overline{D}_{\widetilde{\kappa}_{\text{sin}}} \to 1. \nonumber
\end{align}
Equipped with the formulation of $E_{(0,l)}(i^{*})$ in equation~\eqref{eqn:formulation_proof_non_linear_first} for any $1 \leq l \leq 2\widetilde{r}_{\text{sin}}$, the following system of limits holds:
\begin{align}
\dfrac{1}{L_{n}}K_{(0,l)}(i^{*}) = \dfrac{\sum \limits_{\gamma_{1},\gamma_{2},\tau} \dfrac{P_{\tau}^{(\gamma_{1})}\parenth{(\theta_{1i^{*}}^{0})^{(1)}}}{2^{\gamma_{2}}}\biggr(\sum \limits_{j=1}^{s_{i}} p_{i^{*}j}^{n}\dfrac{\biggr\{(\Delta \theta_{1i^{*}j}^{n})^{(1)}\biggr\}^{\gamma_{1}}(\Delta \theta_{2i^{*}j}^{n})^{\gamma_{2}}}{\gamma_{1}!\gamma_{2}!} \biggr)}{\sum \limits_{j=1}^{s_{i^{*}}}{p_{i^{*}j}^{n}\biggr(\biggr|(\Delta \theta_{1i^{*}j}^{n})^{(1)}\biggr|^{\rtil_{\text{sin}}}+\biggr|\Delta \theta_{2i^{*}j}^{n}\biggr|^{\ceil{\rtil_{\text{sin}}/2}} \biggr)}} \to 0, \label{eqn:proof_non_linear_system}
\end{align}
where the summation with respect to $\gamma_{1},\gamma_{2},\tau$ in the numerator satisfies $\gamma_{1}/2+\tau+2\gamma_{2}=l$, $\tau \leq \gamma_{1}/2$ when $\gamma_{1}$ is an even number while $(\gamma_{1}+1)/2+\tau+2\gamma_{2} = l$, $\tau \leq (\gamma_{1}-1)/2$ when $\gamma_{1}$ is an odd number. Additionally, $\gamma_{1}+\gamma_{2} \leq \rtil_{\text{sin}}$. 

Recall that $\rtil_{\text{sin}} = \rtil((\theta_{1i_{\text{max}}}^{0})^{(1)},k-k_{0}+1)$ where $i_{\text{max}} = \mathop{ \arg \max} \limits_{i \in \mathcal{A}} \rtil((\theta_{1i}^{0})^{(1)},k-k_{0}+1)$. Therefore, $\rtil_{\text{sin}} \geq \rtil((\theta_{1i^{*}}^{0})^{(1)},k-k_{0}+1) \geq \rtil((\theta_{1i^{*}}^{0})^{(1)},s_{i^{*}})$ as $s_{i^{*}} \leq k - k_{0} + 1$. From the definition of $\rtil((\theta_{1i^{*}}^{0})^{(1)},s_{i^{*}})$ in Definition~\ref{definition:smallest_number}, the system of polynomial limit~\eqref{eqn:proof_non_linear_system} does not hold given the values of $\rtil((\theta_{1i^{*}}^{0})^{(1)},s_{i^{*}})$. Therefore, it does not happen under $\rtil_{\text{sin}}$. As a consequence, setting 1 that $i^{*} \in \mathcal{A}$ will not hold.
\paragraph{Setting 2 - $i^{*} \in \mathcal{A}^{c}$:} Since $\rtil_{\text{sin}} \geq 3$, it is clear that $(2,2,2) \prec (\rtil_{\text{sin}},2,\ceil{\rtil_{\text{sin}}/2})$. It implies that 
\begin{align}
\overline{D}_{\widetilde{\kappa}_{\text{sin}}}(G_{n},G_{0}) \lesssim \widetilde{D}(G_{n},G_{0}) : = \sum \limits_{j=1}^{s_{i^{*}}}{p_{i^{*}j}^{n}\biggr(\biggr|(\Delta \theta_{1i^{*}j}^{n})^{(1)}\biggr|^{2}+\biggr|(\Delta \theta_{1i^{*}j}^{n})^{(2)}\biggr|^{2}+\biggr|\Delta \theta_{2i^{*}j}^{n}\biggr|^{2} \biggr)}. \nonumber
\end{align}
Since $\overline{E}_{(l_{1},l_{2})}(i^{*}) \to 0$ for all $(l_{1},l_{2}) \in \overline{\mathcal{B}}_{\text{sub}}$, it leads to
\begin{align}
\overline{F}_{(l_{1},l_{2})}(i^{*}) : = \frac{\overline{D}_{\widetilde{\kappa}_{\text{sin}}}(G_{n},G_{0})}{\widetilde{D}(G_{n},G_{0})}\overline{E}_{(l_{1},l_{2})}(i^{*}) \to 0, \nonumber
\end{align} 
for all $(l_{1},l_{2}) \in \overline{\mathcal{B}}_{\text{sub}}$. We can check that the vanishing of $\overline{F}_{(0,4)}(i^{*})$ to 0 leads to
\begin{align}
\parenth{\sum \limits_{j=1}^{s_{i^{*}}} p_{i^{*}j}^{n} \biggr|(\Delta \theta_{1i^{*}j}^{n})^{(2)}\biggr|^{2}}\bigg/\widetilde{D}(G_{n},G_{0}) \to 0. \label{eqn:proof_non_linear_second_case_first_limit}
\end{align}
Furthermore, the vanishings of $\overline{F}_{(l_{1},l_{2})}(i^{*})$ to 0 for $(l_{1},l_{2}) \in \{(0,1),(1,2), (2,1), (2,2), (3,2), (4,2)\}$ lead to the system of polynomial limits similar to system of polynomial limits~\eqref{eqn:simple_system_polynomial_limits} where the index $i$ in this system is replaced by $i^{*}$ and the distance $D_{\kappa}(G_{n},G_{0})$ is replaced by $\widetilde{D}(G_{n},G_{0})$. Due to the fact that $i^{*} \in \mathcal{A}^{c}$, following the argument after system of limits~\eqref{eqn:simple_system_polynomial_limits}, we obtain that
\begin{align}
\sum \limits_{j=1}^{s_{i^{*}}}{p_{i^{*}j}^{n}\biggr(\biggr|(\Delta \theta_{1i^{*}j}^{n})^{(1)}\biggr|^{2}+\biggr|\Delta \theta_{2i^{*}j}^{n}\biggr|^{2} \biggr)}\bigg/\widetilde{D}(G_{n},G_{0}) \to 0. \label{eqn:proof_non_linear_second_case_second_limit}
\end{align}
Invoking the results from equations~\eqref{eqn:proof_non_linear_second_case_first_limit} and~\eqref{eqn:proof_non_linear_second_case_second_limit} leads to
\begin{align}
1 = \sum \limits_{j=1}^{s_{i^{*}}}{p_{i^{*}j}^{n}\biggr(\biggr|(\Delta \theta_{1i^{*}j}^{n})^{(1)}\biggr|^{2}+\biggr|(\Delta \theta_{1i^{*}j}^{n})^{(2)}\biggr|^{2}+\biggr|\Delta \theta_{2i^{*}j}^{n}\biggr|^{2} \biggr)}\bigg/\widetilde{D}(G_{n},G_{0}) \to 0, \nonumber
\end{align}
which is a contradiction. Therefore, setting 2 that $i^{*} \in \mathcal{A}^{c}$ will not hold. 

As a consequence, not all the coefficients of $A_{n}(i)/D_{\widetilde{\kappa}_{\text{sin}}}(G_{n},G_{0})$ and $B_{n}/D_{\widetilde{\kappa}_{\text{sin}}}(G_{n},G_{0})$ go to 0 as $n \to \infty$ for all $i \in [k_{0}]$. From here, by means of the Fatou's argument as that of the previous proofs, we achieve the conclusion regarding the convergence rate of MLE under non-linearity setting II of $G_{0}$.
\paragraph{Proof of claim $\mathcal{B}_{\text{sub}} \subset \mathcal{B}(i)$ for any $i \in \mathcal{A}$:} First of all, we demonstrate that the elements $X^{l_{1}}\dfrac{\partial^{l_{2}}{f}}{\partial{h_{1}^{l_{2}}}}(Y|h_{1}(X,\theta_{1i}^{0}),h_{2}(X,\theta_{2i}^{0}))\overline{f}(X)$ where $(l_{1},l_{2}) \in \mathcal{B}_{\text{sub}}$ are originated from some partial derivatives of $f$ with respect to $\theta_{1}$ and $\theta_{2}$. In fact, by means of Lemma~\ref{lemma:representation_partial_derivative}, the pairs of indices $(l_{1},l_{2}) = (0,l) \in \mathcal{B}_{\text{sub}}$ for $1 \leq l \leq 2\rtil_{\text{sin}}$ correspond to the elements coming from the partial derivatives $\dfrac{\partial^{|\gamma|}f}{\partial{(\theta_{1}^{(1)})^{\gamma_{1}}\partial{\theta_{2}^{\gamma_{2}}}}}(Y|h_{1}(X,\theta_{1}^{0}),h_{2}(X,\theta_{2}^{0}))$ for $1 \leq \abss{\gamma} \leq \rtil_{\text{sin}}$. Additionally, the pair $(2,0) \in \mathcal{B}_{\text{sub}}$ is associated with element from the derivation of $\dfrac{\partial^{2}f}{\partial{(\theta_{1}^{(2)}})^{2}}(Y|h_{1}(X,\theta_{1}^{0}),h_{2}(X,\theta_{2}^{0}))$. 

Furthermore, it is not hard to verify that the collection of $X^{l_{1}}\dfrac{\partial^{l_{2}}{f}}{\partial{h_{1}^{l_{2}}}}(Y|h_{1}(X,\theta_{1i}^{0}),h_{2}(X,\theta_{2i}^{0}))$\\$\overline{f}(X)$ for $(l_{1},l_{2}) \in \mathcal{B}_{\text{sub}}$ is linearly independent with respect to $X$ and $Y$. Therefore, we achieve the conclusion that $\mathcal{B}_{\text{sub}} \subset \mathcal{B}(i)$. 
\subsection{Proof of Theorem~\ref{theorem:lower_bound_Gaussian_family_without_offset_O2}}
\label{subsec:proof:theorem:lower_bound_Gaussian_family_without_offset_O2}
Similar to the previous proofs, it is sufficient to demonstrate the following results:
\begin{eqnarray}
\lim \limits_{\epsilon \to 0} \inf \limits_{G \in \Ocal_{k,\overline{c}_{0}} (\Omega): \widetilde{W}_{\kappa}(G,G_{0}) \leq \epsilon} V(p_{G},p_{G_{0}})/\widetilde{W}_{\kappa}^{\|\kappa\|_{\infty}}(G,G_{0}) > 0, \label{eqn:proof_Gaussian_family_without_offset_O2_first} \\
\inf \limits_{G \in \Ocal_{k} (\Omega)} h(p_{G},p_{G_{0}})/\widetilde{W}_{\kappa'}^{\|\kappa'\|_{\infty}}(G,G_{0}) = 0, \label{eqn:proof_Gaussian_family_without_offset_O2_second}
\end{eqnarray}
for any $\kappa' \prec \kappa$ where $\kappa=(2,2,2,2)$. Proof of inequality~\eqref{eqn:proof_Gaussian_family_without_offset_O2_first} is in Appendix~\ref{subsec:proof_Gaussian_family_without_offset_O2_first} while proof of equality~\eqref{eqn:proof_Gaussian_family_without_offset_O2_second} is in Appendix~\ref{subsec:proof_Gaussian_family_without_offset_O2_second}.
\subsubsection{Proof for inequality~\eqref{eqn:proof_Gaussian_family_without_offset_O2_first}}
\label{subsec:proof_Gaussian_family_without_offset_O2_first}
We assume that the conclusion of inequality~\eqref{eqn:proof_Gaussian_family_without_offset_O2_first} does not hold. It indicates that we can find a sequence $G_{n}$ that has representation~\eqref{eqn:proof_Gaussian_firstcase_notation} such that $V(p_{G_{n}},p_{G_{0}})/\widetilde{W}_{\kappa}^{\|\kappa\|_{\infty}}(G_{n},G_{0}) \to 0$ and $\widetilde{W}_{\kappa}(G_{n},G_{0}) \to 0$. As $\theta_{2ij}^{n}$ has three dimensions, in this proof we denote $\Delta \theta_{2ij}^{n} = ((\Delta \theta_{2ij}^{n})^{(1)}, (\Delta \theta_{2ij}^{n})^{(2)}, (\Delta \theta_{2ij}^{n})^{(3)})$ for all $1 \leq i \leq k_{0}$ and $1 \leq j \leq s_{i}$. According to Lemma~\ref{lemma:generalized_Wasserstein_distance_polynomials}, we have
\begin{eqnarray}
\widetilde{W}_{\kappa}^{\|\kappa\|_{\infty}}(G_{n},G_{0}) \precsim \sum \limits_{i=1}^{k_{0}}\sum \limits_{j=1}^{s_{i}}{p_{ij}^{n}\biggr(\biggr|\Delta \theta_{1ij}^{n}\biggr|^{2}+\biggr|(\Delta \theta_{2ij}^{n})^{(1)}\biggr|^{2}+\biggr|(\Delta \theta_{2ij}^{n})^{(2)}\biggr|^{2} + \biggr|(\Delta \theta_{2ij}^{n})^{(3)}\biggr|^{2}\biggr)} \nonumber \\
+ \sum \limits_{i=1}^{k_{0}} |\sum \limits_{j=1}^{s_{i}}p_{ij}^{n} - \pi_{i}^{0}| := D_{\kappa}(G_{n},G_{0}). \nonumber
\end{eqnarray}
By means of Taylor expansion up to the second order, we have
\begin{eqnarray}
p_{G_{n}}(X,Y)-p_{G_{0}}(X,Y) & = & \nonumber \\
& & \hspace{- 5 em} \sum \limits_{i=1}^{k_{0}}\sum \limits_{j=1}^{s_{i}}p_{ij}^{n}\sum \limits_{1 \leq |\alpha| \leq 2} \dfrac{1}{\alpha!}\biggr\{\Delta \theta_{1ij}^{n}\biggr\}^{\alpha_{1}}\biggr\{(\Delta \theta_{2ij}^{n})^{(1)}\biggr\}^{\alpha_{2}}\biggr\{(\Delta \theta_{2ij}^{n})^{(2)}\biggr\}^{\alpha_{3}} \biggr\{(\Delta \theta_{2ij}^{n})^{(3)}\biggr\}^{\alpha_{4}} \nonumber \\
& & \hspace{- 5 em}  \times \dfrac{\partial^{|\alpha|}{f}}{\partial{\theta_{1}^{\alpha_{1}}}\partial{(\theta_{2}^{(1)})^{\alpha_{2}}}\partial{(\theta_{2}^{(2)})^{\alpha_{3}}}\partial{(\theta_{2}^{(3)})^{\alpha_{4}}}}\parenth{Y|h_{1}(X,\theta_{1i}^{0}),h_{2}(X,\theta_{2i}^{0})}\overline{f}(X) \nonumber \\
& & \hspace{- 5 em} + \sum \limits_{i=1}^{k_{0}}\biggr(\sum \limits_{j=1}^{s_{i}}p_{ij}^{n} - \pi_{i}^{0}\biggr)f(Y|h_{1}(X,\theta_{1i}^{0}),h_{2}(X,\theta_{2i}^{0}))\overline{f}(X) + R(X,Y)
\nonumber \\
& & \hspace{- 5 em} : = A_{n} + B_{n} + R(X,Y), \nonumber
\end{eqnarray}
where the Taylor remainder $R(X,Y)$ is such that $R(X,Y)/D_{\kappa}(G_{n},G_{0}) \to 0$ as $n \to \infty$. From the formulations of expert functions $h_{1}, h_{2}$, we find that
\begin{eqnarray}
\dfrac{\partial^{|\alpha|}{f}}{\partial{\theta_{1}^{\alpha_{1}}}\partial{(\theta_{2}^{(1)})^{\alpha_{2}}}\partial{(\theta_{2}^{(2)})^{\alpha_{3}}}\partial{(\theta_{2}^{(3)})^{\alpha_{4}}}}\parenth{Y|h_{1}(X,\theta_{1}),h_{2}(X,\theta_{2})} & = & \nonumber \\
& & \hspace{- 14 em} \dfrac{X^{2\alpha_{1}+\alpha_{3} + 2\alpha_{4}}}{2^{\alpha_{2} + \alpha_{3} + \alpha_{4}}}\dfrac{\partial^{\alpha_{1}+2\alpha_{2}+2\alpha_{3} + 2\alpha_{4}}{f}}{\partial{h_{1}^{\alpha_{1}+2\alpha_{2}+2\alpha_{3} + 2\alpha_{4}}}}(Y|h_{1}(X|\theta_{1}),h_{2}(X|\theta_{2})), \nonumber
\end{eqnarray}
for any $\alpha = (\alpha_{1}, \alpha_{2}, \alpha_{3}, \alpha_{4}) \in \mathbb{N}^{4}$. Based on the above equation, we can express $A_{n}$ as follows:
\begin{eqnarray}
A_{n} & = & \sum_{i = 1}^{k_{0}} \sum \limits_{j=1}^{s_{i}}p_{ij}^{n}\sum \limits_{1 \leq |\alpha| \leq 2} \dfrac{1}{\alpha!}\biggr\{\Delta \theta_{1ij}^{n}\biggr\}^{\alpha_{1}}\biggr\{(\Delta \theta_{2ij}^{n})^{(1)}\biggr\}^{\alpha_{2}}\biggr\{(\Delta \theta_{2ij}^{n})^{(2)}\biggr\}^{\alpha_{3}} \biggr\{(\Delta \theta_{2ij}^{n})^{(3)}\biggr\}^{\alpha_{4}} \nonumber \\ \nonumber \\
& & \hspace{4 em} \times  \dfrac{X^{2\alpha_{1}+\alpha_{3} + 2\alpha_{4}}}{2^{\alpha_{2} + \alpha_{3} + \alpha_{4}}}\dfrac{\partial^{\alpha_{1}+2\alpha_{2}+2\alpha_{3} + 2\alpha_{4}}{f}}{\partial{h_{1}^{\alpha_{1}+2\alpha_{2}+2\alpha_{3} + 2\alpha_{4}}}}(Y|h_{1}(X|\theta_{1i}^{0}),h_{2}(X|\theta_{2i}^{0}))\overline{f}(X). \label{eqn:proof_Gaussian_family_without_offset_O2_third}
\end{eqnarray}
If we define
\begin{eqnarray}
\mathcal{F} = \biggr\{X^{l_{1}}\dfrac{\partial^{l_{2}}{f}}{\partial{h_{1}^{l_{2}}}}(Y|h_{1}(X|\theta_{1i}^{0}),h_{2}(X|\theta_{2i}^{0}))\overline{f}(X): & & \nonumber \\
& & \hspace{-15 em} l_{1} = 2 \alpha_{1} + \alpha_{3} + 2 \alpha_{4}, \ l_{2} = \alpha_{1} + 2 \alpha_{2} + 2\alpha_{3} + 2\alpha_{4}, \ 0 \leq |\alpha| \leq 2, \ 1 \leq i \leq k_{0} \biggr\}, \nonumber
\end{eqnarray}
then the elements of $\mathcal{F}$ are linearly independent with respect to $X$ and $Y$. The proof argument of this claim is similar to that in equation~\eqref{eqn:lemma_linear_independence_Gaussian_firstcase_first} in Section~\ref{Section:proof_weakly_identifiable_covariate_independent}. Therefore, we can treat $A_{n}/D_{\kappa}(G_{n},G_{0})$, $B_{n}/D_{\kappa}(G_{n},G_{0})$ as a linear combination of elements of $\mathcal{F}$. 

Similar to the proof of Theorem~\ref{theorem:lower_bound_Gaussian_family} in Section~\ref{Section:proof_weakly_identifiable_covariate_independent}, we denote $F_{l_{1},l_{2}}(\theta_{1i}^{0},\theta_{2i}^{0})$ as the coefficient of $X^{l_{1}}\dfrac{\partial^{l_{2}}{f}}{\partial{h_{1}^{l_{2}}}}(Y|h_{1}(X|\theta_{1i}^{0}),h_{2}(X|\theta_{2i}^{0}))\overline{f}(X)$ in $A_{n}$ and $B_{n}$ for any $l_{1} = 2 \alpha_{1} + \alpha_{3} + 2 \alpha_{4}, \ l_{2} = \alpha_{1} + 2 \alpha_{2} + 2\alpha_{3} + 2\alpha_{4}, \ 0 \leq |\alpha| \leq 2$ and $1 \leq i \leq k_{0}$. Then, we can check that the coefficients of $X^{l_{1}}\dfrac{\partial^{l_{2}}{f}}{\partial{h_{1}^{l_{2}}}}(Y|h_{1}(X|\theta_{1i}^{0}),h_{2}(X|\theta_{2i}^{0}))\overline{f}(X)$ in $A_{n}/D_{\kappa}(G_{n},G_{0})$ and $B_{n}/D_{\kappa}(G_{n},G_{0})$ will be $F_{l_{1},l_{2}}(\theta_{1i}^{0},\theta_{2i}^{0})/D_{\kappa}(G_{n},G_{0})$. 

Assume that all of these coefficients go to 0 as $n \to \infty$. By taking the summation of $|F_{0,0}(\theta_{1i}^{0},\theta_{2i}^{0})/D_{\kappa}(G_{n},G_{0})|$ for all $1 \leq i \leq k_{0}$, we obtain that
\begin{eqnarray}
\biggr(\sum \limits_{i=1}^{k_{0}} |\sum \limits_{j=1}^{s_{i}}{p_{ij}^{n}} - \pi_{i}^{0}|\biggr)/D_{\kappa}(G_{n},G_{0}) \to 0. \label{eqn:proof_Gaussian_family_without_offset_O2_fourth}
\end{eqnarray}
When $l_{1} = 4$ and $l_{2} = 2$, the only $\alpha$ that satisfies these equations is $\alpha = (2, 0, 0, 0)$. It indicates that the summation of $|F_{4,2}(\theta_{1i}^{0},\theta_{2i}^{0})/D_{\kappa}(G_{n},G_{0})|$ for all $1 \leq i \leq k_{0}$ leads to
\begin{align}
    \biggr(\sum \limits_{i=1}^{k_{0}}\sum \limits_{j=1}^{s_{i}}{p_{ij}^{n} \biggr|\Delta \theta_{1ij}^{n}\biggr|^{2}}\biggr)/D_{\kappa}(G_{n},G_{0}) \to 0. \label{eqn:proof_Gaussian_family_without_offset_O2_fifth}
\end{align}
When $l_{1} = 0$ and $l_{2} = 4$, only $\alpha = (0, 2, 0, 0)$ satisfies these equations. By summing all the coefficients $|F_{0,4}(\theta_{1i}^{0},\theta_{2i}^{0})/D_{\kappa}(G_{n},G_{0})|$ for all $1 \leq i \leq k_{0}$, we find that
\begin{align}
    \biggr(\sum \limits_{i=1}^{k_{0}}\sum \limits_{j=1}^{s_{i}}{p_{ij}^{n} \biggr|(\Delta \theta_{2ij}^{n})^{(1)}\biggr|^{2}}\biggr)/D_{\kappa}(G_{n},G_{0}) \to 0. \label{eqn:proof_Gaussian_family_without_offset_O2_sixth}
\end{align}
Similarly, when $(l_{1}, l_{2}) \equiv (2,4)$, we have $\alpha = (0, 0, 2, 0)$ or when $(l_{1}, l_{2}) \equiv (4, 4)$, we have $\alpha = (0, 0, 0, 2)$. By considering the summation of the coefficients of $|F_{2,4}(\theta_{1i}^{0},\theta_{2i}^{0})/D_{\kappa}(G_{n},G_{0})|$ or $|F_{4,4}(\theta_{1i}^{0},\theta_{2i}^{0})/D_{\kappa}(G_{n},G_{0})|$ for all $1 \leq i \leq k_{0}$, we arrive at
\begin{align}
    \biggr(\sum \limits_{i=1}^{k_{0}}\sum \limits_{j=1}^{s_{i}}{p_{ij}^{n} \biggr|(\Delta \theta_{2ij}^{n})^{(2)}\biggr|^{2}}\biggr)/D_{\kappa}(G_{n},G_{0}) \to 0, \nonumber \\
    \biggr(\sum \limits_{i=1}^{k_{0}}\sum \limits_{j=1}^{s_{i}}{p_{ij}^{n} \biggr|(\Delta \theta_{2ij}^{n})^{(3)}\biggr|^{2}}\biggr)/D_{\kappa}(G_{n},G_{0}) \to 0. \label{eqn:proof_Gaussian_family_without_offset_O2_seventh}
\end{align}
Combining the results from equations~\eqref{eqn:proof_Gaussian_family_without_offset_O2_fourth}-\eqref{eqn:proof_Gaussian_family_without_offset_O2_seventh} leads to
\begin{align*}
    1 = \frac{\sum \limits_{i=1}^{k_{0}}\sum \limits_{j=1}^{s_{i}}{p_{ij}^{n}\biggr(\biggr|\Delta \theta_{1ij}^{n}\biggr|^{2}+\biggr|(\Delta \theta_{2ij}^{n})^{(1)}\biggr|^{2}+\biggr|(\Delta \theta_{2ij}^{n})^{(2)}\biggr|^{2} + \biggr|(\Delta \theta_{2ij}^{n})^{(3)}\biggr|^{2}\biggr)}
+ \sum \limits_{i=1}^{k_{0}} |\sum \limits_{j=1}^{s_{i}}p_{ij}^{n} - \pi_{i}^{0}|}{D_{\kappa}(G_{n},G_{0})} \to 0,
\end{align*}
which is a contradiction. As a consequence, not all the coefficients in the linear combinations of $A_{n}/D_{\kappa}(G_{n},G_{0})$ and $B_{n}/D_{\kappa}(G_{n},G_{0})$ go to 0 as $n \to \infty$. From here, using the Fatou's argument in Step 4 of the proof of Theorem~\ref{theorem:total_variation_bound_over-fitted_MECFG} and the fact that the elements of $\mathcal{F}$ are linearly independent with respect to $X$ and $Y$, we achieve the conclusion of claim~\eqref{eqn:proof_Gaussian_family_without_offset_O2_first}.
\subsubsection{Proof for equality~\eqref{eqn:proof_Gaussian_family_without_offset_O2_second}}
\label{subsec:proof_Gaussian_family_without_offset_O2_second}
Our construction of a sequence $G_{n} \in \Ocal_{k}(\Omega)$ to satisfy equality~\eqref{eqn:proof_Gaussian_family_without_offset_O2_second} will be similar to that of equality~\eqref{eqn:general_overfit_second} in the proof of Theorem~\ref{theorem:total_variation_bound_over-fitted_MECFG}. Here, we briefly sketch the proof for equality~\eqref{eqn:proof_Gaussian_family_without_offset_O2_second} to avoid unnecessary repetition. For any $\kappa' \prec \kappa = (2, 2, 2, 2)$, we have $\min \limits_{1 \leq i \leq 4} (\kappa')^{(i)} < 2$. Now, we construct a sequence of mixing measures, $G_{n} = \sum_{i=1}^{k_{0}+1} \pi_{i}^{n}\delta_{\parenth{\theta_{1i}^{n},\theta_{2i}^{n}}}$, with $k_{0}+1$ components as follows: $(\pi_{i}^{n}, \theta_{1i}^{n},\theta_{2i}^{n}) \equiv (\pi_{i-1}^{0}, \theta_{1(i-1)}^{0}, \theta_{2(i-1)}^{0})$ for $3 \leq i \leq k_{0} + 1$. Additionally, $\pi_{1}^{n} = \pi_{2}^{n} = 1/2$, $(\theta_{11}^{n},\theta_{21}^{n}) \equiv (\theta_{11}^{0} - 1/n, \theta_{21}^{0} - \vec{1}_{3}/n)$, and $(\theta_{12}^{n},\theta_{22}^{n}) \equiv (\theta_{11}^{0} + 1/n, \theta_{21}^{0} + \vec{1}_{3}/n)$. From here, by performing Taylor expansion up to the first order around $\theta_{11}^{0}$ and $\theta_{21}^{0}$ as in the proof of equality~\eqref{eqn:general_overfit_second}, we find that
\begin{align*}
    p_{G_{n}}(X,Y) - p_{G_{0}}(X,Y) = \overline{R}(X,Y),
\end{align*}
where $\overline{R}(X,Y)$ is a Taylor remainder such that the following limit holds
\begin{align}
\int \frac{\overline{R}^{2}(X,Y)}{p_{G_{0}}(X,Y)\widetilde{W}_{\kappa'}^{2\|\kappa'\|_{\infty}}(G_{n},G_{0})}d(X,Y) \precsim \frac{\mathcal{O}(n^{-4})}{n^{-2 \min_{1 \leq i \leq 4} \{(\kappa')^{(i)}\}}} \to 0. \nonumber
\end{align}
As a consequence, we eventually achieve that
\begin{align}
h^{2}(p_{G_{n}}, p_{G_{0}})\bigg/\widetilde{W}_{\kappa'}^{2\|\kappa'\|_{\infty}}(G_{n},G_{0}) \to 0 \nonumber
\end{align}
as $n \to \infty$, which leads to the conclusion of equality~\eqref{eqn:proof_Gaussian_family_without_offset_O2_second}.
\subsection{Proof of Theorem~\ref{theorem:lower_bound_Gaussian_family_O3}}
\label{subsec:proof:theorem:lower_bound_Gaussian_family_O3}
It is sufficient to demonstrate the following results:
\begin{eqnarray}
\lim \limits_{\epsilon \to 0} \inf \limits_{G \in \Ocal_{k,\overline{c}_{0}} (\Omega): \widetilde{W}_{\kappa}(G,G_{0}) \leq \epsilon} V(p_{G},p_{G_{0}})/\widetilde{W}_{\kappa}^{\|\kappa\|_{\infty}}(G,G_{0}) > 0, \label{eqn:proof_Gaussian_family_O3_first} \\
\inf \limits_{G \in \Ocal_{k} (\Omega)} h(p_{G},p_{G_{0}})/\widetilde{W}_{\kappa'}^{\|\kappa'\|_{\infty}}(G,G_{0}) = 0, \label{eqn:proof_Gaussian_family_O3_second}
\end{eqnarray}
for any $\kappa' \prec \kappa$ where $\kappa=(2,2,2,2)$. Proof of inequality~\eqref{eqn:proof_Gaussian_family_O3_first} is in Appendix~\ref{subsec:proof_Gaussian_family_O3_first} while proof of equality~\eqref{eqn:proof_Gaussian_family_O3_second} is in Appendix~\ref{subsec:proof_Gaussian_family_O3_second}.
\subsubsection{Proof for inequality~\eqref{eqn:proof_Gaussian_family_O3_first}}
\label{subsec:proof_Gaussian_family_O3_first}
We assume that the conclusion of inequality~\eqref{eqn:proof_Gaussian_family_O3_first} does not hold. It indicates that we can find a sequence $G_{n}$ that has representation~\eqref{eqn:proof_Gaussian_firstcase_notation} such that $V(p_{G_{n}},p_{G_{0}})/\widetilde{W}_{\kappa}^{\|\kappa\|_{\infty}}(G_{n},G_{0}) \to 0$ and $\widetilde{W}_{\kappa}(G_{n},G_{0}) \to 0$. As $\theta_{1ij}^{n}$ and $\theta_{2ij}$ both have two dimensions, in this proof we denote $\Delta \theta_{1ij}^{n} = ((\Delta \theta_{1ij}^{n})^{(1)}, (\Delta \theta_{1ij}^{n})^{(2)})$ and $\Delta \theta_{2ij}^{n} = ((\Delta \theta_{2ij}^{n})^{(1)}, (\Delta \theta_{2ij}^{n})^{(2)})$ for all $1 \leq i \leq k_{0}$ and $1 \leq j \leq s_{i}$. From Lemma~\ref{lemma:generalized_Wasserstein_distance_polynomials}, we obtain that
\begin{eqnarray}
\widetilde{W}_{\kappa}^{\|\kappa\|_{\infty}}(G_{n},G_{0}) \precsim \sum \limits_{i=1}^{k_{0}}\sum \limits_{j=1}^{s_{i}}{p_{ij}^{n}\biggr(\biggr|(\Delta \theta_{1ij}^{n})^{(1)}\biggr|^{2}+\biggr|(\Delta \theta_{1ij}^{n})^{(2)}\biggr|^{2}+\biggr|(\Delta \theta_{2ij}^{n})^{(1)}\biggr|^{2} + \biggr|(\Delta \theta_{2ij}^{n})^{(2)}\biggr|^{2}\biggr)} \nonumber \\
+ \sum \limits_{i=1}^{k_{0}} |\sum \limits_{j=1}^{s_{i}}p_{ij}^{n} - \pi_{i}^{0}| := D_{\kappa}(G_{n},G_{0}). \nonumber
\end{eqnarray}
By means of Taylor expansion up to the second order, we have
\begin{eqnarray}
p_{G_{n}}(X,Y)-p_{G_{0}}(X,Y) & = & \nonumber \\
& & \hspace{- 6 em} \sum \limits_{i=1}^{k_{0}}\sum \limits_{j=1}^{s_{i}}p_{ij}^{n}\sum \limits_{1 \leq |\alpha| \leq 2} \dfrac{1}{\alpha!}\biggr\{(\Delta \theta_{1ij}^{n})^{(1)}\biggr\}^{\alpha_{1}}\biggr\{(\Delta \theta_{1ij}^{n})^{(2)}\biggr\}^{\alpha_{2}}\biggr\{(\Delta \theta_{2ij}^{n})^{(1)}\biggr\}^{\alpha_{3}} \biggr\{(\Delta \theta_{2ij}^{n})^{(2)}\biggr\}^{\alpha_{4}} \nonumber \\
& & \hspace{- 6 em}  \times \dfrac{\partial^{|\alpha|}{f}}{\partial{(\theta_{1}^{(1)})^{\alpha_{1}}}\partial{(\theta_{1}^{(2)})^{\alpha_{2}}}\partial{(\theta_{2}^{(1)})^{\alpha_{3}}}\partial{(\theta_{2}^{(2)})^{\alpha_{4}}}}\parenth{Y|h_{1}(X,\theta_{1i}^{0}),h_{2}(X,\theta_{2i}^{0})}\overline{f}(X) \nonumber \\
& & \hspace{- 6 em} + \sum \limits_{i=1}^{k_{0}}\biggr(\sum \limits_{j=1}^{s_{i}}p_{ij}^{n} - \pi_{i}^{0}\biggr)f(Y|h_{1}(X,\theta_{1i}^{0}),h_{2}(X,\theta_{2i}^{0}))\overline{f}(X) + R(X,Y)
\nonumber \\
& & \hspace{- 6 em} : = A_{n} + B_{n} + R(X,Y), \nonumber
\end{eqnarray}
where the Taylor remainder $R(X,Y)$ is such that $R(X,Y)/D_{\kappa}(G_{n},G_{0}) \to 0$ as $n \to \infty$. From the formulations of expert functions $h_{1}, h_{2}$, we find that
\begin{eqnarray}
\dfrac{\partial^{|\alpha|}{f}}{\partial{(\theta_{1}^{(1)})^{\alpha_{1}}}\partial{(\theta_{1}^{(2)})^{\alpha_{2}}}\partial{(\theta_{2}^{(1)})^{\alpha_{3}}}\partial{(\theta_{2}^{(2)})^{\alpha_{4}}}}\parenth{Y|h_{1}(X,\theta_{1}),h_{2}(X,\theta_{2})} & = & \nonumber \\
& & \hspace{- 14 em} \dfrac{X^{2\alpha_{2}+\alpha_{3} + 3\alpha_{4}}}{2^{\alpha_{3} + \alpha_{4}}}\dfrac{\partial^{\alpha_{1}+\alpha_{2}+2\alpha_{3} + 2\alpha_{4}}{f}}{\partial{h_{1}^{\alpha_{1}+\alpha_{2}+2\alpha_{3} + 2\alpha_{4}}}}(Y|h_{1}(X|\theta_{1}),h_{2}(X|\theta_{2})), \nonumber
\end{eqnarray}
for any $\alpha = (\alpha_{1}, \alpha_{2}, \alpha_{3}, \alpha_{4}) \in \mathbb{N}^{4}$. Based on the above equation, we can express $A_{n}$ as follows:
\begin{eqnarray}
A_{n} & = & \sum \limits_{i=1}^{k_{0}}\sum \limits_{j=1}^{s_{i}}p_{ij}^{n}\sum \limits_{1 \leq |\alpha| \leq 2} \dfrac{1}{\alpha!}\biggr\{(\Delta \theta_{1ij}^{n})^{(1)}\biggr\}^{\alpha_{1}}\biggr\{(\Delta \theta_{1ij}^{n})^{(2)}\biggr\}^{\alpha_{2}}\biggr\{(\Delta \theta_{2ij}^{n})^{(1)}\biggr\}^{\alpha_{3}} \biggr\{(\Delta \theta_{2ij}^{n})^{(2)}\biggr\}^{\alpha_{4}} \nonumber \\ \nonumber \\
& & \hspace{4 em} \times \dfrac{X^{2\alpha_{2}+\alpha_{3} + 3\alpha_{4}}}{2^{\alpha_{3} + \alpha_{4}}}\dfrac{\partial^{\alpha_{1}+\alpha_{2}+2\alpha_{3} + 2\alpha_{4}}{f}}{\partial{h_{1}^{\alpha_{1}+\alpha_{2}+2\alpha_{3} + 2\alpha_{4}}}}(Y|h_{1}(X|\theta_{1i}^{0}),h_{2}(X|\theta_{2i}^{0}))\overline{f}(X). \label{eqn:proof_Gaussian_family_O3_third}
\end{eqnarray}
If we define
\begin{eqnarray}
\mathcal{F} = \biggr\{X^{l_{1}}\dfrac{\partial^{l_{2}}{f}}{\partial{h_{1}^{l_{2}}}}(Y|h_{1}(X|\theta_{1i}^{0}),h_{2}(X|\theta_{2i}^{0}))\overline{f}(X): & & \nonumber \\
& & \hspace{-15 em} l_{1} = 2 \alpha_{2} + \alpha_{3} + 3 \alpha_{4}, \ l_{2} = \alpha_{1} + \alpha_{2} + 2\alpha_{3} + 2\alpha_{4}, \ 0 \leq |\alpha| \leq 2, \ 1 \leq i \leq k_{0} \biggr\}, \nonumber
\end{eqnarray}
then we can check that the elements of $\mathcal{F}$ are linearly independent with respect to $X$ and $Y$. Therefore, we can treat $A_{n}/D_{\kappa}(G_{n},G_{0})$, $B_{n}/D_{\kappa}(G_{n},G_{0})$ as a linear combination of elements of $\mathcal{F}$. We denote $F_{l_{1},l_{2}}(\theta_{1i}^{0},\theta_{2i}^{0})$ as the coefficient of $X^{l_{1}}\dfrac{\partial^{l_{2}}{f}}{\partial{h_{1}^{l_{2}}}}(Y|h_{1}(X|\theta_{1i}^{0}),h_{2}(X|\theta_{2i}^{0}))\overline{f}(X)$ in $A_{n}$ and $B_{n}$ for any $l_{1} = 2 \alpha_{2} + \alpha_{3} + 3 \alpha_{4}, \ l_{2} = \alpha_{1} + \alpha_{2} + 2\alpha_{3} + 2\alpha_{4}, \ 0 \leq |\alpha| \leq 2$ and $1 \leq i \leq k_{0}$. Then, we can check that the coefficients of $X^{l_{1}}\dfrac{\partial^{l_{2}}{f}}{\partial{h_{1}^{l_{2}}}}(Y|h_{1}(X|\theta_{1i}^{0}),h_{2}(X|\theta_{2i}^{0}))\overline{f}(X)$ in $A_{n}/D_{\kappa}(G_{n},G_{0})$ and $B_{n}/D_{\kappa}(G_{n},G_{0})$ will be $F_{l_{1},l_{2}}(\theta_{1i}^{0},\theta_{2i}^{0})/D_{\kappa}(G_{n},G_{0})$. 

Assume that all of these coefficients go to 0 as $n \to \infty$. By taking the summation of $|F_{0,0}(\theta_{1i}^{0},\theta_{2i}^{0})/D_{\kappa}(G_{n},G_{0})|$ for all $1 \leq i \leq k_{0}$, we obtain that
\begin{eqnarray}
\biggr(\sum \limits_{i=1}^{k_{0}} |\sum \limits_{j=1}^{s_{i}}{p_{ij}^{n}} - \pi_{i}^{0}|\biggr)/D_{\kappa}(G_{n},G_{0}) \to 0. \label{eqn:proof_Gaussian_family_O3_fourth}
\end{eqnarray}
When $l_{1} = 0$ and $l_{2} = 2$, the only $\alpha$ that satisfies these equations is $\alpha = (2, 0, 0, 0)$. Given that result, the summation of $|F_{0,2}(\theta_{1i}^{0},\theta_{2i}^{0})/D_{\kappa}(G_{n},G_{0})|$ for all $1 \leq i \leq k_{0}$ leads to
\begin{align}
    \biggr(\sum \limits_{i=1}^{k_{0}}\sum \limits_{j=1}^{s_{i}}{p_{ij}^{n} \biggr|(\Delta \theta_{1ij}^{n})^{(1)}\biggr|^{2}}\biggr)/D_{\kappa}(G_{n},G_{0}) \to 0. \label{eqn:proof_Gaussian_family_O3_fifth}
\end{align}
When $l_{1} = 4$ and $l_{2} = 2$, only $\alpha = (0, 2, 0, 0)$ satisfies these equations. Summing all the coefficients $|F_{4,2}(\theta_{1i}^{0},\theta_{2i}^{0})/D_{\kappa}(G_{n},G_{0})|$ for all $1 \leq i \leq k_{0}$ leads to
\begin{align}
    \biggr(\sum \limits_{i=1}^{k_{0}}\sum \limits_{j=1}^{s_{i}}{p_{ij}^{n} \biggr|(\Delta \theta_{1ij}^{n})^{(2)}\biggr|^{2}}\biggr)/D_{\kappa}(G_{n},G_{0}) \to 0. \label{eqn:proof_Gaussian_family_O3_sixth}
\end{align}
With similar arguments, when $(l_{1}, l_{2}) \equiv (2,4)$, we have $\alpha = (0, 0, 2, 0)$ or when $(l_{1}, l_{2}) \equiv (6, 4)$, we have $\alpha = (0, 0, 0, 2)$. By considering respectively the summation of the coefficients of $|F_{2,4}(\theta_{1i}^{0},\theta_{2i}^{0})/D_{\kappa}(G_{n},G_{0})|$ or $|F_{6,4}(\theta_{1i}^{0},\theta_{2i}^{0})/D_{\kappa}(G_{n},G_{0})|$ for all $1 \leq i \leq k_{0}$, we find that
\begin{align}
    \biggr(\sum \limits_{i=1}^{k_{0}}\sum \limits_{j=1}^{s_{i}}{p_{ij}^{n} \biggr|(\Delta \theta_{2ij}^{n})^{(1)}\biggr|^{2}}\biggr)/D_{\kappa}(G_{n},G_{0}) \to 0, \nonumber \\
    \biggr(\sum \limits_{i=1}^{k_{0}}\sum \limits_{j=1}^{s_{i}}{p_{ij}^{n} \biggr|(\Delta \theta_{2ij}^{n})^{(2)}\biggr|^{2}}\biggr)/D_{\kappa}(G_{n},G_{0}) \to 0. \label{eqn:proof_Gaussian_family_O3_seventh}
\end{align}
Combining the results from equations~\eqref{eqn:proof_Gaussian_family_O3_fourth}-\eqref{eqn:proof_Gaussian_family_O3_seventh}, we obtain
\begin{align*}
    1 = \frac{\sum \limits_{i=1}^{k_{0}}\sum \limits_{j=1}^{s_{i}}{p_{ij}^{n}\biggr(\biggr|(\Delta \theta_{1ij}^{n})^{(1)}\biggr|^{2}+\biggr|(\Delta \theta_{1ij}^{n})^{(2)}\biggr|^{2}+\biggr|(\Delta \theta_{2ij}^{n})^{(1)}\biggr|^{2} + \biggr|(\Delta \theta_{2ij}^{n})^{(2)}\biggr|^{2}\biggr)}
+ \sum \limits_{i=1}^{k_{0}} |\sum \limits_{j=1}^{s_{i}}p_{ij}^{n} - \pi_{i}^{0}|}{D_{\kappa}(G_{n},G_{0})} \to 0,
\end{align*}
which is a contradiction. As a consequence, not all the coefficients in the linear combinations of $A_{n}/D_{\kappa}(G_{n},G_{0})$ and $B_{n}/D_{\kappa}(G_{n},G_{0})$ go to 0 as $n \to \infty$. From here, using the Fatou's argument in Step 4 of the proof of Theorem~\ref{theorem:total_variation_bound_over-fitted_MECFG} and the fact that the elements of $\mathcal{F}$ are linearly independent with respect to $X$ and $Y$, we achieve the conclusion of inequality~\eqref{eqn:proof_Gaussian_family_O3_first}.
\subsubsection{Proof for equality~\eqref{eqn:proof_Gaussian_family_O3_second}}
\label{subsec:proof_Gaussian_family_O3_second}
The proof of equality~\eqref{eqn:proof_Gaussian_family_O3_second} is similar in spirit to that of equality~\eqref{eqn:proof_Gaussian_family_without_offset_O2_second}; therefore, we only provide a proof sketch for this equality. For any $\kappa' \prec \kappa = (2, 2, 2, 2)$, we construct a sequence $G_{n} \in \Ocal_{k}(\Omega)$ such that $G_{n} = \sum_{i=1}^{k_{0}+1} \pi_{i}^{n}\delta_{\parenth{\theta_{1i}^{n},\theta_{2i}^{n}}}$, with $k_{0}+1$ components as follows: $(\pi_{i}^{n}, \theta_{1i}^{n},\theta_{2i}^{n}) \equiv (\pi_{i-1}^{0}, \theta_{1(i-1)}^{0}, \theta_{2(i-1)}^{0})$ for $3 \leq i \leq k_{0} + 1$. Additionally, $\pi_{1}^{n} = \pi_{2}^{n} = 1/2$, $(\theta_{11}^{n},\theta_{21}^{n}) \equiv (\theta_{11}^{0} - \vec{1}_{2}/n, \theta_{21}^{0} - \vec{1}_{2}/n)$, and $(\theta_{12}^{n},\theta_{22}^{n}) \equiv (\theta_{11}^{0} + \vec{1}_{2}/n, \theta_{21}^{0} + \vec{1}_{2}/n)$. From here, by performing Taylor expansion up to the first order around $\theta_{11}^{0}$ and $\theta_{21}^{0}$ as in the proof of equality~\eqref{eqn:general_overfit_second}, we find that
\begin{align*}
    p_{G_{n}}(X,Y) - p_{G_{0}}(X,Y) = \overline{R}(X,Y),
\end{align*}
where $\overline{R}(X,Y)$ is a Taylor remainder such that the following limit holds
\begin{align}
\int \frac{\overline{R}^{2}(X,Y)}{p_{G_{0}}(X,Y)\widetilde{W}_{\kappa'}^{2\|\kappa'\|_{\infty}}(G_{n},G_{0})}d(X,Y) \precsim \frac{\mathcal{O}(n^{-4})}{n^{-2 \min_{1 \leq i \leq 4} \{(\kappa')^{(i)}\}}} \to 0. \nonumber
\end{align}
As a consequence, we eventually achieve that
\begin{align}
h^{2}(p_{G_{n}}, p_{G_{0}})\bigg/\widetilde{W}_{\kappa'}^{2\|\kappa'\|_{\infty}}(G_{n},G_{0}) \to 0 \nonumber
\end{align}
as $n \to \infty$, which leads to the conclusion of equality~\eqref{eqn:proof_Gaussian_family_O3_second}. 
\subsection{Proof of Theorem~\ref{theorem:lower_bound_Gaussian_family_O3_O1}}
\label{subsec:proof:theorem:lower_bound_Gaussian_family_O3_O1}
We will demonstrate that
\begin{eqnarray}
\lim \limits_{\epsilon \to 0} \inf \limits_{G \in \Ocal_{k,\overline{c}_{0}} (\Omega): \widetilde{W}_{\kappa}(G,G_{0}) \leq \epsilon} V(p_{G},p_{G_{0}})/\widetilde{W}_{\kappa}^{\|\kappa\|_{\infty}}(G,G_{0}) > 0, \label{eqn:proof_Gaussian_family_O3_O1_first} \\
\inf \limits_{G \in \Ocal_{k} (\Omega)} h(p_{G},p_{G_{0}})/\widetilde{W}_{\kappa'}^{\|\kappa'\|_{\infty}}(G,G_{0}) = 0, \label{eqn:proof_Gaussian_family_O3_O1_second}
\end{eqnarray}
for any $\kappa' \prec \kappa$ where $\kappa = (\overline{r}, 2, \ceil{\overline{r}/2}, 2)$. Without loss of generality, we assume that $\bar{r}$ is an even number. The proof when $\bar{r}$ is an odd number is similar. Proof of inequality~\eqref{eqn:proof_Gaussian_family_O3_O1_first} is in Appendix~\ref{subsec:proof_Gaussian_family_O3_O1_first} while proof of equality~\eqref{eqn:proof_Gaussian_family_O3_O1_second} is in Appendix~\ref{subsec:proof_Gaussian_family_O3_O1_second}.
\subsubsection{Proof for inequality~\eqref{eqn:proof_Gaussian_family_O3_O1_first}}
\label{subsec:proof_Gaussian_family_O3_O1_first}
Assume that the conclusion of inequality~\eqref{eqn:proof_Gaussian_family_O3_O1_first} does not hold. Therefore, we can find a sequence $G_{n}$ that has representation~\eqref{eqn:proof_Gaussian_firstcase_notation} such that $V(p_{G_{n}},p_{G_{0}})/\widetilde{W}_{\kappa}^{\|\kappa\|_{\infty}}(G_{n},G_{0}) \to 0$ and $\widetilde{W}_{\kappa}(G_{n},G_{0}) \to 0$. As $\theta_{1ij}$ has two dimensions and $\theta_{2ij}^{n}$ has two dimensions, in this proof we denote $\Delta \theta_{1ij}^{n} = ((\Delta \theta_{1ij}^{n})^{(1)}, (\Delta \theta_{1ij}^{n})^{(2)})$ and $\Delta \theta_{2ij}^{n} = ((\Delta \theta_{2ij}^{n})^{(1)}, (\Delta \theta_{2ij}^{n})^{(2)})$ for all $1 \leq i \leq k_{0}$ and $1 \leq j \leq s_{i}$. From Lemma~\ref{lemma:generalized_Wasserstein_distance_polynomials}, we have
\begin{eqnarray}
\widetilde{W}_{\kappa}^{\|\kappa\|_{\infty}}(G_{n},G_{0}) \precsim \sum \limits_{i=1}^{k_{0}}\sum \limits_{j=1}^{s_{i}} p_{ij}^{n}\biggr(\biggr|(\Delta \theta_{1ij}^{n})^{(1)}\biggr|^{\bar{r}}+\biggr|(\Delta \theta_{1ij}^{n})^{(2)}\biggr|^{2} + \biggr|(\Delta \theta_{2ij}^{n})^{(1)}\biggr|^{\bar{r}/2}+\biggr|(\Delta \theta_{2ij}^{n})^{(2)}\biggr|^{2}\biggr) \nonumber \\
+ \sum \limits_{i=1}^{k_{0}} |\sum \limits_{j=1}^{s_{i}}p_{ij}^{n} - \pi_{i}^{0}| := D_{\kappa}(G_{n},G_{0}). \nonumber
\end{eqnarray}
Invoking Taylor expansion up to the $\bar{r}$-th order, we find that
\begin{eqnarray}
p_{G_{n}}(X,Y)-p_{G_{0}}(X,Y) & = & \nonumber \\
& & \hspace{- 10 em} \sum \limits_{i=1}^{k_{0}}\sum \limits_{j=1}^{s_{i}}p_{ij}^{n}\sum \limits_{1 \leq |\alpha| \leq \bar{r}} \dfrac{1}{\alpha!}\biggr\{(\Delta \theta_{1ij}^{n})^{(1)}\biggr\}^{\alpha_{1}} \biggr\{(\Delta \theta_{1ij}^{n})^{(2)}\biggr\}^{\alpha_{2}} \biggr\{(\Delta \theta_{2ij}^{n})^{(1)}\biggr\}^{\alpha_{3}}\biggr\{(\Delta \theta_{2ij}^{n})^{(2)}\biggr\}^{\alpha_{4}} \nonumber \\
& & \hspace{- 10 em}  \times \dfrac{\partial^{|\alpha|}{f}}{\partial{(\theta_{1}^{(1)})^{\alpha_{1}}}\partial{(\theta_{1}^{(2)})^{\alpha_{2}}}\partial{(\theta_{2}^{(1)})^{\alpha_{3}}}\partial{(\theta_{2}^{(2)})^{\alpha_{4}}}}\parenth{Y|h_{1}(X,\theta_{1i}^{0}),h_{2}(X,\theta_{2i}^{0})}\overline{f}(X) \nonumber \\
& & \hspace{- 10 em} + \sum \limits_{i=1}^{k_{0}}\biggr(\sum \limits_{j=1}^{s_{i}}p_{ij}^{n} - \pi_{i}^{0}\biggr)f(Y|h_{1}(X,\theta_{1i}^{0}),h_{2}(X,\theta_{2i}^{0}))\overline{f}(X) + R(X,Y)
\nonumber \\
& & \hspace{- 10 em} : = A_{n} + B_{n} + R(X,Y), \nonumber
\end{eqnarray}
where the Taylor remainder $R(X,Y)$ is such that $R(X,Y)/D_{\kappa}(G_{n},G_{0}) \to 0$ as $n \to \infty$. Since $h_{1}(X, \theta_{1}) = \theta_{1}^{(1)} + \theta_{1}^{(2)} X^2$ and $h_{2}(X, \theta_{2}) = \theta_{2}^{(1)} + \theta_{2}^{(2)} X^2$, we can verify that
\begin{eqnarray}
\dfrac{\partial^{|\alpha|}{f}}{\partial{(\theta_{1}^{(1)})^{\alpha_{1}}}\partial{(\theta_{1}^{(2)})^{\alpha_{2}}}\partial{(\theta_{2}^{(1)})^{\alpha_{3}}}\partial{(\theta_{2}^{(2)})^{\alpha_{4}}}}\parenth{Y|h_{1}(X,\theta_{1}),h_{2}(X,\theta_{2})} & = & \nonumber \\
& & \hspace{- 17 em} \dfrac{X^{2\alpha_{2}+2\alpha_{4}}}{2^{\alpha_{3} + \alpha_{4}}}\dfrac{\partial^{\alpha_{1}+ \alpha_{2} + 2\alpha_{3}+2\alpha_{4}}{f}}{\partial{h_{1}^{\alpha_{1}+ \alpha_{2} + 2\alpha_{3}+2\alpha_{4}}}}(Y|h_{1}(X|\theta_{1}),h_{2}(X|\theta_{2})), \nonumber
\end{eqnarray}
for any $\alpha = (\alpha_{1}, \alpha_{2}, \alpha_{3}, \alpha_{4}) \in \mathbb{N}^{4}$. Given the above equation, we can rewrite $A_{n}$ as follows:
\begin{eqnarray}
A_{n} & = &  \sum \limits_{i=1}^{k_{0}}\sum \limits_{j=1}^{s_{i}}p_{ij}^{n}\sum \limits_{1 \leq |\alpha| \leq \bar{r}} \dfrac{1}{\alpha!}\biggr\{(\Delta \theta_{1ij}^{n})^{(1)}\biggr\}^{\alpha_{1}} \biggr\{(\Delta \theta_{1ij}^{n})^{(2)}\biggr\}^{\alpha_{2}} \biggr\{(\Delta \theta_{2ij}^{n})^{(1)}\biggr\}^{\alpha_{3}}\biggr\{(\Delta \theta_{2ij}^{n})^{(2)}\biggr\}^{\alpha_{4}} \nonumber \\
& & \times  \dfrac{X^{2\alpha_{2}+2\alpha_{4}}}{2^{\alpha_{3} + \alpha_{4}}}\dfrac{\partial^{\alpha_{1}+ \alpha_{2} + 2\alpha_{3}+2\alpha_{4}}{f}}{\partial{h_{1}^{\alpha_{1}+ \alpha_{2} + 2\alpha_{3}+2\alpha_{4}}}}(Y|h_{1}(X|\theta_{1i}^{0}),h_{2}(X|\theta_{2i}^{0}))\overline{f}(X). \label{eqn:proof_Gaussian_family_O3_O1_third}
\end{eqnarray}
Similar to the previous proofs, we define
\begin{eqnarray}
\mathcal{F} = \biggr\{X^{l_{1}}\dfrac{\partial^{l_{2}}{f}}{\partial{h_{1}^{l_{2}}}}(Y|h_{1}(X|\theta_{1i}^{0}),h_{2}(X|\theta_{2i}^{0}))\overline{f}(X): & & \nonumber \\
& & \hspace{-15 em} l_{1} = 2 \alpha_{2} + 2 \alpha_{4}, \ l_{2} = \alpha_{1} + \alpha_{2} + 2 \alpha_{3} + 2\alpha_{4}, \ 0 \leq |\alpha| \leq \bar{r}, \ 1 \leq i \leq k_{0} \biggr\}, \nonumber
\end{eqnarray}
Based on the proof argument similar to that of the claim in equation~\eqref{eqn:lemma_linear_independence_Gaussian_firstcase_first}, we can demonstrate that the elements of $\mathcal{F}$ are linearly independent with respect to $X$ and $Y$. Therefore, we can treat $A_{n}/D_{\kappa}(G_{n},G_{0})$, $B_{n}/D_{\kappa}(G_{n},G_{0})$ as a linear combination of elements of $\mathcal{F}$. We denote $F_{l_{1},l_{2}}(\theta_{1i}^{0},\theta_{2i}^{0})$ as the coefficient of $X^{l_{1}}\dfrac{\partial^{l_{2}}{f}}{\partial{h_{1}^{l_{2}}}}(Y|h_{1}(X|\theta_{1i}^{0}),h_{2}(X|\theta_{2i}^{0}))\overline{f}(X)$ in $A_{n}$ and $B_{n}$ for any $l_{1} = 2 \alpha_{2} + 2 \alpha_{4}$, $l_{2} = \alpha_{1} + \alpha_{2} + 2 \alpha_{3} + 2\alpha_{4}$, $0 \leq |\alpha| \leq \bar{r}$ and $1 \leq i \leq k_{0}$. Then, we can check that the coefficients of $X^{l_{1}}\dfrac{\partial^{l_{2}}{f}}{\partial{h_{1}^{l_{2}}}}(Y|h_{1}(X|\theta_{1i}^{0}),h_{2}(X|\theta_{2i}^{0}))\overline{f}(X)$ in $A_{n}/D_{\kappa}(G_{n},G_{0})$ and $B_{n}/D_{\kappa}(G_{n},G_{0})$ will be $F_{l_{1},l_{2}}(\theta_{1i}^{0},\theta_{2i}^{0})/D_{\kappa}(G_{n},G_{0})$. 

Assume that all of these coefficients go to 0 as $n \to \infty$. By taking the summation of $|F_{0,0}(\theta_{1i}^{0},\theta_{2i}^{0})/D_{\kappa}(G_{n},G_{0})|$ for all $1 \leq i \leq k_{0}$, we obtain that
\begin{eqnarray}
\biggr(\sum \limits_{i=1}^{k_{0}} |\sum \limits_{j=1}^{s_{i}}{p_{ij}^{n}} - \pi_{i}^{0}|\biggr)/D_{\kappa}(G_{n},G_{0}) \to 0. \label{eqn:proof_Gaussian_family_O3_O1_fourth}
\end{eqnarray}
When $l_{1} = 4$ and $l_{2} = 2$, we can check that $\alpha = (0, 2, 0, 0)$ is the only solution to these equations. Similarly, when $l_{1} = 4$ and $l_{2} = 4$, the only solution to these equations is $\alpha = (0, 0, 0, 2)$. Therefore, the summation of $|F_{4,2}(\theta_{1i}^{0},\theta_{2i}^{0})/D_{\kappa}(G_{n},G_{0})|$ and $|F_{4,4}(\theta_{1i}^{0},\theta_{2i}^{0})/D_{\kappa}(G_{n},G_{0})|$ for all $1 \leq i \leq k_{0}$ leads to
\begin{align}
    \biggr(\sum \limits_{i=1}^{k_{0}}\sum \limits_{j=1}^{s_{i}}{p_{ij}^{n} \biggr|(\Delta \theta_{1ij}^{n})^{(2)}\biggr|^{2}}\biggr)/D_{\kappa}(G_{n},G_{0}) \to 0, \nonumber \\
    \biggr(\sum \limits_{i=1}^{k_{0}}\sum \limits_{j=1}^{s_{i}}{p_{ij}^{n} \biggr|(\Delta \theta_{2ij}^{n})^{(2)}\biggr|^{2}}\biggr)/D_{\kappa}(G_{n},G_{0}) \to 0. \label{eqn:proof_Gaussian_family_O3_O1_fifth}
\end{align}
Combining the results from equations~\eqref{eqn:proof_Gaussian_family_O3_O1_fourth}-\eqref{eqn:proof_Gaussian_family_O3_O1_fifth}, we find that
\begin{align*}
    \frac{\sum \limits_{i=1}^{k_{0}}\sum \limits_{j=1}^{s_{i}}{p_{ij}^{n}\biggr(\biggr|(\Delta \theta_{1ij}^{n})^{(2)}\biggr|^{2}+\biggr|(\Delta \theta_{2ij}^{n})^{(2)}\biggr|^{2}\biggr)}
+ \sum \limits_{i=1}^{k_{0}} |\sum \limits_{j=1}^{s_{i}}p_{ij}^{n} - \pi_{i}^{0}|}{D_{\kappa}(G_{n},G_{0})} \to 0.
\end{align*}
The above result indicates that
\begin{align*}
    \frac{\sum \limits_{i=1}^{k_{0}}\sum \limits_{j=1}^{s_{i}}{p_{ij}^{n}\biggr(\biggr|(\Delta \theta_{1ij}^{n})^{(1)}\biggr|^{\bar{r}}+\biggr|(\Delta \theta_{2ij}^{n})^{(1)}\biggr|^{\bar{r}/2} \biggr)}}{D_{\kappa}(G_{n},G_{0})} \to 1.
\end{align*}
Hence, we can find an index $i^{*} \in \{1, 2, \ldots, k_{0}\}$ such that 
\begin{align*}
    L = \frac{\sum \limits_{j=1}^{s_{i^{*}}}{p_{i^{*}j}^{n}\biggr(\biggr|(\Delta \theta_{1i^{*}j}^{n})^{(1)}\biggr|^{\bar{r}}+\biggr|(\Delta \theta_{2i^{*}j}^{n})^{(1)}\biggr|^{\bar{r}/2} \biggr)}}{D_{\kappa}(G_{n},G_{0})} \not \to 0
\end{align*}
as $n \to \infty$. By denoting $M_{l_{2}}(\theta_{1i}^{0}, \theta_{2i}^{0}) = F_{0, l_{2}}(\theta_{1i}^{0}, \theta_{2i}^{0})/ \sum \limits_{j=1}^{s_{i^{*}}}{p_{i^{*}j}^{n}\biggr(\biggr|(\Delta \theta_{1i^{*}j}^{n})^{(1)}\biggr|^{\bar{r}}+\biggr|(\Delta \theta_{2i^{*}j}^{n})^{(1)}\biggr|^{\bar{r}/2} \biggr)}$ for all $1 \leq i \leq k_{0}$, we obtain that
\begin{align*}
    M_{l_{2}}(\theta_{1i}^{0}, \theta_{2i}^{0}) = \frac{1}{L} \frac{F_{0, l_{2}}(\theta_{1i}^{0}, \theta_{2i}^{0})}{D_{\kappa}(G_{n},G_{0})} \to 0
\end{align*}
for any $1 \leq i \leq k_{0}$ and $0 \leq l_{2} \leq 2 \bar{r}$. From the formulation of $F_{0, l_{2}}(\theta_{1i}^{0}, \theta_{2i}^{0})$, the above limits with $M_{l_{2}}(\theta_{1i^{*}}^{0}, \theta_{2i^{*}}^{0})$ can be rewritten as:
\begin{eqnarray}
\dfrac{\sum \limits_{j=1}^{s_{i^{*}}}p_{i^{*}j}^{n}\sum \limits_{\substack{\alpha_{1}+2\alpha_{3}=l_{2} \\ \alpha_{1}+\alpha_{3} \leq \overline{r}}}{\dfrac{\biggr\{(\Delta \theta_{1i^{*}j}^{n})^{(1)}\biggr\}^{\alpha_{1}}\biggr\{(\Delta \theta_{2i^{*}j}^{n})^{(1)}\biggr\}^{\alpha_{3}}}{2^{\alpha_{3}}\alpha_{1}!\alpha_{3}!}}}{\sum \limits_{j=1}^{s_{i^{*}}}{p_{i^{*}j}^{n}\biggr(\biggr|(\Delta \theta_{1i^{*}j}^{n})^{(1)}\biggr|^{\bar{r}}+\biggr|(\Delta \theta_{2i^{*}j}^{n})^{(1)}\biggr|^{\bar{r}/2} \biggr)}} \to 0 \nonumber
\end{eqnarray}
for any $0 \leq \l_{2} \leq 2 \bar{r}$. According to the argument in Step 3 of the proof of Theorem~\ref{theorem:lower_bound_Gaussian_family}, that system of limits cannot happen. As a consequence, not all the coefficients in the linear combinations of $A_{n}/D_{\kappa}(G_{n},G_{0})$ and $B_{n}/D_{\kappa}(G_{n},G_{0})$ go to 0 as $n \to \infty$. From here, using the Fatou's argument in Step 4 of the proof of Theorem~\ref{theorem:total_variation_bound_over-fitted_MECFG} and the fact that the elements of $\mathcal{F}$ are linearly independent with respect to $X$ and $Y$, we achieve the conclusion of inequality~\eqref{eqn:proof_Gaussian_family_O3_O1_first}.
\subsubsection{Proof for equality~\eqref{eqn:proof_Gaussian_family_O3_O1_second}}
\label{subsec:proof_Gaussian_family_O3_O1_second}
The proof of equality~\eqref{eqn:proof_Gaussian_family_O3_O1_second} is similar to that of equality~\eqref{eqn:proof_Gaussian_firstcase_second}. Hence, we only provide the proof sketch. In this proof, we consider two settings of $\kappa' \prec \kappa = (\overline{r}, 2, \overline{r}/2, 2)$.
\paragraph{Case 1:} $\kappa' = \parenth{\kappa'^{(1)}, \kappa'^{(2)}, \kappa'^{(3)}, \kappa'^{(4)}}$ when at least one of $\kappa'^{(2)}, \kappa'^{(4)} < 2$. Under this setting, we construct $G_{n} = \sum_{i = 1}^{k_{0}+1} \pi_{i}^{n}\delta_{(\theta_{1i}^{n},\theta_{2i}^{n})}$ such that $(\pi_{i}^{n},\theta_{1i}^{n}, \theta_{2i}^{n}) \equiv (\pi_{i-1}^{0}, \theta_{1(i-1)}^{0}, \theta_{2(i-1)}^{0})$ for $3 \leq i \leq k_{0}+1$. Additionally, $\pi_{1}^{n} = \pi_{2}^{n} = \pi_{1}^{0}/2$, $\parenth{(\theta_{1i}^{n})^{(1)},(\theta_{2i}^{n})^{(1)}} = \parenth{(\theta_{11}^{0})^{(1)},(\theta_{21}^{0})^{(1)}}$ for $1 \leq i \leq 2$, and $(\theta_{11}^{n})^{(2)} = (\theta_{11}^{0})^{(2)} - 1/n$, $(\theta_{21}^{n})^{(2)} = (\theta_{21}^{0})^{(2)} - 1/n$, $(\theta_{12}^{n})^{(2)} = (\theta_{11}^{0})^{(2)} + 1/n$, $(\theta_{22}^{n})^{(2)} = (\theta_{21}^{0})^{(2)} + 1/n$. From this construction of $G_{n}$, we can verify that $\widetilde{W}_{\kappa'}^{\|\kappa'\|_{\infty}}(G_{n},G_{0}) \asymp n^{-\min \{\kappa'^{(2)}, \kappa'^{(4)}\}} = o(n^{-2})$. Based on Taylor expansion up to the first order around $\theta_{11}^{0}, \theta_{21}^{0}$, we have
\begin{align}
p_{G_{n}}(X,Y) - p_{G_{0}}(X,Y) = \overline{R}_{1}(X,Y), \nonumber
\end{align}
where $\overline{R}_{1}(X,Y)$ is Taylor remainder such that
\begin{align}
\frac{h^{2}(p_{G_{n}}, p_{G_{0}})}{\widetilde{W}_{\kappa'}^{2\|\kappa'\|_{\infty}}(G_{n},G_{0})} \precsim \int \frac{\overline{R}_{1}^{2}(X,Y)}{p_{G_{0}}(X,Y)\widetilde{W}_{\kappa'}^{2\|\kappa'\|_{\infty}}(G_{n},G_{0})}d(X,Y) \precsim \frac{\mathcal{O}(n^{-4})}{n^{-2 \min \{\kappa'^{(2)}, \kappa'^{(4)}\}}} \to 0 \nonumber
\end{align}
as $n \to \infty$. Therefore, we achieve the conclusion of equality~\eqref{eqn:proof_Gaussian_family_O3_O1_second} under Case 1. 
\paragraph{Case 2:} $\kappa' = \parenth{\kappa'^{(1)}, 2, \kappa'^{(3)}, 2}$ when $(\kappa_{1}'^{(1)},\kappa'^{(3)}) \prec (\overline{r},\overline{r}/2)$. Under this setting of $\kappa'$, we construct $G_{n} = \sum_{i=1}^{k}\pi_{i}^{n}\delta_{(\theta_{1i}^{n},\theta_{2i}^{n})}$ such that $(\pi_{i+k-k_{0}}^{n}, \theta_{1(i+k-k_{0})}^{n}, \theta_{2(i+k-k_{0})}^{n}) = (\pi_{i}^{0},\theta_{1i}^{0},\theta_{2i}^{0})$ for $2 \leq i \leq k_{0}$. For $1 \leq j \leq k - k_{0} + 1$, we choose $(\theta_{1j}^{n})^{(2)} = (\theta_{11}^{0})^{(2)}, (\theta_{2j}^{n})^{(2)} = (\theta_{21}^{0})^{(2)}$, and
\begin{align}
(\theta_{1j}^{n})^{(1)} = (\theta_{11}^{0})^{(1)} + \frac{a_{j}^{*}}{n}, \ (\theta_{2j}^{n})^{(1)} = (\theta_{21}^{0})^{(1)} + \frac{2b_{j}^{*}}{n^{2}}, \ \pi_{j}^{n} = \frac{\pi_{1}^{0}(c_{j}^{*})^{2}}{\sum_{i=1}^{k-k_{0}+1} (c_{j}^{*})^{2}}, \nonumber
\end{align}
where $(c_{i}^{*},a_{i}^{*},b_{i}^{*})_{i=1}^{k-k_{0}+1}$ are the non-trivial solution of system of polynomial equations~\eqref{eqn:system_polynomial_Gaussian_first} when $r = \overline{r} - 1$. From here, by performing Taylor expansion around $(\theta_{11}^{0},\theta_{21}^{0})$ and arguing similarly as Case 2 in the proof of equality~\eqref{eqn:proof_Gaussian_firstcase_second} in Section~\ref{subsection:key_equality_algebraic_dependent_linear_first}, we obtain that
\begin{align}
p_{G_{n}}(X,Y) - p_{G_{0}}(X,Y) = \sum \limits_{l = \overline{r}}^{2\overline{r}-2} \mathcal{O}(n^{-\overline{r}})\dfrac{\partial^{l}{f}}{\partial{h_{1}^{l}}}(Y|h_{1}(X,\theta_{1i}^{0}),h_{2}(X,\theta_{2i}^{0}))\overline{f}(X) + \overline{R}_{2}(X,Y), \nonumber 
\end{align}
where $\overline{R}_{2}(X,Y)$ is Taylor remainder such that the following limit holds
\begin{align}
\int \frac{\overline{R}_{2}^{2}(X,Y)}{p_{G_{0}}(X,Y)\widetilde{W}_{\kappa'}^{2\|\kappa'\|_{\infty}}(G_{n},G_{0})}d(X,Y) \precsim \frac{\mathcal{O}(n^{-2\overline{r}})}{n^{-2\min\{\kappa'^{(1)},\kappa'^{(3)}\}}} \to 0. \nonumber
\end{align}
Therefore, we achieve that
\begin{align}
h^{2}(p_{G_{n}}, p_{G_{0}})\bigg/\widetilde{W}_{\kappa'}^{2\|\kappa'\|_{\infty}}(G_{n},G_{0}) \to 0 \nonumber
\end{align}
as $n \to \infty$. As a consequence, we reach the conclusion of equality~\eqref{eqn:proof_Gaussian_family_O3_O1_second} under Case 2.
\subsection{Proof of Theorem~\ref{theorem:lower_bound_Gaussian_family_O3_O2_O1}}
\label{subsection:proof_theorem:lower_bound_Gaussian_family_O3_O2_O1}
We will demonstrate that
\begin{eqnarray}
\lim \limits_{\epsilon \to 0} \inf \limits_{G \in \Ocal_{k,\overline{c}_{0}} (\Omega): \widetilde{W}_{\kappa}(G,G_{0}) \leq \epsilon} V(p_{G},p_{G_{0}})/\widetilde{W}_{\kappa}^{\|\kappa\|_{\infty}}(G,G_{0}) > 0, \label{eqn:proof_Gaussian_family_O3_O2_O1_first} \\
\inf \limits_{G \in \Ocal_{k} (\Omega)} h(p_{G},p_{G_{0}})/\widetilde{W}_{\kappa'}^{\|\kappa'\|_{\infty}}(G,G_{0}) = 0, \label{eqn:proof_Gaussian_family_O3_O2_O1_second}
\end{eqnarray}
for any $\kappa' \prec \kappa$ where $\kappa = (\overline{r}, 2, \ceil{\overline{r}/2}, 2, 2)$. Without loss of generality, we assume that $\bar{r}$ is an even number. The proof when $\bar{r}$ is an odd number is similar. Proof of inequality~\eqref{eqn:proof_Gaussian_family_O3_O2_O1_first} is in Appendix~\ref{subsec:proof_Gaussian_family_O3_O2_O1_first} while proof of equality~\eqref{eqn:proof_Gaussian_family_O3_O2_O1_second} is in Appendix~\ref{subsec:proof_Gaussian_family_O3_O2_O1_second}.
\subsubsection{Proof for inequality~\eqref{eqn:proof_Gaussian_family_O3_O2_O1_first}}
\label{subsec:proof_Gaussian_family_O3_O2_O1_first}
Assume that the conclusion of inequality~\eqref{eqn:proof_Gaussian_family_O3_O2_O1_first} does not hold. It suggests that we can find a sequence $G_{n}$ that has representation~\eqref{eqn:proof_Gaussian_firstcase_notation} such that $V(p_{G_{n}},p_{G_{0}})/\widetilde{W}_{\kappa}^{\|\kappa\|_{\infty}}(G_{n},G_{0}) \to 0$ and $\widetilde{W}_{\kappa}(G_{n},G_{0}) \to 0$. In this proof we denote $\Delta \theta_{1ij}^{n} = ((\Delta \theta_{1ij}^{n})^{(1)}, (\Delta \theta_{1ij}^{n})^{(2)})$ and $\Delta \theta_{2ij}^{n} = ((\Delta \theta_{2ij}^{n})^{(1)}, (\Delta \theta_{2ij}^{n})^{(2)}, (\Delta \theta_{2ij}^{n})^{(3)})$ for all $1 \leq i \leq k_{0}$ and $1 \leq j \leq s_{i}$. From Lemma~\ref{lemma:generalized_Wasserstein_distance_polynomials}, we have
\begin{eqnarray}
\widetilde{W}_{\kappa}^{\|\kappa\|_{\infty}}(G_{n},G_{0}) \precsim \sum \limits_{i=1}^{k_{0}}\sum \limits_{j=1}^{s_{i}} p_{ij}^{n}\biggr(\biggr|(\Delta \theta_{1ij}^{n})^{(1)}\biggr|^{\bar{r}}+\biggr|(\Delta \theta_{1ij}^{n})^{(2)}\biggr|^{2} + \biggr|(\Delta \theta_{2ij}^{n})^{(1)}\biggr|^{\bar{r}/2}+\biggr|(\Delta \theta_{2ij}^{n})^{(2)}\biggr|^{2} \nonumber \\
+ \biggr|(\Delta \theta_{2ij}^{n})^{(3)}\biggr|^{2}\biggr)
+ \sum \limits_{i=1}^{k_{0}} |\sum \limits_{j=1}^{s_{i}}p_{ij}^{n} - \pi_{i}^{0}| := D_{\kappa}(G_{n},G_{0}). \nonumber
\end{eqnarray}
An application of Taylor expansion up to the $\bar{r}$-th order leads to
\begin{eqnarray}
p_{G_{n}}(X,Y)-p_{G_{0}}(X,Y) & = & \nonumber \\
& & \hspace{- 13 em} \sum \limits_{i=1}^{k_{0}}\sum \limits_{j=1}^{s_{i}}p_{ij}^{n}\sum \limits_{1 \leq |\alpha| \leq \bar{r}} \dfrac{1}{\alpha!}\biggr\{(\Delta \theta_{1ij}^{n})^{(1)}\biggr\}^{\alpha_{1}} \biggr\{(\Delta \theta_{1ij}^{n})^{(2)}\biggr\}^{\alpha_{2}} \biggr\{(\Delta \theta_{2ij}^{n})^{(1)}\biggr\}^{\alpha_{3}}\biggr\{(\Delta \theta_{2ij}^{n})^{(2)}\biggr\}^{\alpha_{4}} \biggr\{(\Delta \theta_{2ij}^{n})^{(3)}\biggr\}^{\alpha_{5}} \nonumber \\
& & \hspace{- 13 em}  \times \dfrac{\partial^{|\alpha|}{f}}{\partial{(\theta_{1}^{(1)})^{\alpha_{1}}}\partial{(\theta_{1}^{(2)})^{\alpha_{2}}}\partial{(\theta_{2}^{(1)})^{\alpha_{3}}}\partial{(\theta_{2}^{(2)})^{\alpha_{4}}}\partial{(\theta_{2}^{(3)})^{\alpha_{5}}}}\parenth{Y|h_{1}(X,\theta_{1i}^{0}),h_{2}(X,\theta_{2i}^{0})}\overline{f}(X) \nonumber \\
& & \hspace{- 13 em} + \sum \limits_{i=1}^{k_{0}}\biggr(\sum \limits_{j=1}^{s_{i}}p_{ij}^{n} - \pi_{i}^{0}\biggr)f(Y|h_{1}(X,\theta_{1i}^{0}),h_{2}(X,\theta_{2i}^{0}))\overline{f}(X) + R(X,Y)
\nonumber \\
& & \hspace{- 13 em} : = A_{n} + B_{n} + R(X,Y), \nonumber
\end{eqnarray}
where the Taylor remainder $R(X,Y)$ is such that $R(X,Y)/D_{\kappa}(G_{n},G_{0}) \to 0$ as $n \to \infty$. Since $h_{1}(X, \theta_{1}) = \theta_{1}^{(1)} + \theta_{1}^{(2)} X^2$ and $h_{2}(X, \theta_{2}) = \theta_{2}^{(1)} + \theta_{2}^{(2)} X + \theta_{2}^{(3)} X^2$, we find that
\begin{eqnarray}
\dfrac{\partial^{|\alpha|}{f}}{\partial{(\theta_{1}^{(1)})^{\alpha_{1}}}\partial{(\theta_{1}^{(2)})^{\alpha_{2}}}\partial{(\theta_{2}^{(1)})^{\alpha_{3}}}\partial{(\theta_{2}^{(2)})^{\alpha_{4}}}\partial{(\theta_{2}^{(3)})^{\alpha_{5}}}}\parenth{Y|h_{1}(X,\theta_{1}),h_{2}(X,\theta_{2})} & = & \nonumber \\
& & \hspace{- 17 em} \dfrac{X^{2\alpha_{2}+\alpha_{4} + 2\alpha_{5}}}{2^{\alpha_{3} + \alpha_{4} + \alpha_{5}}}\dfrac{\partial^{\alpha_{1}+ \alpha_{2} + 2\alpha_{3}+2\alpha_{4} + 2\alpha_{5}}{f}}{\partial{h_{1}^{\alpha_{1}+ \alpha_{2} + 2\alpha_{3}+2\alpha_{4} + 2\alpha_{5}}}}(Y|h_{1}(X|\theta_{1}),h_{2}(X|\theta_{2})), \nonumber
\end{eqnarray}
for any $\alpha = (\alpha_{1}, \alpha_{2}, \alpha_{3}, \alpha_{4}, \alpha_{5}) \in \mathbb{N}^{5}$. Therefore, we can rewrite $A_{n}$ as follows:
\begin{eqnarray}
A_{n} & = & \nonumber \\
& &  \hspace{-4 em} \sum \limits_{i=1}^{k_{0}}\sum \limits_{j=1}^{s_{i}}p_{ij}^{n}\sum \limits_{1 \leq |\alpha| \leq \bar{r}} \dfrac{1}{\alpha!}\biggr\{(\Delta \theta_{1ij}^{n})^{(1)}\biggr\}^{\alpha_{1}} \biggr\{(\Delta \theta_{1ij}^{n})^{(2)}\biggr\}^{\alpha_{2}} \biggr\{(\Delta \theta_{2ij}^{n})^{(1)}\biggr\}^{\alpha_{3}}\biggr\{(\Delta \theta_{2ij}^{n})^{(2)}\biggr\}^{\alpha_{4}} \biggr\{(\Delta \theta_{2ij}^{n})^{(3)}\biggr\}^{\alpha_{5}} \nonumber \\
& & \times  \dfrac{X^{2\alpha_{2}+\alpha_{4} + 2\alpha_{5}}}{2^{\alpha_{3} + \alpha_{4} + \alpha_{5}}}\dfrac{\partial^{\alpha_{1}+ \alpha_{2} + 2\alpha_{3}+2\alpha_{4} + 2\alpha_{5}}{f}}{\partial{h_{1}^{\alpha_{1}+ \alpha_{2} + 2\alpha_{3}+2\alpha_{4} + 2\alpha_{5}}}}(Y|h_{1}(X|\theta_{1i}^{0}),h_{2}(X|\theta_{2i}^{0}))\overline{f}(X). \label{eqn:proof_Gaussian_family_with_offset_O2_third}
\end{eqnarray}
Similar to the proof of Theorem~\ref{theorem:lower_bound_Gaussian_family_without_offset_O2}, we define
\begin{eqnarray}
\mathcal{F} = \biggr\{X^{l_{1}}\dfrac{\partial^{l_{2}}{f}}{\partial{h_{1}^{l_{2}}}}(Y|h_{1}(X|\theta_{1i}^{0}),h_{2}(X|\theta_{2i}^{0}))\overline{f}(X): & & \nonumber \\
& & \hspace{-15 em} l_{1} = 2 \alpha_{2} + \alpha_{4} + 2 \alpha_{5}, \ l_{2} = \alpha_{1} + \alpha_{2} + 2 \alpha_{3} + 2\alpha_{4} + 2\alpha_{5}, \ 0 \leq |\alpha| \leq \bar{r}, \ 1 \leq i \leq k_{0} \biggr\}, \nonumber
\end{eqnarray}
Based on the proof argument similar to that of the claim in equation~\eqref{eqn:lemma_linear_independence_Gaussian_firstcase_first}, we can demonstrate that the elements of $\mathcal{F}$ are linearly independent with respect to $X$ and $Y$. Therefore, we can treat $A_{n}/D_{\kappa}(G_{n},G_{0})$, $B_{n}/D_{\kappa}(G_{n},G_{0})$ as a linear combination of elements of $\mathcal{F}$. We denote $F_{l_{1},l_{2}}(\theta_{1i}^{0},\theta_{2i}^{0})$ as the coefficient of $X^{l_{1}}\dfrac{\partial^{l_{2}}{f}}{\partial{h_{1}^{l_{2}}}}(Y|h_{1}(X|\theta_{1i}^{0}),h_{2}(X|\theta_{2i}^{0}))\overline{f}(X)$ in $A_{n}$ and $B_{n}$ for any $l_{1} = 2 \alpha_{2} + \alpha_{4} + 2 \alpha_{5}$, $l_{2} = \alpha_{1} + \alpha_{2} + 2 \alpha_{3} + 2\alpha_{4} + 2\alpha_{5}$, $0 \leq |\alpha| \leq \bar{r}$ and $1 \leq i \leq k_{0}$. Then, we can check that the coefficients of $X^{l_{1}}\dfrac{\partial^{l_{2}}{f}}{\partial{h_{1}^{l_{2}}}}(Y|h_{1}(X|\theta_{1i}^{0}),h_{2}(X|\theta_{2i}^{0}))\overline{f}(X)$ in $A_{n}/D_{\kappa}(G_{n},G_{0})$ and $B_{n}/D_{\kappa}(G_{n},G_{0})$ will be $F_{l_{1},l_{2}}(\theta_{1i}^{0},\theta_{2i}^{0})/D_{\kappa}(G_{n},G_{0})$. 

Assume that all of these coefficients go to 0 as $n \to \infty$. By taking the summation of $|F_{0,0}(\theta_{1i}^{0},\theta_{2i}^{0})/D_{\kappa}(G_{n},G_{0})|$ for all $1 \leq i \leq k_{0}$, we obtain that
\begin{eqnarray}
\biggr(\sum \limits_{i=1}^{k_{0}} |\sum \limits_{j=1}^{s_{i}}{p_{ij}^{n}} - \pi_{i}^{0}|\biggr)/D_{\kappa}(G_{n},G_{0}) \to 0. \label{eqn:proof_Gaussian_family_with_offset_O2_fourth}
\end{eqnarray}
When $l_{1} = 4$ and $l_{2} = 2$, we can check that $\alpha = (0, 2, 0, 0, 0)$ is the only solution to these equations. Therefore, the summation of $|F_{4,2}(\theta_{1i}^{0},\theta_{2i}^{0})/D_{\kappa}(G_{n},G_{0})|$ for all $1 \leq i \leq k_{0}$ leads to
\begin{align}
    \biggr(\sum \limits_{i=1}^{k_{0}}\sum \limits_{j=1}^{s_{i}}{p_{ij}^{n} \biggr|(\Delta \theta_{1ij}^{n})^{(2)}\biggr|^{2}}\biggr)/D_{\kappa}(G_{n},G_{0}) \to 0. \label{eqn:proof_Gaussian_family_with_offset_O2_fifth}
\end{align}
When $l_{1} = 2$ and $l_{2} = 4$, only $\alpha = (0, 0, 0, 2, 0)$ satisfies these equations. By taking into account all the coefficients $|F_{2,4}(\theta_{1i}^{0},\theta_{2i}^{0})/D_{\kappa}(G_{n},G_{0})|$ for all $1 \leq i \leq k_{0}$, we find that
\begin{align}
    \biggr(\sum \limits_{i=1}^{k_{0}}\sum \limits_{j=1}^{s_{i}}{p_{ij}^{n} \biggr|(\Delta \theta_{2ij}^{n})^{(2)}\biggr|^{2}}\biggr)/D_{\kappa}(G_{n},G_{0}) \to 0. \label{eqn:proof_Gaussian_family_with_offset_O2_sixth}
\end{align}
Similarly, when $l_{1} = 4$ and $l_{2} = 4$, we have $\alpha = (0, 0, 0, 0, 2)$ as the unique solution to these equations. By considering the summation of the coefficients $|F_{4,4}(\theta_{1i}^{0},\theta_{2i}^{0})/D_{\kappa}(G_{n},G_{0})|$ for all $1 \leq i \leq k_{0}$, we arrive at
\begin{align}
    \biggr(\sum \limits_{i=1}^{k_{0}}\sum \limits_{j=1}^{s_{i}}{p_{ij}^{n} \biggr|(\Delta \theta_{2ij}^{n})^{(3)}\biggr|^{2}}\biggr)/D_{\kappa}(G_{n},G_{0}) \to 0. \label{eqn:proof_Gaussian_family_with_offset_O2_seventh}
\end{align}
Combining the results from equations~\eqref{eqn:proof_Gaussian_family_with_offset_O2_fourth}-\eqref{eqn:proof_Gaussian_family_with_offset_O2_seventh}, we find that
\begin{align*}
    \frac{\sum \limits_{i=1}^{k_{0}}\sum \limits_{j=1}^{s_{i}}{p_{ij}^{n}\biggr(\biggr|(\Delta \theta_{1ij}^{n})^{(2)}\biggr|^{2}+\biggr|(\Delta \theta_{2ij}^{n})^{(2)}\biggr|^{2}+\biggr|(\Delta \theta_{2ij}^{n})^{(3)}\biggr|^{2} \biggr)}
+ \sum \limits_{i=1}^{k_{0}} |\sum \limits_{j=1}^{s_{i}}p_{ij}^{n} - \pi_{i}^{0}|}{D_{\kappa}(G_{n},G_{0})} \to 0.
\end{align*}
It indicates that
\begin{align*}
    \frac{\sum \limits_{i=1}^{k_{0}}\sum \limits_{j=1}^{s_{i}}{p_{ij}^{n}\biggr(\biggr|(\Delta \theta_{1ij}^{n})^{(1)}\biggr|^{\bar{r}}+\biggr|(\Delta \theta_{2ij}^{n})^{(1)}\biggr|^{\bar{r}/2} \biggr)}}{D_{\kappa}(G_{n},G_{0})} \to 1.
\end{align*}
Hence, we can find an index $i^{*} \in \{1, 2, \ldots, k_{0}\}$ such that 
\begin{align*}
    L = \frac{\sum \limits_{j=1}^{s_{i^{*}}}{p_{i^{*}j}^{n}\biggr(\biggr|(\Delta \theta_{1i^{*}j}^{n})^{(1)}\biggr|^{\bar{r}}+\biggr|(\Delta \theta_{2i^{*}j}^{n})^{(1)}\biggr|^{\bar{r}/2} \biggr)}}{D_{\kappa}(G_{n},G_{0})} \not \to 0
\end{align*}
as $n \to \infty$. By denoting $M_{l_{2}}(\theta_{1i}^{0}, \theta_{2i}^{0}) = F_{0, l_{2}}(\theta_{1i}^{0}, \theta_{2i}^{0})/ \sum \limits_{j=1}^{s_{i^{*}}}{p_{i^{*}j}^{n}\biggr(\biggr|(\Delta \theta_{1i^{*}j}^{n})^{(1)}\biggr|^{\bar{r}}+\biggr|(\Delta \theta_{2i^{*}j}^{n})^{(1)}\biggr|^{\bar{r}/2} \biggr)}$ for all $1 \leq i \leq k_{0}$, we obtain that
\begin{align*}
    M_{l_{2}}(\theta_{1i}^{0}, \theta_{2i}^{0}) = \frac{1}{L} \frac{F_{0, l_{2}}(\theta_{1i}^{0}, \theta_{2i}^{0})}{D_{\kappa}(G_{n},G_{0})} \to 0
\end{align*}
for any $1 \leq i \leq k_{0}$ and $0 \leq l_{2} \leq 2 \bar{r}$. From the formulation of $F_{0, l_{2}}(\theta_{1i}^{0}, \theta_{2i}^{0})$, the above limits with $M_{l_{2}}(\theta_{1i^{*}}^{0}, \theta_{2i^{*}}^{0})$ can be rewritten as:
\begin{eqnarray}
\dfrac{\sum \limits_{j=1}^{s_{i^{*}}}p_{i^{*}j}^{n}\sum \limits_{\substack{\alpha_{1}+2\alpha_{3}=l_{2} \\ \alpha_{1}+\alpha_{3} \leq \overline{r}}}{\dfrac{\biggr\{(\Delta \theta_{1i^{*}j}^{n})^{(1)}\biggr\}^{\alpha_{1}}\biggr\{(\Delta \theta_{2i^{*}j}^{n})^{(1)}\biggr\}^{\alpha_{3}}}{2^{\alpha_{3}}\alpha_{1}!\alpha_{3}!}}}{\sum \limits_{j=1}^{s_{i^{*}}}{p_{i^{*}j}^{n}\biggr(\biggr|(\Delta \theta_{1i^{*}j}^{n})^{(1)}\biggr|^{\bar{r}}+\biggr|(\Delta \theta_{2i^{*}j}^{n})^{(1)}\biggr|^{\bar{r}/2} \biggr)}} \to 0 \nonumber
\end{eqnarray}
for any $0 \leq \l_{2} \leq 2 \bar{r}$. According to the argument in Step 3 of the proof of Theorem~\ref{theorem:lower_bound_Gaussian_family}, that system of limits cannot happen. As a consequence, not all the coefficients in the linear combinations of $A_{n}/D_{\kappa}(G_{n},G_{0})$ and $B_{n}/D_{\kappa}(G_{n},G_{0})$ go to 0 as $n \to \infty$. From here, using the Fatou's argument in Step 4 of the proof of Theorem~\ref{theorem:total_variation_bound_over-fitted_MECFG} and the fact that the elements of $\mathcal{F}$ are linearly independent with respect to $X$ and $Y$, we achieve the conclusion of claim~\eqref{eqn:proof_Gaussian_family_O3_O2_O1_first}.
\subsubsection{Proof for equality~\eqref{eqn:proof_Gaussian_family_O3_O2_O1_second}}
\label{subsec:proof_Gaussian_family_O3_O2_O1_second}
Due to the similarity of this proof to the previous proofs, we will only provide a proof sketch of equality~\eqref{eqn:proof_Gaussian_family_O3_O2_O1_second}. We divide the proof into two settings of $\kappa' \prec \kappa = (\overline{r}, 2, \overline{r}/2, 2, 2)$.
\paragraph{Case 1:} $\kappa' = \parenth{\kappa'^{(1)}, \kappa'^{(2)}, \kappa'^{(3)}, \kappa'^{(4)}, \kappa'^{(5)}}$ when at least one of $\kappa'^{(2)}, \kappa'^{(4)}, \kappa'^{(5)} < 2$. Under this setting, we construct $G_{n} = \sum_{i = 1}^{k_{0}+1} \pi_{i}^{n}\delta_{(\theta_{1i}^{n},\theta_{2i}^{n})}$ such that $(\pi_{i}^{n},\theta_{1i}^{n}, \theta_{2i}^{n}) \equiv (\pi_{i-1}^{0}, \theta_{1(i-1)}^{0}, \theta_{2(i-1)}^{0})$ for $3 \leq i \leq k_{0}+1$. Additionally, $\pi_{1}^{n} = \pi_{2}^{n} = \pi_{1}^{0}/2$, $\parenth{(\theta_{1i}^{n})^{(1)},(\theta_{2i}^{n})^{(1)}} = \parenth{(\theta_{11}^{0})^{(1)},(\theta_{21}^{0})^{(1)}}$ for $1 \leq i \leq 2$, and $(\theta_{11}^{n})^{(2)} = (\theta_{11}^{0})^{(2)} - 1/n$, $(\theta_{21}^{n})^{(2)} = (\theta_{21}^{0})^{(2)} - 1/n$, $(\theta_{21}^{n})^{(3)} = (\theta_{21}^{0})^{(3)} - 1/n$, $(\theta_{12}^{n})^{(2)} = (\theta_{11}^{0})^{(2)} + 1/n$, $(\theta_{22}^{n})^{(2)} = (\theta_{21}^{0})^{(2)} + 1/n$, $(\theta_{22}^{n})^{(3)} = (\theta_{21}^{0})^{(3)} + 1/n$. From this construction of $G_{n}$, we can verify that $\widetilde{W}_{\kappa'}^{\|\kappa'\|_{\infty}}(G_{n},G_{0}) \asymp n^{-\min \{\kappa'^{(2)}, \kappa'^{(4)}, \kappa'^{(5)}\}} = o(n^{-2})$. Given that formulation of $G_{n}$, when we perform Taylor expansion up to the first order around $\theta_{11}^{0}, \theta_{21}^{0}$, the following equation holds
\begin{align}
p_{G_{n}}(X,Y) - p_{G_{0}}(X,Y) = \overline{R}_{1}(X,Y), \nonumber
\end{align}
where $\overline{R}_{1}(X,Y)$ is Taylor remainder such that
\begin{align}
\frac{h^{2}(p_{G_{n}}, p_{G_{0}})}{\widetilde{W}_{\kappa'}^{2\|\kappa'\|_{\infty}}(G_{n},G_{0})} \precsim \int \frac{\overline{R}_{1}^{2}(X,Y)}{p_{G_{0}}(X,Y)\widetilde{W}_{\kappa'}^{2\|\kappa'\|_{\infty}}(G_{n},G_{0})}d(X,Y) \precsim \frac{\mathcal{O}(n^{-4})}{n^{-2 \min \{\kappa'^{(2)}, \kappa'^{(4)}, \kappa'^{(5)}\}}} \to 0 \nonumber
\end{align}
as $n \to \infty$. Therefore, we achieve the conclusion of equality~\eqref{eqn:proof_Gaussian_family_O3_O2_O1_second} under Case 1. 
\paragraph{Case 2:} $\kappa' = \parenth{\kappa'^{(1)}, 2, \kappa'^{(3)}, 2, 2}$ when $(\kappa_{1}'^{(1)},\kappa'^{(3)}) \prec (\overline{r},\overline{r}/2)$. Under this setting of $\kappa'$, we construct $G_{n} = \sum_{i=1}^{k}\pi_{i}^{n}\delta_{(\theta_{1i}^{n},\theta_{2i}^{n})}$ such that $(\pi_{i+k-k_{0}}^{n}, \theta_{1(i+k-k_{0})}^{n}, \theta_{2(i+k-k_{0})}^{n}) = (\pi_{i}^{0},\theta_{1i}^{0},\theta_{2i}^{0})$ for $2 \leq i \leq k_{0}$. For $1 \leq j \leq k - k_{0} + 1$, we choose $(\theta_{1j}^{n})^{(2)} = (\theta_{11}^{0})^{(2)}, (\theta_{2j}^{n})^{(2)} = (\theta_{21}^{0})^{(2)}, (\theta_{2j}^{n})^{(3)} = (\theta_{21}^{0})^{(3)}$ and
\begin{align}
(\theta_{1j}^{n})^{(1)} = (\theta_{11}^{0})^{(1)} + \frac{a_{j}^{*}}{n}, \ (\theta_{2j}^{n})^{(1)} = (\theta_{21}^{0})^{(1)} + \frac{2b_{j}^{*}}{n^{2}}, \ \pi_{j}^{n} = \frac{\pi_{1}^{0}(c_{j}^{*})^{2}}{\sum_{i=1}^{k-k_{0}+1} (c_{j}^{*})^{2}}, \nonumber
\end{align}
where $(c_{i}^{*},a_{i}^{*},b_{i}^{*})_{i=1}^{k-k_{0}+1}$ are the non-trivial solution of system of polynomial equations~\eqref{eqn:system_polynomial_Gaussian_first} when $r = \overline{r} - 1$. From here, by performing Taylor expansion around $(\theta_{11}^{0},\theta_{21}^{0})$ and arguing similarly as Case 2 in the proof of equality~\eqref{eqn:proof_Gaussian_firstcase_second} in Section~\ref{subsection:key_equality_algebraic_dependent_linear_first}, we obtain that
\begin{align}
p_{G_{n}}(X,Y) - p_{G_{0}}(X,Y) = \sum \limits_{l = \overline{r}}^{2\overline{r}-2} \mathcal{O}(n^{-\overline{r}})\dfrac{\partial^{l}{f}}{\partial{h_{1}^{l}}}(Y|h_{1}(X,\theta_{1i}^{0}),h_{2}(X,\theta_{2i}^{0}))\overline{f}(X) + \overline{R}_{2}(X,Y), \nonumber 
\end{align}
where $\overline{R}_{2}(X,Y)$ is Taylor remainder such that the following limit holds
\begin{align}
\int \frac{\overline{R}_{2}^{2}(X,Y)}{p_{G_{0}}(X,Y)\widetilde{W}_{\kappa'}^{2\|\kappa'\|_{\infty}}(G_{n},G_{0})}d(X,Y) \precsim \frac{\mathcal{O}(n^{-2\overline{r}})}{n^{-2\min\{\kappa'^{(1)},\kappa'^{(3)}\}}} \to 0. \nonumber
\end{align}
Therefore, we achieve that
\begin{align}
h^{2}(p_{G_{n}}, p_{G_{0}})\bigg/\widetilde{W}_{\kappa'}^{2\|\kappa'\|_{\infty}}(G_{n},G_{0}) \to 0 \nonumber
\end{align}
as $n \to \infty$. As a consequence, we reach the conclusion of equality~\eqref{eqn:proof_Gaussian_family_O3_O2_O1_second} under Case 2. 
\section{Auxiliary results} \label{Section:Appendix_B}
\label{sec:auxi_result}
In this appendix, we provide two lemmas for the whole results in the paper. To streamline the discussion, we recall that $G_{0} = \sum \limits_{i=1}^{k_{0}}{\pi_{i}^{0}\delta_{(\theta_{1i}^{0},\theta_{2i}^{0})}}$ is the true mixing measure with exactly $k_{0}$ components such that $\theta_{ji}^{0} \in \Omega_{j}$ for all $1 \leq j \leq 2$ and $1 \leq i \leq k_{0}$ where $\Omega_{j} \subset \Rspace^{q_{j}}$ are compact sets for some given $q_{j} \geq 1$ as $1 \leq j \leq 2$. Furthermore, $\Omega = \Omega_{1} \times \Omega_{2}$. 
\begin{lemma} \label{lemma:relabel_sequence} Assume that $\kappa \in \mathbb{N}^{q_{1}+q_{2}}$ is a given vector order of generalized transportation distance and $k > k_{0}$. For any sequence $G_{n} \in \Ocal_{k}(\Omega)$ such that $\widetilde{W}_{\kappa}(G_{n},G_{0}) \to 0$ as $n \to \infty$, we can find a subsequence of $G_{n}$ (by which we replace by the whole sequence $G_{n}$ for the simplicity of presentation) that has the following properties:
\begin{itemize}
\item[(a)] (Fixed number of components) $G_{n}$ has exactly $\overline{k}$ number of components where $k_{0}+1 \leq \overline{k} \leq k$.
\item[(b)] (Universal representation) $G_{n}$ can be represented as:
\begin{align}
G_{n} = \sum \limits_{i=1}^{k_{0}+\overline{l}} \sum \limits_{j=1}^{s_{i}} p_{ij}^{n} \delta_{(\theta_{1ij}^{n},\theta_{2ij}^{n})}, \nonumber
\end{align}
where $\overline{l} \geq 0$ is some non-negative integer number and $s_{i} \geq 1$ for $1 \leq i \leq k_{0}+\overline{l}$ such that $\sum \limits_{i=1}^{k_{0}+\overline{l}} s_{i} = \overline{k}$. Furthermore, $(\theta_{1ij}^{n},\theta_{2ij}^{n}) \to (\theta_{1i}^{0},\theta_{2i}^{0})$ and $\sum \limits_{j=1}^{s_{i}} p_{ij}^{n} \to \pi_{i}^{0}$ for all $1 \leq i \leq k_{0}+\overline{l}$. Here, $\pi_{i}^{0} = 0$ as $k_{0} +1 \leq i \leq \overline{k}$ while $(\theta_{1i}^{0},\theta_{2i}^{0})$ are extra limit points from the convergence of components of $G_{n}$ as $k_{0} +1 \leq i \leq \overline{k}$.  
\end{itemize}
\end{lemma}
\begin{lemma} \label{lemma:generalized_Wasserstein_distance_polynomials}
Given the assumptions with $G_{0}$ and $G_{n}$ as those in Lemma~\ref{lemma:relabel_sequence}, we denote $\eta_{i}^{0} = \parenth{\theta_{1i}^{0},\theta_{2i}^{0}}$ and $\eta_{ij}^{n} = \parenth{\theta_{1ij}^{n},\theta_{2ij}^{n}}$ for $1 \leq i \leq k_{0}$ and $1 \leq j \leq s_{i}$. For any $\kappa \in \mathbb{N}^{q_{1}+q_{2}}$, we define the following distance:
\begin{align}
D_{\kappa}(G_{n},G_{0}) := \sum \limits_{i=1}^{k_{0} + \overline{l}}\sum \limits_{j=1}^{s_{i}} p_{ij}^{n} d_{\kappa}^{\|\kappa\|_{\infty}}\parenth{\eta_{ij}^{n},\eta_{i}^{0}} + \sum \limits_{i=1}^{k_{0}+ \overline{l}}\abss{\sum \limits_{j=1}^{s_{i}}p_{ij}^{n} - \pi_{i}^{0}} \nonumber
\end{align}
where the pseudo-metric $d_{\kappa}(.,.)$ is defined as in Section~\ref{subsection:generalized_Wasserstein_metric}. Then, the following holds:
\begin{align}
\widetilde{W}_{\kappa}^{\|\kappa\|_{\infty}} (G_{n},G_{0}) \precsim D_{\kappa}(G_{n},G_{0}). \nonumber
\end{align}
\end{lemma}
The proofs of the above lemmas are similar to those in~\citep{Ho-Nguyen-SIAM-18}; therefore, they are omitted.
\section{Convergence rate of density estimation} 
\label{Section:Appendix_C}
In this appendix, we provide a proof for convergence rate of density estimation of over-specified GMCF in Proposition~\ref{proposition:convergence_rates_density_estimation_MECFG}. Our proof technique follows standard result on density estimation for M-estimators in~\citep{Vandegeer-2000}. To ease the presentation, we adapt several notion from the empirical process theory into the setting of over-specified GMCF. 
\subsection{Key notation and results}
We denote $\mathcal{P}_{k}(\Omega) : = \{p_{G}(X,Y): \ G \in \Ocal_{k}(\Omega)\}$. Additionally, we define $N\parenth{\epsilon, \mathcal{P}_{k}(\Omega),\|.\|_{\infty}}$ as the covering number of metric space $\parenth{\mathcal{P}_{k}(\Omega),\|.\|_{\infty}}$ and $H_{B}(\epsilon, \mathcal{P}_{k}(\Omega), h)$ as the bracketing entropy of $\mathcal{P}_{k}(\Omega)$ under Hellinger distance $h$. We start with the following result regarding the upper bounds of these terms.
\begin{lemma} \label{lemma:entropy_control}
Suppose that $\Omega_{1}$ and $\Omega_{2}$ are respectively two bounded subsets of $\mathbb{R}^{q_{1}}$ and $\mathbb{R}^{q_{2}}$. Then, for any $0 < \epsilon < 1/2$, the following results hold
\begin{align}
\log N\parenth{\epsilon, \mathcal{P}_{k}(\Omega),\|.\|_{\infty}} & \precsim \log (1/\epsilon), \label{eqn:covering_control} \\
H_{B}(\epsilon, \mathcal{P}_{k}(\Omega), h) & \precsim \log(1/\epsilon). \label{eqn:bracketing_control}
\end{align}
\end{lemma}  
The detail proof of Lemma~\ref{lemma:entropy_control} is deferred to Appendix~\ref{subsection:proof_lemma_entropy_control}. To utilize the above bounds with covering number and bracketing entropy of $\mathcal{P}_{k}(\Omega)$, we will resort to Theorem 7.4 of~\cite{Vandegeer-2000} for density estimation with MLE. In particular, we denote the following key notation:
\begin{align}
\overline{\mathcal{P}}_{k}(\Omega) : = \{p_{(G+G_{0})/2}(X,Y): G \in \Ocal_{k}(\Omega)\}, \ \overline{\mathcal{P}}_{k}^{1/2}(\Omega) : = \{p_{(G+G_{0})/2}^{1/2}(X,Y): G \in \Ocal_{k}(\Omega)\}. \nonumber
\end{align}
For any $\delta > 0$, we define the Hellinger ball centered around $p_{G_{0}}(X,Y)$ and intersected with $\overline{\mathcal{P}}_{k}^{1/2}(\Omega)$ as follows:
\begin{align}
\overline{\mathcal{P}}_{k}^{1/2}(\Omega, \delta) : = \{f^{1/2} \in \overline{\mathcal{P}}_{k}^{1/2}(\Omega): h(f,p_{G_{0}}) \leq \delta \}. \nonumber
\end{align}
Furthermore, the size of this set can be captured by the following integral:
\begin{align}
\mathcal{J}_{B}\parenth{\delta, \overline{\mathcal{P}}_{k}^{1/2}(\Omega, \delta)} : = \int \limits_{\delta^2/2^{13}}^{\delta} H_{B}^{1/2}\parenth{u, \overline{\mathcal{P}}_{k}^{1/2}(\Omega, u), \|.\|_{2}}du \vee \delta. \nonumber
\end{align}
Equipped with the above notation, the results from Theorem 7.4 of~\cite{Vandegeer-2000} regarding convergence rates of density estimation from MLE can be formulated as follows.
\begin{theorem} \label{theorem:density_estimation_Vandegeer}
Take $\Psi(\delta) \geq \mathcal{J}_{B}\parenth{\delta, \overline{\mathcal{P}}_{k}^{1/2}(\Omega, \delta)}$ in such a way that $\Psi(\delta)/\delta^{2}$ is a non-increasing function of $\delta$. Then, for a universal constant $c$ and for
\begin{align}
\sqrt{n}\delta_{n}^{2} \geq c\Psi(\delta_{n}), \nonumber
\end{align}
we have for all $\delta \geq \delta_{n}$ that
\begin{align}
\Prob\parenth{h(p_{\widehat{G}_{n}}, p_{G_{0}}) > \delta} \leq c\exp\parenth{-\frac{n\delta^2}{c^2}}. \nonumber
\end{align}
\end{theorem}
\subsection{Proof for Proposition~\ref{proposition:convergence_rates_density_estimation_MECFG}}
Given Theorem~\ref{theorem:density_estimation_Vandegeer}, we are ready to finish the proof of Proposition~\ref{proposition:convergence_rates_density_estimation_MECFG}. In fact, we have 
\begin{align}
H_{B} \parenth{u, \overline{\mathcal{P}}_{k}^{1/2}(\Omega, u), \|.\|_{2}} \leq H_{B} \parenth{u, \mathcal{P}(\Omega, u), h},
\end{align}
for any $u > 0$. The above inequality leads to
\begin{align}
\mathcal{J}_{B}\parenth{\delta, \overline{\mathcal{P}}_{k}^{1/2}(\Omega, \delta)} & \leq \int \limits_{\delta^2/2^{13}}^{\delta} H_{B}^{1/2}\parenth{u, \mathcal{P}_{k}(\Omega, u), h}du \vee \delta \nonumber \\
& \precsim \int \limits_{\delta^2/2^{13}}^{\delta} \log(1/u)du \vee \delta, \nonumber
\end{align}
where the second inequality is due to the inequality~\eqref{eqn:bracketing_control} in Lemma~\ref{lemma:entropy_control}. Therefore, we can choose $\Psi(\delta) = \delta\parenth{\log(1/\delta}^{1/2}$ such that $\Psi(\delta) \geq \mathcal{J}_{B}\parenth{\delta, \overline{\mathcal{P}}_{k}^{1/2}(\Omega, \delta)}$. From here, with $\delta_{n} = \mathcal{O}\parenth{\brackets{\log n/ n}^{1/2}}$, the result of Theorem~\ref{theorem:density_estimation_Vandegeer} indicates that
\begin{align}
\Prob(h(p_{\widehat{G}_{n}},p_{G_{0}}) > C(\log n/n)^{1/2}) \precsim \exp(-c \log n) \nonumber
\end{align}
for some universal positive constants $C$ and $c$ that depend only on $\Omega$. As a consequence, we reach the conclusion of Proposition~\ref{proposition:convergence_rates_density_estimation_MECFG}. 
\subsection{Proof for Lemma~\ref{lemma:entropy_control}}
\label{subsection:proof_lemma_entropy_control}
The proof of the lemma follows the argument of Theorem 3.1 in~\citep{Ghosal-2001}. To facilitate the proof argument, our proof is divided into two parts.
\paragraph{Proof for covering number bound~\eqref{eqn:covering_control}:}
For any set $\mathcal{E}$, we denote $\mathcal{E}_{\epsilon}$ an $\epsilon$-net of $\mathcal{E}$ if each element of $\mathcal{E}$ is within $\epsilon$ distance from some elements of $\mathcal{E}_{\epsilon}$. Since $\Omega_{1}$ and $\Omega_{2}$ are two bounded subsets of $\Rspace^{q_{1}}$ and $\Rspace^{q_{2}}$ respectively, there exist corresponding $\epsilon$-nets $\overline{\Omega}_{1}(\epsilon)$ and $\overline{\Omega}_{2}(\epsilon)$ of these sets with $M_{1}$ and $M_{2}$ elements. We can validate that
\begin{align}
M_{1} \leq c_{1}(q_{1},k,\Omega_{1})\parenth{\frac{1}{\epsilon}}^{q_{1}k}, \ M_{2} \leq c_{2}(q_{2},k,\Omega_{2})\parenth{\frac{1}{\epsilon}}^{q_{2}k}, \nonumber
\end{align}
where $c_{i}(q_{i},k,\Omega_{i})$ are universal constants depending only on $q_{i}, k, \Omega_{i}$ for $1 \leq i \leq 2$. Furthermore, we denote $\Delta(\epsilon)$ an $\epsilon$-net for $k$-dimensional simplex. It is known that the cardinality of $\Delta(\epsilon)$ is upper bounded by $(5/\epsilon)^{k}$. We denote
\begin{align}
\mathcal{S} : = \{p_{G} \in \mathcal{P}_{k}(\Omega): \ \text{weights and components of} \ G \ \text{are on} \ \Delta(\epsilon) \times \overline{\Omega}_{1}(\epsilon) \times \overline{\Omega}_{2}(\epsilon)\}. \nonumber 
\end{align}
For each $p_{G} \in \mathcal{P}_{k}(\Omega)$ where $G = \sum_{i=1}^{k'} \pi_{i}\delta_{(\theta_{1i},\theta_{2i})}$ such that $k' \leq k$, we denote $\overline{G} = \sum_{i = 1}^{k'} \pi_{i}\delta_{(\theta_{1i}^{*},\theta_{2i}^{*})}$ such that $(\theta_{1i}^{*},\theta_{2i}^{*}) \in \overline{\Omega}_{1}(\epsilon) \times \overline{\Omega}_{2}(\epsilon)$ and $(\theta_{1i}^{*},\theta_{2i}^{*})$ are the closest points to $(\theta_{1i},\theta_{2i})$ in this set for $1 \leq i \leq k'$. Additionally, we denote $G^{*} = \sum_{i=1}^{k'} \pi_{i}^{*}\delta_{(\theta_{1i}^{*},\theta_{2i}^{*})}$ where $\pi_{i}^{*} \in \Delta(\epsilon)$ and $\pi^{*}$ are the closest points to $\pi_{i}$ in this set for $1 \leq i \leq k'$. From the formulation of $G^{*}$, it is clear that $p_{G^{*}} \in \mathcal{S}$. Invoking triangle inequality with sup-norm, the following inequality holds:
\begin{align}
\|p_{G}(X,Y) - p_{G_{*}}(X,Y)\|_{\infty} \leq \|p_{G}(X,Y) - p_{\overline{G}}(X,Y)\|_{\infty} + \|p_{\overline{G}}(X,Y) - p_{G_{*}}(X,Y)\|_{\infty}. \nonumber
\end{align}
According to the definition of $\overline{G}$ and $G^{*}$, direct computation leads to
\begin{align}
\|p_{\overline{G}}(X,Y) - p_{G_{*}}(X,Y)\|_{\infty} \leq \sum_{i=1}^{k'} \abss{\pi_{i}^{*} - \pi_{i}}\|f\parenth{Y|h_{1}(X,\theta_{1i}^{*}),h_{2}(X,\theta_{2i}^{*})}\overline{f}(X)\|_{\infty} \precsim \epsilon.
\end{align}
Furthermore, given the formulation of $\overline{G}$, we obtain that
\begin{align}
\|p_{G}(X,Y) - p_{\overline{G}}(X,Y)\|_{\infty} & \leq \sum \limits_{i=1}^{k'} \pi_{i}\|\overline{f}(X)\big[f\parenth{Y|h_{1}(X,\theta_{1i}^{*}),h_{2}(X,\theta_{2i}^{*})} \nonumber \\
& \hspace{6 em} - f\parenth{Y|h_{1}(X,\theta_{1i}),h_{2}(X,\theta_{2i})}\big]\|_{\infty} \nonumber \\
& \precsim \sum \limits_{i=1}^{k'} \pi_{i}\parenth{\enorm{\theta_{1i}^{*} - \theta_{1i}} + \enorm{\theta_{2i}^{*} - \theta_{2i}}} \precsim \epsilon, \nonumber
\end{align}
where the second inequality is due to the fact that the expert functions $h_{1}$ and $h_{2}$ are twice differentiable with respect to their parameters $\theta_{1}$ and $\theta_{2}$ and the space $\mathcal{X}$ is a bounded set. This inequality implies that the covering number for metric space $(\mathcal{P}_{k}(\Omega),\|.\|_{\infty})$ will be upper bounded by the cardinality of $\mathcal{S}$. More precisely, we obtain the following bound
\begin{align}
 N\parenth{\epsilon, \mathcal{P}_{k}(\Omega),\|.\|_{\infty}} \leq c_{1}(q_{1},k,\Omega_{1})c_{2}(q_{2},k,\Omega_{2})\parenth{\frac{5}{\epsilon}}^{k}\parenth{\frac{1}{\epsilon}}^{(q_{1}+q_{2})k}. \nonumber
\end{align}
Putting the above results together, we reach to the conclusion of the bound with covering number~\eqref{eqn:covering_control}. 
\paragraph{Proof for bracketing entropy control~\eqref{eqn:bracketing_control}:} Recall that, from the assumption with expert functions $h_{1}$ and $h_{2}$, we have $h_{1}(X,\theta_{1}) \in [-a, a]$ and $h_{2}(X,\theta_{2}) \in [\underline{\gamma},\overline{\gamma}]$ for all $X \in \mathcal{X}$, $\theta_{1} \in \Omega_{1}$, and $\theta_{2} \in \Omega_{2}$ where $a$ is some positive constant depending only on $\mathcal{X}$ and $\Omega_{1}$.

Now, let $\eta \leq \epsilon$ to be some positive number that we will chose later. From the formulation of univariate location-scale Gaussian distribution, we can check that
\begin{align}
f(Y|h_{1}(X,\theta_{1}),h_{2}(X,\theta_{2})) \leq \frac{1}{\sqrt{2\pi}\underline{\gamma}}\exp\parenth{-Y^2/(8\overline{\gamma}^2)}, \nonumber
\end{align}
for any $\abss{Y} \geq 2a$ and $X \in \mathcal{X}$. Therefore, if we define
\begin{align}
H(X,Y) = \begin{cases} \dfrac{1}{\sqrt{2\pi}\underline{\gamma}}\exp\parenth{-Y^2/(8\overline{\gamma}^2)}\overline{f}(X), & \text{for} \ \abss{Y} \geq 2a \\ \dfrac{1}{\sqrt{2\pi}\underline{\gamma}}\overline{f}(X), & \text{for} \ \abss{Y} < 2a, \end{cases}
\end{align}
then we can verify that that $H(X,Y)$ is an envelope of $\mathcal{P}_{k}(\Omega)$. We denote $g_{1},\ldots,g_{N}$ an $\eta$-net over $\mathcal{P}_{k}(\Omega)$. Then, we construct the brackets $[p_{i}^{L}(X,Y),p_{i}^{U}(X,Y)]$ as follows:
\begin{align}
p_{i}^{L}(X,Y) : = \max \{g_{i}(X,Y) - \eta, 0\}, \ p_{i}^{U}(X,Y) : = \max \{g_{i}(X,Y) + \eta, H(X,Y) \} \nonumber 
\end{align}
for $1 \leq i \leq N$. We can verify that $\mathcal{P}_{k}(\Omega) \subset \cup_{i=1}^{N} [p_{i}^{L}(X,Y), p_{i}^{U}(X,Y)]$ and $p_{i}^{U}(X,Y) - p_{i}^{L}(X,Y) \leq \min \{2\eta, H(X,Y)\}$. Direct computations lead to
\begin{align}
\int \parenth{p_{i}^{U}(X,Y) - p_{i}^{L}(X,Y)} d(X,Y) & \leq \int \limits_{\abss{Y} < 2a} \parenth{p_{i}^{U}(X,Y) - p_{i}^{L}(X,Y)} d(X,Y) \nonumber \\
& + \int \limits_{\abss{Y} \geq 2a} \parenth{p_{i}^{U}(X,Y) - p_{i}^{L}(X,Y)} d(X,Y) \nonumber \\
& \leq \overline{C} \eta + \exp\parenth{-\overline{C}^{2}/(2\overline{\gamma}^2)} \leq c\eta, \nonumber
\end{align}
where $\overline{C} = \max\{2a,\sqrt{8}\overline{\gamma}\}\log(1/\eta)$ and $c$ is some positive universal constant. The above bound leads to
\begin{align}
H_{B}(c\eta, \mathcal{P}_{k}(\Omega), \|.\|_{1}) \leq N \precsim \log(1/\eta). \nonumber
\end{align}
By choosing $\eta = \epsilon/c$, we have
\begin{align}
H_{B}(\epsilon, \mathcal{P}_{k}(\Omega), \|.\|_{1}) \precsim \log(1/\epsilon). \nonumber
\end{align}
Due to the inequality $h^{2} \leq \|.\|_{1}$ between Hellinger distance and total variational distance, we reach the conclusion of bracketing entropy bound~\eqref{eqn:bracketing_control}.
\bibliography{Nhat,NPB,Nguyen}
\end{document}